\documentclass[a4paper,10pt]{article}
\usepackage{amsmath,amssymb,amsthm}
\usepackage{dsfont}
\usepackage{mathrsfs}
\usepackage[mathscr]{euscript}

\newtheorem{theo}{Theorem}[section] 
\newtheorem{defi}[theo]{Definition}
\newtheorem{lemm}[theo]{Lemma} 
\newtheorem{prop}[theo]{Proposition}

\newcommand{\Na}{\mathbb N}                   

\newcommand{\Za}{\mathbb Z}                   

\newcommand{\Ra}{\mathbb R}                   

\newcommand{\Ca}{\mathbb C}                   

\newcommand{\scal}[1]{\langle #1 \rangle}

\newcommand{\finpreuve}{\hfill $\Box$}

\newcommand{\name}{$\underline{\qquad \qquad}$}

\begin{document}

\title{\sc  Global in time Strichartz inequalities on  asymptotically flat manifolds with temperate trapping}

\author{Jean-Marc Bouclet    \and Haruya Mizutani  }


\date{ }

\maketitle

{\abstract We prove global Strichartz inequalities for the Schr\"odinger equation on a large class of asymptotically conical manifolds. Letting $ P $ be the nonnegative Laplace operator and $ f_0 \in C_0^{\infty}(\Ra) $ be a smooth cutoff equal to $1$ near zero, we show first that the low frequency part of any solution $ e^{-itP} u_0 $, i.e. $ f_0 (P) e^{-itP} u_0 $, enjoys the same global Strichartz estimates as on $ \Ra^n $ in dimension $ n \geq 3 $. We also show that the high energy part $ (1-f_0)(P) e^{-itP} u_0$ also satisfies global Strichartz estimates without loss of derivatives outside a compact set, even if the manifold has trapped geodesics but in a temperate sense. We then show that the full solution $ e^{-itP}u_0 $ satisfies global space-time Strichartz estimates if the trapped set is empty or sufficiently filamentary, and we derive a scattering theory for the $ L^2 $ critical nonlinear Schr\"odinger equation in this geometric framework.}




\section{Introduction and main results}
\setcounter{equation}{0}

In the past ten or fifteen years, a lot of activity has been devoted to study Strichartz inequalities on manifolds. We recall that these inequalities were stated first on $ \Ra^n $ for the wave  equation \cite{Stri} and then  the Schr\"odinger one  \cite{GiVe}; for the Schr\"odinger equation and a pair $ (p,q) \in [2,\infty] \times [2,\infty] $, they read
$$ || u ||_{L^p (\Ra,L^q)} \lesssim || u_0 ||_{L^2} , \qquad u (t) = e^{it \Delta} u_0 , \qquad \mbox{if} \ \  \frac{2}{p} + \frac{n}{q}  = \frac{n}{2}, \ \ (n,p,q) \ne (2,2,\infty) . $$
(A pair $ (p,q) $ satisfying the last two conditions is called Schr\"odinger admissible.) The strong interest on Strichartz inequalities is mainly related to their key role  in the study of nonlinear dispersive equations (see {\it e.g.} \cite{Cazenavebook,Taobook}). 

  On compact manifolds these estimates may be different as those on $  \Ra^n$, either due to the strong confinment leading to derivative losses  for the Sch\"odinger equation \cite{BGT} (the $ L^2 $ norm of initial data is replaced by some Sobolev norm) or to the absence of global in time estimates (if initial data are eigenfunctions the solutions are periodic in time). 
  
  One may ask to which extent the estimates  on $ \Ra^n $ still hold on noncompact manifolds,  at least in the  class of asymptotically flat ones.  For the Schr\"odinger equation, the only one considered from now on, this problem was considered in several articles for local in time estimates \cite{StTa,RZ,HTW,BoucletTzvetkov1,MizutaniCPDE}. From the geometrical point of view, those papers consider stronger and stronger perturbations, namely from compactly supported perturbations of the flat metric on $ \Ra^n $ to long range perturbations of conical metrics on manifolds. We refer to Definition \ref{defasympconique} for a description of long range asymptotically conical metrics but point out here that long range perturbations are natural in  that it is the only type of decay which is invariant under a change of radial coordinates (see \cite{BoucletOsaka}).

  

  Global in time estimates for long range perturbations are considerably more delicate to obtain and have been considered in fewer papers \cite{Tataru1,MarzMetcTata,HassellZhang} (see also \cite{BoucletTzvetkov2} with a low frequency cutoff).
  

To prove global Strichartz inequalities on curved backgrounds, one has to face two difficulties. The first one, which does not happen on $ \Ra^n $, is the possible occuring of trapped geodesics (geodesics not escaping to infinity, in the future or in the past). This trapping  is only sensitive at high frequencies and may affect the estimates by a loss of derivatives. However, if it is sufficiently weak, one can still expect Strichartz estimates without loss as shown in \cite{BGH} locally in time.  Trapping is already a problem for local in time estimates hence a fortiori for global in time ones. 

The second difficulty stems in  the analysis of low frequencies. Indeed, except in a few model situations such as $ \Ra^n $ or flat cones \cite{Ford}   where the fundamental solution of the Schr\"odinger equation can be computed explicitly, the only robust  strategy accessible so far is  to  localize  the solution in frequency, {\it e.g.} by mean of a Littlewood-Paley decomposition, and then to prove Strichartz estimates for the spectrally localized components by using microlocal techniques to derive appropriate dispersive estimates. Due to the uncertainty principle, low frequency data cannot be studied purely by microlocal techniques and thus  require additional non trivial estimates.   On $ \Ra^n $ (or a pure cone), one may use a global scaling argument to reduce the analysis of low frequency blocks to the study at frequency one, but this is in general impossible on manifolds.

  The first breakthrought on global in time Strichartz estimates was done by Tataru in \cite{Tataru1} where he considered long range and globally small perturbations of the Euclidean metric, with $ C^2 $ and time dependent coefficients. In  this framework, no trapping could occur. The results were then improved in \cite{MarzMetcTata} by allowing more general perturbations in a compact set, including some weak trapping. Recently, Hassell-Zhang \cite{HassellZhang}  partially extended those results by  considering the general geometric framework of asymptotically conic manifolds and including very short range potentials, but using a non trapping condition.  
  
  In the present paper, we improve on those references in the following directions. On one hand, we  consider a class of asymptotically conic manifolds which is larger than the one of Hassell-Zhang, and contains all usual smooth long range perturbations of the Euclidean metric. More importantly, we allow the possibility to have trapped trajectories  and, assuming this trapping to be temperate (assumption (\ref{aprioriresolvente})), show that the solutions to the linear Schr\"odinger equation enjoy the same global in time estimates without loss as on $ \Ra^n $ outside a large enough compact set. This fact is a priori not clear at all since, by the infinite speed of propagation of the Schr\"odinger equation,  one may fear that the geometry and the form of the initial datum inside a compact set has an influence on the solution all the way to spatial infinity.  This question was  considered first in \cite{BoucletTzvetkov1} locally in time and then in \cite{MarzMetcTata} globally in time case but our approach in this paper allows to deal with much stronger types of trapping than in this last reference (see the discussion after Theorem \ref{theohighinfty}).

  As a byproduct of this analysis, we derive global space-time Strichartz estimates without loss if there is no trapping (thus recovering the results of Hassell-Zhang for a larger class of manifolds, when there is no potential) or if the trapping is filamentary in the sense of \cite{NZ,BGH}. In particular, we extend to the global in time case one of the results of \cite{BGH}.
  
    Then, we apply these estimates to the  scattering theory of the $L^2$ critical nonlinear Schr\"odinger equation with small data on a manifold with filamentary (or empty) trapped set (Theorem \ref{theoremenonlineaire}).

    From the technical point of view, an important part of our paper is  devoted to construct tools adapted to the analysis of low frequencies. In particular, along the way, we develop  a new version of the Isozaki-Kitada parametrix for long range metrics. Recall that the Isozaki-Kitada parametrix was introduced on $ \Ra^n $ to study the scattering theory of Schr\"odinger operators with long range potentials \cite{IsKi}. One of the new features of our parametrix is the treatment of low frequencies which, to our knowledge, does not seem to have been much considered before, up to the reference \cite{DeSk} in the context of scattering by potentials on $ \Ra^n $ which is very different from ours (especially at  low energy). We derive related $ L^2 $ propagation estimates which are needed in the present paper but can be of interest for other questions of scattering theory, such as the study of scattering matrices at low energy. In a more directly oriented PDE perspective, the methods developed in this paper also allow  to handle other dispersive models like the wave or Klein-Gordon equations \cite{Duong}.

Let us now state our results more precisely.

Let $ ({\mathcal M},G) $ be an asymptotically conic manifold, possibly with a boundary, {\it i.e.} a manifold diffeomorphic away from a compact set to a product $ (R_{\mathcal M},+\infty) \times {\mathcal S} $, for some closed Riemannian manifold $ ({\mathcal S}, \bar{g} ) $, such that $ G $ is a long range perturbation of the exact conical metric $ dr^2 + r^2 \bar{g} $. To state a precise definition, we  denote by $ \Gamma (T_q^p {\mathcal S}) $ the space of $ (p,q) $ tensors on $ {\mathcal S} $, {\it i.e.} sections of $ (\otimes^{p} T {\mathcal S} ) \otimes  (\otimes^q T^* {\mathcal S}) $, and for a given  smooth map $ e = e (r) $ defined on  $ (R_{\mathcal M}, + \infty) $  with values in $   \Gamma (T_q^p {\mathcal S}) $, we will note
$$ e \in S^{- \nu} \qquad \Longleftrightarrow \qquad N_{pq} \big( \partial_r^j e (r) \big) \lesssim \scal{r}^{-\nu-j}  \qquad \mbox{ for each semi-norm }  N_{pq}  \mbox{ of }  \Gamma (T_q^p {\mathcal S}) \mbox{ and } j \geq 0 . $$ 
If $ (\theta_1 , \ldots , \theta_{n-1}) $ are local coordinates on $ {\mathcal S} $, this means equivalently that $e$ is a linear combination of terms of the form
$  e_{i_1 \cdots i_q}^{j_1 \cdots j_p} (r,\theta) d \theta_{i_1} \otimes \cdots \otimes d \theta_{i_q} \otimes \partial_{\theta_{j_1}} \otimes \cdots \otimes \partial_{\theta_{j_p}} $ such that, for each $j$ and $ \alpha $, we have an estimate $| \partial_r^j \partial_{\theta}^{\alpha} e_{i_1 \cdots i_q}^{j_1 \cdots j_p} (r,\theta) |\lesssim \scal{r}^{-\nu-j}  $ locally uniformly in $ \theta $. Here  $ \scal{\cdot} $ is the standard japanese bracket.
 

\begin{defi} \label{defasympconique} A Riemannian manifold $ ({\mathcal M},G) $ is {\bf asymptotically conic} if there exists a continuous  and proper function $ r : {\mathcal M} \rightarrow [ 0 , + \infty ) $, a compact subset $ {\mathcal K} \Subset {\mathcal M} $ and a closed Riemannian manifold $ ({\mathcal S}, \bar{g} ) $ such that for some $ R_{\mathcal M} > 0 $ there is a diffeomorphism 
$$ \Omega : {\mathcal M} \setminus {\mathcal K} \ni m \mapsto  \big( r (m) , \omega (m) \big) \in (R_{\mathcal M} , + \infty) \times {\mathcal S}  $$
through which
$$ G = \Omega^* \big( A(r)dr^2 + 2 r B (r) dr + r^2 g (r) \big) $$
where $ A (r) \in \Gamma (T_0^0 {\mathcal S}) $, $ B (r) \in \Gamma (T_1^0 {\mathcal S}) $ and $ g (r) \in \Gamma (T_2^0{\mathcal S} )$ is a Riemannian metric on $ {\mathcal S} $ 
such that, for some $ \nu > 0 $,
\begin{eqnarray}
 A - 1 \in S^{-\nu}, \qquad B \in S^{-\nu}, \qquad g(\cdot)- \bar{g} \in S^{-\nu} .   \label{metriqueconiquegenerale}
\end{eqnarray}
If $ A \equiv 1 $ and $ B \equiv 0 $, one says the metric $ G $ is in {\bf normal form}.
\end{defi}
Without loss of generality, we will assume that $ G $ is in normal form (see Appendix \ref{appendiceformenormale}). This plays no role in the present introduction but will be useful in later sections.

 Everywhere in the sequel, we denote by $ L^q ({\mathcal M}) $ or just $ L^q $ the Lebesgue spaces associated to the  Riemannian measure on $ {\mathcal M} $. We let $ P $ be the Friedrichs extension of $ - \Delta_G $ on $ L^2 ({\mathcal M}) $, namely the unique selfadjoint realization if $ {\mathcal M} $ has no boundary or the Dirichlet one if $ \partial {\mathcal M} $ is not empty. One interest of our geometric framework is that, if $ n \geq 3 $, we have a Sobolev estimate
 \begin{eqnarray}
 ||v||_{L^{2^*}({\mathcal M})} \leq C || P^{1/2} v ||_{L^2 ({\mathcal M})}, \qquad 2^* = \frac{2n}{n-2} , \label{Sobolevestimate}
 \end{eqnarray}
 for all $ v $ in the domain of $ P^{1/2} $ (see Appendix \ref{sectionSobolevhomogene} for a proof).

For $ u_0 \in L^2 ({\mathcal M}) $, we let $  u (t) := e^{-itP}u_0 $ be the solution to the Schr\"odinger equation
$$ i \partial_t u - P u = 0, \qquad u_{| t= 0} = u_0 . $$
Let $ f_0 \in C_0^{\infty}(\Ra) $ be such that $ f_0 \equiv 1 $ on $ [-1,1] $ and split $ u (t) =  u_{\rm low} (t) + u_{\rm high} (t) $ 
according to  low and high frequencies, {\it i.e.}
\begin{eqnarray}
 u_{\rm low} (t) := f_0 (P) e^{-itP}u_0, \qquad u_{\rm high}(t) = (1-f_0)(P) e^{-itP} u_0 . \label{defhighlow}
\end{eqnarray}
\begin{theo} \label{theoremlow}[Global space-time low frequency estimates]  Assume that   $ n \geq 3 $ and let $ (p,q) $ be a Schr\"odinger admissible pair. Then there exists $ C > 0 $ such that, for all $ u_0 \in L^2 ({\mathcal M}) $,
\begin{eqnarray}
 \big| \big| u_{\rm low} \big| \big|_{L^p (\Ra;L^q({\mathcal M}))} \leq C ||u_0 ||_{L^2 ({\mathcal M})} . \label{highfrequencyestimate}
\end{eqnarray} 

\end{theo}

Notice that in this theorem $ \partial {\mathcal M} $ may be empty or not.

\bigskip

\noindent {\it Proof.} Paragraph \ref{proofTheoremlow}.

\bigskip

\begin{theo} \label{theohighinfty}[Global in time high frequency estimates at spatial infinity] Assume that $ n \geq 2 $ and that for some $  M > 0 $ large enough, we have for all $ \chi \in C_c^{\infty} ({\mathcal M}) $
\begin{eqnarray}
\big| \big|  \chi (P- \lambda \pm i 0)^{-1} \chi \big| \big|_{L^2({\mathcal M})\rightarrow L^2 ({\mathcal M})} \lesssim_{\chi} \lambda^{M}, \qquad \lambda \geq  1 . \label{aprioriresolvente}
\end{eqnarray}
Then there exists $ R \gg 1 $  such that for any Schr\"odinger admissible pair $ (p,q) $ there exists $ C> 0 $ such that
\begin{eqnarray}
 \big| \big| {\mathds 1}_{\{r> R \}} u_{\rm high} \big| \big|_{L^p (\Ra;L^q({\mathcal M}))} \leq C ||u_0 ||_{L^2 ({\mathcal M})} , \label{lowfrequencyestimate}
\end{eqnarray} 
for all $ u_0 \in L^2 ({\mathcal M}) $.
\end{theo}

If we recast the global in time estimates at spatial infinity of \cite[Theorem 1.5]{MarzMetcTata} in our framework, these authors show that 
$$ || {\mathds 1}_{\{r> R \}} u_{\rm high} ||_{L^p (\Ra;L^q)} \leq C ||u_0 ||_{L^2 } +  || {\mathds 1}_{\{r < R \}} u_{\rm high} ||_{L^2 (\Ra;L^2)}   $$
where the last term can be  controlled by $ || u_0 ||_{L^2} $ thanks to (\ref{aprioriresolvente})  if $ M \leq  0$ (the usual non trapping case is $ M = -1/2 $) but not clearly otherwise. In our result,  the right hand side of (\ref{lowfrequencyestimate}) does not involve any corrective term depending on $u$ and holds for any $ M $.

 Note that examples of situations where  bounds of the form (\ref{aprioriresolvente})  hold include  \cite{NZ,ChWu} in some  trapping geometries and, of course,  the nontrapping case \cite{Vodev}. We point out that, without any dynamical assumption, the upper bounds on the high energy resolvent are in general of order $ O (\exp (C \lambda^{1/2})) $ (see \cite{CardosoVodev} and the references therein); it would be interesting to know wether (\ref{lowfrequencyestimate}) persists or fails in this most general case. Our method allows to deal with any polynomial growth in $ \lambda^{1/2} $ (hence the name of {\it temperate} trapping), but not clearly with an exponential one.


We also remark that, as in Theorem \ref{theoremlow}, the boundary of $ {\mathcal M} $ does not need to be empty but  this observation is less relevant here for we consider  estimates near infinity.

 Theorems \ref{theoremlow} and \ref{theohighinfty} reduce the proof of Strichartz estimates on $u$ to estimates on  $  {\mathds 1}_{\{r \leq R \}} u_{\rm high} $. This leads to the following result.

\begin{theo}[Global spacetime estimates without loss] \label{globalestimates} Assume that $ n \geq 3 $ and $ \partial {\mathcal M} $ is empty. If either
\begin{itemize}
\item{ the geodesic flow is non trapping and $(p,q)$ is any Schr\"odinger admissible pair,}
\item{the trapped set satisfies the assumptions of \cite{BGH} and $(p,q)$ is any non endpoint Schr\"odinger admissible pair,}
\end{itemize}
  then there exists $ C > 0 $ such that 
\begin{eqnarray}
 \big| \big| u \big| \big|_{L^p (\Ra;L^q({\mathcal M}))} \leq C ||u_0 ||_{L^2 ({\mathcal M})} , \label{estimatesansperte}
\end{eqnarray} 
for all $ u_0 \in L^2 ({\mathcal M}) $.
\end{theo}

This theorem improves on the result of \cite{HassellZhang} in two directions: Hassell-Zhang only consider the nontrapping case and, even in the nontrapping situation, we consider more general types of ends.  It also provides a global in time version of the estimates of \cite{BGH} in the asymptotically conic case.

We state this result in the boundaryless case in order to give complete proofs or references. We emphasize however that using the techniques of \cite{Ivanovici} it can certainly be extended to the case when $ {\mathcal M} $ has a stricly geodesically concave boundary and is non trapping for the associated billiard flow. 

We  recall finally the well known fact  that inhomogeneous Strichartz estimates, for non endpoint pairs, can be derived from the homogeneous ones (\ref{estimatesansperte})  by using the Christ-Kiselev Lemma \cite{ChKi}; this is sufficient for the applications to the nonlinear equations studied in Section \ref{sectionnonlineaire}.

\medskip

Here is the plan of our paper. In Section \ref{sectionNotation}, we record notation about charts, partitions of unity, scaling operators, etc. that will be used in further sections. In Section \ref{sectionpseudo}, we describe the pseudo-differential calculus adapted to our framework, including a rescaled one for low frequency estimates which is not quite standard. In Section \ref{sectionLittlewoodPaley}, we prove Littlewood-Paley decompositions at low and high frequencies. In Sections \ref{sectionscatteringclassique} and \ref{sectionIsozakiKitada}, we construct an Isozaki-Kitada parametrix for the microlocalized Schr\"odinger group, both at high and low frequencies. We use it in Section \ref{sectionpropagation} to derive some $ L^2 $ propagation estimates to be used in Section \ref{SectionStrichartz} where the theorems stated in this introduction are proved. Finally, in Section \ref{sectionnonlineaire}, we give nonlinear applications of our Strichartz estimates.

\medskip

\noindent {\bf Acknowledgments.} JMB is partially supported by ANR Grant  GeRaSic, ANR-13-BS01-0007-01. HM is partially supported by JSPS Wakate (B) 25800083. 

\section{Notation} \label{sectionNotation}
\setcounter{equation}{0}

In this  section, we collect some notation or definitions that will be used throughout this paper.

\medskip 

\noindent {\bf Coordinate charts.} If $ \kappa : U_{\kappa} \subset {\mathcal S} \rightarrow V_{\kappa} \subset \Ra^{n-1} $  is a coordinate chart on $ {\mathcal S} $ then, upon the identification of $ (R_{\mathcal M},+\infty) \times U_{\kappa} $ with a subset of $ {\mathcal M} $, $ (r,\omega)  \mapsto (r,\kappa(\omega))  $ defines a coordinate chart on $ {\mathcal M} $. 
We define $ \Pi_{\kappa} $ and $ \Pi_{\kappa}^{-1} $ respectively as the pullback and pushforward operators associated to this chart on $ {\mathcal M} $, {\it i.e.}
\begin{eqnarray}
\big(\Pi_{\kappa} v \big)(r,\omega) = v \big(r, \kappa(\omega) \big), \qquad   \big( \Pi_{\kappa}^{-1} u \big) (r,\theta) = u \big( r , \kappa^{-1} (\theta) \big)  . \label{notationpullpush}
\end{eqnarray}
If $ \tau : V_1 \rightarrow V_2  $ is a diffeomorphism between  open subsets of $ \Ra^{n-1}$ (typically a transition map between charts of $ {\mathcal S} $), we also define $ \Pi_{\tau} $ and $ \Pi_{\tau}^{-1} $ as above for the  diffeomorphism $ (r,\theta) \mapsto (r,\tau(\theta))  $ between $  \Ra \times V_1 \rightarrow \Ra \times V_2 $. With such a definition, if $  \kappa_j :  U_j \rightarrow V_j $, $ j = 1,2$,  are two coordinates charts on $ {\mathcal S} $, it follows that
\begin{eqnarray}
 \Pi_{\kappa_2 }^{-1} \Pi_{\kappa_1} = \Pi_{\tau_{12} }^{-1}, \qquad \tau_{12} := \kappa_2 \circ \kappa^{-1}_1 : \kappa_1 (U_1 \cap U_2) \rightarrow \kappa_2 (U_1 \cap U_2) .  \label{notationtransition}
\end{eqnarray} 
We choose a finite atlas on  $ {\mathcal S} $ composed of charts  with the property that  $ \kappa^* \bar{g} =: \bar{g}_{lm}(\theta) d \theta_l d \theta_m $ satisfies the following uniform estimates on each $ V_{\kappa} $: 
\begin{eqnarray}
C^{-1}_0 I_{n-1} \leq \big(  \bar{g}_{l m}(\theta) \big) \leq C_0 I_{n-1} , \label{hjkellipticity} \\
 \big| \partial^{\alpha} \bar{g}_{l m}(\theta) \big| \leq C_{\alpha}  .  \label{hjksmooth}
\end{eqnarray}
We will also use the matrices $ \bar{g} (\theta) := (\bar{g}_{lm}(\theta)) $, $ (\bar{g}^{lm}(\theta)) := \bar{g} (\theta)^{-1} $ as well as the function $ |\bar{g}(\theta)|:= \mbox{det} \bar{g}(\theta)^{1/2} $.

\medskip

\noindent {\bf Partitions of unity.}  We pick  a partition of unity $
1 = \sum_{\kappa} \varphi_{\kappa} (\omega) $  on $ {\mathcal S} $,  with $ \varphi_{\kappa} \in C_0^{\infty} \big( U_{\kappa} \big) $ and where  the sum  over $ \kappa $, as well as all similar sums below, is taken over the finite atlas we chose above.  For each $\kappa$, we also pick   $ \tilde{\varphi}_{\kappa} , \tilde{\tilde{\varphi}}_{\kappa} \in C_0^{\infty}  \big( U_{\kappa} \big) $ such that $ \tilde{\varphi}_{\kappa} \equiv 1 $ near $  \mbox{supp} \big( \varphi_{\kappa} \big) $ and $ \tilde{\tilde{\varphi}}_{\kappa} \equiv 1 $ near $ \mbox{supp}(\tilde{\varphi}_{\kappa}) $. We then pick $ \zeta , \tilde{\zeta} , \tilde{\tilde{\zeta}} \in C^{\infty} (\Ra) $ supported in $ (R_{\mathcal M}, \infty) $, equal to $1$ near infinity and such that $ \tilde{\zeta} \equiv 1 $ near the support of $ \zeta $, $ \tilde{\tilde{\zeta}} \equiv 1 $ near the support of $ \tilde{\zeta} $ and define
\begin{eqnarray}
 \psi_{\kappa} (r,\omega) : =  \zeta(r) \varphi_{\kappa}(\omega), \qquad \tilde{\psi}_{\kappa} (r,\omega) : =  \tilde{\zeta}(r) \tilde{\varphi}_{\kappa}(\omega) , \qquad
 \tilde{\tilde{\psi}}_{\kappa} (r,\omega) : =  \tilde{\tilde{\zeta}}(r) \tilde{\tilde{\varphi}}_{\kappa}(\omega) .  \label{notationpartition}
\end{eqnarray} 
Their interest is that they are supported on coordinate patches of $ {\mathcal M} $ and that
\begin{eqnarray}
 \sum_{\kappa} \psi_{\kappa} = \zeta (r) \equiv 1 \ \ \mbox{near infinity}, \qquad \tilde{\psi}_{\kappa} \equiv 1 \ \ \mbox{near} \ \ \mbox{supp}(\psi_{\kappa}) , \qquad \tilde{\tilde{\psi}}_{\kappa} \equiv 1 \ \ \mbox{near} \ \ \mbox{supp}(\tilde{\psi}_{\kappa}) .  \label{partitionfixee}
 \end{eqnarray}
 They will be useful to globalize pseudo-differential operators on $ {\mathcal M} $.

\medskip

\noindent {\bf Rescaling operators at infinity.} For $ \epsilon \in (0,1] $, we will use the operators $ {\mathscr D}_{\epsilon} $ defined by
\begin{eqnarray}
 {\mathscr D}_{\epsilon} v (r,\omega) = \epsilon^{\frac{n}{2}} v (\epsilon r , \omega) , \qquad \mbox{if supp}(v) \subset \{ r > R_{\mathcal M} \}   .  \label{definitrescaling}
\end{eqnarray}
Here $v$ is a function on $ {\mathcal M} $ but we will also freely use $ {\mathscr D}_{\epsilon} $ for functions on $ \Ra^n $ supported in $ (R_{\mathcal M}, \infty) \times V $, for any $ V \subset \Ra^{n-1} $. Note that $ {\mathscr D}_{\epsilon} v $ is supported in $  \{ r > \epsilon^{-1} R_{\mathcal M} \} $. The normalization factor $ \epsilon^{n/2} $ ensures that
$$ || {\mathscr D}_{\epsilon} v ||_{L^2({\mathcal M})} \approx ||v||_{L^2({\mathcal M})} $$
({\it i.e.} their quotient is bounded from above and below uniformly in $ \epsilon $)
since  the measure in $ \{ r > R_{\mathcal M} \} $ is comparable to the exact conic measure $ r^{n-1} |\bar{g}(\theta)| dr d \theta $.
We define similarly
\begin{eqnarray}
 {\mathscr D}_{\epsilon}^{-1} v (r,\omega) = \epsilon^{-\frac{n}{2}} v (\epsilon^{-1} r , \omega) , \qquad \mbox{if supp}(v) \subset \{ r > \epsilon^{-1} R_{\mathcal M} \}   .  \label{definitrescalinginverse}
\end{eqnarray}
Of course we have also the equivalence $  || {\mathscr D}_{\epsilon}^{-1} v ||_{L^2({\mathcal M})} \approx ||v||_{L^2({\mathcal M})}  $.

\medskip

\noindent {\bf Modified japanese bracket.}
Everywhere in the text, we will replace the usual japanese bracket $ \scal{r} = (1+r^2)^{1/2} $ by another positive function still denoted by $ \scal{r} $ and such that
\begin{eqnarray}
 \scal{r} = \begin{cases} 1 & \mbox{on a large enough compact set}  \\ r & \mbox{for } r \gg 1 . \end{cases}  \label{scalrM}
 \end{eqnarray}
 By large enough compact set, we mean that $ \scal{r} = 1 $ in a neighborhood of the region where $ \zeta (r) \ne 1 $ (see {\it e.g.} (\ref{partitionfixee}) for $ \zeta $). The interest is that commutators of powers of $ \scal{r}$ with differential operators will be automatically supported in a region where $ \zeta (r) = 1 $. More generally, commutators with powers of $ \scal{\epsilon r} $ will be supported where $ \zeta (\epsilon r) = 1 $.

\medskip

\noindent {\bf Laplacian.}  With the metric  in  normal form, the operator $ - P = \Delta_G $ reads in local coordinates near infinity
\begin{eqnarray}
 \Delta_G = \partial_r^2  +  r^{-2} g^{jk}(r,\theta)  \partial_{\theta_j \theta_k}^2 + (n-1)r^{-1} \partial_r + w (r,\theta) \partial_r + w_k (r,\theta) \partial_{\theta_k}   \label{formuleLaplaciensansrescaling}
\end{eqnarray}  
where $ ( g^{jk}(r,\theta) ) = ( g_{jk}(r,\theta) )^{-1} $ if $ g (r) =  g_{jk}(r,\theta) d \theta_j d \theta_k $. The lower order coefficients are
\begin{eqnarray}
w (r,\theta) =  \frac{\partial_r |g(r,\theta)|}{|g(r,\theta)|} \in S^{-1-\nu},
\end{eqnarray}
since $ |g(r,\theta)| := \mbox{det} (g_{jk}(r,\theta))^{1/2} = |\bar{g}(\theta)| + S^{-\nu} $, and
\begin{eqnarray}
 w_k (r,\theta) = \frac{1}{r^2}  \frac{1}{|g(r,\theta)|} \partial_{\theta_j} \big( g^{jk}(r,\theta) |g(r,\theta)|  \big) \in S^{-2} .
\end{eqnarray}
%
The description of the first order terms will be particularly useful to solve transport equations (see Proposition \ref{prop-transport}). 
   It is also useful to observe that, using the rescaled variable $ \breve{r} = \epsilon r $,
\begin{eqnarray}
    \frac{\Delta_G}{\epsilon^2}  =  {\mathscr D}_{\epsilon} \Delta_{G_{\epsilon}}  {\mathscr D}_{\epsilon}^{-1}, \qquad G_{\epsilon} = d \breve{r}^2 + \breve{r}^2 g ( \breve{r} / \epsilon ) ,  \label{Laplacienrescaleref}
\end{eqnarray}    
that is
\begin{eqnarray}
\Delta_{G_{\epsilon}} =   \partial_{\breve{r}}^2  +  \breve{r}^{-2} g^{jk}(\breve{r}/\epsilon,\theta)  \partial_{\theta_j \theta_k}^2 + (n-1)\breve{r}^{-1} \partial_{\breve{r}} + \epsilon^{-1} w (\breve{r}/\epsilon ,\theta) \partial_{\breve{r}} + \epsilon^{-2}w_k (\breve{r}/ \epsilon ,\theta) \partial_{\theta_k} .  \nonumber
\end{eqnarray}
We will see in Lemma \ref{lemmedebonrescaling} that the negative powers of $ \epsilon  $ in front of $ w (\breve{r}/\epsilon ,\theta) $ and $ w_k (\breve{r}/ \epsilon ,\theta) $ are harmless in $  \{ \breve{r} \gtrsim 1 \} $, {\it i.e.} in the region $ \{  \epsilon r \gtrsim 1\} $.

To distinguish clearly between what is globally defined and what is defined in a chart, we will use the notation
$$ P_{\kappa} =  \Pi_{\kappa}^{-1} P \Pi_{\kappa} ,$$
for the expression of $ P $ in local coordinates (that is minus the right hand side of (\ref{formuleLaplaciensansrescaling})) and 
\begin{eqnarray}
 P_{\epsilon , \kappa} =  {\mathscr D}_{\epsilon}^{-1} \frac{P_{\kappa}}{\epsilon^2 } {\mathscr D}_{\epsilon} . \label{refoprescale}
 \end{eqnarray}
for its rescaled version (that is minus the above expression of $ \Delta_{G_{\epsilon}} $). We denote respectively by
\begin{eqnarray}
  p_{\kappa}  = \rho^2 +  r^{-2} g^{jk}(r,\theta) \eta_j \eta_k \qquad \mbox{and} \qquad  p_{\epsilon,\kappa} =  \breve{\rho}^2 +  \breve{r}^{-2} g^{jk}(\breve{r}/\epsilon,\theta) \eta_j \eta_k \label{symbolesprincipaux}
\end{eqnarray}
the principal symbols of $ P_{\kappa} $ and $ P_{\epsilon,\kappa} $ in local coordinates near infinity.


\section{Pseudodifferential calculus} \label{sectionpseudo}
\setcounter{equation}{0}
\subsection{ Operators on $ \Ra^n $.} 
We shall use  symbols in 
the classes $ \widetilde{S}^{m,\mu} $  which are defined as follows. For $ m, \mu \in \Ra $,  $ \widetilde{S}^{m,\mu} $ is the set of symbols on $ \Ra^{2n} $ such that
\begin{eqnarray}
 \big| \partial_r^j \partial_{\theta}^{\alpha} \partial_{\rho}^k \partial_{\eta}^{\beta} a (r,\theta,\rho,\eta) \big| \leq C \scal{r}^{m-j-|\beta|} \left( \scal{\rho} + \frac{\scal{\eta}}{\scal{r}} \right)^{\mu-k-|\beta|}  \label{typedecroissancetilde}
\end{eqnarray} 
for all $ r, \rho \in \Ra $ and $ \theta , \eta \in \Ra^{n-1} $. As usual, the best constants $C$ are semi-norms which define the topology of $ \widetilde{S}^{m,\mu} $. We also set $ \widetilde{S}^{-\infty,\mu} := \cap_m \widetilde{S}^{m,\mu} $. We use the semiclassical quantization
$$ O \! p^h (a) = a (r,\theta,hD_r , h D_{\theta}) , $$
with $ h \in (0,1] $. Note that we put $h$ in exponent in this notation to distinguish it with the one of rescaled pseudo-differential operators introduced in Definition \ref{defrescaledpseudo} below; high frequencies are raised, while low frequencies will be lowered! 

We need to consider admissible  symbols, {\it i.e.} $h$ dependent families of symbols with an asymptotic expansion in $h$ in the following usual sense
$$ a_h \sim \sum_{j \geq 0} h^j a_j \  \mbox{in} \ \widetilde{S}^{m,\mu} \ \ \stackrel{\rm def}{\Longleftrightarrow} \ \ \mbox{for all} \ N, \ \  h^{-N} \left( a_h - \sum_{k < N} h^j a_j \right)  \ \mbox{is bounded in} \  \widetilde{S}^{m-N,\mu-N}  . $$
Note that this implies in particular that each $a_j$ belongs to $ \widetilde{S}^{m-j,\mu-j}  $.
We call the symbol in the right hand side the remainder of order $ N $. When $ m = - \infty $, the above expansion means that it holds for every finite $m$.

The pseudo-differential calculus in the classes $ \widetilde{S}^{m,\mu} $ enjoys the usual symbolic properties since the weight $  \scal{\rho} + \frac{\scal{\eta}}{\scal{r}}  $ is temperate, for it is easily seen that 
$$  \left( \scal{\rho + \delta_{\rho}} + \frac{\scal{\eta + \delta_{\eta}}}{\scal{r + \delta_{r}}} \right) \lesssim \left( \scal{\rho} + \frac{\scal{\eta}}{\scal{r}} \right) (1+ |\delta_r| +|\delta_{\rho}|+ |\delta_{\eta}|)^2 , $$
for all $ r , \delta_r , \rho , \delta_{\rho} \in \Ra $ and $ \eta , \delta_{\eta} \in \Ra^{n-1} $. In particular, we have the following rules.

\begin{prop}[Symbolic calculus in $ \widetilde{S}^{m,\mu} $] \label{bonpseudodiff} Let $ m , m^{\prime}, \mu , \mu^ {\prime} \in \Ra $.
\begin{itemize}
\item{{\bf Adjoint\footnote{for clarity, we will denote by $ \dag $ the adjoints w.r.t. to the Lebesgue measure and keep the notation $ * $ for adjoints w.r.t. the Riemannian measure}:} for every $ a \in \widetilde{S}^{m,\mu} $, one has
$$ O \! p^h (a)^{\dag} = O \! p^h \big(a^{\dag}_h \big) , \  \ \  a^{\dag}_h \sim \sum_{j \geq 0} h^j \left( \sum_{k+|\alpha|=j} \frac{D_r^k D_{\theta}^{\alpha} \partial_{\rho}^k \partial_{\eta}^{\alpha} \bar{a}}{k! \alpha !}  \right) \ \ \mbox{in} \ \widetilde{S}^{m,\mu}. $$}
\item{{\bf Composition:} for every $ a \in \widetilde{S}^{m,\mu} $ and $ b \in \widetilde{S}^{m^{\prime},\mu^{\prime}}  $, one has
$$ O \! p^h (a)  O \! p^h (b)  = O \! p^h \big((a\# b)_h \big) , \  \ \  (a\#b)_h \sim \sum_{j \geq 0} h^j \left( \sum_{k+|\alpha|=j} \frac{\partial_{\rho}^k \partial_{\eta}^{\alpha} a D_r^k D_{\theta}^{\alpha} b}{k! \alpha !} \right) \ \ \mbox{in} \  \widetilde{S}^{m+m^{\prime},\mu+\mu^{\prime}} . $$}
\item{{\bf Invariance by angular diffeomorphisms:}  let $ \tau : V_1 \rightarrow V_2  $ be a diffeomorphism between two open subsets of $ \Ra^{n-1} $. 
 For  all $ a \in \widetilde{S}^{m,\mu} $ such that
 \begin{eqnarray}
 \emph{supp}(a) \subset \Ra \times K \times \Ra^n  \qquad \mbox{for some} \ K \Subset V_1 ,
 \end{eqnarray}
 and for all $ \varphi  \in C_0^{\infty}(V_1) $, one has 
$$   \Pi_{\tau  }^{-1}  O \! p^h (a) \varphi (\theta) \Pi_{\tau } = O \! p^h \big (a^{\tau}(h) \big), \qquad a^{\tau}(h) \sim \sum_{j \geq 0} h^j a_j^{\tau} \ \ \mbox{in} \ \ \widetilde{S}^{m,\mu}  , $$
with symbols $ a_j^{\tau} $ such that
\begin{eqnarray}
\emph{supp} (a_j^{\tau}) \subset  \left\{  \big(r , \tau (\theta),\rho,  ( d \tau (\theta)^T)^{-1 } \eta \big)  \ | \ (r,\theta,\rho,\eta) \in \emph{supp}(a) \right\} \subset \Ra \times V_2 \times \Ra^n . \label{supportgeneral}
\end{eqnarray}
}
\item{{\bf $ L^2 $ boundedness:} There exists a constant $ C (a) $ depending on a finite number of semi-norms of $a \in \widetilde{S}^{0,0}$  such that, for all such $a$ and all $h \in (0,1] $,
\begin{eqnarray}
\big| \big|   O \! p^h (a)   \big| \big|_{L^2 (\scal{r}^{n-1}dr d \theta) \rightarrow L^2(\scal{r}^{n-1}dr d \theta)}  \leq C(a) .
\label{CalderonVaillancourtapoids}
\end{eqnarray}
Here and below, $ L^2(\scal{r}^{n-1}dr d \theta) $ is a shorthand for $ L^2(\Ra^n,\scal{r}^{n-1}dr d \theta) $.}
\end{itemize}
\end{prop}
We point out that all terms of the expansions as well as the remainders depend equicontinuously  on $a$ (or $ (a,b) $ in the second item). In the fourth item, we consider the measure $ \scal{r}^{n-1} d r d \theta $ for this is of course the good model near infinity for the Riemannian measure of $G$. The $ L^2 $ boundedness is a consequence of the usual Calder\'on-Vaillancourt Theorem since
$$ \big| \big|   O \! p^h (a)   \big| \big|_{L^2 (\scal{r}^{n-1}dr d \theta) \rightarrow L^2(\scal{r}^{n-1}dr d \theta)} =  \big| \big|   \scal{r}^{\frac{n-1}{2}} O \! p^h (a)  \scal{r}^{\frac{1-n}{2}} \big| \big|_{L^2 (drd\theta) \rightarrow L^2 (drd\theta)}  , $$
where, by the second item of Proposition \ref{bonpseudodiff}, $ \scal{r}^{\frac{n-1}{2}} O \! p^h (a)  \scal{r}^{\frac{1-n}{2}} =  O \! p^h (a(h))   $ for some admissible family $ a(h) \in \widetilde{S}^{0,0} $.

\bigskip

We next introduce the convenient definition of rescaled pseudo-differential operators. 

\begin{defi}[Rescaled pseudo-differential operators] \label{defrescaledpseudo} If $a \in \widetilde{S}^{m,\mu} (\Ra_{\breve{r}} \times \Ra^{n-1}_{\theta} \times \Ra_{\breve{\rho}} \times \Ra^{n-1}_{\eta}) $ for some $ m , \mu \in \Ra $, we set
\begin{eqnarray*}
 O \! p_{\epsilon} (a) & := & {\mathscr D}_{\epsilon} O \! p^1 (a) {\mathscr D}_{\epsilon}^{-1} .
 \end{eqnarray*}
\end{defi}
More explicitly, 
$$ O \! p_{\epsilon} (a) = a \left( \epsilon r , \theta , \frac{D_r}{\epsilon} , D_{\theta} \right)  . $$
To clarify the presentation, we distinguish the variables $ (r,\rho)  $ and $ ( \breve{r} , \breve{\rho} ) $ which have to be thought  as
$$ \epsilon r = \breve{r}, \qquad \frac{\rho}{\epsilon} = \breve{\rho} . $$
In the typical situation we shall encounter, we will consider $ a (\breve{r},\theta,\breve{\rho},\eta) = b \left( \breve{r} , \theta , \breve{\rho} ,  \breve{r}^{-1} \eta \right) $ for which 
\begin{eqnarray*}
 O \! p_{\epsilon} (a)  =  b \left( \epsilon r , \theta , \frac{D_r}{\epsilon} , \frac{1}{r} \frac{D_{\theta}}{\epsilon}  \right) .
\end{eqnarray*} 
If $b$ is compactly supported in momentum, this corresponds  to a low frequency localization.

Let us comment a little bit more on Definition \ref{defrescaledpseudo}. Rescaled pseudo-differential operators will be used to approximate low frequency localization of $P$, {\it i.e.} operators of the form $ f (P/ \epsilon^2) $ with $ f \in C_0^{\infty} (\Ra_+) $. By the uncertainty principle, one can only expect to get such an approximation where $r$ is large, typically $ r \gtrsim \epsilon^{-1} $, which corresponds to considering symbols $a$ (or $b$ as above)  supported in $  \breve{r} \gtrsim 1 $. This is consistent with the following easily verified property.


\begin{lemm} \label{lemmedebonrescaling} Let $ a \in S^{\mu} (\Ra_r \times \Ra^{n-1}_{\theta}) $ with $ \mu \in \Ra $. Let 
$$ a_{\epsilon} (\breve{r},\theta) := \epsilon^{\mu} a \left( \breve{r} /\epsilon , \theta \right) . $$
Then $ (a_{\epsilon} )_{\epsilon \in (0,1]} $ belongs to a bounded subset of $ S^{\mu} \big( (1,\infty)_{\breve{r}} \times \Ra^{n-1}_{\theta} \big) $, {\it i.e.}
$$ \big| \partial_{\breve{r}}^j \partial_{\theta}^{\alpha} a_{\epsilon} (\breve{r},\theta) \big| \leq C_{j \alpha} \breve{r}^{\mu - j}, \qquad \breve{r} \geq 1 , \ \theta \in \Ra^{n-1} , \ \epsilon \in (0,1] . $$
\end{lemm}
The meaning of this lemma is that $ a_{\epsilon} $ is only singular for $ \breve{r} $ close to $0$ (the threshold $ \breve{r} \geq 1 $ could be replaced by $ \breve{r} \geq c $  for any $ c > 0 $ positive). In other words, as long as one works in the region $ \epsilon r \gtrsim 1 $, rescaling does not produce singular symbols. 

 We further illustrate the interest of rescaled pseudo-differential operators by keeping in mind the example of (\ref{Laplacienrescaleref}). For  $ k + |\beta| \leq 2 $ and $ a \in S^{ k + |\beta| - 2 - \nu} $ ($ \nu \geq 0 $), we will consider in pratice operators of the form
\begin{eqnarray*}
 \frac{1}{\epsilon^2} a (r,\theta) ( r^{-1} D_{\theta})^{\beta} D_r^{k} & = & {\mathscr D}_{\epsilon} \left( \frac{1}{\epsilon^{2-k-|\alpha|}} a \left( \breve{r} /\epsilon , \theta \right) ( \breve{r}^{-1} D_{\theta})^{\alpha} D_{\breve{r}}^k \right) {\mathscr D}_{\epsilon}^{-1} , \\
 & = & {\mathscr D}_{\epsilon} \left(  \epsilon^{\nu} a_{\epsilon} \left( \breve{r} /\epsilon , \theta \right) ( \breve{r}^{-1} D_{\theta})^{\alpha} D_{\breve{r}}^k \right) {\mathscr D}_{\epsilon}^{-1}
\end{eqnarray*} 
with
$$ a_{\epsilon} (\breve{r},\theta) = \epsilon^{k + |\beta| - 2 - \nu} a (\breve{r}/\epsilon, \theta) . $$ 
Studying such operators in  $ \{ \epsilon r \gtrsim 1 \} $ corresponds to study  $ a_{\epsilon} (\breve{r},\theta) ( \breve{r}^{-1} D_{\theta})^{\alpha} D_{\breve{r}}^k   $ in $ \{ \breve{r} \gtrsim 1 \} $;  by Lemma \ref{lemmedebonrescaling},    $a_{\epsilon} $ is bounded in $ S^{k + |\beta|-2-\nu}( (1,\infty)_{\breve{r}} \times \Ra^{n-1}_{\theta} ) $, and allows to use pseudo-differential calculus in the variables $ (\breve{r},\theta,\breve{\rho},\eta) $. Typically, to construct a parametrix for $ \chi (\epsilon r) (P/\epsilon^2  + i)^{-1} $ in $ \{ \epsilon r \geq R \} $, we will consider symbols of the form
$$ \chi (\breve{r}) \frac{1}{\breve{\rho}^2 + \breve{r}^{-2} g^{jk} (\breve{r}/\epsilon,\theta) \eta_j \eta_k + i } $$
with $ \chi $ supported in $ (R,+\infty) $. By Lemma \ref{lemmedebonrescaling}, this $\epsilon$-dependent symbol belongs to a bounded subset of $ \widetilde{S}^{-2,0}  $, allowing to perform the usual iterative parametrix construction (see paragraph \ref{subsectioncalculfonc}). 

\subsection{ Operators on $ {\mathcal M} $.} 
Let us define the space $ {\mathscr S}({\mathcal M}) $ by
\begin{eqnarray}
 u \in {\mathscr S} ({\mathcal M}) \qquad \Longleftrightarrow \qquad u \in \cap_{m > 0} \mbox{Dom}(P^m) \ \ \mbox{and} \ \ \ r^j \partial_r^k \partial_{\theta}^{\alpha} u \in L^2 \ \ \mbox{for all } j,k,\alpha,   \label{definitiondeSdeM}
 \end{eqnarray}
the second condition in the right hand side being a condition at infinity (it is invariant by change of coordinates on $ {\mathcal S} $). It is the natural Schwartz space on $ {\mathcal M} $ and will be  convenient for our purposes.

Using the charts introduced in Section \ref{sectionNotation}, we will  note everywhere in this paper
\begin{eqnarray}
  O \! p^h_{\kappa} (a) := \Pi_{\kappa} O \! p^h (a)  \Pi_{\kappa}^{-1} .  \label{definitionnotationdiffeo}
\end{eqnarray}
If nothing is specified about $ a \in \widetilde{S}^{m,\mu}(\Ra^{2n}) $, such operators are defined from $ C_0^{\infty}((R_{\mathcal M},\infty) \times U_{\kappa}) $ to $ C^{\infty}( (R_{\mathcal M},\infty) \times U_{\kappa} ) $. If in addition $ \mbox{supp} (a) \subset (R_{\mathcal M},\infty ) \times V_{\kappa} \times \Ra^n $, which will always be the case in this paper, they map  $ C_0^{\infty}((R_{\mathcal M},\infty) \times U_{\kappa}) $  to $ C^{\infty}({\mathcal M}) $. In practice, we will only consider globally defined operators of the form
\begin{eqnarray}
 O \! p^h_{\kappa} (a) \tilde{\psi}_{\kappa} = O \! p^h_{\kappa} (a) \tilde{\psi}_{\kappa} (r,\omega)  \label{pseudoglobalhf}
\end{eqnarray} 
where the cutoff $ \tilde{\psi}_{\kappa} $ localizes inside $ (R_{\mathcal M},\infty) \times U_{\kappa}  $ (see (\ref{notationpartition})) and where we will use symbols spatially supported in $  (R_{\mathcal M},\infty) \times U_{\kappa} $ ({\it e.g.} in the support of $ \psi_{\kappa} (r,\kappa^{-1}(\theta))$ - see again (\ref{notationpartition})). We point out that such operators are localized near infinity, where we will focus essentially all our analysis. 
Note also that since pseudo-differential operators on $ \Ra^n $ with symbols in $ {\widetilde S}^{m,\mu} $ map the Schwartz space (on $ \Ra^n $) into itself, we have  
$$  O \! p^h_{\kappa} (a) \tilde{\psi}_{\kappa} : {\mathscr S}({\mathcal M}) \rightarrow {\mathscr S}({\mathcal M})  .  $$
We define analogously rescaled pseudo-differential operators on $ {\mathcal M} $ by
$$     O \! p_{\epsilon,\kappa} (a) :=  \Pi_{\kappa} O \! p_{\epsilon} (a)  \Pi_{\kappa}^{-1}    $$
and will consider, for symbols supported in $ (R_{\mathcal M},\infty ) \times V_{\kappa} \times \Ra^n $,
\begin{eqnarray}
 O \! p_{\epsilon, \kappa} (a) \tilde{\psi}_{\kappa} (\epsilon r ) = O \! p_{\epsilon,\kappa} (a) \tilde{\psi}_{\kappa} (\epsilon r,\omega)    \label{pseudolocalbf}
\end{eqnarray}
(we will often drop the dependence on $ \omega $ from the notation, though $ \tilde{\psi}_{\kappa} (\epsilon r ) $ really depends also on $ \omega \in {\mathcal S}$). It is important to note that if $ a $ is spatially localized in $ (R_{\mathcal M},\infty)_{\breve{r}} \times V_{\kappa} $ then the range of $  O \! p_{\epsilon, \kappa} (a) \tilde{\psi}_{\kappa} (\epsilon r ) $ contains only functions supported in $ (\epsilon^{-1} R_{\mathcal M},\infty)_r  \times U_{\kappa} $; in other words, such operators are localized in $\{ \epsilon r > R_{\mathcal M} \}$ and will be used as microlocalization in this region only. We finally note that we will often use $ \epsilon $ dependent symbols, similar to those considered in Lemma \ref{lemmedebonrescaling}.

For further use and to illustrate that such definitions fit the usual expected properties of a pseudo-differential calculus, we compute adjoints with respect to the Riemannian measure  $ r^{n-1}|g(r,\theta)| d r d \theta $ (see Section \ref{sectionNotation} for $ |g(r,\theta)| $).
 Let  $a = a(r,\theta,\rho,\eta)$ be a symbol spatially supported inside $ (R_{\mathcal M},+\infty) \times V_{\kappa} $, i.e. with support in $ (r,\theta) $ contained in $ (R,\infty) \times K $ for some $ R > R_{\mathcal M} $ and $ K \Subset V_{\kappa} $. Then, using Proposition \ref{bonpseudodiff} and elementary computations, we find
\begin{eqnarray}
 \left( O \! p^h_{\kappa} (a) \tilde{\psi}_{\kappa} \right)^* & = &  \tilde{\psi}_{\kappa} \Pi_{\kappa} \left( \frac{1}{r^{n-1} |g(r,\theta)|} O \! p^h (a)^{\dag} r^{n-1} |g(r,\theta)| \right) \Pi_{\kappa}^{-1} \nonumber  \\
 & = &  \tilde{\psi}_{\kappa} O \! p^h_{\kappa} (b(h)) \psi_{1,\kappa}   \label{adjointglobalhf}
 \end{eqnarray}
for some admissible symbol $ b(h) $ in the same class as $a$ and $ \psi_{1,\kappa} $ supported in $ (R_{\mathcal M},+\infty) \times U_{\kappa} $. Similarly
\begin{eqnarray}
 \left( O \! p_{\epsilon, \kappa} (a) \tilde{\psi}_{\kappa}(\epsilon r) \right)^* & = &  \tilde{\psi}_{\kappa} (\epsilon r) \Pi_{\kappa} {\mathscr D}_{\epsilon} \left( \frac{1}{\breve{r}^{n-1} |g(\breve{r}/\epsilon,\theta)|} O \! p^1 (a)^{\dag} \breve{r}^{n-1} |g(\breve{r}/\epsilon,\theta)| \right) {\mathscr D}_{\epsilon}^{-1} \Pi_{\kappa}^{-1} \nonumber \\
 & = &  \tilde{\psi}_{\kappa} (\epsilon r) O \! p_{\epsilon, \kappa} (b_{\epsilon}) \psi_{1,\kappa}(\epsilon r)  \label{adjointgloballf}
 \end{eqnarray}
with $ (b_{\epsilon} )_{\epsilon \in (0,1]} $ bounded in the same class as $ a $, also using here Lemma \ref{lemmedebonrescaling} to handle  $ |g(\breve{r}/\epsilon,\theta)|^{\pm 1} $. 

To get $ L^2 $ or $ L^q $ estimates, we will use the following proposition.
\begin{prop} \label{borneL1Linfiniepsilon} Let $ \psi  $ be bounded and supported in $ (R_{\mathcal M},+\infty) \times U_{\kappa} $ and $ q \in [1,\infty] $. Then
\begin{eqnarray}
 \big| \big| \psi (\epsilon r , \omega) \Pi_{\kappa} {\mathscr D}_{\epsilon} \big| \big|_{L^{q} (\scal{r}^{n-1}dr d\theta) \rightarrow L^{q} ({\mathcal M})} & \lesssim & \epsilon^{\frac{n}{2} - \frac{n}{q}} \label{pourqprime} \\
  \big| \big|  {\mathscr D}_{\epsilon}^{-1} \Pi_{\kappa}^{-1} \psi (\epsilon r , \omega)  \big| \big|_{L^q ({\mathcal M}) \rightarrow L^q (\scal{r}^{n-1}dr d\theta)}  & \lesssim & \epsilon^{\frac{n}{q} - \frac{n}{2}} \label{pourqegal2}
 \end{eqnarray}
 for $ \epsilon \in (0,1] $.
\end{prop}

\noindent {\it Proof.} It follows from an elementary change of variable together with the observation that, on the support of $ \psi \big(\epsilon r , \kappa^{-1}(\theta) \big) $,
$$  r^{n-1} |g(r,\theta)| / C \leq  \scal{r}^{n-1} \leq C r^{n-1} |g(r,\theta)| $$
for some $ C > 1 $. \finpreuve

\bigskip

We note in particular that, when $q=2$, Proposition \ref{borneL1Linfiniepsilon} together with (\ref{CalderonVaillancourtapoids})  imply  that
%
\begin{eqnarray}
\big| \big|    O \! p_{\epsilon,\kappa} (a) \tilde{\psi}_{\kappa} (\epsilon r)   \big| \big|_{L^2 \rightarrow L^2} \leq C (a), \qquad \epsilon \in (0,1] , \label{CalderonVaillancourtrescale}
\end{eqnarray}
with $ C (a) $ bounded as long as $a$ belongs to a bounded subset of $ \widetilde{S}^{0,0} $ ($a$ being spatially supported in $(R_{\mathcal M},\infty) \times V_{\kappa}$). For completeness, we also recall that at high frequency, under the same assumptions on $a$,
\begin{eqnarray}
\big| \big|    O \! p^h _{\kappa} (a) \tilde{\psi}_{\kappa}    \big| \big|_{L^2 \rightarrow L^2} \leq C (a), \qquad h \in (0,1] , \label{CalderonVaillancourtsemiclassique}
\end{eqnarray}
which is more standard (and does not use Proposition \ref{borneL1Linfiniepsilon}).

We will also need $ L^q $ estimates on pseudo-differential operators.
\begin{prop} \label{bornepseudoLq} Let $ a \in \widetilde{S}^{-\infty,0} $ be spatially supported in $ (R_{\mathcal M},+\infty) \times V_{\kappa}  $. Let $1 \leq  q_1 \leq q_2 \leq \infty  $. Then
\begin{eqnarray}
\big| \big| O \! p_{\kappa}^h (a) \tilde{\psi}_{\kappa}  \big| \big|_{L^{q_1} \rightarrow L^{q_2} } & \leq & C h^{ \frac{n}{q_2} - \frac{n}{q_1} } , \\
\big| \big| O \! p_{\epsilon, \kappa} (a) \tilde{\psi}_{\kappa} (\epsilon r) \big| \big|_{L^{q_1} \rightarrow L^{q_2}} & \leq & C \epsilon^{\frac{n}{q_1} - \frac{n}{q_2}} ,
\end{eqnarray}
The constant is bounded as long as $a$ belongs to a bounded subset of $ \widetilde{S}^{-\infty,0} $.
\end{prop}

\noindent {\it Proof.} Write $ a (r,\theta,\rho,\eta) = b (r,\theta,\rho,\eta/r) $ so that $b$ becomes a Schwartz function in the momentum variables, uniformly in $ (r,\theta) $. The estimate in the semiclassical case follows from the similar estimate for $ O \! p^h (a) $ from $ L^{q_1} (\scal{r}^{n-1} dr d \theta) $ to $ L^{q_2} (\scal{r}^{n-1}dr d \theta) $ obtained from the usual Schur test and interpolation argument, by exploiting that its kernel with respect to $ \scal{r}^{n-1} dr d \theta $ reads
$$ (2\pi h)^{-n} r^{n-1} \hat{b} \left(r,\theta,  \frac{r^{\prime}-r}{h}, r \frac{\theta^{\prime} - \theta}{h}\right) \scal{r^{\prime}}^{1-n} $$
where $ \hat{ \ } $ is the Fourier transform in the momentum variables. The low frequency case follows from the above one with $ h = 1 $ together with Proposition \ref{borneL1Linfiniepsilon}. \finpreuve

\subsection{Functional calculus} \label{subsectioncalculfonc}
   We will use operators of the form (\ref{pseudoglobalhf}) or (\ref{pseudolocalbf})  to describe functions of $ P $. In the semiclassical or high frequency regime, this is mostly standard, see {\it e.g.}  \cite{BoucletFourier,MizutaniCPDE},  though we will need a sharper description of the remainders than in those references. We will also consider the  low frequency regime, which is less standard but can be easily handled by considering appropriate spatial localizations and rescaled operators. 
    The first and main step is  to construct a parametrix for  $ (P / \epsilon^2 - z )^{-1} $ in the region $ \{ \epsilon r > R_{\mathcal M} \}  $. To do so, we need  basically to use that
   \begin{eqnarray}
    \left( \frac{P}{\epsilon^2}  - z \right)= \Pi_{\kappa} {\mathscr D}_{\epsilon} \big(  P_{\epsilon,\kappa} - z \big) {\mathscr D}_{\epsilon}^{-1} \Pi_{\kappa}^{-1}  \label{rescalingdiffexplicite}
  \end{eqnarray}
(see (\ref{refoprescale}))  namely that $ P / \epsilon^2  $ is a rescaled (pseudo-)differential operator whose symbol is not singular w.r.t. $ \epsilon $ in the region $ \{ \breve{r} > R_{\mathcal M}\} $ thanks to Lemma \ref{lemmedebonrescaling}. One can then apply the usual elliptic parametrix scheme to $ P_{\epsilon,\kappa} - z $ to construct an approximate inverse. Taking into account the composition rules of Proposition \ref{bonpseudodiff}, we obtain the following technical result. 



\begin{prop} \label{leplusaffreux} Let $ \psi , \tilde{\psi}, \tilde{\tilde{\psi}} $ be smooth functions supported in a patch $ (R,\infty) \times U_{\kappa} $ with $ R > R_{\mathcal M} $, all belonging to $ S^0 $ and such that
$$ \tilde{\psi} \equiv 1 \ \ \mbox{near} \ \ \emph{supp}(\psi) , \qquad \tilde{\tilde{\psi}} \equiv 1 \ \ \mbox{near} \ \ \emph{supp}(\tilde{\psi}). $$
Then  for $ j,N \in \Na $ and $ z \in \Ca \setminus [0,+\infty) $, one has
\begin{itemize}
\item{{\bf High frequency parametrix:} for $ h \in (0,1] $,
$$  \psi (r,\omega)  (h^2 P - z)^{-j}  =  \sum_{l = 0}^{N-1 } h^l {\psi}  O \! p^h_{\kappa} (q_{l}(z)) \tilde{\psi}  + h^N {\mathcal R}_{\rm high} (z,h)  $$
where each $ q_{l}(z) \in \widetilde{S}^{-2 j -l,-l}  $ is a  linear combination of $  a_k ( p_{\kappa} - z )^{-j-k} $ for some symbol $ a_k \in \widetilde{S}^{2k-l,-l} $ independent of $z$, and with
\begin{eqnarray*}
{\mathcal R}_{\rm high}(z,h) & = &     \psi   O \! p^h_{\kappa} (r(z,h))  \tilde{\tilde{\psi}} (h^2 P -z)^{-j} 
\end{eqnarray*}
where
 $ r(z,h) \in \widetilde{S}^{-N,-N} $ with seminorms growing polynomially in $ 1 / {\rm dist}(z,\Ra_+) $ uniformly in $h$ as long as $z$ belongs to a bounded set of $ \Ca \setminus [0,+\infty ) $.}
\item{{\bf Low frequency parametrix:} for $ \epsilon \in (0,1] $,
$$  \psi (\epsilon r,\omega)  ( P/ \epsilon^2 - z)^{-j}   =  \sum_{l = 0}^{N-1 }     {\psi} (\epsilon r, \omega)  O \! p_{\epsilon,\kappa} ( q_{\epsilon,l}(z))  \tilde{\psi} (\epsilon r , \omega)     +    {\mathcal R}_{\rm low} (z,\epsilon) $$
where each $ q_{\epsilon,l}(z) \in \widetilde{S}^{-2j-l,-l} $ is a  linear combinations of $  a_{\epsilon,k} ( p_{\epsilon,\kappa} - z )^{-j-k} $ with symbols $ a_{\epsilon,k} \in \widetilde{S}^{2k-l,-l} $ bounded w.r.t. $ \epsilon $, and 
\begin{eqnarray*}
{\mathcal R}_{\rm low}(z,\epsilon) & = &      \psi (\epsilon r , \omega)  O \! p_{\epsilon, \kappa} ( r_{\epsilon}(z))  \tilde{\tilde{\psi}} (\epsilon r, \omega)   (P/ \epsilon^2 -z)^{-j}
\end{eqnarray*}
where $ r_{\epsilon}(z) \in \widetilde{S}^{-N,-N} $ with seminorms growing polynomially in $ 1 / {\rm dist}(z,\Ra_+) $ uniformly in $ \epsilon $ as long as $z$ belongs to a bounded set of $ \Ca \setminus [0,+\infty) $.}
\end{itemize}
\end{prop}
 We refer to (\ref{symbolesprincipaux}) for the definitions of
 $ p_{\kappa} $ and $ p_{\epsilon,\kappa} $.

Note that the spatial localizations are different at high and low frequency. We also point out that the low frequency parametrix is not an asymptotic expansion in $ \epsilon $, but it only says that $   ( P/  \epsilon^2 - z)^{-j}  \psi (\epsilon r,\omega)  $ is a sum of rescaled pseudo-differential operators and of a remainder which is smoothing and spatially decaying like $ \scal{ \epsilon r }^{-N} $. We finally  remark that a similar proposition  holds for $  (h^2 P - z)^{-j} \psi (r,\omega)  $ and $  ( P/ \epsilon^2 - z)^{-j} \psi (\epsilon r,\omega)  $ (this follows by taking the adjoints and using  (\ref{adjointglobalhf})-(\ref{adjointgloballf})). We will use this occasionally.
 



\bigskip

As a first application, we record the following result where we use the function $ \zeta $ introduced in (\ref{notationpartition})-(\ref{partitionfixee}).

\begin{prop} \label{injectionsdeSobolevfinales} If $ j > n/4 $, then
\begin{eqnarray}
 \big| \big| \zeta (r) (h^2 P + 1)^{-j} \big| \big|_{L^2 \rightarrow L^{\infty}} \lesssim h^{-\frac{n}{2}}, \qquad h \in (0,1] , \label{bytakingtheadjointestimate}
\end{eqnarray} 
and
$$ \big| \big| \zeta (\epsilon r) ( P/ \epsilon^2 + 1)^{-j} \big| \big|_{L^2 \rightarrow L^{\infty}} \lesssim \epsilon^{\frac{n}{2}}, \qquad \epsilon \in (0,1] . $$
\end{prop}

Recall that for simplicity we have set $ L^q = L^q ({\mathcal M}) $ (see after Definition \ref{defasympconique}).

\bigskip

\noindent {\it Proof.} We prove only the second estimate, the first one being standard (see {\it e.g.} \cite{BoucletFourier}). We use  Proposition \ref{leplusaffreux} with  $ \psi $ replaced by  $ \psi_{\kappa} $, $ \tilde{\psi} $  by $ \tilde{\psi}_{\kappa} $   etc. (see (\ref{notationpartition})), and with $ N > n/2 $.  Then $ \zeta (\epsilon r)  ( P/ \epsilon^2 + 1)^{-j} $ is a sum over $ \kappa $ of parametrices as in Proposition \ref{leplusaffreux}. For each $ \kappa $, consider the first term 
$$  {\psi}_{\kappa} (\epsilon r)  O \! p_{\epsilon,\kappa} ( q_{\epsilon,0}(-1)) \tilde{\psi}_{\kappa} (\epsilon r) = \left(
 \psi_{\kappa} (\epsilon r , \omega) \Pi_{\kappa} {\mathscr D}_{\epsilon} \right) \left( O \! p^1 \big( q_{\epsilon,0}(-1) \big) \right)  
 \left( {\mathscr D}_{\epsilon}^{-1} \Pi_{\kappa}^{-1} \tilde{\psi}_{\kappa} (\epsilon r , \omega)  \right) $$
where $ q_{\epsilon,0}(-1) $ belongs to (a bounded set of) $ \widetilde{S}^{-2j,0} $.
The result is a consequence of the fact that $    O \! p^1 ( q_{\epsilon,0}(-1))     $ maps $ L^2 (\scal{r}^{n-1}dr d \theta) $ into $ L^{\infty} (\scal{r}^{n-1}dr d \theta) $ since  $ 2j > n/2 $ (see \cite[Lemma 2.4]{BoucletFourier})), together with the estimates (\ref{pourqprime}) (with $q^{\prime}= \infty$) and (\ref{pourqegal2}) (with $q=2$). The other terms  are treated analogously, as well as the remainder $ {\mathcal R}_{\rm low}(-1,\epsilon) $ by using additionally that $  ||  (P/\epsilon^2 + 1)^{-j}||_{L^2 \rightarrow L^2} \leq 1$ for the remainder.  \finpreuve


\bigskip

To describe the remainders  that will be involved in the different parametrices we are going to construct, it us useful introduce the following norms
\begin{eqnarray}
| | u  | |_{{\mathscr H}^{2j}_{\mu}}  =  \big| \big| \scal{r}^{\mu} (h^2P+1)^j u \big| \big|_{L^2} , \qquad   || u ||_{{\mathscr L}^{2j}_{\mu}}  = \left| \left| 
\scal{\epsilon r}^{\mu}  (P/\epsilon^2 + 1)^{j} u \right| \right|_{L^2} , \label{normshighlow}
\end{eqnarray}
for $ \mu \in \Ra $, $ j \in \Za $ and $ u \in {\mathscr S} (\mathcal M) $. The first one is a standard  weighted semiclassical  Sobolev norm, which will be used at high frequency, and the second one will be used at low frequency. We will only consider these norms  on   $ {\mathscr S}({\mathcal M}) $ for this space is  stable by the resolvent of $ P $ (this is fairly standard or can be checked by using the parametrix of Proposition \ref{leplusaffreux} for $ \epsilon = h = 1 $) so that the norms (\ref{normshighlow}) make clearly sense.  We also point out that we do not define the spaces $ {\mathscr H}^{2j}_{\mu} $ nor $ {\mathscr L}^{2j}_{\mu} $ (which should be the closures of $ {\mathscr S}({\mathcal M}) $ for the corresponding norms) and will only use their norms on $ {\mathscr S}({\mathcal M}) $. The interest of using such norms is to state estimates which are uniform in $ \epsilon $ or $h$. It is also worth recalling that the japanese bracket used in (\ref{normshighlow}) is the modified one chosen in (\ref{scalrM}).

   Given a family of operators $ A_{\epsilon} $ preserving $ {\mathscr S} ({\mathcal M}) $, we will write
   $$ A_{\epsilon} = O_{{\mathscr L}_{\mu_1}^{2j_1} \rightarrow {\mathscr L}_{\mu_2}^{2j_2} } (1) \qquad \Longleftrightarrow \qquad ||A_{\epsilon} u  ||_{{\mathscr L}_{\mu_2}^{2j_2}} \leq C  || u ||_{{\mathscr L}_{\mu_1}^{2j_1}}  \ \ \mbox{for all} \ \epsilon \in (0,1], \ u \in {\mathscr S}({\mathcal M}) , $$
  the point being that the constant is independent of $ \epsilon $. The notation $ A_{h} = O_{{\mathscr H}_{\mu_1}^{2j_1} \rightarrow {\mathscr H}_{\mu_2}^{2j_2} } (1) $
  is  defined similarly.


\begin{prop} \label{propositionpourestes} For all $ j , j^{\prime} \in \Za $ and $ \mu , \mu^{\prime} \in \Ra $, we have
\begin{itemize}
\item{{\bf Global estimates:}
\begin{eqnarray}
(P/\epsilon^2 + 1)^{j^{\prime}} = O_{{\mathscr L}_{\mu}^{2j} \rightarrow {\mathscr L}_{\mu }^{2(j-j^{\prime})}} (1) , 
\qquad (h^2P + 1)^{j^{\prime}} = O_{{\mathscr H}_{\mu}^{2j} \rightarrow {\mathscr H}_{\mu }^{2(j - j^{\prime})}} (1) \label{firstbound}
\end{eqnarray}
and, as multiplication operators,
\begin{eqnarray}
\scal{\epsilon r}^{\mu^{\prime}} = O_{{\mathscr L}_{\mu}^{2j} \rightarrow {\mathscr L}_{\mu - \mu^{\prime}}^{2j}} (1) , 
\qquad \scal{ r}^{\mu^{\prime}} = O_{{\mathscr H}_{\mu}^{2j} \rightarrow {\mathscr H}_{\mu - \mu^{\prime}}^{2j}} (1) . \label{secondbound}
\end{eqnarray}}
\item{{\bf Embeddings estimates:}
\begin{eqnarray}
\mu^{\prime} \leq \mu \qquad \mbox{and} \qquad j^{\prime} \leq j \qquad \Longrightarrow \qquad  I = O_{{\mathscr L}_{\mu}^{2j} \rightarrow {\mathscr L}_{\mu^{\prime}}^{2j^{\prime}}} (1) , 
\qquad I = O_{{\mathscr H}_{\mu}^{2j} \rightarrow {\mathscr H}_{ \mu^{\prime}}^{2{j^{\prime}}}} (1) . \label{thirdbound}
\end{eqnarray}}
\item{{\bf Action of pseudo-differential operators:} Let $ \tilde{\psi} \in S^0 $ be a smooth function supported in the patch $ (R_{\mathcal M},\infty) \times U_{\kappa} $ and $ a \in \widetilde{S}^{2j^{\prime},\mu^{\prime}} $ be spatially supported in $ (R_{\mathcal M},\infty) \times V_{\kappa} $. Then
\begin{eqnarray}
 O \! p_{\epsilon, \kappa}(a) \tilde{\psi} (\epsilon r)  = O_{{\mathscr L}_{\mu}^{2j} \rightarrow {\mathscr L}_{\mu - \mu^{\prime} }^{2(j - j^{\prime})}} (1) , \qquad
 O \! p^h_{\kappa} (a) \tilde{\psi} = O_{{\mathscr H}_{\mu}^{2j} \rightarrow {\mathscr H}_{\mu - \mu^{\prime} }^{2(j - j^{\prime})}} (1) . \label{fourthbound}
\end{eqnarray}
These uniform bounds remain valid as long as $a$ belongs to a bounded subset of $ \widetilde{S}^{2j^{\prime},\mu^{\prime}}  $. }
\end{itemize}
We recall that in \emph{(\ref{fourthbound})} $ \tilde{\psi}(\epsilon r) $ and $ \tilde{\psi} $ are respectively shortands for $ \tilde{\psi}(\epsilon r , \omega) $ and $ \tilde{\psi}(r,\omega) $.
\end{prop}

\noindent {\it Proof.} In all cases, we consider only the low frequency estimates, the semiclassical ones being similar and more standard. (\ref{firstbound}) is an immediate consequence of the definitions of the norms (\ref{normshighlow}). We next prove the first estimate of (\ref{secondbound}). We observe first that for any $ j \in \Za $ and $ \mu \in \Ra $, there exists $ C > 0 $ such that
\begin{eqnarray}
 C^{-1} | | u  | |_{{\mathscr L}^{2j}_{\mu}} \leq \big| \big|  ( P/ \epsilon^2 +1)^j \scal{\epsilon r}^{\mu} u \big| \big|_{L^2}  \leq  C | | u  | |_{{\mathscr L}^{2j}_{\mu}}  ,  \label{epsilonL}
\end{eqnarray} 
for all $u \in {\mathscr S} ({\mathcal M})$ and $ \epsilon  \in (0,1] $. Indeed, let us write 
$$ \scal{\epsilon r}^{\mu}  (P/\epsilon^2 + 1)^{j}  = \left( \scal{\epsilon r}^{\mu} (P/\epsilon^2 + 1)^{j} \scal{\epsilon r}^{-\mu} (P/\epsilon^2 + 1)^{-j} \right) (P/\epsilon^2 + 1)^j \scal{\epsilon r}^{\mu}  . $$
The lower bound in (\ref{epsilonL}) would then follow from the uniform $ L^2 \rightarrow L^2 $ of the parenthesis. Assume for instance that $ j \geq 0 $. Then the parenthesis in the right hand side is the sum of the identity and
\begin{eqnarray}
  \scal{\epsilon r}^{\mu} \left[  (P/\epsilon^2 + 1)^{j} , \scal{\epsilon r}^{-\mu} \right] \zeta (\epsilon r) (P/\epsilon^2 + 1)^{-j}  \label{opabove}
\end{eqnarray}  
where one can insert the cutoff $ \zeta (\epsilon r) $ of the partition of unity (\ref{partitionfixee}) since the commutator is supported in the region where $ \zeta (\epsilon r) = 1 $ by (\ref{scalrM}). The operator (\ref{opabove}) is uniformly bounded on $ L^2 $ since the composition of
$$ \scal{\epsilon r}^{\mu} \left[  (P/\epsilon^2 + 1)^{j} , \scal{\epsilon r}^{-\mu} \right]  = \sum_{\kappa}
 \Pi_{\kappa} {\mathscr D}_{\epsilon}  \scal{ r}^{\mu} \left[  (P_{\kappa,\epsilon} + 1)^{j} , \scal{ r}^{-\mu} \right]  {\mathscr D}_{\epsilon}^{-1} \Pi_{\kappa}^{-1} \psi_{\kappa} (\epsilon r)
 $$
 (see (\ref{rescalingdiffexplicite}))
 with the low energy parametrix for $ \zeta (\epsilon r) (P/\epsilon^2+1)^{-j} $ (derived from Proposition \ref{leplusaffreux} and the partition of unity (\ref{partitionfixee})) is uniformly bounded on $ L^2 $. This follows by using the composition rules of Proposition \ref{bonpseudodiff} together with (\ref{CalderonVaillancourtrescale}) and the bound $ || (P/\epsilon^2 +1)^{-j} ||_{L^2 \rightarrow L^2} \leq 1 $. The case $ j < 0 $ and the upper bound are proved similarly (using possibly the parametrix of $ (P/\epsilon^2 + 1)^{-j} \zeta (\epsilon r) $). Now, with (\ref{epsilonL}) at hand,   the first estimate of (\ref{secondbound}) follows from
 $$ \big| \big| \scal{\epsilon r}^{\mu^{\prime}} u \big| \big|_{{\mathcal L}_{\mu-\mu^{\prime}}^{2j}} \leq C \big| \big| (P/\epsilon^2 + 1)^j \scal{\epsilon r}^{\mu - \mu^{\prime} + \mu^{\prime}} u \big| \big|_{L^2} \leq C^2  \big| \big|  u \big| \big|_{{\mathcal L}_{\mu}^{2j}}  . $$
 Similarly the first estimate of (\ref{thirdbound}) follows from (\ref{epsilonL}) since
 \begin{eqnarray*}
   || \scal{\epsilon r}^{\mu^{\prime}} (P/\epsilon^2 + 1)^{j^{\prime}} u ||_{L^2}   & \leq &   || \scal{\epsilon r}^{\mu} (P/\epsilon^2 + 1)^{j^{\prime}} u  ||_{L^2} \\
   & \leq & C  || (P/\epsilon^2 + 1)^{j^{\prime}} \scal{\epsilon r}^{\mu} u ||_{L^2}
  \leq C  || (P/\epsilon^2 + 1)^{j} \scal{\epsilon r}^{\mu} u ||_{L^2}  . 
 \end{eqnarray*} 
  We finally consider (\ref{fourthbound}). By using the equivalence of norms (\ref{epsilonL}), the result follows from the uniform $ L^2 $ boundedness of
 $$ (P/\epsilon^2 + 1 )^{j-j^{\prime}} \scal{\epsilon r}^{\mu - \mu^{\prime}}     O \! p_{\epsilon,\kappa} (a) \tilde{\psi} (\epsilon r)  \scal{\epsilon r}^{-\mu} (P/\epsilon^2 + 1)^{-j} . $$
By the composition rule of Proposition \ref{bonpseudodiff}, we may assume that $ \mu = \mu^{\prime} = 0 $ up to the replacement of $a$ by $  \tilde{a}$ such that $ O \! p_1 (\tilde{a}) =
 \scal{ r}^{\mu - \mu^{\prime}}   O \! p_1 (a)  \scal{ r}^{-\mu} $. Then if both $ j - j^{\prime} $ and $ - j $ are non negative, the result follows by using (\ref{rescalingdiffexplicite}), the composition rule and the $ L^2$ bound (\ref{CalderonVaillancourtrescale}). Otherwise we expand the negative powers of  $ P/\epsilon^2 + 1  $ by mean of Proposition \ref{leplusaffreux}  so that we can compose rescaled operators supported in the same patch and conclude again with (\ref{CalderonVaillancourtrescale}). \finpreuve

 \begin{theo} \label{theoremcaclulfonctionnel} For all $ f \in C_0^{\infty} (\Ra) $ and all given $N$,
$$ \zeta (r) f (h^2 P) = \sum_{l = 0}^{ N-1} \sum_{  \kappa} h^l \psi_{\kappa} O \! p^h_{\kappa} (a_{\kappa,l}) \tilde{\psi}_{\kappa}  + h^N {\mathscr R}_{\rm high} (f,h)  $$
where $ a_{\kappa,l} \in \widetilde{S}^{-\infty,-l}$ with $ \emph{supp}(a_{\kappa,l}) \subset \emph{supp}(f \circ p_{\kappa}) $ 
and, for any $ M > 0 $ and $ \mu \in \Ra $,
$$ {\mathscr R}_{\rm high} (f,h) =  O_{{\mathscr H}_{\mu}^{-2M} \rightarrow {\mathscr H}_{\mu + N }^{2M}} (1) .  $$
Also 
\begin{eqnarray}
 \zeta (\epsilon r) f (P/\epsilon^2) =  \sum_{l = 0}^{ N-1} \sum_{  \kappa}   \psi_{\kappa} (\epsilon r ) O \! p_{\epsilon,\kappa} (a_{\epsilon,\kappa,l}) \tilde{\psi}_{\kappa} (\epsilon r)  +  {\mathscr R}_{\rm low} (f,\epsilon)  
 \label{lowfrequencyparametrix}
 \end{eqnarray}
where $ ( a_{\epsilon,\kappa,l} )_{\epsilon \in (0,1]}$ belongs to a bounded subset of $ \widetilde{S}^{-\infty,-l}$ with $ \emph{supp}(a_{\epsilon,\kappa,l}) \subset \emph{supp}(f \circ p_{\epsilon, \kappa}) $
  and, for any $ M > 0 $ and $ \mu \in \Ra $
\begin{eqnarray}  
 {\mathscr R}_{\rm low} (f,\epsilon) = O_{{\mathscr L}_{\mu}^{-2M} \rightarrow {\mathscr L}_{\mu + N }^{2M}} (1)  .  
 \label{discussionreste}
\end{eqnarray} 
\end{theo}

\noindent {\it Proof.} We consider only the proof of the low frequency parametrix (\ref{lowfrequencyparametrix}), the high frequency one being similar and more standard (see {\it e.g.} \cite{BoucletTzvetkov2} in the asymptotically Euclidean case). Note first that the $l$-th term in the sum (\ref{lowfrequencyparametrix}) is, for any $M$,  $ O_{{\mathscr L}_{\mu}^{-2M} \rightarrow {\mathscr L}_{\mu + l }^{2M}} (1) $ by (\ref{fourthbound}). Therefore, up to putting additional terms of the expansion in the remainder, it suffices to prove (\ref{lowfrequencyparametrix}) with a remainder satisfying, instead of (\ref{discussionreste}),
\begin{eqnarray}
  {\mathscr R}_{\rm low} (f,\epsilon) = O_{{\mathscr L}_{\mu}^{-2M_N} \rightarrow {\mathscr L}_{\mu_N }^{2M_N}} (1) , \qquad \mbox{with} \ \  M_N , \mu_N \rightarrow \infty \ \mbox{as} \ \ N \rightarrow \infty .  \label{conclusionreduite}
\end{eqnarray}
 Using  the Helffer-Sj\"{o}strand formula $ f (P/ \epsilon^2 ) = \int_{\Ca} \bar{\partial} \tilde{f}(z) (P/\epsilon^2-z)^{-1} L (dz) $ ($\tilde{f} \in C_0^{\infty}(\Ca)$ being an almost analytic extension of $ f $, see \cite{DimassiSjostrand}) together with Proposition \ref{leplusaffreux}, we get (\ref{lowfrequencyparametrix}) with a remainder which is a sum over $ \kappa $ of integrals of the form
$$ {\mathcal R}_{{\rm low} , \kappa}(f,\epsilon) =  \int_{\Ca} \bar{\partial} \tilde{f}(z)   \psi_{\kappa} (\epsilon r)  O \! p_{\epsilon,\kappa} ( r_{\epsilon,\kappa}(z))  \tilde{\tilde{\psi}}_{\kappa} (\epsilon r)    (P/ \epsilon^2 -z)^{-1} L (dz) $$
where $ r_{\epsilon,\kappa}(z) \in \widetilde{S}^{-N,-N} $  has semi-norms growing polynomially in $ |\mbox{Im}(z)|^{-1} $ (which is harmless since $ \bar{\partial} \tilde{f}(z) = O (|\mbox{Im}(z)|^{\infty}) $). In the above integral, we write 
$$  (P/ \epsilon^2 -z)^{-1} =  (P/ \epsilon^2 -z)^{-1}\big(  1 -  \zeta (\epsilon r) \big) +   (P/ \epsilon^2 -z)^{-1} \zeta (\epsilon r) . $$
Using Proposition \ref{propositionpourestes}, we observe that, for any $M$,  $  1 -  \zeta (\epsilon r) =  O_{{\mathscr L}_{\mu}^{-2M} \rightarrow {\mathscr L}_{0 }^{-2M}} (1) $, for it is compactly supported in $ \epsilon r $. We also have  $  (P/ \epsilon^2 -z)^{-1} =  O_{{\mathscr L}_{0}^{-2M} \rightarrow {\mathscr L}_{0 }^{-2(M-1)}} (|\mbox{Im}(z)|^{-1})  $ thanks to the spectral theorem. By Proposition \ref{propositionpourestes}, we also get that, for some $ \sigma = \sigma (M,N) $,
\begin{eqnarray}
   \psi_{\kappa} (\epsilon r)  O \! p_{\epsilon,\kappa} ( r_{\epsilon,\kappa}(z))  \tilde{\tilde{\psi}}_{\kappa} (\epsilon r)   =  O_{{\mathscr L}_{0}^{-2(M-1)} \rightarrow {\mathscr L}_{N }^{N-2(M-1)}} (|\mbox{Im}(z)|^{- \sigma }) .  \label{restexelow}
\end{eqnarray}  
All this implies that, for any given $ \mu  $ and $ N $,
\begin{eqnarray}
  {\mathcal R}_{{\rm low} , \kappa}(f,\epsilon) (1-\zeta(\epsilon r)) =   O_{{\mathscr L}_{\mu}^{-2 M} \rightarrow {\mathscr L}_{N }^{N-2(M-1)}} (1) .  \label{restelocallow}
\end{eqnarray}  
To analyse  $  {\mathcal R}_{{\rm low} , \kappa}(f,\epsilon)\zeta(\epsilon r) $, we use a parametrix for $ (P/\epsilon^2-z)^{-1} \zeta (\epsilon r) $ obtained analogously to the one of Proposition  \ref{leplusaffreux}: for any $ N^{\prime} \in \Na $,  $ (P/\epsilon^2-z)^{-1} \zeta (\epsilon r) $ is a sum of rescaled pseudo-differential operators with symbols in $ \widetilde{S}^{-2,0} $ and a remainder which is $ (P/\epsilon^2-z)^{-1} $ composed (to the right) with a sum of rescaled pseudo-differential operators with symbols in $ \widetilde{S}^{-N^{\prime},-N^{\prime}} $. This implies that, for any $ \mu $ and $ M$, and by choosing $ N^{\prime} > |\mu| $,
$  (P/\epsilon^2-z )^{-1} \zeta (\epsilon r) $ is of the form 
$$   O_{{\mathscr L}_{\mu}^{-2 M} \rightarrow {\mathscr L}_{\mu }^{-2(M-1)}} (|\mbox{Im}(z)|^{-\sigma^{\prime}}) + (P/\epsilon^2 - z)^{-1}  O_{{\mathscr L}_{\mu}^{-2 M} \rightarrow {\mathscr L}_{0 }^{N^{\prime}-2M}} (|\mbox{Im}(z)|^{-\sigma^{\prime}}) $$
for some $ \sigma^{\prime} = \sigma^{\prime} (M,N^{\prime}) > 0 $. Using an estimate similar to
 (\ref{restexelow}) together with the fact that $   (P/ \epsilon^2 -z)^{-1} =  O_{{\mathscr L}_{0}^{N^{\prime}-2M} \rightarrow {\mathscr L}_{0 }^{N^{\prime}-2(M-1)}} (|\mbox{Im}(z)|^{-1})   $, we get
\begin{eqnarray}
  {\mathcal R}_{{\rm low} , \kappa}(f,\epsilon) \zeta(\epsilon r) =   O_{{\mathscr L}_{\mu}^{-2 M} \rightarrow {\mathscr L}_{N }^{N-2(M-1)}} (1) . 
  \nonumber
\end{eqnarray}  
Together with (\ref{restelocallow}), this yields (\ref{conclusionreduite}) by choosing $ M = M_N = N/4 $  for instance. 
\finpreuve

\bigskip

As a first consequence of Theorem \ref{theoremcaclulfonctionnel}, we have the following estimates.

\begin{prop}[$ L^{\infty} \rightarrow L^{\infty} $ boundedness at spatial infinity] \label{clarificationLinfini} For all $ f \in C_0^{\infty} (\Ra) $,
$$ || \zeta (r) f (h^2 P) ||_{L^{\infty} \rightarrow L^{\infty}} \lesssim 1, \qquad h \in (0,1] $$
and
$$ || \zeta (\epsilon r) f (P/\epsilon^2) ||_{L^{\infty} \rightarrow L^{\infty}} \lesssim 1, \qquad \epsilon \in (0,1] . $$ 
\end{prop}

\noindent {\it Proof.} We consider only the low frequency case. The high frequency one is essentially standard, and can be proved e.g. as in \cite{BoucletFourier}. We thus consider $ \zeta (\epsilon r) f (P/\epsilon^2) $ which we expand using (\ref{lowfrequencyparametrix}). The (rescaled) pseudo-differential terms are bounded uniformly on $ L^{\infty} $ by Proposition \ref{bornepseudoLq}. Choosing $ M = N $ and $ \mu = - N $ in (\ref{discussionreste}), the remainder can be written
$$ {\mathscr R}_{\rm low} (f,\epsilon) = \tilde{\zeta} (\epsilon r)(P/\epsilon^2+1)^{-N} B_{\epsilon} \scal{\epsilon r}^{-N} $$
with $ || B_{\epsilon} ||_{L^2 \rightarrow L^2} \lesssim 1 $. This follows from Proposition \ref{propositionpourestes} and  that  $ \tilde{\zeta} (\epsilon r)  \zeta (\epsilon r) = \zeta (\epsilon r)$.
If $ N > n/2$,  we have $ || \scal{\epsilon r}^{-N} ||_{L^{\infty} \rightarrow L^2} \lesssim \epsilon^{-n/2} $ so, using the second estimate of  Proposition \ref{injectionsdeSobolevfinales} with $ \tilde{\zeta} $ instead of $ \zeta $, we get
$$ || \tilde{\zeta} (\epsilon r)(P/\epsilon^2+1)^{-N} B_{\epsilon} \scal{\epsilon r}^{-N} ||_{L^{\infty} \rightarrow L^{\infty}} \lesssim \epsilon^{n/2} \epsilon^{-n/2} \lesssim 1 $$
 which yields the result.
 \finpreuve

\bigskip

To illustrate another application of Proposition \ref{propositionpourestes}, we record  some rough a priori estimates on the propagator $ e^{-itP} $ which will be useful in Section \ref{sectionpropagation}. 
For $ k \geq 0 $ integer, we define $ \gamma (k) $ by $ \gamma (0) = 0 $ and $ \gamma (k+1) = 2 \gamma(k) + 1 $ ({\it i.e.} $ \gamma (k) = 2^k - 1 $).

\begin{prop}[Rough propagation estimates] \label{proproughprop}  For $ \mu \in \Ra $ denote by $ \lceil \mu \rceil $ the smallest integer $ \geq |\mu| $.  Then for all $ j \in \Za $,
\begin{eqnarray}
 e^{-itP} = O_{{\mathcal H}_{\mu}^{2j} \rightarrow {\mathcal H}_{\mu}^{2j- 2\gamma (\lceil \mu \rceil)}} \left( \scal{t/h}^{\gamma(\lceil \mu \rceil)} \right) , \label{scalingnaturel}
\end{eqnarray}
meaning that $ \scal{t/h}^{- \gamma(\lceil \mu \rceil)} e^{-itP} = O_{{\mathcal H}_{\mu}^{2j} \rightarrow {\mathcal H}_{\mu}^{2j-2\gamma(\lceil \mu \rceil)}}(1) $ uniformly in $t \in \Ra$.
Similarly
\begin{eqnarray}
 e^{-itP} =  O_{{\mathcal L}_{\mu}^{2j} \rightarrow {\mathcal L}_{\mu}^{2j- 2\gamma(\lceil \mu \rceil)}} \left( \scal{\epsilon^2 t}^{\gamma(\lceil \mu \rceil)} \right) . \label{scalingexplique}
\end{eqnarray}
\end{prop}

This proposition will be very useful to handle the remainders of some microlocal propagation estimates. The knowledge of the power $ \gamma(\lceil \mu \rceil) $ is not very important, the main interest being only the polynomial growth w.r.t to $ \scal{t/h} $ and $ \scal{\epsilon^2 t} $. We rather comment on the different scalings in $h$ and $ \epsilon $. The estimate (\ref{scalingnaturel}) reflects roughly that waves localized at frequency $ 1/h $ move at speed $ 1/h $. Based on this intuition, one could expect  to get a bound in term of $ \scal{ \epsilon t} $ in (\ref{scalingexplique}) for waves localized at frequency $ \epsilon $. The reason why we have  bounds in term of $ \scal{\epsilon^2 t} $ is that we use the rescaled spatial weights $ \scal{\epsilon r}^{\mu} $.  Another way to see that the scalings are natural is to consider the flat Laplacian on $ \Ra^n $ and to observe that for every symbol $a$ one has $ e^{i t \Delta} a (x,D) = a (x-2tD, D) e^{it\Delta} $
we see easily that
$$ e^{i t \Delta} a (x, h D) = a \left(x-2 \frac{t}{h} h D, hD\right) e^{it\Delta}   $$
and that
$$ e^{i t \Delta} a \left(\epsilon x,  \frac{D}{\epsilon} \right) = a \left( \epsilon x-2 t \epsilon^2 \frac{D}{\epsilon} ,  \frac{D}{\epsilon} \right) e^{it\Delta}  , $$
where the power $ \epsilon^2 $ on $t$ follows  both from writing $ D = \epsilon (D/\epsilon) $ and from the scaling in $x$.

We finally note that Proposition \ref{proproughprop} uses implicitly that $ {\mathscr S} ({\mathcal M}) $ is preserved by $ e^{-itP} $ (recall our convention to consider the $ {\mathscr H}_{\mu}^{2j} $ and $ {\mathscr L}_{\mu}^{2j} $ norms only on $ {\mathscr S}({\mathcal M}) $). This fact can be checked by routine arguments using exactly the commutator techniques involved in the proof below, but we omit this aspect  and focus only on the estimates in time.

\bigskip

\noindent {\it Proof of Proposition \ref{proproughprop}.} Let us show (\ref{scalingnaturel}). By (\ref{secondbound}), 
it suffices to show that $ \scal{r}^{\mu} e^{-itP} \scal{r}^{-\mu} $ satisfies the expected bound between $ {\mathcal H}_{0}^{2j} $ and $ {\mathcal H}_{0}^{2(j- \gamma(\lceil \mu \rceil))} $. If $ \mu = 0 $, this is a straightforward consequence of
$$ || (h^2 P + 1)^{j} e^{-itP} (h^2 P + 1)^{-j} ||_{L^2 \rightarrow L^2} = 1 . $$
Assume next that $ \lceil \mu \rceil = 1 $ and compute first the commutator
$$ \big[ \scal{r}^{|\mu|} , e^{-it P} \big] =  i \int_0^t e^{-i(t-s) P} [ P , \scal{r}^{|\mu|} ] e^{-isP} ds . $$
Using that  $ [ P , \scal{r}^{|\mu|} ] $ is $ h^{-1} $ times a sum of semiclassical differential operators with symbols in $ \widetilde{S}^{1,|\mu|-1} \subset \widetilde{S}^{2,0}  $ as in (\ref{fourthbound}) (they are supported in $ r \gg 1$ by (\ref{scalrM})),  we can write the commutator
$  \big[ \scal{r}^{|\mu|} , e^{-it P} \big] = O_{{\mathcal H}_{0}^{2j} \rightarrow {\mathcal H}_{0}^{2(j-1)}}  \left( |t|/ h \right) $. Thus, using that
$$   \scal{r}^{\mu} e^{-itP} \scal{r}^{-\mu}  = e^{-itP} + \begin{cases} \scal{r}^{-|\mu|} \big[ e^{-it P} , \scal{r}^{|\mu|}  \big] & \mbox{if} \ \mu < 0 \\ \mbox{} \\
  \big[ \scal{r}^{|\mu|} , e^{-it P} \big] \scal{r}^{-|\mu|} & \mbox{if} \ \mu \geq 0 \end{cases}  $$
 we get the result since $ \scal{r}^{-|\mu|} $ is bounded on each $ {\mathscr H}_0^{2k} $ by Proposition \ref{propositionpourestes}. If $ \lceil \mu \rceil > 1 $
we proceed by induction by writing, {\it e.g.} if $ \mu > 0 $,
$$ \scal{r}^{\mu} e^{-itP} \scal{r}^{-\mu} =  \scal{r}^{\mu-1} \left( e^{-itP} +  i \int_0^t e^{-i(t-s) P} [ P , \scal{r} ] e^{-isP} ds \scal{r}^{-1} \right)  \scal{r}^{1-\mu}  $$
The induction assumption and Proposition \ref{propositionpourestes} then show that the right hand side is of order
$$ O \left( \scal{t/h}^{\gamma (\lceil \mu \rceil - 1)} \right) + \int_0^t O \left(  \scal{(t-s)/h}^{\gamma (\lceil \mu \rceil - 1)} \right) O (h^{-1}) O    \left(  \scal{s/h}^{\gamma (\lceil \mu \rceil - 1)} \right) ds  $$ 
as an operator from $ {\mathcal H}_{\mu}^{2j} $ to $ {\mathcal H}_{\mu}^{2j- 2\gamma (\lceil \mu \rceil)} $. Using the definition of $ \gamma (.) $, we get (\ref{scalingnaturel}). The proof of (\ref{scalingexplique}) is similar, the gain in $ \epsilon^2 $ following  from the fact that
$$  [ P , \scal{\epsilon r } ]  = \epsilon^2 [ P/\epsilon^2 , \scal{\epsilon r} ]   = O_{{\mathcal L}_{\mu}^{2j} \rightarrow {\mathcal L}_{\mu}^{2j- 2}}  (\epsilon^2) $$ for all $ \mu $ and $j$ since the commutator in the middle is a linear combination of rescaled pseudo-differential operators as in (\ref{fourthbound}).  \finpreuve

\bigskip

\section{Spectral localizations} \label{sectionLittlewoodPaley}
\setcounter{equation}{0}
The purpose of this section is to prove  Theorems \ref{thlowfreqLq} and \ref{thhighfreqLq} which provide Littlewood-Paley type estimates, at low and high  frequencies respectively. Specific comments are given after each theorem. We only point out here that we adopt a pragmatic point of view, in the sense that we do not try to mimic exactly the usual form of Littlewood-Paley estimates on  $ \Ra^n $ ({\it e.g.} by using non trivial heat kernel bounds
) but rather provide robust and spatially localized versions of such estimates which seem naturally adapted to the proof of Strichartz estimates. In particular, the form of the decompositions are not the same at high and low frequencies; this is related to the fact that we use different types of estimates to treat the remainder terms.


   We use the function $ f_0 $ introduced in (\ref{defhighlow}) and consider  $ f (\lambda) = f_0 (\lambda) - f_0 (2 \lambda) $ so that $ f \in C_0^{\infty} (\Ra \setminus 0) $ and,  for all $ \lambda \in \Ra $, 
$$ \sum_{\ell = 0}^{\infty} f (2^{\ell} \lambda) =  {\mathds 1}_{\Ra \setminus 0} (\lambda) f_0 (\lambda) , \qquad  \sum_{\ell=1}^{\infty} f (2^{-\ell} \lambda ) = 1 - f_0 ( \lambda)  . $$
The spectral theorem then implies that, in the strong sense on $ L^2 ({\mathcal M}) $,
\begin{eqnarray}
 f_0 (P) = \sum_{\ell \geq  0} f (2^{\ell} P ) , \qquad (1 - f_0 ) (P) = \sum_{\ell \geq 1} f (2^{- \ell} P ) ,  \label{pour1ersomme}
\end{eqnarray} 
using in the first sum that  $0$ is not an eigenvalue of $ P $.

\subsection{Low frequencies}
In this paragraph we prove the following result.
\begin{theo} \label{thlowfreqLq} Assume that $ n \geq 3 $. Let $ \chi \in C_0^{\infty}(\Ra) $ be equal to $1$ on a  large enough interval so that $ (1-\chi)  = (1-\chi) \zeta  $ (see (\ref{partitionfixee})). Then
$$ ||f_0(P) v ||_{L^{2^*}} \lesssim \left(\sum_{\epsilon^2 = 2^{-\ell}} \big| \big| (1-\chi)(\epsilon r) f (P/\epsilon^{2}) v \big| \big|_{L^{2^*}}^2 +  \big| \big| \scal{r}^{-1} f (P/\epsilon^{2}) v \big| \big|_{L^{2}}^2\right)^{1/2} , $$
for all $v \in L^2 $. In the sum $\ell $ belongs to $ \Za_+ $.
\end{theo}
Let us comment that this Littlewood-Paley estimate holds for the exponent $ 2^* $ (and presumably for exponents between $ 2 $ and $  2^* $)
 which is sufficient and somewhat natural for applications to Strichartz estimates. Indeed, the first half of the sum is appropriately localized to use microlocal techniques  while the second one can be treated in a straightforward fashion by using the $ L^2 $ estimates (\ref{remplacerrefvide0})-(\ref{remplacerefvide}).

Theorem \ref{thlowfreqLq} is a consequence of the next two propositions in which we pick $ \tilde{f} \in C_0^{\infty} (\Ra \setminus 0;\Ra) $  such that $ \tilde{f} = 1 $ on $ \mbox{supp}(f) $ and let
\begin{eqnarray}
 \widetilde{\mathscr Q} (\epsilon) =  (1-\chi)(\epsilon r) \sum_{\kappa}  \psi_{\kappa} (\epsilon r) O \! p_{\epsilon,\kappa} \big( \tilde{f} \circ p_{\epsilon, \kappa} \big) \tilde{\psi}_{\kappa} (\epsilon r) , \label{formedeQepsilon}
\end{eqnarray} 
that is the first term of the parametrix of $ (1 - \chi)(\epsilon r) \tilde{f}(P / \epsilon^{2} )  $
according to Theorem \ref{theoremcaclulfonctionnel}.
Here and everywhere in this paragraph, we set $ \epsilon^2 = 2^{-\ell} $. 
\begin{prop} \label{resteL2LP} If $ n \geq 3 $, then 
\begin{eqnarray}
|| f_0 (P) v ||_{L^{2^*}} \lesssim  \sup_M \left| \left| \sum_{ \ell = 0 }^{ M} \widetilde{\mathscr Q} (\epsilon) (1-\chi)(\epsilon r) f (P / \epsilon^{2}  ) v \right| \right|_{L^{2^*}}  +  \left( \sum_{\ell \geq 0} \big| \big|  \scal{ r}^{-1} f ( P / \epsilon^{2} ) v \big| \big|_{L^2 }^2 \right)^{1/2}
\nonumber 
\end{eqnarray}
for all $ v \in L^2  $.  
\end{prop}
Up to the homogeneous Sobolev inequality (\ref{Sobolevestimate}) this proposition rests on purely $ L^2 \rightarrow L^2 $ estimates. In particular, we feel it is quite robust and could be used in other contexts ({\it e.g.} allow more general compactly supported perturbations).

To state the second proposition, we need to define the family of square functions
$$ \widetilde{S}_M w := \left( \sum_{\ell = 0}^M \big| \widetilde{\mathscr Q} (\epsilon)^* w \big|^2 \right)^{1/2} , \qquad M \geq 0 , $$
 where the adjoint is taken with respect to the Riemannian measure.

\begin{prop} \label{singularintegral} For all $ q_1 \in (1,2] $ there exists $ C > 0 $ such that 
$$ \big| \big| \widetilde{S}_M w \big| \big|_{L^{q_1} } \leq C || w ||_{L^{q_1} } , $$
for all $ M \geq 0 $ and all $ w \in C_0^{\infty}({\mathcal M}) $.
\end{prop}
This proposition is a consequence of fairly standard singular integral estimates, by exploiting the explicit form of the Schwartz kernel of $ \widetilde{\mathscr Q}(\epsilon) $. Note  that we do not need to assume  $ n \geq 3 $ here.

Before proving these two technical results, we prove Theorem \ref{thlowfreqLq}.

\bigskip

\noindent {\bf Proof of Theorem \ref{thlowfreqLq}.} Let us set $ S_M v := \left( \sum_{\ell = 0}^M \big| (1-\chi)(\epsilon r) f (  P / \epsilon^2 ) v \big|^2 \right)^{1/2}  $. Then, by the usual trick {\it i.e.} the Cauchy-Schwarz inequality in $\ell$ and the H\"older inequality in space, we have
$$ \left| \left( w , \sum_{\ell=0}^M  \widetilde{\mathscr Q} (\epsilon) (1-\chi)(\epsilon r) f ( P / \epsilon^2 ) v \right) \right| \leq || \widetilde{S}_M w ||_{L^{2_*}}  || S_M v ||_{L^{2^*}}  $$
so using Proposition \ref{singularintegral} for $ q_1 = 2_* $, we obtain
$$ \left| \left|\sum_{\ell=0}^M  \widetilde{\mathscr Q} (\epsilon) (1-\chi)(\epsilon r) f ( P / \epsilon^2 ) v \right| \right|_{L^{2^*}} \leq C   || S_M v ||_{L^{2^*}} . $$
We conclude by  using
$   || S_M v ||_{L^{2^*}} \leq \left( \sum_{\ell \geq 0} || (1-\chi)(\epsilon r) f ( P/ \epsilon^2 ) \psi ||_{L^{2^*}}^2 \right)^{1/2} $, which follows from the Minkowski inequality
since $ 2^* \geq 2 $, together with  Proposition \ref{resteL2LP}. \finpreuve

\bigskip

To prove Proposition \ref{resteL2LP}, we recall first for clarity the following well known results.
\begin{prop} \label{TTSchur} Let $( T_\ell )_{\ell} $ be a sequence of linear operators on a Hilbert space $ {\mathcal H} $. 
\begin{enumerate}
\item{(Discrete Schur estimate) If   $||T_j^* T_\ell||_{{\mathcal H}  \rightarrow {\mathcal H }} \lesssim 2^ {-|\ell-j|/2}$, then there is $ C $ such that
$$ \big| \big| \sum T_\ell v_\ell \big| \big|_{\mathcal H} \leq C \left( \sum ||v_\ell||_{\mathcal H}^2 \right)^{1/2} , $$
 for all sequence $ (v_\ell)  $ of $ {\mathcal H} $.}
\item{(Cotlar-Stein estimate)   If   $||T_j^* T_\ell||_{{\mathcal H}  \rightarrow {\mathcal H }} + ||T_j T_\ell^*||_{{\mathcal H}  \rightarrow {\mathcal H }} \lesssim 2^ {-|\ell-j|/2}$, then there is $ C $ such that  
$$ \big| \big| \sum T_\ell v \big| \big|_{\mathcal H} \leq C ||v||_{\mathcal H} , $$
for all $ v \in {\mathcal H} $.}
\end{enumerate}
\end{prop}

 
 We will apply the Schur estimate to two types of operators. The first one is very elementary: if we let 
 $$ T_\ell = 2^{\ell/2} P^{1/2} (2^\ell P + 1)^{-1} $$ then, assuming for instance $ \ell \geq j $ so that $ \frac{j+\ell}{2} = - \frac{|j-\ell|}{2} + \ell $, we have
\begin{eqnarray}
 || T_j^* T_\ell ||_{L^2 \rightarrow L^2} = 2^{-\frac{|j-\ell|}{2}} \big| \big| (2^j P + 1)^{-1} 2^\ell P (2^\ell P + 1)^{-1} \big| \big|_{L^2 \rightarrow L^2} \leq 2^{-\frac{|j-\ell|}{2}}  \label{TTstarelementaire}
\end{eqnarray}
by using the spectral sheorem. The second type of operators requires a lemma.
\begin{lemm} \label{lemmepseudoorthogonauxrescales} Let $ \kappa : U_{\kappa} \rightarrow V_{\kappa} $ be a chart on $ {\mathcal S} $. Let  $ \psi $ be smooth on $ {\mathcal M} $,  supported in $ (R_{\mathcal M},\infty ) \times U_{\kappa} $ and belonging to $ S^0 $.  For $ s =0 $ or $ 1 $, denote
$$ T_\ell = \big( P^{1/2}/\epsilon \big)^s      O\! p_{\epsilon,\kappa} (a_{\epsilon} + b_{\epsilon})   \psi ( \epsilon r,\omega)  , $$
where, for some given $ J \Subset (0,+\infty) $ (independent of $ \epsilon $),
 $$ (a_{\epsilon})_{\epsilon} \ \mbox{is bounded in} \ \widetilde{S}^{-\infty,0}, \qquad \emph{supp}(a_{\epsilon}) \subset p_{\epsilon , \kappa}^{-1}(J)
 ,  \qquad (b_{\epsilon})_{\epsilon} \ \mbox{is bounded in} \ \widetilde{S}^{-\infty,-1} , $$
 are all spatially supported in $ (R_{\mathcal M},\infty) \times V_{\kappa} $.
Then, if $ s = 0, 1 $,
\begin{eqnarray}
 \big| \big| T_j^* T_\ell \big| \big|_{L^2 \rightarrow L^2} \leq C 2^{-\frac{|j-\ell|}{2}}  \label{pourSchur}
\end{eqnarray} 
and, if $ s= 0 $,
\begin{eqnarray}
    \big| \big| T_j T_\ell^* \big| \big|_{L^2 \rightarrow L^2} \leq C 2^{-\frac{|j-\ell|}{2}}  . \label{pourCotlar}
\end{eqnarray}
\end{lemm}

\bigskip

\noindent {\it Proof.}  We start with two preliminary remarks. First, it suffices to prove both estimates when $ \ell \geq j $ (otherwise take the adjoint). The second one is   that, if $ s = 0 $,    $ T_\ell^* $ is of the same form as $ T_\ell $ (see (\ref{adjointgloballf})) up to perhaps changing the function $ \psi $. In particular, proving (\ref{pourSchur}) is sufficient. Let us prove (\ref{pourSchur}) when $s=1$. For simplicity, we set $ \psi (\breve{r}) = \psi (\breve{r},\kappa^{-1}(\theta)) $. Using Proposition \ref{borneL1Linfiniepsilon} with $q=2$, it suffices to show that
\begin{eqnarray}
 \left| \left|  O \! p^1 (c_{\epsilon_j}) \psi (\breve{r}) {\mathscr D}_{\epsilon_j}^{-1} \frac{P_{\kappa}}{\epsilon_j \epsilon_\ell} {\mathscr D}_{\epsilon_\ell} O\!p^1 (d_{\epsilon_\ell})  \right| \right|_{L^2 (\scal{\breve{r}}^{n-1} d \breve{r} d\theta) \rightarrow L^2 (\scal{\breve{r}}^{n-1} d\breve{r} d\theta) } \lesssim 2^{-\frac{|j-\ell|}{2}}  \label{reductionH1}
\end{eqnarray} 
with $ \epsilon_\ell = 2^{-\ell/2} $,  and $ (c_{\epsilon})_{\epsilon}, (d_{\epsilon})_{\epsilon} $ bounded families of $ \widetilde{S}^{-\infty,0} $ supported in $ (R_{\mathcal M},\infty) \times V_{\kappa} $ with respect to $ (\breve{r},\theta) $.      Using that $ \ell \geq j $ and (\ref{refoprescale}), we write 
$$ \frac{P_{\kappa}}{\epsilon_j \epsilon_\ell} {\mathscr D}_{\epsilon_\ell} = 2^{\frac{j-\ell}{2}} \frac{P_{\kappa}}{ \epsilon_\ell^2} {\mathscr D}_{\epsilon_\ell} =   2^{-\frac{|j-\ell|}{2}}   {\mathscr D}_{\epsilon_\ell}  P_{\epsilon_\ell, \kappa} . $$
 Since $ P_{\epsilon_\ell,\kappa} O\!p_1 (d_{\epsilon_\ell}) = O\!p_1 (e_{\epsilon_\ell}) $ for some bounded family $ (e_{\epsilon})_{\epsilon}  $ of $ \widetilde{S}^{-\infty,0} $ (with support contained in the one of $ d_{\epsilon} $), (\ref{reductionH1}) follows from (\ref{CalderonVaillancourtapoids}) together with 
 \begin{eqnarray}
   \left| \left| \psi (\breve{r})  {\mathscr D}_{\epsilon_j}^{-1}  {\mathscr D}_{\epsilon_\ell}  \right| \right|_{L^2 ((R_{\mathcal M},\infty) \times \Ra^{n-1}, \scal{\breve{r}}^{n-1} d\breve{r} d\theta) \rightarrow L^2 (\Ra^n,\scal{\breve{r}}^{n-1} d\breve{r} d\theta) } \lesssim 1 , \label{pseudounitaritecruciale} 
\end{eqnarray} 
 which follows for instance from the unitarity of $ {\mathscr D}_{\epsilon}^{\pm 1} $ on $ L^2 ( (0,\infty) \times \Ra^{n-1} , \breve{r}^{n-1}d \breve{r} d \theta ) $.
We next prove (\ref{pourSchur}) when $s=0$. It suffices to show that
\begin{eqnarray}
 \left| \left|  O \! p_1 (\tilde{b}_{\epsilon_j}) \psi (\breve{r}) {\mathscr D}_{\epsilon_j}^{-1}  {\mathscr D}_{\epsilon_\ell} O\!p_1 (d_{\epsilon_\ell}) \psi \right| \right|_{L^2 (\scal{\breve{r}}^{n-1} d \breve{r} d\theta) \rightarrow L^2 (\scal{\breve{r}}^{n-1} d\breve{r} d\theta) } \lesssim 2^{-\frac{|j-\ell|}{2}}  \label{reductionCS1}
\end{eqnarray}
and
\begin{eqnarray}
 \left| \left|  O \! p^1 (\bar{a}_{\epsilon_j}) \psi (\breve{r}) {\mathscr D}_{\epsilon_j}^{-1}  {\mathscr D}_{\epsilon_\ell} O\!p^1 (d_{\epsilon_\ell})  \right| \right|_{L^2 (\scal{\breve{r}}^{n-1} d \breve{r} d\theta) \rightarrow L^2 (\scal{\breve{r}}^{n-1} d \breve{r} d\theta) } \lesssim 2^{-\frac{|j-\ell|}{2}}  \label{reductionCS0}
\end{eqnarray}
 whenever $ (\tilde{b}_{\epsilon})_{\epsilon} \in \widetilde{S}^{-\infty,-1} $ and $ (d_{\epsilon})_{\epsilon} \in \widetilde{S}^{-\infty,0} $ are spatially supported in $ (R_{\mathcal M},\infty) \times V_{\kappa} $.
 To prove (\ref{reductionCS1}), we use 
 $$ \breve{r}^{-1} {\mathscr D}_{\epsilon_j}^{-1} {\mathscr D}_{\epsilon_\ell} = \epsilon_\ell \epsilon_j^{-1} {\mathscr D}_{\epsilon_j}^{-1} {\mathscr D}_{\epsilon_\ell} \breve{r}^{-1} $$ to write
 $$ O \! p^1 (\tilde{b}_{\epsilon_j}) \psi (\breve{r}) {\mathscr D}_{\epsilon_j}^{-1}  {\mathscr D}_{\epsilon_\ell} = O \! p^1 (\tilde{b}_{\epsilon_j}) \psi (\breve{r}) \breve{r} \breve{r}^{-1} {\mathscr D}_{\epsilon_j}^{-1}  {\mathscr D}_{\epsilon_\ell} = 2^{-\frac{|j-\ell|}{2}} O \! p^1 (\tilde{b}_{\epsilon_j}) \breve{r} \psi {\mathscr D}_{\epsilon_j}^{-1}  {\mathscr D}_{\epsilon_\ell}   \breve{r}^{-1} $$
 and conclude again from the $ L^2 (\scal{\breve{r}}^{n-1}dr d \theta) $  boundedness of $ O \! p^1 (\tilde{b}_{\epsilon}) \breve{r}  $ and
 $ \breve{r}^{-1}  O\!p^1 (d_{\epsilon}) $ (there is no singularity at $ \breve{r} = 0 $ since $ d_{\epsilon} $ is supported in $ \{  r \geq R_{\mathcal M}\} $) together with (\ref{pseudounitaritecruciale}).
We finally prove (\ref{reductionCS0}).  The support assumption on $ a_{\epsilon} $ implies that $ a_{\epsilon} / p_{\epsilon,\kappa} $ is a smooth symbol in $ \widetilde{S}^{-\infty,0} $ so, using in addition that $ \psi \in S^0 $, we can write by symbolic calculus
$$ O \! p^1 (\bar{a}_{\epsilon}) \psi (\breve{r}) = O \! p^1 (\bar{a}_{\epsilon}/ p_{\epsilon,\kappa}) \psi P_{\epsilon,\kappa} + O \! p^1 \big( \tilde{\tilde{b}}_{\epsilon} \big) \tilde{\psi} (\breve{r}) $$
with $ ( \tilde{\tilde{b}}_{\epsilon} )_{\epsilon}  $ bounded in $ \widetilde{S}^{-\infty,-1} $ and some cutoff $ \widetilde{\psi} (\breve{r},\theta) \in S^0 $,  both supported in $ (R_{\mathcal M}, \infty ) \times V $ with respect to $ (\breve{r},\theta) $. The contribution of the second term in the right hand side follows from
  (\ref{reductionCS1}). For the first term, one can use (\ref{reductionH1}) once observed that
$$ O \! p^1 (\bar{a}_{\epsilon_j}/ p_{\epsilon_j,\kappa}) \psi (\breve{r}) P_{\epsilon_j,\kappa} {\mathscr D}_{\epsilon_j}^{-1}  {\mathscr D}_{\epsilon_\ell} = O \! p^1 (\bar{a}_{\epsilon_j}/ p_{\epsilon_j,\kappa} ) \psi (\breve{r}) {\mathscr D}_{\epsilon_j}^{-1}  \frac{P_{\kappa}}{\epsilon_j^2} {\mathscr D}_{\epsilon_\ell}  $$
and that $ \epsilon_j^{-2} = 2^{\frac{j-\ell}{2}}( \epsilon_j \epsilon_\ell )^{-1} $ (so that we actually get an estimate of order $ 2^{-|j-\ell|} $ for this term).  This completes the proof. \finpreuve

\bigskip 

\noindent {\bf Proof of Proposition \ref{resteL2LP}.} Let us write $ (1-\chi)(\epsilon r) \tilde{f}(P/\epsilon^2)  = \widetilde{\mathscr Q} (\epsilon) + \widetilde{\mathscr R} (\epsilon)  $.
Using that $ f \tilde{f} = f $ and that $ 1 = \chi (\epsilon r) + (1-\chi)(\epsilon r) $, we have
$$ f (P / \epsilon^2) = \widetilde{\mathscr Q}(\epsilon) (1-\chi)(\epsilon r) f (P/\epsilon^2) + \epsilon T (\epsilon) \scal{ \epsilon r}^{-1} f (P/\epsilon^2) $$
with
\begin{eqnarray}
 T (\epsilon) = \epsilon^{-1} \left( \chi (\epsilon r) \tilde{f}(P/\epsilon^2) + (1-\chi)(\epsilon r) \tilde{f}(P/ \epsilon^2) \chi (\epsilon r) +   \widetilde{\mathscr R}(\epsilon) \right) \scal{\epsilon r} .  \label{lessobvious}
\end{eqnarray} 
Using the first sum in (\ref{pour1ersomme}) (the strong convergence also holds in $ L^{2^*} $ by Sobolev embedding) and the homogeneous Sobolev estimate (\ref{Sobolevestimate}), we have
$$ ||f_0(P)v||_{L^{2^*}}  \lesssim   \sup_M \left(  \left| \left| \sum_{ \ell = 0 }^{ M} \widetilde{\mathscr Q} (\epsilon) (1-\chi)(\epsilon r) f (P / \epsilon^{2}  ) v \right| \right|_{L^{2^*}}  +  \left| \left| \sum_{\ell = 0}^M   P^{1/2} T (\epsilon)  \epsilon \scal{\epsilon r}^{-1} f ( P / \epsilon^{2} ) v \right| \right|_{L^2 } \right) $$
where it suffices to estimate the second norm. Using Theorem \ref{theoremcaclulfonctionnel}, one can write
\begin{eqnarray}
 P^{1/2} T (\epsilon) =   \epsilon^{-1} P^{1/2} (P/\epsilon^2 + 1)^{-1} B (\epsilon) ,   \label{withthisestimateathand}
 \end{eqnarray}
with $ B (\epsilon) $ bounded on $ L^2  $ uniformly in $ \epsilon $. The less obvious contribution of terms of (\ref{lessobvious}) is the uniform $ L^2 $ boundedness  of  $ (P/\epsilon^2 +1) \chi (\epsilon r) \tilde{f}(P/ \epsilon^2) \scal{\epsilon r} $. One can analyze it as follows. On one hand, the commutator $ [ (P/\epsilon^2 +1) , \chi (\epsilon r)  ] $ being a sum of rescaled (pseudo-)differential operators vanishing outside the support of  $ \zeta (\epsilon r)  $, one can use Theorem \ref{theoremcaclulfonctionnel} to get a parametrix for  $  [ (P/\epsilon^2 +1) , \chi (\epsilon r)  ]  \tilde{f}(P/ \epsilon^2) \scal{\epsilon r}  $ from which the uniform $ L^2 $ boundedness follows. On the other hand, $ \chi (\epsilon r)  (P/\epsilon^2 +1)  \tilde{f}(P/ \epsilon^2) \scal{\epsilon r}  = \chi (\epsilon r)  \tilde{f}_1(P/ \epsilon^2) \scal{\epsilon r}   $ with $ \tilde{f}_1 \in C_0^{\infty} $. We then write $  \scal{\epsilon r}  = \chi (\epsilon r) \scal{\epsilon r}  +  (1 - \chi)(
 \epsilon r) \scal{\epsilon r}  $ whose first term is obviously uniformly bounded on $L^2$ while one can use the parametrix for $ \tilde{f}_1(P/ \epsilon^2)(1- \chi) (\epsilon r)  $ to see that $  \chi (\epsilon r)  \tilde{f}_1(P/ \epsilon^2) (1-\chi)(\epsilon r) \scal{\epsilon r} $ is uniformly bounded on $L^2$.
 Now with (\ref{withthisestimateathand}) at hand, by using (\ref{TTstarelementaire}) and Lemma \ref{lemmepseudoorthogonauxrescales} with $s=1$  together with the Schur estimate of Proposition \ref{TTSchur}, we have
\begin{eqnarray}
  \sup_M   \left| \left| \sum_{\ell = 0}^M   P^{1/2} T (\epsilon)  \epsilon \scal{\epsilon r}^{-1} f ( P / \epsilon^{2} ) v \right| \right|_{L^2 } \leq C \left( \sum_{k \geq 0} || \epsilon \scal{\epsilon r}^{-1} f (P/\epsilon^2)v ||_{L^2}^2 \right)^{1/2} . 
  \nonumber
 \end{eqnarray}
In the right hand side of this inequality, we finally use that 
$$ \epsilon \scal{\epsilon r}^{-1} \lesssim \scal{r}^{-1} $$
 and we get the result. \finpreuve

\bigskip

We now consider the proof of Proposition \ref{singularintegral}.

\bigskip

\noindent {\bf Proof of Proposition \ref{singularintegral}.} It follows the same line as the one for the standard Littlewood-Paley decomposition (see {\it e.g.} \cite{MuSc}).  Let $ (\varrho_\ell )_{\ell \geq 0} $ be the usual Rademacher sequence (realized as functions of $ t \in [0,1] $). By the Khintchine inequality, it suffices to show that
$$ \left| \left|  \sum_{\ell =0}^{ M} \varrho_\ell(t) \widetilde{\mathscr Q}(\epsilon)^* \right| \right|_{L^{q_1} ({\mathcal M}) \rightarrow L^{q_1} ({\mathcal M})} \leq C , \qquad t \in [0,1], \ M \geq 0 . $$
This in turn follows from the Marcinkiewicz interpolation theorem provided we prove the above estimate for $ q_1 = 2 $ as well as weak type $ (1,1) $ estimates for $ \sum_{\ell \leq M} \varrho_\ell(t) \widetilde{\mathscr Q}(\epsilon)^* $  uniformly in $ t$ and $ M $. Using the form of $ \widetilde{Q}(\epsilon) $ given by (\ref{formedeQepsilon}), the uniform $ L^2 \rightarrow L^2 $ bound follows from the Cotlar-Stein estimate of Proposition \ref{TTSchur} together with the estimates (\ref{pourSchur}) and (\ref{pourCotlar}) (with $s=0$) of Lemma \ref{lemmepseudoorthogonauxrescales}. The weak type $ (1,1)$ estimate follows from standard estimates on Calder\'on-Zygmund operators as explained in Section \ref{sectionCalderonZygmund}.  \finpreuve

\subsection{High frequencies}
The purpose of this paragraph is to prove the following result.
\begin{theo} \label{thhighfreqLq}  Let $ N \geq 0 $ and $ \chi \in C_0^{\infty}(\Ra) $ be equal to $1$ on a large enough set  so that $ \zeta \equiv 1 $ near the support of $ 1 - \chi $. Let $ q \in [2,\infty) $. Then
$$ ||(1-\chi)(r) (1 - f_0)(P) v ||_{L^{q}} \lesssim \left(\sum_{h^2 = 2^{-\ell}} \big| \big| (1-\chi)( r) f (h^2 P) v \big| \big|_{L^{q}}^2 +  h^N \big| \big| \scal{r}^{-N} f (h^2 P) v \big| \big|_{L^{2}}^2\right)^{1/2} , $$
for all $v  \in {\mathscr S}({\mathcal M}) $. In the sum $\ell $ belongs to $ \Na $.
\end{theo}
This  is a spatially localized Littlewood-Paley decomposition similar to the one of \cite{BoucletSMF}. The improvement  here is that the nonlocal $ L^2 $ correction involves the weight $ \scal{r}^{-N} $ which will allow us to use the resolvent estimates (\ref{aprioriresolvente}) and their time dependent counterparts (see paragraph \ref{soussectionresolvante} and Section \ref{SectionStrichartz}).

To prove this theorem, we pick again $ \tilde{f} \in C_0^{\infty} (\Ra \setminus 0;\Ra) $  such that $ \tilde{f} = 1 $ on $ \mbox{supp}(f) $. 
We define the square functions
$$ \Sigma_M v = \left( \sum_{\ell =1}^M |(1-\chi)(r)f(h^2 P)v|^2 \right)^{1/2} $$
and 
$$ \widetilde{\Sigma}_M w =   \left( \sum_{\ell =1}^M | \zeta (r)\tilde{f}(h^2 P)w|^2 \right)^{1/2} . $$
Here and throughout this paragraph, we set $ h^2 = 2^{-\ell} $. 

\bigskip

\noindent {\bf Proof of Theorem \ref{thhighfreqLq}.} It is very close to that of Theorem \ref{thlowfreqLq}. We only explain the new arguments. Using the second sum in (\ref{pour1ersomme}), we write $ \big( w , (1-\chi)(r) (1-f_0)(P) v \big) $ as the limit as $ M \rightarrow \infty $ of $ \sum_{\ell=1}^M (w,(1-\chi)(r) f (h^2 P) v) $.  By standard semiclassical estimates based on Theorem \ref{theoremcaclulfonctionnel} and the fact that $ (1- \tilde{f}) $ vanishes near the support of $f$, we see that, for any $ N $,
$$ ( 1 - \tilde{f})(h^2 P) (1-\chi)(r) f (h^2 P) = h^N B_N (h) \scal{r}^{-N} f (h^2P ) , $$
with
$$  || B_N (h) ||_{L^2 \rightarrow L^q} \leq C, \qquad h \in (0,1] . $$
Therefore, using additionally that $ \zeta (r) (1-\chi)(r) = (1-\chi)(r) $, we have
\begin{eqnarray}
 \big|  \big( w , (1-\chi)(r) (1-f_0)(P) v \big) \big| & \leq & \sup_M \left| \sum_{\ell=1}^M \big(\zeta (r)\tilde{f}(h^2P )w,(1-\chi)(r) f (h^2 P) v \big) \right| + \nonumber \\
 & &\ \ C ||w||_{L^{q^{\prime}}} \sum_{\ell \geq1}  h^N || \scal{r}^{-N} f (h^2P ) v ||_{L^2} .  \nonumber
\end{eqnarray} 
By proceeding as in the proof of Theorem \ref{thlowfreqLq}, in particular by using that the supremum above is bounded by $ \sup_M || \widetilde{\Sigma}_M w ||_{L^{q^{\prime}}} ||\Sigma_M v ||_{L^q} $, Theorem \ref{thhighfreqLq} follows from Proposition \ref{squarehautefrequence} below. \finpreuve

\bigskip

\begin{prop} \label{squarehautefrequence} For all $ q_1 \in (1,2] $, there exists $ C > 0 $ such that 
$$ \big| \big| \widetilde{\Sigma}_M w \big| \big|_{L^{q_1}} \leq C ||w||_{L^{q_1}} $$
for all $M \geq 1 $ and all $ w \in {\mathscr S} ({\mathcal M}) $.
\end{prop}

\noindent {\it Proof.} As in the proof of Proposition \ref{singularintegral}, it suffices to show that $ \sum_{\ell = 1}^M \varrho_\ell (t) \zeta (r) \tilde{f}(h^2P) $ is bounded on $ L^2 $ and satisfies weak type $ (1,1) $ estimates, uniformly in $t$ and $M$. The uniform boundedness on $ L^2 $ follows from the spectral theorem and the fact that the  functions
$$  \lambda \mapsto   \sum_{\ell = 1}^M \varrho_\ell (t)  \tilde{f}(2^{-\ell} \lambda)  $$
belong to $ L^{\infty}(\Ra) $ uniformly in $t,M$ since at most a finite number ($ \lambda , M,t $ independent) of  terms of the sum do not vanish. To prove the weak type $(1,1) $ estimate, we use Theorem \ref{theoremcaclulfonctionnel} to decompose
$$ \zeta (r) \tilde{f}(h^2P) = {\mathscr Q}_{\rm high} (h)  + h {\mathcal R}_{\rm high}(h)    $$
with $ {\mathcal R}_{\rm high}(h) $ uniformly (in $h$) bounded on $ L^1 $ and $ L^2 $. The uniform boundedness on $ L^2 $ is obvious. To see the uniform boundedness on $  L^1 $, one uses an expansion of $ \zeta (r) \tilde{f}(h^2P) $ to a sufficiently high order $ N_0 + 1 $ so that one can write   $ h {\mathcal R}_{\rm high}(h)  = h^{1+N_0} \scal{r}^{-N_0} B (h) (h^2 P + 1)^{-N_0}  $ with $ B (h) $ uniformly bounded on $ L^2 $. Then using on one hand that $  (h^2 P + 1)^{-N_0} : O_{L^1 \rightarrow L^2} (h^{-n/2})  $ (by taking the adjoint estimate of (\ref{bytakingtheadjointestimate}) near infinity and using a standard elliptic regularity estimate on any compact set -- including near the boundary if any) and on the other hand that $ \scal{r}^{-N_0} : L^2 \rightarrow L^1 $, one get the desired $ L^1 \rightarrow L^1 $ estimate. In particular, we have
\begin{eqnarray}
 \sum_{\ell \geq 1} h \big| \big|  {\mathcal R}_{\rm high}(h) \big| \big|_{L^1 \rightarrow L^1} < \infty . \label{noneasy}
 \end{eqnarray}
Then, it suffices to prove the uniform weak type $ (1,1) $ estimates for $ \sum_{\ell = 1}^M \varrho_\ell (t)  {\mathscr Q}_{\rm high} (h)  $
and this follows again from standard arguments on Calder\'on-Zygmund operators (see Appendix \ref{sectionCalderonZygmund}).
\finpreuve

\bigskip

\noindent {\bf Remark.} In more general situations, {\it e.g.} with non smooth coefficients in a compact set, it may be not easy to prove (\ref{noneasy}). Actually, it would suffice to have $ \sum_{\ell \geq 1} h \big| \big|  {\mathcal R}_{\rm high}(h) \big| \big|_{L^{2_*} \rightarrow L^{2_*}} < \infty  $ for our purpose. It would restrict the range of exponents in Proposition \ref{squarehautefrequence} to $ [2_*,2] $, and thus those of Theorem \ref{thhighfreqLq} to $ [2,2^*] $, but this would be sufficient for Strichartz estimates.

\section{Classical scattering} \label{sectionscatteringclassique}   
\setcounter{equation}{0}

In this section, we construct real phase functions solutions to  Hamilton-Jacobi equations that will be used to construct  Isozaki-Kitada type parametrices. The transport equations associated to such parametrices  are also studied.

Everywhere in this section, we work in a single  chart at infinity  $ (R_{\mathcal M},\infty) \times V_{\kappa} $.
Since we want to consider both  high and low frequencies parametrices, we have to analyze the Hamiltonian flow of $ p_{\kappa} $ and $ p_{\epsilon,\kappa} $. Observing that $ p_{\kappa} = p_{\kappa,1} $, we will state the main results only for   $ p_{\epsilon,\kappa} $ for $ 0 < \epsilon \leq 1 $. 

We  let
$  \phi^s_{\epsilon,\kappa}  $ be the Hamiltonian flow of $ p_{\epsilon,\kappa} $ as well as
$$ \phi_0^s (r,\vartheta, \varrho,\eta) := \big( r + 2 s \varrho , \vartheta , \varrho , \eta \big) $$
be the Hamiltonian flow of $ \varrho^2 $. We denote the time by $s$ here since  it will be interpreted as a rescaled version of $t$ in the applications (either $ s = t/h $ or $ s = t \epsilon^2 $).

  For $ R \gg 1 $, $ V \subset V_{\kappa} $,  and $ \varepsilon > 0 $, we define the subset of $ \Ra \times (\Ra^{n-1})^2 $
\begin{eqnarray}  
 \Theta (R,V,\varepsilon) = \left \{ (r,  \theta  ,  \vartheta ) \ | \ r > R, \ \theta \in V , \ |\theta - \vartheta| < \varepsilon  \right\} .  
\label{Thetadef}
\end{eqnarray}
To describe the asymptotic behaviour of our phases, and to  take into account the dependence on $ \epsilon $ of the functions we are going to consider  ({\it e.g.}  the components of the flow $ \phi^s_{\epsilon,\kappa} $), the following definition will be useful.
\begin{defi} \label{defi-asymp-symbol} Let $ R > 0 $, $ V \subset V_{\kappa} $ and $ \varepsilon  > 0 $. For $ \mu \in \Ra $,
\begin{enumerate}
\item{$ S_{\mu} $ is the set of ($\epsilon$ dependent families of) functions on $ \Theta (R,V,\varepsilon) $ such that
$$  \big| \partial_r^j \partial_{\theta}^{\alpha}  \partial_{\vartheta}^{\beta}  a_{\epsilon} (r,\theta,\vartheta) \big| \leq C r^{\mu - j},  $$
for all $  (r,\theta,\vartheta) \in \Theta(R,V,\varepsilon)  $ and all $ \epsilon \in (0,1] $ (the constant is independent of $ \epsilon $).}
\item{For any integer $ m \geq 0 $,  we denote by $ S_{\mu} (\theta - \vartheta)^m $ the set of all functions of the form
$$ \sum_{|\gamma| = m} a_{\epsilon,\gamma} (r,\theta,\vartheta) (\theta - \vartheta)^{\gamma} ,  $$
with $ a_{\epsilon,\gamma} \in S_{\mu} $.}
\item{Given 
real numbers $ \mu_1 , \mu_2 $ and integers $ m_1,m_2 $, the equality
$$ a_{\epsilon} = b_{\epsilon} + S_{\mu_1} (\vartheta - \theta)^{m_1} +  S_{\mu_2} (\vartheta - \theta)^{m_2} $$
means that $ a_{\epsilon} - b_{\epsilon} $ is the sum of an element of $ S_{\mu_1} (\vartheta - \theta)^{m_1} $ and a one of $ S_{\mu_2} (\vartheta - \theta)^{m_2} $.}
\end{enumerate}
\end{defi}

The main result of this section is the following theorem.

\begin{theo}[Eikonal equation] \label{eikonal-equation}  Fix an open subset $ V \Subset V_{\kappa} $. Assume that $ V $ is convex. Then we can find $ R \gg 1 $ and $ 0 < \varepsilon \ll 1 $ such that for all $ \epsilon \in (0,1] $, there exists a  smooth function 
$$ \psi_{\epsilon} :  \Theta (R,V,\varepsilon) \rightarrow \Ra $$
such that the function $ \varphi_{\epsilon} (r,\theta,\varrho, \vartheta) : = \varrho \psi_{\epsilon} (r,\theta,\vartheta) $ satisfies the following properties:
\begin{enumerate}
\item{It solves the equation 
\begin{eqnarray}
 p_{\epsilon,\kappa}  \big(r,\theta, \partial_r \varphi_{\epsilon} , \partial_{\theta} \varphi_{\epsilon} \big) = \varrho^2  .  \label{eikonaladeriver}
\end{eqnarray}}
\item{The range of 
$$ (r,\theta,\varrho,\vartheta) \mapsto (r,\theta,\partial_r \varphi_{\epsilon} , \partial_{\theta} \varphi_{\epsilon}), \qquad (r,\theta,\vartheta) \in \Theta (R,V,\varepsilon), \ \pm \varrho > 0 , $$ is contained in a set  $ \widetilde{\Gamma}^{\pm}_{\rm st} $ where $ (\phi^s_{\epsilon,\kappa})_{\pm s \geq 0} $ and the limit $ \lim_{s \rightarrow \pm \infty} \phi_0^{-s} \circ \phi^s_{\epsilon,\kappa}  =: F^{\pm}_{\epsilon, \kappa} $ are defined. Furthermore,  one has
\begin{eqnarray}
 F^{\pm}_{\epsilon,\kappa} \big(r,\theta, \partial_r \varphi_{\epsilon} , \partial_{\theta} \varphi_{\epsilon}  \big) = \big( \partial_{\varrho} \varphi_{\epsilon}  , \vartheta , \varrho , - \partial_{\vartheta} \varphi_{\epsilon} \big) . \label{fonctioncaracteristiqueparametres}
 \end{eqnarray}}
\item{One has the expansions
\begin{eqnarray}
 \psi_{\epsilon} & = & r  + S_{1-\nu} (\vartheta - \theta) + S_1 (\vartheta - \theta)^2 .  \label{expansioneikonal4} \\
\partial_r \psi_{\epsilon} & = & 1 + S_{-\nu} (\vartheta - \theta)   + S_0 (\vartheta - \theta)^2    \label{expansioneikonal1} \\
\partial_{\theta} \psi_{\epsilon} & = &   r    \bar{g} (\theta) ( \vartheta - \theta ) +  S_{1-\nu}  (\vartheta - \theta) +  S_1 (\vartheta - \theta)^2   \label{expansioneikonal2} \\
\partial_{\vartheta} \psi_{\epsilon} & = &  - r    \bar{g} (\theta) ( \vartheta - \theta ) +  S_{1-\nu}  (\vartheta - \theta) +  S_1 (\vartheta - \theta)^2   \label{expansioneikonal3}  
\end{eqnarray}}
\end{enumerate}
\end{theo}

\noindent {\bf Remark.} Be careful not to mistake $ \epsilon $ (the low frequency parameter) for $ \varepsilon $ which is a small enough but fixed number defining $ \Theta (R,V,\varepsilon) $.

\bigskip

The purpose of the next proposition is to solve transport equations associated to $ \varphi_{\epsilon} $ and which will  be used in Section \ref{sectionIsozakiKitada}. We consider equations of the form
\begin{eqnarray}
( \partial_{\rho,\eta} p_{\epsilon,\kappa}) \big(r,\theta,\partial_r \varphi_{\epsilon}, \partial_{\theta} \varphi_{\epsilon} \big) \cdot \partial_{r,\theta} u + b_{\epsilon} (r,\theta,\vartheta,\varrho)  u = f_{\epsilon} (r,\theta,\varrho,\vartheta) , \label{transport}
\end{eqnarray}
where  $f_{\epsilon}$ is a given short range symbol (see condition (\ref{shortrangesymbol})) and
\begin{eqnarray}
 b_{\epsilon} := - P_{\epsilon,\kappa} \varphi_{\epsilon} . 
\label{bepsilon}
\end{eqnarray}
In practice, we will study these equations only locally in $ \varrho $, namely on sets of the form
\begin{eqnarray}
 \Theta^{\pm} (R,V,I,\varepsilon) := \{ (r,\theta,\varrho,\vartheta) \ | \ (r,\theta,\vartheta) \in \Theta (R,V,\varepsilon) , \  \varrho^2 \in I , \ \pm  \varrho  > 0  \} 
 \label{Thetaenergie}
\end{eqnarray}
where $ I \Subset (0,+\infty) $ is a given relatively compact interval.
 The natural domains to work on are actually  the  larger sets (of trajectories starting in $  \Theta^{\pm} (R,V,I,\varepsilon) $)
$$ {\mathcal T}^{\pm}_{\epsilon} (R,V,I,\varepsilon) := \big\{  \big( \big( \bar{r}^s_{\epsilon} , \bar{\vartheta}^s_{\epsilon} \big) ( r , \theta , \partial_{r,\theta} \varphi_{\epsilon} (r,\theta,\varrho,\vartheta) ) , \varrho, \vartheta \big)  \  | \ (r,\theta,\varrho,\vartheta) \in \Theta^{\pm}(R,V,I,\varepsilon), \ \pm s \geq 0 \big\} , $$
where
\begin{eqnarray}
 \big( \bar{r}^s_{\epsilon} , \bar{\vartheta}^s_{\epsilon} , \bar{\varrho}^s_{\epsilon}, \bar{\eta}^s_{\epsilon} \big) = \ \mbox{components of} \ \phi^s_{\epsilon,\kappa} . \label{componentsofphikappa}
\end{eqnarray}
It will follow from the proof below that   $ \big( r , \theta , \partial_{r,\theta} \varphi_{\epsilon} (r,\theta,\varrho,\vartheta) \big) $ belongs to
a set where the flow $ \phi^s_{\epsilon, \kappa} $ is well defined for all $ \pm s \geq 0 $ (if $  \pm \varrho  > 0 $) so that the sets $ {\mathcal T}^{\pm}_{\epsilon} (R,V,I,\varepsilon)  $ are well defined.
\begin{prop}[Transport equations]  \label{prop-transport} Let $ \Theta (R,V,\varepsilon) $ be as in Theorem \ref{eikonal-equation} and $ I \Subset (0,+\infty) $.
\begin{enumerate}
\item{{\bf Form of characteristics:} For all $ (r,\theta,\varrho,\vartheta) \in  \Theta^{\pm}(R,V,I,\varepsilon)  $,  $ \pm s \geq 0 $ and $ \epsilon \in (0,1] $ define
$$ \big( \check{r}^s_{\epsilon} , \check{\theta}^s_{\epsilon} , \check{\rho}^s_{\epsilon} , \check{\eta}^s_{\epsilon} \big) := \phi^s_{\epsilon,\kappa} \big( r , \theta , \partial_{r,\theta} \varphi_{\epsilon} (r,\theta,\varrho,\vartheta) \big). $$
Then
$$  \big( \check{\rho}^s_{\epsilon} , \check{\eta}^s_{\epsilon} \big) =  \big( \partial_{r,\theta} \varphi_{\epsilon} \big)  \big( \check{r}^s_{\epsilon} , \check{\theta}^s_{\epsilon}  , \varrho , \vartheta  \big) . $$
In particular,
$$ ( \partial_{\rho,\eta} p_{\epsilon,\kappa}) \big( \check{r}^s_{\epsilon} , \check{\theta}^s_{\epsilon} , ( \partial_{r,\theta}  \varphi_{\epsilon} ) (\check{r}^s_{\epsilon} , \check{\theta}^s_{\epsilon} , \varrho , \vartheta ) \big) =  \big( \dot{\check{r}}^s_{\epsilon} , \dot{\check{\theta}}^s_{\epsilon} \big)  . $$}
\item{{\bf Time integrability of} $ b_{\epsilon} $ {\bf along characteristics:} For all $ j,\alpha,k,\beta $, there exists $ C $ independent of $ \epsilon \in (0,1] $ such that
\begin{eqnarray}
 \left| \partial_r^j \partial_{\theta}^{\alpha} \partial_{\varrho}^k \partial_{\vartheta}^{\beta} \left(  b_{\epsilon} (\check{r}^s_{\epsilon} , \check{\theta}^s_{\epsilon} , \varrho, \vartheta) \right) \right| \leq  C \scal{s/r}^{-1-\nu} r^{-1-\nu-j}  + C \scal{s/r}^{-2} r^{-1-j} \label{passhortrange}
\end{eqnarray}
for $ \pm s \geq 0 $ and $ (r,\theta,\varrho,\vartheta) \in  \Theta^{\pm}(R,V,I,\varepsilon)  $. }
\item{{\bf Form of solutions:} Assume that $f_{\epsilon}$ belongs to  $ S_{-1-\mu} := S_{-1-\mu} \big(   {\mathcal T}^{\pm}_{\epsilon}(R,V,I,\varepsilon) \big) $ for some $ \mu > 0 $, {\it i.e.} on $   {\mathcal T}^{\pm}_{\epsilon}(R,V,I,\varepsilon)  $
\begin{eqnarray}
  | \partial_r^j \partial_{\theta}^{\alpha} \partial_{\varrho}^k \partial_{\vartheta}^{\beta}  f_{\epsilon}(r,\theta,\vartheta,\varrho)| \lesssim \scal{r}^{-1-\mu - j} ,  \label{shortrangesymbol}  
\end{eqnarray}
uniformly in $ \epsilon $.
Then, given a constant $ C $, the solution to (\ref{transport}) going to $  C $ as $ r \rightarrow \infty $ is given by
$$  C \exp \left( \int_0^{\pm \infty} b_{\epsilon} ( \check{r}^s_{\epsilon} , \check{\theta}^s_{\epsilon},\varrho,\vartheta) ds \right) - \int_0^{\pm \infty} f_{\epsilon}  ( \check{r}^s_{\epsilon} , \check{\theta}^s_{\epsilon},\varrho,\vartheta)  \exp \left( \int_0^s  b ( \check{r}^{s_1}_{\epsilon} , \check{\theta}^{s_1}_{\epsilon} , \varrho,\vartheta) ds_1 \right) ds . $$
This solution is still defined on ${\mathcal T}^{\pm}_{\epsilon} (R,V,I,\varepsilon)$ and, if $ C = 0 $,  it belongs to $ S_{-\mu} $.
}
\end{enumerate}
\end{prop}

\noindent {\bf Remark.} In the asymptotically Euclidean case with global coordinates, the usual construction of the Isozaki-Kitada phase  shows that  $ b_{\epsilon} $  is a short range symbol, which implies easily its integrability in time when evaluated  along a trajectory. Here, it only follows from the asymptotics of Theorem \ref{eikonal-equation} that
$$ b_{\epsilon} = S_{-1-\nu} + S_{-1} (\theta - \vartheta) $$
which in general fails to be short range because of the second term. However, when evaluated along a trajectory, we will  recover the integrability in time  (\ref{passhortrange}) by exploiting the decay in time of $  \check{\theta}^s_{\epsilon} - \vartheta $ (see (\ref{integrabilityintimedetheta})).

\bigskip

  We will prove Theorem \ref{eikonal-equation}, and Proposition \ref{prop-transport} likewise, only in the case $ \epsilon =1  $. Indeed, by Lemma \ref{lemmedebonrescaling},  if we define $ v (r,\theta) := (v^{jk}(r,\theta)) $ by
 \begin{eqnarray}
 g^{jk} (r,\theta) = \bar{g}^{jk}(\theta) + v^{jk} (r,\theta) ,  \label{definitionpotentielmetrique}
\end{eqnarray} 
 we have 
  $ p_{\epsilon,\kappa} =  \rho^2 + r^{-2} \bar{g}^{jk} (\theta) \eta_j \eta_k +r^{-2} v^{jk}(r/\epsilon,\theta) \eta_j \eta_k $ where $  v (r/\epsilon,\theta) $  is bounded in $ S^{-\nu} $ as $ \epsilon \in (0,1] $ (it is actually $ O (\epsilon^{\nu}) $).  The analysis below for $ \epsilon = 1 $ still applies uniformly for $ \epsilon \in (0,1] $, but only at the expense of heavier statments and notation\footnote{in the same spirit, since we don't need to use the distinction between $ \breve{r} $ and $r$ in this section; we use the simpler notation $ r $ though $ p_{\epsilon,\kappa} $ must be though as  a function of $ \breve{r} $}.  Thus, for simplicity, we will drop $ \epsilon $ and $ \kappa $  from the notation (except on $ V_{\kappa} $) everywhere below.
  
 We let   $ p = p_{\kappa,1}(r,\theta,\rho,\eta) $ and $ (\bar{r}^s , \bar{\vartheta}^s , \bar{\varrho}^s , \bar{\eta}^s) := \phi^s $ be the components of $ \phi^s (= \phi^s_{1,\kappa})$, namely the solution to
\begin{eqnarray}
 \dot{\bar{r}}^s = (\partial_{\rho} p) (\phi^s), \qquad  \dot{\bar{\vartheta}}^s = (\partial_{\eta} p) (\phi^s) , \qquad     \dot{\bar{\varrho}}^s = - (\partial_{r} p) (\phi^s) \qquad  \dot{\bar{\eta}}^s = - (\partial_{\theta} p) (\phi^s) , \label{pourcaracteristiques}
\end{eqnarray}
with initial condition
$$  (\bar{r}^s , \bar{\vartheta}^s , \bar{\varrho}^s , \bar{\eta}^s)_{|s=0} = (r,\theta,\rho,\eta) . $$
We will see it exists for $ \pm s \geq 0 $ on strongly outgoing (+)/ incoming (-) areas defined, for $ V \subset V_{\kappa} $, $ R > R_{\mathcal M} $ and $ 0 <\varepsilon <1 $, by
$$ \widetilde{\Gamma}^{\pm}_{\rm st} (R,V,\varepsilon) = \big\{  (r,\theta,\rho,\eta) \ | \ r > R , \ \theta \in V, \ \pm \rho > (1-\varepsilon^{2}) p^{1/2}  \big\} $$
where $ p = p (r,\theta,\rho,\eta) $. Note that the square on $ \varepsilon $ ensures that the condition $ \pm \rho > (1-\varepsilon^{2}) p^{1/2} $ is equivalent to $ |\eta|/r \lesssim \varepsilon $ and $ \pm \rho > 0 $. 

These sets are conical  ({\it i.e.} invariant under $ (\rho,\eta) \mapsto ( \lambda \rho, \lambda \eta ) $ for any $ \lambda > 0 $) and symmetric w.r.t. eachother, {\it i.e.} 
$$ (r,\theta,\rho,\eta) \in  \widetilde{\Gamma}^{+}_{\rm st} (R,V,\varepsilon) \qquad  \Longleftrightarrow \qquad   (r,\theta,- \rho,-\eta) \in  \widetilde{\Gamma}^{-}_{\rm st} (R,V,\varepsilon) .  $$
This symmetry together with the  property that,  for any $ \lambda \in \Ra $ and as long as the flow exists, 
\begin{eqnarray}
\begin{gathered}
  \big(\bar{r}^s , \bar{\vartheta}^s\big) (r,\theta, \lambda \rho, \lambda \eta)  =   \big(\bar{r}^{\lambda s} , \bar{\vartheta}^{\lambda s} \big) (r,\theta,\rho,\eta)  , \\
   \big(\bar{\varrho}^s , \bar{\eta}^s \big) (r,\theta, \lambda \rho, \lambda \eta)  =  \lambda  \big(  \bar{\varrho}^{\lambda s} ,  \bar{\eta}^{\lambda s} \big) (r,\theta,\rho,\eta) ,
   \end{gathered}
   \label{symmetriehomogeneite}
\end{eqnarray}  
will allow us to restrict the analysis to strongly outgoing regions and times $ s \geq 0 $.
  The same homogeneity properties hold for  $ \phi_0^s $ which in turns implies they also hold for $ F^{\pm} $.

The reason for denoting  the angular position by $ \bar{\vartheta}^s$ rather than $ \bar{\theta}^s $ and the radial momentum by $ \bar{\varrho}^s $
rather than $ \bar{\rho}^s $ is the following one. Let introduce 
\begin{eqnarray}
 \big( \bar{r},\bar{\vartheta} , \bar{\varrho} , \bar{\eta} \big) := \lim_{s \rightarrow + \infty}   \big( \bar{r}^s - 2 s \bar{\varrho}^s,\bar{\vartheta}^s , \bar{\varrho}^s , \bar{\eta}^s \big) = \lim_{s \rightarrow + \infty} \phi_0^{-s} \circ \phi^{s}(r,\theta,\rho,\eta),  \label{ref-class-scat}
\end{eqnarray}
which will be shown to exist for $ (r,\theta,\rho,\eta) $ in a strongly outgoing area $  \widetilde{\Gamma}^{+}_{\rm st}  $ (the parameters of which we omit here).
The item 2 of Theorem \ref{eikonal-equation} means that $ \varphi $ is a generating function of the Lagragian submanifold
\begin{eqnarray}
 \Lambda^{+} = \big\{ \big( (r,\theta,\rho,\eta), (   \bar{r},\bar{\vartheta} , \bar{\varrho} , \bar{\eta} \big) \big)  \ | \ (r,\theta,\rho,\eta) \in \widetilde{\Gamma}^{+}_{\rm st}  \big\} ,  \label{lagrangienne}
\end{eqnarray} 
   {\it i.e.}  the graph of the symplectic map $ F^+ $. The existence of $ \varphi $ rests on the fact that $ \Lambda^{+} $ can be parametrized by $ (r,\theta) $, the initial positions, and by  $ (\varrho,\vartheta ) $, the final  radial momentum and angular position. In particular, it is crucial to distinguish between the variables $ \theta $ and $ \vartheta $ which motivates our choice of notation.



Before starting the proof of Theorem \ref{eikonal-equation} which will come after several preparatory results, we introduce one more notation, for $ I \Subset (0,\infty) $, 
\begin{eqnarray}
  \widetilde{\Gamma}^{\pm}_{\rm st} (R,V,I,\varepsilon) = \big\{  (r,\theta,\rho,\eta) \in \widetilde{\Gamma}_{\rm st}^{\pm} (R,V,\varepsilon)   \ | \ \  p(r,\theta,\rho,\eta) \in I   \big\}  . \label{referencestrongarea}
\end{eqnarray}  
It allows to localize flow estimates in the energy shell $ p^{-1}(I) $, without loss of generality by the above homogeneity properties. Occasionally, we will also use $  \Gamma^{+}_{\rm st} (R,V,I,\varepsilon) $ defined by
\begin{eqnarray}
 (r,\theta,\rho,\xi) \in \Gamma^{+}_{\rm st} (R,V,I,\varepsilon)  \qquad \Longleftrightarrow \qquad   (r,\theta,\rho,r \xi) \in \widetilde{\Gamma}^{+}_{\rm st} (R,V,I,\varepsilon) .  
\label{refpourdiffeo}
\end{eqnarray}
Note that $ \rho^2 + g^{jk}(r,\theta) \xi_j \xi_k \in I $ on $\Gamma^{+}_{\rm st} (R,V,I,\varepsilon)$ so 
 $  \rho , \xi $ (and $ \theta $) are bounded there. In particular, all symbolic estimates on functions defined on 
 $\Gamma^{+}_{\rm st} (R,V,I,\varepsilon)$ will be only with respect to $r$. 

\medskip

To start the proof we recall a result from  \cite{MizutaniCPDE}.

\begin{prop}[Long time geodesic flow estimates] \label{blackboxdynamic} Let $ V_0 \Subset V_{\kappa} $ and $ I_0 \subset (0,\infty) $. One can choose $ R_0 \gg 1 $ large enough and $  0 <  \varepsilon_0 \ll 1 $ such that
\begin{enumerate}
\item{ for  all $ (r,\theta,\rho,\eta) \in  \widetilde{\Gamma}^{+}_{\rm st} (R_0,V_0,I_0, \varepsilon_0)  $, $ \phi^s (r,\theta,\rho,\eta) $ is defined for all $  s \geq 0 $ and
$$ (\bar{r}^s , \bar{\vartheta}^s ) \in (R_0,\infty) \times V_{\kappa} , \qquad  s \geq 0 , $$}
\item{for all $ ( j , \alpha , k , \beta ) \in {\mathbb Z}_+^{2n} $, there exists $ C > 0 $ such that for all  $ (r,\theta,\rho,\eta) \in \widetilde{\Gamma}^{+}_{\rm st} (R_0,V_0,I_0, \varepsilon_0)  $and all $ s > 0 $, 
 $$ \left| \partial_r^j \partial_{\theta}^{\alpha} \partial_{\rho}^k \partial_{\eta}^{\beta} \left( \frac{\bar{r}^s - r - 2 s \rho}{s}, \bar{\vartheta}^s  , \bar{\varrho}^s  , \frac{\bar{\eta}^s }{r} \right) \right| \leq C r^{-j-|\beta|} . $$
Moreover, there exists $ C >0  $  such that
\begin{eqnarray}
  (r + s ) / C \leq \bar{r}^s \leq C (r + s) ,  \label{Mizdecay}
\end{eqnarray}
for $  s \geq 0 $  and $ (r,\theta,\rho,\eta) \in   \widetilde{\Gamma}^{+}_{\rm st} (R_0,V_0,I_0,\varepsilon_0)  $.}
\end{enumerate}
\end{prop}


In the rest of the section, we choose $ R_0,V_0,I_0 $ and $ \varepsilon_0 $  as in Proposition \ref{blackboxdynamic}.

\bigskip

To study (\ref{ref-class-scat}) it will be convenient to use asymptotics in suitable symbol classes, in the spirit of those of Definition \ref{defi-asymp-symbol}. Given functions $a$ and $ b $ on $ \widetilde{\Gamma}^{+}_{\rm st} (R_0,V_0,I_0,\varepsilon_0) $, a real number $ \mu $ and an integer $ m $,  we define
$$ a = b + \widetilde{S}_{\mu} (\eta/r)^m  \ \ \ \stackrel{\rm def}{\Longleftrightarrow}  \ \ \ a - b = \sum_{|\gamma| = m} c_{\gamma}\left( r , \theta , \rho , \frac{\eta}{r} \right) \frac{\eta^{\gamma}}{r^{|\gamma|}}  \ \ \mbox{with} \ c_{\gamma} \in S^{\mu} , $$
where $ S^{\mu} =  S^{\mu} (\Gamma^{+}_{\rm st} (R_0,V_0,I_0,\varepsilon_0)) $.
 The relation
$  a = b + \widetilde{S}_{\mu_1} (\eta/r)^{m_1}  + \widetilde{S}_{\mu_2} (\eta/r)^{m_2}   $, with $ m_1, m_2 $  integers and $ \mu_1, \mu_2 \in \Ra $, is defined analogously.

   It is useful to record the following characterization of symbols of the form $ c (r,\theta,\rho, \eta / r) $.
\begin{lemm} \label{lemmecharacterization}  A function $ a : \widetilde{\Gamma}_{\rm st}^{+} (R_0,V_0,I_0,\varepsilon_0) \rightarrow \Ca $ is of the form
$$ a (r,\theta,\rho,\eta) = c \left( r , \theta , \rho , \frac{\eta}{r} \right) $$
for some $ c $ in  $ S^{\mu} \big(\Gamma^{+}_{\rm st} (R_0,V_0,I_0,\varepsilon_0) \big) $ if and only if, for all $ (j,\alpha,k,\beta) $,
\begin{eqnarray}
 \left| \partial_r^j \partial_{\theta}^{\alpha} \partial_{\rho}^k \partial_{\eta}^{\beta} a (r,\theta,\rho,\eta) \right| \leq C_{j \alpha k \beta} r^{\mu - j - |\beta|} , \label{condition-tilde}
\end{eqnarray}
for all $ (r,\theta,\rho,\eta) \in \widetilde{\Gamma}_{\rm st}^{+} (R_0,V_0,I_0,\varepsilon_0) $.
 In particular, if $a$ satisfies (\ref{condition-tilde}), then
\begin{eqnarray}
 a (r,\theta,\rho,\eta) = a (r,\theta,\rho,0) + (r \nabla_{\eta} a )(r,\theta,\rho,0) \cdot \frac{\eta}{r} + \widetilde{S}_{\mu} (\eta/r)^2 .  \label{symbolTaylor}
\end{eqnarray}
\end{lemm}

\noindent {\it Proof.} Follows from routine computations by considering $ c (r,\theta,\rho,\xi) : = a (r,\theta,\rho,r \xi) $.  \finpreuve

\bigskip


\begin{prop}[Asymptotics for $ F^+ $] \label{prop-inter-flow} For all  $ (r,\theta,\rho,\eta) \in \widetilde{\Gamma}^{+}_{\rm st} (R_0,V_0,I_0, \varepsilon_0)$, the limit (\ref{ref-class-scat}) exists. Furthermore, we have the expansions
\begin{eqnarray}
 \bar{r}  & = & r + \widetilde{S}_1 (\eta/r)^2  \label{proof1}  \\
  \bar{\varrho} &= & \rho + \widetilde{S}_0 (\eta/r)^2 \label{proof01} \\
 \bar{\eta}  & = & \eta +  \widetilde{S}_1 (\eta/r)^2 \label{proof02}  
\end{eqnarray}
and
\begin{eqnarray}
 \bar{\vartheta} & = & \theta +  \bar{g} (\theta)^{-1} \frac{\eta}{r \rho} + \widetilde{S}_{-\nu} (\eta/r) + \widetilde{S}_0 (\eta/r)^2 .
\label{proof2}
\end{eqnarray}
\end{prop}

Notice that   $  \rho $ is positive on $ \widetilde{\Gamma}^{+}_{\rm st} (R_0,V_0,I_0,\varepsilon_0) $ so the second term is the right hand side of (\ref{proof2}) is well defined. To prove this proposition, we will use the easily verified fact that for $ s \geq 0 $ and $ (r,\theta,\rho,0) \in \widetilde{\Gamma}^{+}_{\rm st}  (R_0,V_0,I_0,\varepsilon_0)  $, we have
\begin{eqnarray}
(\bar{r}^s , \bar{\vartheta}^s , \bar{\varrho}^s , \bar{\eta}^s)_{|\eta = 0} & = & (r + 2 s \rho , \theta,\rho,0), \label{Taylorexplicite0} \\
\big( \partial_{\eta} \bar{r}^s , \partial_{\eta} \bar{\vartheta}^s , \partial_{\eta} \bar{\varrho}^s , \partial_{\eta} \bar{\eta}^s \big)_{|\eta = 0} & = & \left(0 , \frac{2 \bar{g} (\theta)^{-1} + 2 v (r + 2 s \rho, \theta) }{(r + 2 s \rho )^2}  , 0 , I_{n-1} \right) , \label{Taylorexplicite1}
\end{eqnarray}
where we recall that $v$ is defined in (\ref{definitionpotentielmetrique}).
 We will also need the following lemma.
\begin{lemm} \label{lemmeaussicaracteristiques} For all $ ( j , \alpha , k , \beta ) \in {\mathbb Z}_+^{2n} $, setting $ \partial^{\gamma} =  \partial_r^j \partial_{\theta}^{\alpha} \partial_{\rho}^k \partial_{\eta}^{\beta} $, there is $ C > 0 $ such that,
\begin{eqnarray}
\left| \frac{\partial^{\gamma}  \bar{r}^s }{\bar{r}^s} \right| & \leq & C  r^{-j-|\beta|} , \label{premiereestimeelemme} \\
\left|  \partial^{\gamma} 
(r/\bar{r}^s)
\right| & \leq &  C (1 + |s/r|)^{-1} r^{-j-|\beta|} , \label{deuxiemeestimeelemme}
\end{eqnarray}
 for all  $ (r,\theta,\rho,\eta) \in \widetilde{\Gamma}_{\rm st}^{+} (R_0,V_0,I_0,\varepsilon_0)  $ and $ s \geq 0 $.
If furthermore $b \in S^{\mu} \big( (R_0,\infty) \times V_{\kappa} \big) $, then
\begin{eqnarray}
  \big|  \partial^{\gamma}  \big( b (\bar{r}^s,\bar{\vartheta}^s) \big) \big| \leq C (\bar{r}^s / r )^{\mu} r^{\mu-j-\beta} ,  \label{troisiemeestimeelemme}
\end{eqnarray}
with a constant bounded as long as $b $ varies in a bounded set.
\end{lemm}

\noindent {\it Proof.} The estimate (\ref{premiereestimeelemme}) is a simple consequence of the item 2 of Proposition \ref{blackboxdynamic} and the fact that $  \partial^{\gamma}   (r + 2 s \rho) = O (r + s)r^{-j-|\beta|}  $.  Next, by  observing that
$$ \partial^{\gamma} (r/\bar{r}^s) = \mbox{linear comb. of} \ \ \frac{ \partial^{\gamma^1}r}{\bar{r}^s} \frac{ \partial^{\gamma^2} \bar{r}^s}{\bar{r}^s} \cdots \frac{\partial^{\gamma^N}\bar{r}^s}{\bar{r}^s} \ \ \mbox{with} \ \ \gamma^1 + \cdots + \gamma^N = \gamma, \ \ N \leq |\gamma| + 1 , $$
we see that (\ref{deuxiemeestimeelemme}) follows from (\ref{Mizdecay}) and (\ref{premiereestimeelemme}).
Finally, the estimate (\ref{troisiemeestimeelemme}) follows from the item 2 of Proposition \ref{blackboxdynamic} for $ \bar{\vartheta}^s $, (\ref{deuxiemeestimeelemme}) and the fact that $  \partial^{\gamma} \big( b (\bar{r}^s,\bar{\vartheta}^s) \big) $ is a linear combination of
$$ ( \partial_r^{\tilde{j}} \partial_{\vartheta}^{\tilde{\alpha}} b \big) (\bar{r}^s,\bar{\vartheta}^s)   \partial^{\gamma^1_1}\bar{r}^s \cdots  \partial^{\gamma^1_{\tilde{j}}} \bar{r}^s \cdots \partial^{\gamma^{n}_1}\bar{\vartheta}^s_{n-1} \cdots  \partial^{\gamma^{n}_{\tilde{\alpha}_{n-1}}}\bar{\vartheta}^s_{n-1}
$$ with $ \gamma_1^1 + \cdots + \gamma^1_{\tilde{j}} + \cdots + \gamma_{1}^{n} + \cdots + \gamma^n_{\tilde{\alpha}_{n-1}} = \gamma $. \finpreuve

\bigskip

\noindent {\bf Proof of Proposition \ref{prop-inter-flow}.} We give the proofs of (\ref{proof1}) and (\ref{proof2}), the ones of (\ref{proof01}) and (\ref{proof02}) being similar (and slightly simpler).  We start with (\ref{proof2}). Writing $ \bar{\vartheta}^T = \theta + \int_0^T \dot{\bar{\vartheta}}^s ds $ and letting $ T \rightarrow + \infty $,  we obtain
$$ \bar{\vartheta} = \theta +2 \int_0^{+ \infty} \left( \bar{g} \big(\bar{\vartheta}^s \big)^{-1} + v \big(\bar{r}^{s} , \bar{\vartheta}^s \big) \right) \frac{\bar{\eta}^s}{(\bar{r}^s)^2} ds $$
where the integral is convergent since, for fixed $ (r,\theta,\rho,\eta) $, $ \bar{\eta}^s $ is bounded while $ \bar{r}^s \gtrsim r + s $ by (\ref{Mizdecay}).
Then, by using (\ref{deuxiemeestimeelemme}), (\ref{troisiemeestimeelemme})  and the item 2 of Proposition \ref{blackboxdynamic} for $ \bar{\eta}^s / r $, we see that
\begin{eqnarray}
  \partial_r^j \partial_{\theta}^{\alpha} \partial_{\rho}^k \partial_{\eta}^{\beta}  \left( \left( \bar{g} \big(\bar{\vartheta}^s \big)^{-1} + v \big(\bar{r}^{s} , \bar{\vartheta}^s \big) \right) \frac{\bar{\eta}^s}{(\bar{r}^s)^2}  \right) = O\big((1+|s/r|)^{-2}r^{-1-j-|\beta|} \big) .  \label{pourintegrand1}
\end{eqnarray}
Integrating this estimate in $s$ and using the characterization of Lemma \ref{lemmecharacterization}, we find  $ \bar{\vartheta} = \theta + \widetilde{S}_0 $ (see after Proposition \ref{blackboxdynamic} for this notation). Using (\ref{symbolTaylor}) together with (\ref{Taylorexplicite0}) and (\ref{Taylorexplicite1}) we get the improved expansion (\ref{proof2}) since
$$ 2 r \int_0^{+ \infty} \bar{g} ( \theta)^{-1} \frac{I_{n-1}}{(r+2s\rho)^2} ds = \frac{\bar{g}(\theta)^{-1}}{\rho}  , \qquad r \int_0^{+ \infty} v (r+ 2 s \rho , \theta) \frac{I_{n-1}}{(r+2s\rho)^2} ds = \widetilde{S}_{-\nu} . $$
 We next prove (\ref{proof1}). We start by writing  $ \bar{r}^T = r + \int_0^T 2 \bar{\varrho}^s ds $ and
\begin{eqnarray}
 \bar{\varrho}^s = \bar{\varrho}^T - \int_s^T \left(\frac{2}{\bar{r}^u} \big( \bar{g}^{lm} (\bar{\vartheta}^u) + v^{lm} (\bar{r}^u , \bar{\vartheta}^u) \big) - (\partial_r v^{lm}) (\bar{r}^u, \bar{\vartheta}^u)   \right) \frac{\bar{\eta}_l^u \bar{\eta}_m^u}{(\bar{r}^u )^2} du . \label{pourintegrand2}
\end{eqnarray}
Similarly to (\ref{pourintegrand1}), the $ \partial_r^j \partial_{\theta}^{\alpha} \partial_{\rho}^k \partial_{\eta}^{\beta}   $ derivative of the integrand in (\ref{pourintegrand2}) is  $O \big( (1+|u/r|)^{-3} r^{-1-j-|\beta|} \big) $. This implies on one hand that the limit of $ \bar{r}^T  - 2 T \bar{\varrho}^T $ exists as $ T \rightarrow + \infty $ and equals
$$ \bar{r} =  r - 2 \int_0^{+ \infty} \left( \int_s^{+ \infty}  \left(\frac{2}{\bar{r}^u} \big( \bar{g}^{lm} (\bar{\vartheta}^u) + v^{lm} (\bar{r}^u , \bar{\vartheta}^u) \big) - (\partial_r v^{lm}) (\bar{r}^u, \bar{\vartheta}^u)   \right) \frac{\bar{\eta}_l^u \bar{\eta}_m^u}{(\bar{r}^u )^2}  d u  \right) d s  $$
and on the other hand that, for any $ (j,\alpha,k,\beta) $,  the   $ \partial_r^j \partial_{\theta}^{\alpha} \partial_{\rho}^k \partial_{\eta}^{\beta}   $ derivative  of the above double integral is $ O (r^{1-j-\beta}) $. This gives  the rough bound $ \bar{r} = r + \widetilde{S}_1 $ which then improves to  (\ref{proof1}) by using the above expression together with (\ref{symbolTaylor}), (\ref{Taylorexplicite0}) and (\ref{Taylorexplicite1}). \finpreuve

\bigskip

The last intermediate result needed to prove Theorem \ref{eikonal-equation} is the following one.

\begin{prop}[projecting the Lagragian] \label{projLagrangian} Let $ I_1 \Subset I_ 0  $ and $ V_1 \Subset V_0  $ with $ V_1 $ convex. Then
 one can find $ R_1 \gg 1 $ and $ C > 1 $ such that for all $ \varepsilon \ll 1 $,
 the map
\begin{eqnarray}
 (r,\theta,\rho,\eta) \mapsto (r,\theta,\bar{\varrho},\bar{\vartheta})  \label{unnameddiffeo}
\end{eqnarray}
is a diffeomorphism from $ \widetilde{\Gamma}^{+}_{\rm s} (R_1 , V_1 , I_0 , C \varepsilon  ) $ onto an open subset containing $ \Theta^{+} (R_1,V_1,I_1,\varepsilon) .$
On $ \Theta^{+} (R_1,V_1,I_1,\varepsilon)  $, the inverse of (\ref{unnameddiffeo}) is of the form
$$ (r,\theta,\varrho,\vartheta) \mapsto \left( r , \theta , \underline{\rho}  , \underline{\eta}  \right)  $$
with
\begin{eqnarray}
 \ \underline{\rho}(r,\theta,\varrho,\vartheta) & = & \varrho + S_{-\nu} (\vartheta - \theta)   + S_0 (\vartheta - \theta)^2  ,
\label{expansionderivee1} \\
 \underline{\eta}(r,\theta,\varrho,\vartheta)  & = & r   \varrho \bar{g} (\theta) ( \vartheta - \theta ) +  S_{1-\nu}  (\vartheta - \theta) +  S_1 (\vartheta - \theta)^2   . \label{expansionderivee2}
\end{eqnarray}
\end{prop}

Recall that the notation $ \Theta^{+} (R,V,I,\varepsilon) $ is defined in (\ref{Thetaenergie}). To understand informally why we can take proportional  parameters $ C \varepsilon $ and $ \varepsilon $, we recall that the condition $ \rho > (1 - C \varepsilon^2) p^{1/2} $ means that $ |\eta | / r \lesssim  \varepsilon   $ (and $ \rho > 0 $) which, by (\ref{expansionderivee2}), is comparable to the condition $ |\theta - \vartheta| \lesssim \varepsilon $.

\bigskip

\noindent {\it Proof of Proposition \ref{projLagrangian}.} Denote by $ \widetilde{H} $ the map (\ref{unnameddiffeo}). Consider the maps $ H $ and $ K $ defined by
$$ H (r,\theta,\rho,\xi) = \big( r,\theta,\rho , \theta + \rho^{-1} \bar{g} (\theta)^{-1} \xi \big) , \qquad K (r,\theta,\varrho, \vartheta) = \big(r,\theta,\varrho, \varrho \bar{g} (\theta)(\vartheta-\theta) \big) , $$
which are inverse to eachother (on appropriate domains  given below). We also set
$$ E (r,\theta,\rho,\xi) = ( r , \theta , \rho , r \xi ) . $$
    It follows from (\ref{proof01}) and (\ref{proof2}) that
$$ \widetilde{H} (r,\theta,\rho,\eta) = H (r,\theta,\rho,\eta/r) + \widetilde{S}_{-\nu} (\eta/r) + \widetilde{S}_0 (\eta/r)^2 . $$
thus, after composition with $ E \circ K $ and using Lemma \ref{lemmecharacterization}, we see that
\begin{eqnarray}
 \widetilde{H} \circ E \circ K = I + S_{-\nu} (\vartheta-\theta) + S_0 (\vartheta-\theta)^2  .  \label{domainediffeopreuve}
\end{eqnarray}
These computations make sense on the following sets. Since $ \widetilde{H} $ is defined on $ \widetilde{\Gamma}_{\rm st}^{+}(R_0,V_0,I_0,\varepsilon_0) $,  it follows from (\ref{refpourdiffeo}) that (\ref{domainediffeopreuve}) holds on any set which is mapped into $  \Gamma_{\rm s}^{+}(R_0,V_0,I_0,\varepsilon_0)  $ by $ K $. Using (\ref{hjkellipticity}) and the fact that $ I_0 $ is relatively compact,  one can find $ C > 1 $ such that
$$  K \left(   \Theta^{+} (R_0,V_0,I_0, \varepsilon) \right) \subset   \Gamma^{+}_{\rm st} \left(R_0,V_0,I_0 , C \varepsilon  \right),  $$
and thus (\ref{domainediffeopreuve}) holds on $   \Theta^{+} (R_0,V_0,I_0,\varepsilon) $ if $ C \varepsilon < \varepsilon_0 $. Since the right hand side of (\ref{domainediffeopreuve}) is a small perturbation of identity where $ r$ is large and $ \theta - \vartheta $ is small, it follows from a routine argument that if $ R $ is large enough and $ \varepsilon $ is small enough, it is a diffeomorphism on $ \Theta^{\pm} (R,V_1,I_0,\varepsilon) $ onto an open set containing $  \Theta^{\pm} (R,V_1,I_1,\varepsilon/8)  $.  Note that $ \Theta^{+} (R,V_1,I_0,\varepsilon) $ is convex which is useful to justify this fact, for instance to prove the injectivity of (\ref{domainediffeopreuve}) by using the mean value theorem. Note also that $ (r,\theta) $ is unchanged by the left hand side of (\ref{domainediffeopreuve}) and that, when $ \vartheta = \theta $, we have
$  (\widetilde{H} \circ E \circ K  )  (r,\theta,\varrho,\theta) =  (r,\theta,\varrho,\theta)  $. This allows to check that the inverse mapping to (\ref{domainediffeopreuve}) is still of the form $ I + S_{-\nu} (\vartheta-\theta) + S_0 (\vartheta-\theta)^2 $. Composing $ E \circ K  $ with this inverse diffeomorphism, we get the existence of $ \big( \underline{\rho},\underline{\eta} \big) $  and the expansions (\ref{expansionderivee1})-(\ref{expansionderivee2}). \finpreuve


\bigskip

\noindent {\bf Proof of Theorem \ref{eikonal-equation}.}   We  choose $  I_1 = I   $ and $  V_1 = V $ in Proposition \ref{projLagrangian} (recall that $ I_0 $ and $ V_0 $ were chosen arbitrarily). We  prove the items 1 and 2 at the same time.
By  Proposition \ref{prop-inter-flow}, $  F^{+}  $ is well defined on $  \widetilde{\Gamma}_{\rm st}^{+}(R_0,V_0,I_0, \varepsilon_0) $.
Since $ \phi^s  $ and $ \phi_0^{-s} $ are symplectic maps so is $ F^{+} $ and its graph  is Lagrangian. Together with Proposition \ref{projLagrangian}, this implies that the differential form
$$ \underline{\rho} (r,\theta,\varrho,\vartheta) d r + \underline{\eta}(r,\theta,\varrho,\vartheta) d \theta+ \bar{r}(r,\theta,\underline{\rho},\underline{\eta}) d \varrho -  \bar{\eta} (r,\theta,\underline{\rho},\underline{\eta}) d \vartheta $$
is closed on $ \Theta (R,V,I,\varepsilon)  $ for $ \varepsilon $ small enough. Since this set is convex, we  get the existence of a function $ \varphi $, unique up to an additive constant, such that
\begin{eqnarray}
\partial_r \varphi = \underline{\rho} \qquad \partial_{\theta} \varphi = \underline{\eta} \qquad \partial_{\varrho} \varphi  = \bar{r} (r,\theta,\underline{\rho},\underline{\eta}) \qquad 
\partial_{\vartheta} \varphi = - \bar{\eta} (r,\theta,\underline{\rho},\underline{\eta}) . \label{formulesplusprecises}
\end{eqnarray}
To fix the constant and to define  $ \varphi $ globally in $ \varrho $, we observe that (\ref{symmetriehomogeneite}) (for $ \lambda > 0 $) implies that $ \underline{\rho} , \underline{\eta} $ are homogeneous of degree $ 1 $ in $ \varrho $ and $ \bar{r}(r,\theta,\underline{\rho},\underline{\eta}), \bar{\eta} (r,\theta,\underline{\rho},\underline{\eta}) $ of degree $0$. We can thus find a unique solution $ \varphi $ defined for $ (r,\theta,\vartheta) \in \Theta (R,V,\varepsilon) $ and $ \varrho > 0 $, which is
homogeneous of degree $1$ in $ \varrho $. Then  (\ref{formulesplusprecises}) and Proposition \ref{projLagrangian} yield the item 2.  It turns out that if one considers $ \varrho \varphi (r,\theta,1,\vartheta) $ with $ \varrho \in \Ra $, we get the expected solution for it is also a generating function of $ F^- $ for $ \varrho < 0 $ by the symmetry (\ref{symmetriehomogeneite}) for $ \lambda = -1 $.
To prove that $ \varphi $ satisfies the eikonal equation it suffices to observe that
$$ p (r,\theta,\rho,\eta) = \bar{\varrho} (r,\theta,\rho,\eta)^2 , $$
which is well known (see {\it e.g.} \cite{MizutaniCPDE}) and easy to get from Proposition \ref{blackboxdynamic} and the conservation of energy.
By evaluating this equality on $ (r,\theta,\underline{\rho},\underline{\eta}) $, we get (\ref{eikonaladeriver}). For the item 3, (\ref{expansioneikonal1}) and (\ref{expansioneikonal2}) are direct consequences of Proposition \ref{projLagrangian} by (\ref{formulesplusprecises}). The expansions (\ref{expansioneikonal3}) and (\ref{expansioneikonal4}) follow from (\ref{formulesplusprecises}) combined with (\ref{expansionderivee1})-(\ref{expansionderivee2}) and Proposition \ref{prop-inter-flow}.
\finpreuve

\bigskip

We end up this section with the proof of Proposition \ref{prop-transport} on transport equations (recall that $ \epsilon , \kappa $ have been dropped from the notation). As before, we only consider the case when $ s \geq 0 $ (and $ \varrho > 0 $).

\bigskip

\noindent {\bf Proof of Proposition \ref{prop-transport}.} The item 1 follows from the well known method of characteristics (see {\it e.g.} \cite{Evans}) and has nothing to do with our specific geometric context so we only give the main lines.   We let $ (\tilde{r}^s , \tilde{\theta}^s ) $ be the maximal solution to the ODE
$$ \big( \dot{\tilde{r}}^s , \dot{\tilde{\theta}}^s \big) = (\partial_{\rho,\eta} p)  \big( \tilde{r}^s , \tilde{\theta}^s , \partial_{r,\theta} \varphi (\tilde{r}^s,\tilde{\theta}^s,\varrho,\vartheta) \big), \qquad \big( \tilde{r}^0 , \tilde{\theta}^0 \big)  = (r,\theta) . $$
We also let $ (\tilde{\rho}^s , \tilde{\eta}^s ) = \partial_{r,\theta} \varphi (\tilde{r}^s,\tilde{\theta}^s,\varrho,\vartheta) $.
By differentiating (\ref{eikonaladeriver}) in $ (r,\theta) $, one has
$$ (\partial_{r,\theta} p ) (r,\theta , \partial_{r,\theta} \varphi ) + \big( D^2_{r,\theta} \varphi  \big) (\partial_{\rho,\eta} p )(r,\theta, \partial_{r,\theta} \varphi ) = 0 $$
where $  D^2_{r,\theta} \varphi $ is the Hessian matrix of $ \varphi $ (seen as a function of $ (r, \theta) $). By evaluating this identity at $  (\tilde{r}^s , \tilde{\theta}^s , \varrho , \vartheta)  $, we obtain
$$  \big( \dot{\tilde{\rho}}^s , \dot{\tilde{\eta}}^s \big) = (\partial_{\rho,\eta} p)  \big( \tilde{r}^t , \tilde{\theta}^t , \tilde{\rho}^t , \tilde{\eta}^t \big)  $$
which, together with the first equation, shows that $ \big( \tilde{r}^s , \tilde{\theta}^s , \tilde{\rho}^s , \tilde{\eta}^s \big)  $ solves the equation (\ref{pourcaracteristiques}) with initial condition $(r,\theta,\partial_{r,\theta} \varphi)$. Thus $    \big( \check{r}^s , \check{\theta}^s , \check{\rho}^s , \check{\eta}^s \big) = \big( \tilde{r}^s , \tilde{\theta}^s , \tilde{\rho}^s , \tilde{\eta}^s \big)   $ satisfies the expected properties of the first item. To prove the second item, the main observation is that
$$ P \varphi  = S_{-1-\nu} + S_{-1} ( \vartheta - \theta ) , $$
which follows from (\ref{formuleLaplaciensansrescaling}), (\ref{expansioneikonal1}) and (\ref{expansioneikonal2}). 
Using (\ref{fonctioncaracteristiqueparametres}) (see also (\ref{ref-class-scat})), we have 
$$ \check{\theta}^s \rightarrow \bar{\vartheta} (r,\theta,\partial_{r,\theta} \varphi) = \vartheta , \qquad s \rightarrow + \infty . $$ Thus, by integrating $ \dot{\check{\theta}}^s $ from $s$ to $ + \infty $ and using the flow estimates of Proposition \ref{blackboxdynamic} and Lemma \ref{lemmeaussicaracteristiques} together with the estimates on $ \varphi $ given in the item 3 of Theorem \ref{eikonal-equation}, we get
\begin{eqnarray}
 | \partial_r^j \partial_{\theta}^{\alpha} \partial_{\varrho}^k \partial_{\vartheta}^{\beta} ( \check{\theta}^s - \vartheta ) | \lesssim \scal{s/r}^{-1} r^{-j} . 
\label{integrabilityintimedetheta}
\end{eqnarray}
By the same techniques we can estimate the derivatives of $ \check{r}^s $
and we get the result by routine calculations.  The third item follows from the usual method of characteristics for linear transport equations. We only record that to prove that the  solution is in $ S_{-\mu} $ (if $C = 0 $ and $ f \in S_{-1- \mu} $),  it suffices to observe that
$$ \left| \partial_r^j \partial_{\theta}^{\alpha} \partial_{\varrho}^k \partial_{\vartheta}^{\beta} \left(  f (\check{r}^s , \check{\theta}^s , \varrho, \vartheta) \right) \right| \lesssim \scal{s/r}^{-1-\mu} r^{-1-\mu-j}  , $$
on ${\mathcal T}^+ (R,V,I,\varepsilon) $.
 \finpreuve



\section{The Isozaki-Kitada parametrix} \label{sectionIsozakiKitada}
\setcounter{equation}{0}
In this section, we construct a new version of the Isozaki-Kitada parametrix compared to the ones introduced in \cite{BoucletAPDE,MizutaniCPDE}. The novelty stems basically from the parametrization of the Lagragian (\ref{lagrangienne}) in term of the final angular position $ \bar{\vartheta} $ rather than the final angular momentum $ \bar{\eta} $; it turns out that it is more accurate to deal with global in time estimates. 

Before displaying the parametrix, we need some notation and preliminary results for operators on $ \Ra^n $. For $ \mu \in \Ra $,  $ S_{\mu}(\Ra^{2n}) $ denotes the space of symbols defined on $ \Ra^{2n} $ such that
$$ |\partial_r^j \partial_{\theta}^{\alpha} \partial_{\varrho}^k \partial_{\eta}^{\beta} a (r,\theta,\varrho,\vartheta) | \leq C \scal{r}^{\mu-j}, \qquad \mbox{on} \ \Ra^{2n} . $$
We equip it with the standard topology. We will also need the space $ S_{\mu}^{\rm min} (\Ra^{3n})  $ of functions satisfying
$$  |\partial_r^j \partial_{\theta}^{\alpha} \partial_{r^{\prime}}^{j^{\prime}} \partial_{\theta^{\prime}}^{\alpha^{\prime}} \partial_{\varrho}^k \partial_{\eta}^{\beta} A (r,\theta,r^{\prime},\theta^{\prime},\varrho,\vartheta) | \leq C \scal{\min(r,r^{\prime})}^{\mu-j-j^{\prime}}, \qquad \mbox{on} \ \Ra^{3n} . $$
Let us consider first the semiclassical version of the operators.
For $ a \in  S_{\mu} (\Ra^{2n}) $  supported in $ \Theta^{\pm} (R,V,I,\varepsilon ) $ (see (\ref{Thetaenergie})), we define
$$  J^h (a) v (r,\theta)  = (2 \pi h)^{- \frac{n+1}{2}}  \int \! \int \! \int e^{\frac{i}{h} \big( \varphi_{1}(r,\theta,\varrho,\vartheta) - x \varrho  \big)}
a (r,\theta,\varrho,\vartheta)  v (x, \vartheta) d x  d \varrho d \vartheta , $$
where $ \varphi_1 $ is the phase constructed in Theorem \ref{eikonal-equation} with $ \epsilon = 1 $.
The operator $ J^h (a) $ is well defined on $  {\mathscr S}(\Ra^n) $ and it is not hard to check that it maps  $ {\mathscr S}(\Ra^n) $ into itself. Its formal adjoint (with respect to the Lebesgue measure) is given by
$$ J^h (a)^{\dag} u (x,\vartheta)  = (2 \pi h)^{- \frac{n+1}{2}}  \int \! \int \! \int e^{\frac{i}{h} \big( x \varrho - \varphi_{1}(r^{\prime},\theta^{\prime},\varrho,\vartheta) \big)}
\overline{a} (r^{\prime},\theta^{\prime},\varrho,\vartheta)  u (r^{\prime},  \theta^{\prime}) d \varrho  d r^{\prime} d \theta^{\prime} $$
and $ J^h (a)^{\dag} $ also maps the Schwartz space into itself. 
The prototype of our parametrix at {\bf high frequency} will be of the form
\begin{eqnarray}
 J^h (a)  e^{- i t D_x^2} J^h (b)^{\dag}  .
 \label{formelhautefrequence}
\end{eqnarray}
For the parametrix at {\bf low frequency}, 
we will rather consider operators of the form 
\begin{eqnarray}
{\mathscr D}_{\epsilon}  J_{\epsilon} (a_{\epsilon})  e^{- i \epsilon^2 t D_x^2} J_{\epsilon} (b_{\epsilon})^{\dag}  {\mathscr D}_{\epsilon}^{-1} \label{formelbassefrequence}
\end{eqnarray} 
 where $ J_{\epsilon} (a_{\epsilon}) $ is defined by
 $$   J_{\epsilon} (a_{\epsilon}) v (r,\theta)  = (2 \pi )^{- \frac{n+1}{2}}  \int \! \int \! \int e^{i\big( \varphi_{\epsilon}(r,\theta,\varrho,\vartheta) - x \varrho  \big)}
a_{\epsilon} (r,\theta,\varrho,\vartheta)  v (x, \vartheta) d x  d \varrho d \vartheta , $$
{\it i.e.} is defined as $ J^h $ with $h=1$ and $ \varphi_1 $ replaced by $ \varphi_{\epsilon} $. In this case, we need to consider  $ \epsilon $ dependent amplitudes $a_{\epsilon}, b_{\epsilon} $ which will be bounded in their classes with respect to $ \epsilon $ and supported in $ \epsilon $ independent areas of the form $ \Theta^{\pm} (R,V,I,\varepsilon) $.
Omitting  the scaling operators $ {\mathscr D}_{\epsilon} $ and $ {\mathscr D}_{\epsilon}^{-1} $ in (\ref{formelbassefrequence}), we can write the Schwartz kernels (with respect to the Lebesgue measure) of  both (\ref{formelhautefrequence}) and (\ref{formelbassefrequence}) under the following single form
\begin{eqnarray}
 (2 \pi h)^{-n} \int \int e^{\frac{i}{h} \left( \varphi_{\epsilon} (r,\theta,\varrho,\vartheta) - \frac{\epsilon^2 t}{h} \varrho^2  - \varphi_{\epsilon}(r^{\prime},\theta^{\prime},\varrho,\vartheta) \right)  } a_{\epsilon}(r,\theta,\varrho,\vartheta) 
\bar{b}_{\epsilon} (r^{\prime},\theta^{\prime},\varrho,\vartheta) d \varrho d \vartheta . 
\label{kernelexplicitdouble}
\end{eqnarray}
Indeed, (\ref{formelhautefrequence})  corresponds to $ \epsilon = 1 $ and $ h \in (0,1] $, while (\ref{formelbassefrequence}) corresponds to $h=1$ and $ \epsilon \in (0,1] $.
The form of this kernel motivates the introduction of oscillatory integrals of the form
\begin{eqnarray}
 I^h_{\epsilon} (A_{\epsilon,s}) = (2 \pi h)^{-n} \int \! \int e^{ \frac{i}{h} \Phi_{\epsilon} (s,r,\theta,r^{\prime},\theta^{\prime},\varrho,\vartheta)} A_{\epsilon,s} (r,\theta,r^{\prime},\theta^{\prime},\varrho,\vartheta) d \varrho d \vartheta , \label{integraleoscillante}
\end{eqnarray} 
where $ \Phi_{\epsilon} = \Phi_{\epsilon} (s,r,\theta,r^{\prime},\theta^{\prime},\varrho,\vartheta) $ is defined as
$$ \Phi_{\epsilon} :=  \varphi_{\epsilon}(r,\theta,\varrho,\vartheta) - s \varrho^2 - \varphi_{\epsilon} (r^{\prime},\theta^{\prime},\varrho,\vartheta) . $$
In the applications we will take either $ s = t/h $ or $ s = \epsilon^2 t $ (and $h=1$) to fit (\ref{kernelexplicitdouble}).  We will consider amplitudes $ A_{\epsilon,s} $ bounded in  $ S_0^{\min} (\Ra^{3n}) $ with respect to $ (\epsilon,s) $ and satisfying the support condition
\begin{eqnarray}
 \mbox{supp}(A_{\epsilon,s}) \subset  \widehat{\Theta}^{\pm} (R,V,R^{\prime},V^{\prime},I,\varepsilon,\varepsilon^{\prime}) \label{conditiondesupportdouble}
\end{eqnarray} 
where $  \widehat{\Theta}^{\pm} (R,V,R^{\prime},V^{\prime},I,\varepsilon,\varepsilon^{\prime}) $ is the set
$$ \{    (r,\theta,r^{\prime},\theta^{\prime},\varrho , \vartheta) \ |  \  (r,\theta,\vartheta) \in \Theta (R,V,\varepsilon), \ (r^{\prime},\theta^{\prime},\vartheta) \in \Theta (R^{\prime},V^{\prime},\varepsilon^{\prime}), \ \varrho^2 \in I, \ \pm \varrho > 0 \} .$$
We refer to (\ref{Thetadef}) for $ \Theta (R,V,\varepsilon) $ and, as in Theorem \ref{eikonal-equation}, we will assume that $V$ is convex. Note that the above amplitudes are compactly supported with respect to  $ (\varrho,\vartheta ) $. In the same spirit, to cover both definitions of $ J^h $ and $ J_{\epsilon} $ in the next section, we will use 
$$ J_{\epsilon}^h (a_{\epsilon}) v (r,\theta) =   (2 \pi h)^{- \frac{n+1}{2}}  \int \! \int \! \int e^{ \frac{i}{h} \big( \varphi_{\epsilon}(r,\theta,\varrho,\vartheta) - x \varrho  \big)}
a_{\epsilon} (r,\theta,\varrho,\vartheta)  v (x, \vartheta) d x  d \varrho d \vartheta , $$
where $a_{\epsilon}$ is allowed to depend on $ \epsilon $ in a bounded fashion.




\subsection{FIO estimates}

In this section, we record properties on operators $ J^h_{\epsilon} (a_{\epsilon}) $ and oscillatory integrals $ I^h_{\epsilon} (A_{\epsilon,s}) $.  All propositions and lemmas are stated in full generality; however, for notational simplicity only, we will prove them in the outgoing case ($ + $ case) and will omit the dependence on $ \epsilon $ in the notation of proofs, similarly to what we did in Section \ref{sectionscatteringclassique}.

\begin{prop}[Non stationary phase estimates] \label{nonstatprop} Let $ I  = (\varrho_{\rm inf}^2, \varrho_{\rm sup}^2 ) $ with $ \varrho_{\rm sup} > \varrho_{\rm inf} > 0  $. 
\begin{enumerate}
\item{ Let $ \delta \in (0,1) $. If $ \varepsilon $ and $ \varepsilon^{\prime} $ are small enough, then
$$ \left. \begin{array}{c} (1- \delta) r  \geq  (r^{\prime} \pm 2 s \varrho_{\rm sup}) \\ \qquad  \mbox{or} \\
r \leq (1 - \delta) (r^{\prime} \pm 2 s \varrho_{\rm inf}) 
 \end{array} \right\} \qquad \Longrightarrow \qquad  I^h_{\epsilon} (A_{\epsilon,s}) = O \big( h^{\infty}  \scal{ s ,r , r^{\prime} }^{ - \infty} \big) ,
 $$
uniformly in $ \epsilon $, provided that   $ \pm s \geq 0 $ and $ (A_{\epsilon,s}) $ belongs to a bounded set of $ S_0^{\min} (\Ra^{3n}) $ such that (\ref{conditiondesupportdouble}) holds.  }
 \item{Let $ c \in (0,1) $. If $ R  $ is large enough and  $ \varepsilon $ is small enough, then 
 $$ \varepsilon^{\prime} \leq \varepsilon^2   \qquad  \mbox{and} \qquad
  |\theta - \vartheta | \geq c \varepsilon \ \  \mbox{on } \ \emph{supp}(A_{\epsilon,s}) \qquad \Longrightarrow \qquad  I^h_{\epsilon} (A_{\epsilon,s}) = O \big( h^{\infty}  \scal{ s ,r , r^{\prime} }^{ - \infty} \big) 
 $$
uniformly in $ \epsilon $, provided that   $ \pm s \geq 0 $ and $ (A_{\epsilon,s}) $ belongs to a bounded set of $ S_0^{\min} (\Ra^{3n}) $ such that (\ref{conditiondesupportdouble}) holds.  }
 \end{enumerate}
\end{prop}

\noindent {\it Proof.}  For both items, we consider only the outgoing case. For the first one, using the expansion (\ref{expansioneikonal4}) we find that, on the support of the amplitude, 
$$  \partial_{\varrho} \Phi = - 2 s \varrho + r \big(1 + O (\varepsilon) \big) - r^{\prime} (1+O(\varepsilon^{\prime})) . $$
Therefore, if $ (1 - \delta) r \geq  r^{\prime} + 2 s \varrho_{\rm sup} $ we see that
$$ \frac{\partial_{\varrho} \Phi }{r}\geq   \delta - C \varepsilon - C \varepsilon^{\prime}  $$
where the right hand side is larger than $ \delta /2 $ if $ \varepsilon,\varepsilon^{\prime} $ are small enough. Then, repeated integrations by part in $ \varrho $ show that $ I^h (A_s) = O (h^{\infty} r^{-\infty}) $ which yields the result since $ r^{\prime} + |s| \lesssim r$  in this regime. On the other hand, if $ r \leq (1-\delta) (r^{\prime} + 2 s \varrho_{\rm inf}) $ then
$$ \frac{\partial_{\varrho} \Phi }{r^{\prime} + 2 s \varrho_{\rm inf}} \leq -   \delta  + C \varepsilon + C \varepsilon^{\prime} .  $$
Then,  as above, integrations by part in $ \varrho $ show that $ I^h (A_s) = O \big(h^{\infty} (r^{\prime} + |s|)^{-\infty} \big) $ which yields the result since $ r \lesssim r^{\prime} + |s| $. For the second item we observe first that by the item 1, we can assume that $ C^{-1} r \leq r^{\prime} + s \leq C r $. Then using the expansion (\ref{expansioneikonal3}), we have
$$ \partial_{\vartheta} \Phi = - r \varrho \big( \bar{g} (\theta) (\vartheta - \theta) + O (R^{-\nu} \varepsilon) + O (\varepsilon^2) \big) + r^{\prime} \varrho O (\varepsilon^2) $$
on the support of the amplitude. Thus,  if $  |\vartheta - \theta| \geq c \varepsilon $, we see that for $ R $ large enough and $ \varepsilon $ small enough,
$$ \frac{| \partial_{\vartheta} \Phi | }{ r } \gtrsim \varepsilon $$
since $ r^{\prime} / r $ is bounded thanks to the assumption $ r^{\prime} + s \leq C r $. Then, integrating by part in $ \vartheta $, we obtain  $  I^h (A_s) = O \big( h^{\infty}  r^{ - \infty} \big) $ which yields the full decay since we also assume that $ r \gtrsim r^{\prime} + s $.
 \finpreuve

\bigskip




We next state an Egorov type theorem. It is a classical result but we quote it explicitly for we are not in a completely standard situation and also consider $ \epsilon $ dependent phases and symbols.

\begin{prop}[Egorov theorem] \label{LemmeEgorovstat} We can choose $ R^{\prime} \gg 1 $ and $ 0 < \varepsilon^{\prime} \ll 1 $ such that for  all bounded families $ (a_{\epsilon}) \in S_{\mu}(\Ra^{2n}) , (b_{\epsilon}) \in S_{\mu^{\prime}} (\Ra^{2n}) $ such that
$$  \emph{supp}(a_{\epsilon}) \subset \Theta^{\pm} (R^{\prime},V,I,\varepsilon^{\prime}), \qquad  \emph{supp}(b_{\epsilon}) \subset \Theta^{\pm} (R^{\prime},V,I,\varepsilon^{\prime}) $$ 
one has
$$ J^h_{\epsilon}(a_{\epsilon}) J^h_{\epsilon}(b_{\epsilon})^{\dag} = O \! p^h (c_{\epsilon}(h)) , $$
for some admissible $ c_{\epsilon} (h) \in \widetilde{S}^{-\infty,\mu + \mu^{\prime} + 1 - n}(\Ra^{2n}) $ depending in a bounded fashion on $\epsilon$ and such that
$$ c_{\epsilon} (h) \sim \sum_{j \geq 0} h^j c_{\epsilon,j}, \qquad c_{\epsilon,0} =  a_{\epsilon}  (r,\theta,\bar{\varrho}_{\epsilon},\bar{\vartheta}_{\epsilon}) \bar{b}_{\epsilon}
(r,\theta,\bar{\varrho}_{\epsilon},\bar{\vartheta}_{\epsilon}) \big|{\rm det} \  d_{\rho,\eta} \big(  \bar{\varrho}_{\epsilon}, \bar{\vartheta}_{\epsilon} \big) \big| ,$$
where we recall that $ ( \bar{\varrho}_{\epsilon},\bar{\vartheta}_{\epsilon}) $ are  components of $ F_{\epsilon,\kappa}^{\pm} $ (see the item 2 of Theorem \ref{eikonal-equation}), namely
$$ ( \bar{\varrho}_{\epsilon},\bar{\vartheta}_{\epsilon}) := \lim_{s \rightarrow \pm \infty} \big( \bar{\varrho}^s_{\epsilon} , \bar{\vartheta}^{s}_{\epsilon} \big) , \qquad \mbox{where} \ \  \big( \bar{r}^s_{\epsilon} , \bar{\vartheta}^s_{\epsilon} , \bar{\varrho}^s_{\epsilon} , \bar{\eta}^s_{\epsilon} \big) =  \phi_{\epsilon,\kappa}^s (r,\theta,\rho,\eta) .  $$
For $j \geq 1$, $ c_{\epsilon,j} $ has its support contained in the support of $ c_{\epsilon,0} $.
\end{prop}



In several proofs below, the following definition will be useful
$$ A = S_{\mu}^{\min} \big( (\vartheta - \theta) + (\vartheta - \theta^{\prime}) \big)^k  \ \ \ \stackrel{\rm def}{\Longleftrightarrow} \ \ \ A = \sum_{|\alpha| + |\alpha^{\prime}| = k} A_{\alpha \alpha^{\prime}} (\vartheta - \theta)^{\alpha} (\vartheta - \theta^{\prime})^{\alpha^{\prime}}, \ \ A_{\alpha \alpha^{\prime}} \in S_{\mu}^{\min} . $$
Such expansions are of course similar to those in Definition \ref{defi-asymp-symbol}. 

\bigskip

\noindent {\it Proof.} 
We study the kernel (\ref{kernelexplicitdouble})
with $t=0$ (and the dependence on $ \epsilon $ omitted).
Consider the function $ ( \hat{\rho}, \hat{\eta} ) $  of  $ (r,\theta,r^{\prime},\theta^{\prime},\varrho,\vartheta) $ defined by
$$ ( \hat{\rho}, \hat{\eta} ) = \int_0^1 \big( \partial_r \varphi , \partial_{\theta} \varphi \big) (r_{\lambda} , \theta_{\lambda},\varrho,\vartheta) d \lambda $$
where $ r_{\lambda} = r^{\prime} + \lambda (r - r^{\prime}) $ and $ \theta_{\lambda} = \theta^{\prime} + \lambda (\theta - \theta^{\prime}) $. Note that by convexity of $V$ (see after (\ref{conditiondesupportdouble})), $ (r_{\lambda},\theta_{\lambda},\varrho,\vartheta) $ belongs to $ \Theta^{+}(R,V,I,\varepsilon) $ if both $ (r,\theta, \varrho , \vartheta) $ and $ (r^{\prime},\theta^{\prime},\varrho,\vartheta) $ do.
Introduce next
  $$ \hat{\xi}  =  \frac{2}{r+r^{\prime}} \frac{\bar{g}(\theta)^{-1}}{\hat{\rho}} \hat{\eta}  $$ 
so that the  phase  becomes
\begin{eqnarray}
 \varphi (r,\theta,\varrho,\vartheta) - \varphi (r^{\prime},\theta^{\prime},\varrho,\vartheta) = \hat{\rho} (r-r^{\prime}) + \hat{\xi} \cdot \frac{r+r^{\prime}}{2} \hat{\rho} \bar{g} (\theta) (\theta - \theta^{\prime}) .  \label{fullphasepseudo}
\end{eqnarray} 
  By using  (\ref{expansioneikonal1}) and (\ref{expansioneikonal2}), we obtain
$$ \hat{\rho} = \varrho \left( 1 + S_{-\nu}^{\rm min} \big((\vartheta - \theta) + (\vartheta - \theta^{\prime}) \big) + S_{0}^{\rm min} \big((\vartheta - \theta) + (\vartheta - \theta^{\prime}) \big)^2 \right), $$
 and
 \begin{eqnarray}
  \hat{\xi}  =   \vartheta - \theta^{\prime} + \ \frac{r^{\prime} (\theta-\theta^{\prime})}{3(r+r^{\prime})} + S_{-\nu}^{\rm min} \big((\vartheta - \theta) + (\vartheta - \theta^{\prime}) \big) + S_{0}^{\rm min} \big((\vartheta - \theta) + (\vartheta - \theta^{\prime}) \big)^2 .
 \end{eqnarray} 
 Both expansions follow from routine computations, using that $ \vartheta - \theta_{\lambda} =(1- \lambda) (\vartheta - \theta^{\prime}) + \lambda (\vartheta - \theta ) $
 and
 $$ \frac{2}{r+r^{\prime}} \int_0^1 r_{\lambda} (\vartheta - \theta_{\lambda}) d \lambda =  \vartheta - \theta^{\prime} + \ \frac{r^{\prime} (\theta-\theta^{\prime})}{3(r+r^{\prime})} .  $$
All this shows that $ \hat{\rho} $ and (the components of) $\hat{\xi} $ belong to $ S_0^{\rm min} $, and also that
 $$ \big| d_{\varrho,\vartheta} (\hat{\rho} , \hat{\xi}) - I_n \big| \lesssim \min(r,r^{\prime})^{ - \nu} + |\vartheta - \theta | +  |\vartheta - \theta|^{\prime} . $$
Thus, if we assume that $ r ,r^{\prime } > R^{\prime} \gg 1 $ and $ \varepsilon^{\prime} \ll 1 $, $ (\varrho , \vartheta ) \mapsto (\hat{\rho} , \hat{\xi} ) $ is a diffeomorphism from $ \{|\vartheta - \theta|< \varepsilon^{\prime}\} \cap \{ |\vartheta - \theta^{\prime}| < \varepsilon^{\prime} \} \cap \{ \varrho^2 \in I, \ \varrho > 0 \}  $ onto its range. If we denote by $ (\rho,\xi) \mapsto ( \check{\varrho} , \check{\vartheta} ) $ the inverse map (which depends also on $ r,\theta,r^{\prime},\theta^{\prime} $), the fact that $ \hat{\rho},\hat{\xi} \in S_0^{\min} $ implies that,
 \begin{eqnarray}
\big| \partial_r^j \partial_{\theta}^{\alpha} \partial_{r^{\prime}}^{j^{\prime}} \partial_{\theta^{\prime}}^{\alpha^{\prime}} \partial_{\rho}^k \partial_{\xi}^{\beta} (\check{\varrho},\check{\vartheta}) \big| \leq C \min (r,r^{\prime})^{-j-j^{\prime}} \label{estimeecruciale}
 \end{eqnarray}
 on its domain of definition, hence on the support of $ a (r,\theta,\check{\varrho},\check{\vartheta}) \bar{b}(r^{\prime},\theta^{\prime}, \check{\varrho},\check{\vartheta}) $. Also, since $ \hat{\rho} - \varrho  $ is small, $ \rho $ must belong to a compact subset of $ (0,+\infty) $ (remember we prove the outgoing case).
 Then, by using successively the changes of variables $ (\varrho , \vartheta ) \mapsto (\hat{\rho} , \hat{\xi} ) $ and  $ \xi \mapsto \eta :=  \frac{r+r^{\prime}}{2} \rho \bar{g} (\theta) \xi $ (recall (\ref{fullphasepseudo})), the kernel of $ J^h (a) J^h (b)^{\dag} $ becomes
 $$ (2 \pi h)^{-n} \int \! \! \int e^{\frac{i}{h}\big( (r-r^{\prime})\rho + (\theta-\theta^{\prime})\cdot \eta \big)}
 a (r,\theta, \tilde{\varrho},\tilde{\vartheta}) \bar{b} (r^{\prime},\theta^{\prime},\tilde{\varrho},\tilde{\vartheta}) | \partial_{\rho,\eta} (\tilde{\varrho} , \tilde{\vartheta}) | d \rho d \theta $$
 where
 $$ (\tilde{\varrho} , \tilde{\vartheta } ) = ( \hat{\varrho} , \hat{\vartheta} ) \left(r , \theta , r^{\prime} , \theta^{\prime} , \rho ,  \frac{2}{r+r^{\prime}} \frac{\bar{g}(\theta)^{-1}}{\rho} \eta  \right)$$
 and where $ | \partial_{\rho,\eta} (\tilde{\varrho} , \tilde{\vartheta}) | $ is the corresponding Jacobian, which satisfies in particular
\begin{eqnarray}
 | \partial_{\rho,\eta} (\tilde{\varrho} , \tilde{\vartheta}) | = O \big( (r+r^{\prime})^{1-n} \big) . \label{explicationexposant}
\end{eqnarray} 
  Note in addition that, restricted to $ r = r^{\prime} $ and $ \theta = \theta^{\prime} $, $  (\tilde{\varrho} , \tilde{\vartheta } ) = (\bar{\varrho} , \bar{\vartheta} ) $
 since it is the inverse of $ (\varrho,\vartheta) \mapsto \partial_{r,\theta} \varphi (r,\theta,\varrho,\vartheta) $.
   One can then rewrite the kernel with an amplitude $c (h) $ independent of $ (r^{\prime},\theta^{\prime}) $ according to the usual procedure (see {\it e.g.} \cite[Theorem 4.20]{Zworski}). That $ c (h) $ belongs to $ \widetilde{S}^{-\infty, \mu + \mu^{\prime} + 1 -n } $ follows from (\ref{estimeecruciale}), (\ref{explicationexposant}) and the fact that $ a \in S_{\mu} $, $ b \in S_{\mu^{\prime}}  $. 
 This concludes the proof. \finpreuve

\bigskip

We next consider two applications of Proposition \ref{LemmeEgorovstat}.

\begin{prop} \label{borneL2utile} If $ (a_{\epsilon})_{\epsilon} $ is a bounded family in $ S_{0} (\Ra^{2n}) $, supported in $ \Theta^{\pm}(R^{\prime},V,I,\varepsilon^{\prime}) $ (with $ \varepsilon^{\prime} $ as in Proposition \ref{LemmeEgorovstat}), then
$$ \big| \big| J^h_{\epsilon} (a_{\epsilon}) \big| \big|_{ L^2 (dx d \vartheta) \rightarrow L^2 (\scal{r}^{n-1} dr d \theta) } \leq C ,  $$
with a constant $C$ independent of $h, \epsilon \in (0,1]$. Similarly, if $ (b_{\epsilon})_{\epsilon} $ is a bounded family of $ S_{n-1} $  supported in $ \Theta^{\pm}(R^{\prime},V,I,\varepsilon^{\prime}) $,
$$  \big| \big| J^h_{\epsilon} (b_{\epsilon})^{\dag} \big| \big|_{   L^2 (\scal{r}^{n-1} dr d \theta)  \rightarrow L^2 (dx d \vartheta) } \leq C ,  $$
with a constant $C$  independent of $h, \epsilon \in (0,1]$. 
\end{prop}


\bigskip

\noindent {\it Proof.} The first estimate is equivalent to the fact that 
$ \scal{r}^{\frac{n-1}{2}} J^h (a) J^h(a)^{\dag} \scal{r}^{\frac{n-1}{2}} $ is bounded on $ L^2 (\Ra^n) $ equipped with the Lebesgue measure.  By Proposition \ref{LemmeEgorovstat}, $  J^h (a) J^h(a)^{\dag}  $
is of the form $ O \! p^h (c(h)) $ for some admissible symbol $ c(h) \in \widetilde{S}^{-\infty,1-n} $. Thus, when composed on both sides with $ \scal{r}^{\frac{n-1}{2}} $, Proposition \ref{bonpseudodiff} shows we get a pseudo-differential operator with admissible symbol in $ \widetilde{S}^{-\infty,0} $.
Since such pseudo-differential operators are bounded on  $L^2(\Ra^n)$, according to the usual Calder\'on-Vaillancourt Theorem, the result follows. The second estimate is equivalent to the boundedness of $ J_{}^h (b)^{\dag} \scal{r}^{\frac{1-n}{2}}  $ on $ L^2 (\Ra^n) $ and thus follows from the first case by taking the adjoint since $ b \scal{r}^{\frac{1-n}{2}} = \scal{r}^{\frac{n-1}{2}}a  $ for some $a \in S_0 $. \finpreuve

\bigskip


In the next proposition, to take into account the dependence on $ \epsilon $, we introduce the sets
\begin{eqnarray}
  \widetilde{\Gamma}_{\epsilon,{\rm st}}^{\pm} (R,V,I,  \varepsilon ) = \{ (r,\theta,\rho,\eta) \ | \ r > R , \ \theta \in V , \ p_{\epsilon,\kappa} \in I, \ \pm \rho > (1 - \varepsilon^2) p_{\epsilon,\kappa}^{1/2} \} . \label{referencestrongarealow}
\end{eqnarray}
This is the convenient replacement of (\ref{referencestrongarea}) at low frequency. It allows to cover the case $ \epsilon =1 $ used for high frequency parametrices (in which case we drop the dependence on $ \epsilon $), while the regime $ \epsilon \in (0,1) $ will be for low frequency parametrices. In this last case, $ \widetilde{\Gamma}_{\epsilon,{\rm st}}^{\pm} (R,V,I,  \varepsilon ) $ has to be understood as a set of $ (\breve{r},\theta,\breve{\rho},\eta) $. We use only  (\ref{referencestrongarealow}) in the intermediate technical statements but, for clarity, we will use both (\ref{referencestrongarea}) and (\ref{referencestrongarealow}) to state the main result of this section (Theorem \ref{Isozaki-Kitada-explicite}). 

\begin{prop}[Factorizing $ \Psi $DO] \label{propfactorization} Assume we are given $N$ bounded families $ (a_{\epsilon,0}) , \ldots , (a_{\epsilon,N}) $  of symbols supported in $ \Theta^{\pm}(R^{\prime},V,I,\varepsilon^{\prime}) $ 
such that, for some $ c > 0 $ independent of $ \epsilon  $,
$$ a_{\epsilon,j} \in S_{-j}(\Ra^{2n}), \qquad a_{\epsilon,0} \geq c > 0 \ \ \mbox{on some} \ \Theta^{\pm} \big(R^{\prime \prime},V^{\prime \prime},I^{\prime },\varepsilon^{\prime \prime} \big) . $$
Let $ I^{\prime \prime} \Subset I^{\prime} $. Then there exists $ C > 0 $ such that, for all $ 0 < \varepsilon \ll 1 $, $ \mu \in \Ra $ and all bounded family $ (f_{\epsilon}) $ of $ \widetilde{S}^{-\infty,\mu} (\Ra^{2n}) $ such that
$$ \emph{supp}(f_{\epsilon}) \subset \widetilde{\Gamma}_{\epsilon,{\rm st}}^{\pm} (R^{\prime \prime},V^{\prime \prime},I^{\prime \prime},  \varepsilon ),  $$
one can write
$$ O \! p^h (f_{\epsilon}) = \sum_{j+k \leq N} h^{j+k}J^h_{\epsilon}(a_{\epsilon,j}) J^h_{\epsilon} (b_{\epsilon,k})^{\dag} + h^N O \! p^h ( f_{\epsilon,N} (h)) , $$
with  $ ( f_{\epsilon,N} (h))_{\epsilon,h \in (0,1]} $ bounded in  $ \widetilde{S}^{-\infty,\mu-N}(\Ra^{2n}) $ and some $ b_{\epsilon,0} , \ldots , b_{\epsilon,N} $ such that
\begin{eqnarray}
 ( b_{\epsilon,k} )_{\epsilon \in (0,1]} \ \ \mbox{bounded in} \ \ S_{\mu+ n-1-k} (\Ra^{2n}) , 
\end{eqnarray}
and
\begin{eqnarray}  
    \emph{supp}(b_{\epsilon,k}) \subset \emph{supp} \big( f_{\epsilon} ( . , . , \partial_{r} \varphi_{\epsilon},\partial_{\theta} \varphi_{\epsilon}) \big) \subset \Theta^{\pm} (R^{\prime \prime}, V^{\prime \prime},I^{\prime},C \varepsilon) .  \label{supportestimate}
\end{eqnarray}  
For $k=0$, we have explicitly  
\begin{eqnarray}
  b_{\epsilon,0} (r,\theta,\varrho,\vartheta) = \bar{f}_{\epsilon} \big(r,\theta,\partial_r \varphi_{\epsilon}, \partial_{\theta} \varphi_{\epsilon} \big) \frac{| \emph{det} \big( \partial_{\varrho,\vartheta}  \partial_{r, \theta} \varphi_{\epsilon}  \big) | }{\bar{a}_{\epsilon,0} (r,\theta,\varrho,\vartheta)}  . \label{exemplepratique}
\end{eqnarray}  
\end{prop}

Notice that when $ \mu = 0 $, this proposition shows in particular that a bounded pseudo-differential operator can be factorized (up to a nice error) as a product $ J_{\epsilon}^h (a) J_{\epsilon}^h (b)^{\dag} $ where, according to Proposition \ref{borneL2utile}, $ J_{\epsilon}^h (a) $ and $ J_{\epsilon}^h (b)^{\dag} $ are bounded respectively from $ L^2 (dx d \vartheta) $ to $ L^2 (\scal{r}^{n-1} dr d \theta) $ and from  $ L^2 (\scal{r}^{n-1} dr d \theta) $ to $ L^2 (dx d \vartheta) $.

\bigskip

\noindent {\it Proof of Proposition \ref{propfactorization}.} The principle is well known. We recall it briefly to emphasize where the support estimate in (\ref{supportestimate}) comes from. To seek which conditions must be fulfilled by the $b_k$'s we compute first
$$  \sum_{j= 0}^{N-1} \sum_{k = 0}^{ N-1} h^{j+k}J^h (a_j) J^h (b_k)^{\dag} = \sum_{j,k,l=0}^{N-1} h^{j+k+l} O \! p^h (c_{j,k,l}) + h^N O \! p^h (r_N(h)) . $$
By Proposition \ref{LemmeEgorovstat}, the first symbol reads $ c_{0,0,0}  = a_0 (r,\theta,\bar{\varrho},\bar{\vartheta}) \bar{b}_0 (r,\theta,\bar{\varrho},\bar{\vartheta})  |{\rm det}(\partial_{\rho,\eta}(\bar{\varrho},\bar{\vartheta}))| $ so the requirement  that $ c_{0,0,0} =  f $ together with Proposition \ref{projLagrangian} (in particular (\ref{formulesplusprecises})) show that $ b_0 $ must equal (\ref{exemplepratique}). This function is well defined since $  \bar{f} \big(r,\theta,\partial_r \varphi, \partial_{\theta} \varphi \big)  $ is supported in the image of $ \mbox{supp}(f) $ by the map (\ref{unnameddiffeo}) hence, using (\ref{proof01}) and (\ref{proof2}), 
in $ \Theta^{+} (R^{\prime \prime}, V^{\prime \prime},I^{\prime},C \varepsilon)  $ if $ \varepsilon $ is small enough; in particular, $a_0$ is bounded below on such a domain. Using then that $ \mbox{det} \big( \partial_{\varrho,\vartheta}  \partial_{r, \theta} \varphi  \big) \in S_{n-1} $, we see that $ b_0 \in \widetilde{S}^{-\infty,\mu+n-1} $. Then, the next symbol in the expansion is $ \sum_{j+k+l = 1}  c_{j,k,l} $ and we require it to be $0$, which yields the equation
$$ a_0 (r,\theta,\bar{\varrho},\bar{\vartheta}) \bar{b}_1 (r,\theta,\bar{\varrho},\bar{\vartheta})  |{\rm det}(\partial_{\rho,\eta}(\bar{\varrho},\bar{\vartheta}))| = -   \sum_{j+k+l = 1, \atop k =0}  c_{j,k,l}  $$
where, by Proposition \ref{LemmeEgorovstat} and the form of $ b_0 $, the right hand side vanishes outside the support of $ b_0 (r,\theta,\bar{\varrho},\bar{\vartheta}) $. One can thus divide by $a_0$ and find $ b_1 $. Higher order terms are obtained by iterating this process.  \finpreuve

\bigskip

In the sequel, we let $ U_0 (s) = e^{- is h D_x^2} $ be the semiclassical Schr\"odinger group on the line $ \Ra_x $.

\begin{prop}[Propagation estimates for the parametrix] \label{propparam} Let $ I \Subset (0,+\infty) $.  If  $ \varepsilon^{\prime} $ is small enough and $  R^{\prime} $ large enough then for all integer  $ N \geq 0 $, all bounded families $ (a_{\epsilon})_{\epsilon} $  of $ S_0 (\Ra^{2n}) $ and  $(b_{\epsilon})_{\epsilon} $ of $ S_{n-1}(\Ra^{2n}) $, both  supported in $ \Theta^{\pm} (R^{\prime},V,I,\varepsilon^{\prime}) $,  there exists $ C > 0 $ such that
\begin{eqnarray}
 \left| \left|  \scal{r}^{-N} J^h_{\epsilon}(a_{\epsilon}) U_0 (s) J^h_{\epsilon}(b_{\epsilon})^{\dag}  (r \pm s)^N  \right| \right|_{L^2(\scal{r}^{n-1}dr d \theta) \rightarrow L^2 (\scal{r}^{n-1}dr d \theta)} \leq C , \label{propagationparametrix}
\end{eqnarray} 
for all $ \pm s \geq 0 $ and all $ h, \epsilon \in (0,1] $.  In particular, we have
$$ \left| \left|  \scal{r}^{-N_1-N_2} J^h_{\epsilon}(a_{\epsilon}) U_0 (s) J^h_{\epsilon}(b_{\epsilon})^{\dag} \scal{r}^{N_1}  \right| \right|_{L^2(\scal{r}^{n-1}dr d \theta) \rightarrow L^2 (\scal{r}^{n-1}dr d \theta)} \leq C \scal{s}^{-N_2},  $$
if $ N_1 , N_2 \geq 0 $ are integers.
\end{prop}

\noindent {\it Proof.}  The main observation is that $ 2 s \varrho + \partial_{\varrho} \varphi (r^{\prime},\theta^{\prime},\varrho,\vartheta) \gtrsim r^{\prime} +s $ by (\ref{expansioneikonal4}). We can then write
$$ r^{\prime} + s = \frac{r^{\prime} + s}{ \partial_{\varrho} \varphi^{\prime} + 2 s \varrho } ( \partial_{\varrho} \varphi^{\prime} + 2 s \varrho )  $$
where the prime on $ \varphi $ is a shortand for the evaluation at $ (r^{\prime},\theta^{\prime},\vartheta,\varrho) $. Here the fraction belongs to $ S_0 $ uniformly with respect to $ s \geq 0 $. Writing next $  \partial_{\varrho} \varphi^{\prime} + 2 s \varrho =   \partial_{\varrho} \varphi(r,\theta,\varrho,\vartheta) - \partial_{\varrho} \Phi  $ and setting $ \tilde{b} = b (r + s )/(\partial_{\varrho} \varphi + 2 s \varrho) $,   integrating by part in $ \varrho $ shows that
$\scal{r}^{-1} J^h(a)U_0(s) J^h(b)^{\dag} (r +s ) $ reads
$$ J^h \big( \scal{r}^{-1} \partial_{\varrho} \varphi a \big) U_0 (s) J^h (\tilde{b})^{\dag} - i h \scal{r}^{-1} \big( J^h (\partial_{\varrho} a)U_0(s) J^h (\tilde{b})^* + J^h (\partial_{\varrho}  a)U_0(s) J^h ( \partial_{\varrho} \tilde{b})^{\dag} \big) $$
which is bounded on $ L^2 (\scal{r}^{n-1}dr d \theta) $ uniformly in $s \geq 0$ by  Proposition \ref{borneL2utile} since
 $ \scal{r}^{-1} \partial_{\varrho} \varphi a $ and $ \tilde{b} $ belong respectively to $ S_0 (\Ra^{2n}) $ and $ S_{n-1} (\Ra^{2n}) $ (uniformly in $s \geq 0$ for $ \tilde{b} $). This proves the estimate (\ref{propagationparametrix}) with $ N = 1 $.  For $ N \geq 2 $, the result is obtained by iteration of this process.  \finpreuve

\bigskip

We next turn to the proof of dispersive estimates for the oscillatory integrals of the form (\ref{integraleoscillante}).

\begin{prop}[Stationary phase estimates] \label{statphaseIK} Let $ I \Subset (0,+\infty) $. If $ \varepsilon $, $ \varepsilon^{\prime} $ are small enough and $ R , R^{\prime} $ large enough, then for all bounded family $ (A_{\epsilon})_{\epsilon} $ of  $ S^{\min}_{0}(\Ra^{3n}) $ satisfying (\ref{conditiondesupportdouble}), one has
$$ \big|  I^h_{\epsilon} (A_{\epsilon})  \big|  \lesssim \min \left( h^{-n} , |hs|^{-n/2} \right)  , \qquad s \in \Ra,  \ \epsilon, h \in (0,1].  $$
\end{prop}

Notice that, unlike the non stationary phase estimates of Proposition \ref{nonstatprop} and the propagation estimates of Proposition \ref{propparam}, we do not need any sign condition on $s$ here. 

To prove Proposition \ref{statphaseIK} (omitting $ \epsilon $ as before), we will rewrite
$$ \varphi (r,\theta,\varrho,\vartheta) - \varphi (r^{\prime},\theta^{\prime},\varrho,\vartheta) = (r-r^{\prime}) \widetilde{\partial}_r \varphi+ (\theta - \theta^{\prime}) \cdot \widetilde{\partial}_{\theta} \varphi $$
where, setting $ r_{\lambda} = r^{\prime} + \lambda (r-r^{\prime}) $ and $ \theta_{\lambda} = \theta^{\prime} + \lambda (\theta - \theta^{\prime}) $, 
$$ \widetilde{\partial}_r \varphi := \int_0^1 \partial_r \varphi (r_{\lambda},\theta^{\prime},\varrho,\vartheta) d \lambda , \qquad  \widetilde{\partial}_{\theta} \varphi := \int_0^1 \partial_{\theta} \varphi (r,\theta_{\lambda},\varrho,\vartheta) d \lambda . $$
\begin{lemm}[Improved asymptotic expansion]
$$ \widetilde{\partial}_r \varphi = \varrho \left( 1 - \frac{1}{2}  ( \vartheta - \theta^{\prime} ) \cdot \bar{g} (\theta^{\prime})(\vartheta - \theta^{\prime})  + S_{-\nu}^{\rm min} ( \vartheta - \theta^{\prime})^2 + S_0^{\rm min} (\vartheta - \theta^{\prime})^{3} \right) .
$$
\end{lemm}

\noindent {\it Proof.}  Using the notation and  estimates of the proof of Proposition \ref{prop-inter-flow}, we have $ \bar{\eta}^s = \eta + \widetilde{S}_1 (\eta/r) $ and $ \bar{g} (\vartheta^s) = \bar{g} (\theta) + \widetilde{S}_0 (\eta/r)$  (where the remainders $ \widetilde{S}_{0  } (\eta/r) $, $ \widetilde{S}_{  1} (\eta/r) $ depend in a bounded fashion on $s$) so by using the motion equations and letting $s$ go to infinity, we get easily
$$ \bar{\varrho} = \rho + \frac{1}{2\rho} \frac{\eta}{r} \cdot \bar{g} (\theta)^{-1} \frac{\eta}{r} + \widetilde{S}_{-\nu} \left( \frac{\eta}{r} \right)^2 + \widetilde{S}_0 \left( \frac{\eta}{r} \right)^3 . $$
Evaluating this identity at $ (\underline{\rho},\underline{\eta}) = (\partial_r \varphi ,\partial_{\theta} \varphi ) $ and using (\ref{expansioneikonal1})-(\ref{expansioneikonal2}), we find
$$ \varrho = \partial_r \varphi + \frac{\varrho}{2} (\vartheta - \theta) \cdot \bar{g} (\theta)^{-1} (\vartheta - \theta ) + S_{-\nu} \left( \vartheta - \theta \right)^2 + S_0 \left( \vartheta - \theta \right)^3  . $$
 This provides an expansion of $ \partial_r \varphi $ which yields the result after evaluation at $(r_{\lambda},\theta^{\prime}, \varrho,\vartheta)$ and integration on $ [0,1]_{\lambda} $. \finpreuve

\bigskip

\begin{lemm} \label{lemmeangulaire}  Let $ \delta \in (0,1) $. If $ \varepsilon $, $ \varepsilon^{\prime} $ are small enough and $ R , R^{\prime} $ large enough, then
$$ r |\theta - \theta^{\prime}| \geq \delta |s| \qquad \mbox{and} \qquad  |s| \geq h  \qquad \Longrightarrow \qquad   I^h (A) = O \big( h^{-n}  (s/h)^{ - \infty} \big) .  $$
\end{lemm}



\noindent {\it Proof of Lemma \ref{lemmeangulaire}.} Let us observe first that
$$ \widetilde{\partial}_{\theta} \varphi = \varrho r \left[ \bar{g}(\theta^{\prime}) \left( \vartheta - \frac{\theta + \theta^{\prime}}{2} \right) + S_{-\nu}^{\rm min} \big( (\vartheta - \theta) + (\vartheta - \theta^{\prime})\big) + S_0^{\rm min} \big( (\vartheta - \theta) + (\vartheta - \theta^{\prime})\big)^2  \right] , $$
which follows from  the expansion (\ref{expansioneikonal2}) and by writing $ \theta_{\lambda} = \theta^{\prime} + \lambda (\theta - \vartheta) + \lambda (\vartheta - \theta^{\prime}) $.
Then,  in the integral (\ref{integraleoscillante}), we use the one dimensional (linear) change of variable $ \varrho \mapsto \widetilde{\partial}_r \varphi $.  Its inverse  is of the form
\begin{eqnarray}
 \tilde{\varrho} \mapsto   \left( 1 + \frac{1}{2}  ( \vartheta - \theta^{\prime} ) \cdot  \bar{g} (\theta^{\prime}) (\vartheta - \theta^{\prime})  + S_{-\nu}^{\rm min} ( \vartheta - \theta^{\prime})^2 + S_0^{\rm min} (\vartheta - \theta^{\prime})^{3} \right) \tilde{\varrho}  .  \label{expressiontilderho} 
\end{eqnarray} 
Letting $ \tilde{\Phi} $ be the expression of $ \Phi $ composed with  this change of variable, we have
\begin{eqnarray}
\widetilde{\Phi} = (r-r^{\prime}) \tilde{\varrho} - s \tilde{\varrho}^2 \left( 1 + (\vartheta - \theta^{\prime}) \cdot \bar{g}(\theta^{\prime}) (\vartheta - \theta^{\prime})  \right) + \tilde{\varrho} r (\theta - \theta^{\prime}) \bar{g}(\theta^{\prime}) \left( \vartheta - \frac{\theta + \theta^{\prime}}{2} \right) + \widetilde{\Omega}
\end{eqnarray}
with a remainder of the form
\begin{eqnarray}
\widetilde{\Omega} & = &  r (\theta - \theta^{\prime}) \left( S_{-\nu}^{\rm min} \big( (\vartheta - \theta) + (\vartheta - \theta^{\prime}) \big) + S_0^{\rm min} \big( (\vartheta - \theta) + (\vartheta - \theta^{\prime}) \big)^2 \right) + \nonumber \\ &  & \ \  \  s \left( S_{-\nu}^{\rm min} (\vartheta - \theta^{\prime})^2 + S_0^{\rm min} (\vartheta - \theta^{\prime})^3 \right) .
\nonumber
\end{eqnarray}
The interest of this change of variable is that the only term involving $ r  - r^{\prime} $, namely $ ( r - r^{\prime} ) \tilde{\varrho} $, is independent of $ \vartheta $. Therefore, using the above expansion, we have
$$ \partial_{\vartheta} \tilde{\Phi} = \tilde{\varrho}  \bar{g} (\theta^{\prime}) r (\theta - \theta^{\prime}) + s  O (\varepsilon^{\prime})  + r (\theta - \theta^{\prime}) \big( O \big(\min (R,R^{\prime})^{-\nu} \big) + O (\varepsilon) + O (\varepsilon^{\prime}) \big) . $$
Hence, by using $ r |\theta - \theta^{\prime}| \geq \delta |s| $ and by taking $ \varepsilon , \varepsilon^{\prime} $ small enough as well as $ R, R^{\prime} $ large enough, we get a lower bound $ |\partial_{\vartheta} \tilde{\Phi} | \gtrsim |s| $ from which the result follows by integrations by part. \finpreuve

\bigskip

\noindent {\bf Proof of Proposition \ref{statphaseIK}.} The estimate is trivial if $ |s| \leq h $. Thus we  assume that $ |s| \geq h $ and, according to Lemma \ref{lemmeangulaire}, that $ r |\theta - \theta^{\prime}| \leq \delta |s| $ for some small enough $ \delta $ to be chosen below,  otherwise we use that $ h^{-n} |s/h|^{-N} \lesssim h^{-n} |s/h|^{-n/2} = |hs|^{-n/2} $ for any integer $ N \geq n/2 $. Using the same change of variable as in Lemma \ref{lemmeangulaire}, we find that the Hessian matrix of $ \tilde{\Phi} $ reads
$$ d^2_{\tilde{\varrho},\vartheta} \tilde{\Phi} = - 2 s \left( \begin{matrix} 1 & 0 \\ 0  & \tilde{\varrho}^2 \bar{g} (\theta^{\prime}) \end{matrix} \right) + s \big( O (\varepsilon^{\prime}) + O (\min(R,R^{\prime})^{-\nu}) \big)+  O (r|\theta - \theta^{\prime}|) . $$
We choose  $ \delta $ small enough so that $ O (r|\theta-\theta^{\prime}|/s) = O (\delta )  $ is sufficiently small with respect to the first matrix on the right hand side (here we use (\ref{hjkellipticity})).  This imposes to consider $ \varepsilon $ and $ \varepsilon^{\prime} $ sufficiently small too and $ R,R^{\prime} $ sufficiently large  to use Lemma \ref{lemmeangulaire}. Then, by possibly decreasing again $ \varepsilon^{\prime} $ and increasing again $ R,R^{\prime}$, we find that $ s^{-1} d^2_{\tilde{\varrho},\vartheta} \tilde{\Phi} $ is a negative definite matrix uniformly with respect to $ r , r^{\prime} , \theta , \theta^{\prime} $ on the support of the amplitude (and such that $ r |\theta - \theta^{\prime}| \leq \delta |s| $). The result then follows from the stationary phase theorem with $s/h$ as a large parameter.  \finpreuve

\subsection{Construction of the parametrix }
In this paragraph, we state the main result of the section which is Theorem \ref{Isozaki-Kitada-explicite} on the construction of an Isozaki-Kitada type parametrix. Given a chart $ \kappa : U_{\kappa} \subset {\mathcal S} \rightarrow V_{\kappa} \subset \Ra^{n-1} $ and $V \subset V_{\kappa}$ as in Theorem \ref{eikonal-equation}, we  introduce the notation
$$ J^h_{\kappa} (a) := \Pi_{\kappa} J^h (a), \qquad  J^h_{\kappa} (b)^{\dag} := J^h (b)^{\dag} \Pi_{\kappa}^{-1} $$
and
\begin{eqnarray}
 J_{\epsilon,\kappa} (a) := \Pi_{\kappa}  {\mathscr D}_{\epsilon} J_{\epsilon} (a) , \qquad J_{\epsilon,\kappa} (b)^{\dag}:= J_{\epsilon} (b)^{\dag} {\mathscr D}_{\epsilon}^{-1} \Pi_{\kappa}^{-1} .  \label{pourmemoireavecrescaling}
\end{eqnarray} 
where, in (\ref{pourmemoireavecrescaling}), the symbols will depend on $ \epsilon $ in the applications.

As a starting point, we observe that the general formula
$$ e^{-itP} B (0) = B (t) - i \int_0^t e^{-i (t - \tau)P} \left( P B (\tau) - i B^{\prime}(\tau)  \right) d \tau , $$
leads respectively to the identities
\begin{eqnarray}
e^{-itP}  J^h_{\kappa} (a) J^h_{\kappa} (b)^{\dag}  & = & J^h_{\kappa} (a) e^{- i t  D_x^2} J^h_{\kappa} (b)^{\dag} 
    - \mbox{R}_{\rm hi}  
\end{eqnarray}
with
\begin{eqnarray}  
\mbox{R}_{hi} = \frac{i}{h^2} \int_0^{t} e^{-i(t-\tau)P}  \left(   \Pi_{\kappa}  \left[ h^2 P_{\kappa} J^h (a) - J^h (a) h^2 D_x^2 \right] e^{-i\tau D_x^2} J^h (b)^{\dag}   \right)  d \tau , \nonumber
\end{eqnarray}
and similarly,
\begin{eqnarray}
e^{-itP}  J_{\epsilon,\kappa} (a_{\epsilon}) J_{\epsilon,\kappa} (b_{\epsilon})^{\dag}  & = &
   J_{\epsilon,\kappa} (a_{\epsilon}) e^{- i t \epsilon^2 D_x^2} J_{\epsilon,\kappa} (b_{\epsilon})^{\dag}   - \mbox{R}_{\rm lo}  
\end{eqnarray}
with
\begin{eqnarray}  
\mbox{R}_{\rm lo} = i \epsilon^2  \int_0^{t} e^{-i(t-\tau)P}  \left(   \Pi_{\kappa}   {\mathscr D}_{\epsilon}  \left[ P_{\epsilon,\kappa} J_{\epsilon} (a_{\epsilon}) - J_{\epsilon} (a_{\epsilon}) D_x^2 \right] e^{-i\tau \epsilon^2 D_x^2} J_{\epsilon,\kappa} (b_{\epsilon})^{\dag}  \right)  d \tau . \nonumber
\end{eqnarray}
Recall from (\ref{refoprescale}) that $ P \Pi_{\kappa} {\mathscr D}_{\epsilon}  = \epsilon^2 \Pi_{\kappa} {\mathscr D}_{\epsilon}  P_{\epsilon,\kappa}$. Note also the scaling in time.

We seek $ a,b $ and $ a_{\epsilon},b_{\epsilon} $ such that $ R_{\rm hi} $ and $ R_{\rm lo} $ are  respectively small and such that $ J^h_{\kappa} (a) J^h_{\kappa} (b)^{\dag}  $ and $ J_{\epsilon,\kappa} (a_{\epsilon}) J_{\epsilon,\kappa} (b_{\epsilon})^{\dag} $ can be prescribed.

We consider in detail the high frequency case. The first step is to find $$ a = a (h) := a_0 + h a_1 + \cdots + h^M a_M , $$ such that $ h^2 P_{\kappa} J^h (a(h)) - J^h (a(h)) h^2 D_x^2  $ is small, in an appropriate sense (here $ M $ is an arbitrary integer order which is fixed).
A simple calculation yields
\begin{eqnarray}
  h^2 P_{\kappa} J_h (a(h)) - J_h (a(h)) h^2 D_x^2 = J_h \left( c_0 + \cdots + h^{M+2} c_{M+2} \right)  , \label{entrelacementapproche}
\end{eqnarray}  
where
\begin{eqnarray}
c_0 & = & E a_0 \\
c_1 & =&  E a_1 - i T a_0 \\
c_j & = & E a_j - i T a_{j-1} + P_{\kappa} a_{j-2}, \qquad 2 \leq j \leq M \\
c_{M+1} & = & - i T a_M + P_{\kappa} a_{M-1}  \\
c_{M+2} & =& P_{\kappa} a_M
\end{eqnarray}
where $ E $ corresponds to the eikonal term and $ T $ to the transport operator, namely
$$ E = p_{\kappa} (r,\theta,\partial_{r,\theta} \varphi) - \varrho^2, \qquad T = (\partial_{\rho,\eta} p) \big(r,\theta,\partial_{r,\theta} \varphi \big) \cdot \partial_{r,\theta} - P_{\kappa} \varphi .  $$
By  Theorem \ref{eikonal-equation},  we can solve the equation $ E  = 0$ on  $ \Theta (R,V,\varepsilon) $  for any given convex subset $ V \Subset V_{\kappa} $ and some $ R \gg 1 $, $ \varepsilon \ll 1 $. Therefore, solving the system of equations 
\begin{eqnarray}
  c_j = 0, \qquad 0 \leq j \leq M+1  , \label{systemenfonctiondec}
\end{eqnarray}
on subsets of $ \Theta (R,V,\varepsilon) $
amounts to solve transport equations of the form (\ref{transport}), which can thus be done by Proposition \ref{prop-transport} (third item). More precisely, given $ I_0 \Subset (0,+\infty) $ and $ V_0 \Subset V $, we can find $ R_0 > R $, $ 0 < \varepsilon_0 < \varepsilon $ and solutions $ \hat{a}_0^{\pm} , \ldots , \hat{a}_M^{\pm}  $ to (\ref{systemenfonctiondec}) such that
$$ \hat{a}_j^{\pm} \in S_{-j} \left( \Theta^{\pm} (R_0,V_0,I_0,\varepsilon_0) \right) $$
 (see (\ref{Thetaenergie}) for the definition of $ \Theta^{\pm} (R,V,I,\varepsilon) $) with the additional condition that, locally uniformly with respect to $ (\theta , \vartheta, \varrho ) $,
 \begin{eqnarray}
 \hat{a}_0^{\pm} (r,\theta,\vartheta,\varrho) \rightarrow 1 , \qquad r \rightarrow \infty . \label{ellipticiteplusmoins}
\end{eqnarray}
We use the notation $ \hat{a}_j^{\pm} $ to make a clear difference between these symbols defined on $ \Theta^{\pm} (R_0,V_0,I_0,\varepsilon_0) $ and the final $ a_{j}^{\pm} $ defined globally on $ \Ra^{2n} $ in (\ref{defglobalFouriercurved}).
We also point out  the technical fact that, to find solutions $ \hat{a}^{\pm}_j $ defined on $ \Theta^{\pm} (R_0,V_0,I_0,\varepsilon_0)  $, we choose $ R_0 $ and $ \varepsilon_0 $  respectively large enough and small enough to ensure that
\begin{eqnarray}
 \Theta^{\pm} (R_0,V_0,I_0,\varepsilon_0) \subset {\mathcal T}^{\pm}  (R_0,V_0,I_0,\varepsilon_0) \subset \Theta^{\pm} (R,V,I_0, \varepsilon) \label{pourdeuxiemeinclusionsurtout}
\end{eqnarray}
(see prior to Proposition \ref{prop-transport} for $  {\mathcal T}^{\pm} (R,V,I,\varepsilon) $). The interest is to guarantee, if  $ (\bar{r}^s , \bar{\vartheta}^s )$ are the spatial components of the Hamiltonian flow of $ p_{\kappa} $,  that $ (\bar{r}^s(r,\theta,\partial_{r,\theta} \varphi ) , \bar{\vartheta}^s (r,\theta,\partial_{r,\theta} \varphi) , \varrho, \vartheta ) $ belongs to the domain of definition of $ \varphi $ for $ \pm s \geq 0 $ (see Proposition \ref{prop-transport}). The first inclusion in (\ref{pourdeuxiemeinclusionsurtout}) is trivial while the second one is a consequence of
$$ \bar{r}^s (r,\theta,\partial_{r,\theta} \varphi) \gtrsim r, \qquad \big|\bar{\vartheta}^s (r,\theta,\partial_{r,\theta} \varphi) - \theta \big| \lesssim \frac{|\partial_{\theta}\varphi|}{r} \lesssim |\vartheta - \theta| $$
which follow from the flow estimates of Proposition \ref{blackboxdynamic}, (\ref{Taylorexplicite0}) for $ \bar{\vartheta}^s $ and the asymptotics of $ \varphi $ in Theorem \ref{eikonal-equation}.

We next globalize the symbols. Given $ R_1 > R_0 $, $ V_1 \Subset V_0 $, $ I_1 \Subset I_0 $ and $ \varepsilon_1 < \varepsilon_0 $, it is easy to construct
$$ \chi_{\pm} \in S_0 (\Ra^{2n}), \qquad \chi_{\pm} \equiv 1 \ \ \mbox{on} \ \ \Theta^{\pm} (R_1,V_1,I_1,\varepsilon_1), \qquad \mbox{supp} (\chi_{\pm}) \subset \Theta^{\pm} (R_0,V_0,I_0,\varepsilon_0) , $$
by choosing it of the form $ \chi_{1} (r) \chi_2 \left(\theta - \vartheta \right) \chi_3 (\vartheta) \chi_4 (\pm \varrho) $ with suitable $ \chi_1,\chi_2,\chi_3 \in C_0^{\infty} $ and $ \chi_1 \equiv 1 $ near $ + \infty $. We  then define
\begin{eqnarray}
 a_j^{\pm} := \chi_{\pm} \hat{a}^{\pm}_j \in S_{-j} (\Ra^{2n}) . \label{defglobalFouriercurved}
\end{eqnarray}
Notice that if we compute (\ref{entrelacementapproche}) with $ a (h) = a (h)^{\pm} := \sum_j a_j^{\pm} $, we also have to take into account the derivatives falling on the cutoff $ \chi_{\pm} $; we summarize the above results in the following proposition,  including the case of low frequencies which is completely similar.
\begin{prop}[Approximate intertwining] Let $ V $ be a convex relatively compact subset of $ V_{\kappa} $. Then for all $ V_1 \Subset V_0 \Subset V $ and $ I_1 \Subset I_0 \Subset (0,+\infty) $, we can find $ R_1 > R_0 \gg 1 $ and $ 0 < \varepsilon_1 < \varepsilon_0 \ll 1 $ such that:
\begin{enumerate}
\item{{\bf at high frequency:} one can find symbols $ a_j^{\pm} \in S_{-j} (\Ra^{2n})$, $ j \geq 0 $,  supported in $ \Theta^{\pm} (R_0,V_0,I_0,\varepsilon_0) $ such that
$$ a_0^{\pm} (r,\theta,\vartheta,\varrho) \geq 1/2, \qquad \mbox{on} \ \ \Theta^{\pm} (R_1, V_1,I_1,\varepsilon_1) $$
and, if one sets $ a^h = a_0^{\pm} +\cdots + h^M a_M^{\pm} $,
$$ h^2 P_{\kappa} J^h \big( a^h \big) - J^h \big( a^h \big) h^2 D_x^2 = h^{M+2} J^h \big( r_M^h \big) + J^h \big( \check{a}^h \big) + J^h \big( a_{\rm c}^h \big) $$
with $ r_M^h \in S_{-M-2} $, $ \check{a}^h,a_{\rm c}^h \in S_0 $, all supported in $  \Theta^{\pm} (R_0,V_0,I_0,\varepsilon_0) $, bounded with respect to $h$ and, mainly, such that
\begin{eqnarray}
 \emph{supp}  \big( \check{a}^h \big) \subset \{ |\theta - \vartheta| \geq \varepsilon_1 \} , \qquad \emph{supp} \big(a_{\rm c}^h \big) \subset \{ r \leq R_1 \} . \label{pourphasenonstathautefreq}
\end{eqnarray}
}
\item{{\bf At low frequency:}  one can find bounded families of symbols $( a_{\epsilon,j}^{\pm})_{\epsilon \in (0,1]} $ in $ S_{-j} (\Ra^{2n})$, $ j \geq 0 $,  supported in $ \Theta^{\pm} (R_0,V_0,I_0,\varepsilon_0) $ such that
$$ a_{\epsilon,0}^{\pm} (r,\theta,\vartheta,\varrho) \geq 1/2, \qquad \mbox{on} \ \ \Theta^{\pm} (R_1, V_1,I_1,\varepsilon_1) $$
and, if one sets $ a_{\epsilon} = a_{\epsilon,0}^{\pm} +\cdots +  a_{\epsilon,M}^{\pm} $,
$$  P_{\epsilon,\kappa} J_{\epsilon} \big( a_{\epsilon} \big) - J_{\epsilon} \big( a_{\epsilon} \big)  D_x^2 =  J_{\epsilon} \big( r_{\epsilon,M}  \big) + J_{\epsilon} \big( \check{a}_{\epsilon} \big) + J_{\epsilon} \big( a_{\epsilon,{\rm c}} \big) $$
with $ r_{\epsilon,M } \in S_{-M-2} $, $ \check{a}_{\epsilon},a_{\epsilon,{\rm c}} \in S_0 $, all supported in $  \Theta^{\pm} (R_0,V_0,I_0,\varepsilon_0) $, bounded with respect to $\epsilon$ and such that
\begin{eqnarray}
 \emph{supp}  \big( \check{a}_{\epsilon} \big) \subset \{ |\theta - \vartheta| \geq \varepsilon_1 \} , \qquad \emph{supp} \big(a_{\epsilon,{\rm c}} \big) \subset \{ r \leq R_1 \} .  \label{pourphasenonstatbassefreq}
\end{eqnarray}
  }
\end{enumerate}
\end{prop}
We point out that the terms $ \check{a}, a_{\rm c} $ are the contributions of derivatives falling on the cutoff $ \chi_{\pm} $. The properties (\ref{pourphasenonstathautefreq}) and (\ref{pourphasenonstatbassefreq}) will be useful to derive non stationary phase estimates from Proposition \ref{nonstatprop}. The ellipticity condition $ a_{0}^{\pm} \geq 1/2 $ (and likewise for $ a_{\epsilon,0}^{\pm} $) is a consequence of (\ref{ellipticiteplusmoins}).

The next step is a direct application of Proposition \ref{propfactorization}. Here again we  only consider  the procedure in the high frequency case but summarize both high and low frequencies parametrices in Theorem \ref{Isozaki-Kitada-explicite}. Given a symbol $ \chi_{{\rm st}}^{\pm} $ supported in strongly outgoing or incoming area (see (\ref{referencestrongarea}) in which we recall that $ p = p_{\kappa} $), we can factorize the corresponding pseudo-differential operator by mean of Proposition \ref{propfactorization}. More precisely, if $ I_2 \Subset I_1 $ and $ V_2 \Subset V_1 $ are given, then for $ R_2  $ large enough,  $ \varepsilon_2  $ small enough and all
$ \chi_{\rm st}^{\pm} \in \widetilde{S}^{-\infty,0} (\Ra^{2n}) $  supported in $ \widetilde{\Gamma}_{\rm st}^{\pm} (R_2,V_2,I_2,\varepsilon_2) $, one can find symbols $ b_k \in S_{n-1-k} $ supported in $ \Theta^{\pm} (R_2,V_2,I_1,C \varepsilon_2) $, such that $ b^h := b_0^{\pm} + \cdots + h^M b_M^{\pm} $ satisfies
$$ J^h_{\kappa} (a^h) J^h_{\kappa} (b^h)^{\dag} = O \! p_{\kappa}^h (\chi_{\rm st}^{\pm}) \tilde{\psi}_{\kappa} + h^M O \! p_{\kappa}^h ( \tilde{r}_M^h) \tilde{\psi}_{\kappa} $$
with  $ \tilde{r}_M^h \in \widetilde{S}^{-\infty,-M} (\Ra^{2n}) $, boundedly in $h$.  Using Proposition \ref{propositionpourestes}, this can also be written
\begin{eqnarray}
J^h_{\kappa} (a^h) J^h_{\kappa} (b^h)^{\dag} = O \! p_{\kappa}^h (\chi_{\rm st}^{\pm}) \tilde{\psi}_{\kappa} + O_{{\mathscr H}_{-M/2}^{-2M} \rightarrow {\mathscr H}_{M/2}^{2M}} ( h^M ) .
\end{eqnarray}



We  synthetize the analysis of this paragraph in the next theorem. Notice that, at low frequency, we consider the $ \epsilon $ dependent areas $ \widetilde{\Gamma}^{\pm}_{\epsilon,{\rm st}} (R,V,I,\varepsilon) $ introduced in (\ref{referencestrongarealow}).

\begin{theo}[Isozaki-Kitada parametrix] \label{Isozaki-Kitada-explicite} Let $ \kappa : U_{\kappa} \rightarrow V_{\kappa} $ be a chart of the atlas of Section \ref{sectionNotation} and $ V \Subset V_{\kappa} $ be a convex open subset. For all given 
$$ V_2  \Subset V_0 \Subset V \qquad \mbox{and} \qquad I_2  \Subset I_1 \Subset I_0 \Subset (0,+\infty) $$
and can choose $ C > 0 $, $ 0  < \varepsilon_1 < \varepsilon_0 $ and $  R_1 > R_0 $ such that for all $ N \geq 0 $ and all $ 0 < \varepsilon_2 \ll 1 $, $ R_2 \gg R_1 $, the following approximations hold.
\begin{enumerate}
\item{{\bf High frequency:}  there are $ a^h , a_{\rm c}^h , \check{a}^h \in S_0 (\Ra^{2n}) $ supported in $ \Theta^{\pm} (R_0,V_0,I_0,\varepsilon_0) $, satisfying
$$ \emph{supp} (a_{\rm c}^h) \subset \{ r \leq R_1 \} , \qquad \emph{supp} (\check{a}^h) \subset \{ |\theta - \vartheta| \geq \varepsilon_1 \}  $$
and $ r^h_N \in S_{-N} (\Ra^{2n}) $ also supported in $ \Theta^{\pm} (R_0,V_0,I_0,\varepsilon_0) $,  such that for all  $ \chi_{\rm st}^{\pm} \in \widetilde{S}^{-\infty,0} $ satisfying 
$$ \emph{supp} (\chi_{\rm st}^{\pm}) \subset \widetilde{\Gamma}_{\rm st}^{\pm} (R_2,V_2,I_2,\varepsilon_2) $$
one can find $ b^h \in S_0 (\Ra^{2n}) $ such that
$$ \emph{supp} (b^h) \subset \Theta^{\pm} (R_2 , V_2 , I_1 , C \varepsilon_2) $$
and
\begin{eqnarray}
 e^{-it P} O \! p_{\kappa}^h (\chi_{\rm st}^{\pm}) \tilde{\psi}_{\kappa} = J_{\kappa}^h (a^h) e^{-itD_x^2} J_{\kappa}^h (b^h)^{\dag} + R^h_N (t) \label{highfreqIKexpl} 
\end{eqnarray}
with
$$ R^h_N (t) = e^{-it P} O_{{\mathscr H}_{-N}^{2N} \rightarrow {\mathscr H}_{N}^{2N}} \big( h^{N} \big) - \frac{i}{h^2} \int_0^t e^{-i(t-\tau) P}  J^h_{\kappa} \big(a_{\rm c}^h + \check{a}^h + h^N r_N^h\big)  e^{-i \tau D_x^2}  J^h_{\kappa} (b^h)^{\dag} d \tau . $$}
\item{{\bf Low frequency:}  there are $ a_{\epsilon} , a_{\epsilon,{\rm c}} , \check{a}_{\epsilon} \in S_0 (\Ra^{2n}) $ supported in $ \Theta^{\pm} (R_0,V_0,I_0,\varepsilon_0) $, satisfying
$$ \emph{supp} (a_{\epsilon,{\rm c}}) \subset \{ r \leq R_1 \} , \qquad \emph{supp} (\check{a}_{\epsilon}) \subset \{ |\theta - \vartheta| \geq \varepsilon_1 \}  $$
and $ r_{\epsilon,N} \in S_{-N} (\Ra^{2n}) $ also supported in $  \Theta^{\pm} (R_0,V_0,I_0,\varepsilon_0)  $,  such that for all bounded family  $ (\chi_{\epsilon,{\rm st}}^{\pm} )_{\epsilon \in (0,1]} $ of $ \widetilde{S}^{-\infty,0} $ satisfying 
$$ \emph{supp} (\chi_{\epsilon,{\rm st}}^{\pm}) \subset \widetilde{\Gamma}_{\epsilon,{\rm st}}^{\pm} (R_2,V_2,I_2,\varepsilon_2) $$
one can find $ b_{\epsilon} \in S_0 (\Ra^{2n}) $ such that
$$ \emph{supp} (b_{\epsilon}) \subset \Theta^{\pm} (R_2 , V_2 , I_1 , C \varepsilon_2) $$
and
\begin{eqnarray}
 e^{-it P} O \! p_{\epsilon,\kappa} (\chi_{{\rm st},\epsilon}^{\pm}) \tilde{\psi}_{\kappa}(\epsilon r) = J_{\epsilon,\kappa} (a_{\epsilon}) e^{-i \epsilon^2 tD_x^2} J_{\epsilon,\kappa} (b_{\epsilon})^{\dag} + R_{\epsilon,N} (t) \label{lowfreqIKexpl}
\end{eqnarray} 
with
$$ R_{\epsilon,N} (t) = e^{-it P} O_{{\mathscr L}_{-N}^{-2N} \rightarrow {\mathscr L}_{N}^{2N}} \big( 1 \big) - i \int_0^{\epsilon^2 t} e^{-i( \epsilon^2 t-s) \frac{P}{\epsilon^2} }  J_{\epsilon,\kappa} \big(a_{\epsilon,{\rm c}} + \check{a}_{\epsilon} +  r_{\epsilon,N}\big)  e^{-i s D_x^2}  J_{\epsilon,\kappa} (b_{\epsilon})^{\dag} d s . $$}
\end{enumerate}
In both cases, the symbols are bounded uniformly in $h$ and $ \epsilon $ respectively.
\end{theo}

We point out that, so far, we have not justified to which extent the remainder terms in (\ref{highfreqIKexpl}) and (\ref{lowfreqIKexpl}) are small.  We will use Theorem \ref{Isozaki-Kitada-explicite} in subsection \ref{soussectionpropagation} to prove $ L^2 $ propagation estimates for $ e^{-it P} $
and  will see there that the remainders decay as $  \pm t \rightarrow \infty $.  In Section \ref{SectionStrichartz}, we will  use Theorem \ref{Isozaki-Kitada-explicite} in association with the (dual) propagation estimates of subsection \ref{soussectionpropagation} to control the remainders $ R_N^h(t) , R_{\epsilon,N}(t) $ in $ L^1 \rightarrow L^{\infty} $ norm. 

\medskip

\noindent {\bf Remark.} For future purposes, we record that, by using (\ref{expansioneikonal1}) and (\ref{expansioneikonal2}),  $ \varepsilon_0 $ and $ R_0 $ can be chosen respectively small and large enough in such a way (depending on $ V_0 $ and $ I_0 $) that we  have
\begin{eqnarray}
\frac{\partial_{r} \varphi}{\varrho} \in (1/2, 2 ) \qquad \mbox{and}  \qquad C^{-1}  r |\varrho| |\theta - \vartheta| \leq |\partial_{\theta} \varphi (r,\theta,\varrho,\vartheta)| \leq C  r |\varrho| |\theta - \vartheta| , \label{pourdesconditionsdesupport}
\end{eqnarray}
 on $ \Theta^{\pm} (R_0,V_0,I_0,\varepsilon_0) $. In particular, $ \partial_r \varphi $ and $ \varrho $ have the same strict sign.

\section{Propagation estimates} \label{sectionpropagation}
\setcounter{equation}{0}
\subsection{Finite time estimates} \label{paragraphetempsfini}


In this paragraph, we  prove propagation estimates over finite times depending on the spatial localization of the symbols and  on the frequency regime.

We introduce first some notation. We are going to work on $T^*( (R_{\mathcal M},\infty) \times {\mathcal S})  $ which is isomorhic to $   T^* (R_{\mathcal M},\infty) \times T^* {\mathcal S} $, so we will write its elements as $ (r , \rho , \varpi) $ with $ (r,\rho) \in (R_{\mathcal M},\infty) \times \Ra $ and $ \varpi \in T^* {\mathcal S} $. We then let  $ p_{\epsilon} = p_{\epsilon} (r,\rho,\varpi) $ be the  principal symbols of  $ - \Delta_{G_{\epsilon}} $ (see (\ref{Laplacienrescaleref})) which is intrinsically defined  on $ T^* \big( (R_{\mathcal M},\infty) \times {\mathcal S} \big) $.  We let  $ \phi^s_{\epsilon} $ be the associated Hamiltonian flow. Notice that, for $ \epsilon = 1 $, $ p_1 $ is the principal symbol of  $-\Delta_G $. Note also that the flow $ \phi^s_{\epsilon} $ is  not complete on $ T^*\big(  (R_{\mathcal M},\infty) \times {\mathcal S} \big) $. 
We then set
$$ \big( \bar{r}^s_{\epsilon} , \bar{\varrho}^s_{\epsilon} \big) := \mbox{component of} \ \phi^s_{\epsilon}(r,\rho,\varpi) \ \mbox{on} \ T^* (R_{\kappa},\infty) . $$
For $ R > R_{\mathcal M} $ and $ - 1 < \sigma < 1 $, we finally consider
\begin{eqnarray}
 \widetilde{\Gamma}^{\pm}_{\epsilon} (R,\sigma) = \{ (r,\rho,\varpi) \in T^*( (R_{\mathcal M},\infty) \times {\mathcal S})  \ | \ r > R , \ \pm \rho > \sigma p_{\epsilon}^{1/2}  \} . \label{notationpasterrible}
\end{eqnarray} 
It is an open conical subset of $ T^*( (R_{\mathcal M},\infty) \times {\mathcal S}) \setminus 0 $ (the strict inequality in (\ref{notationpasterrible}) prevents $ (\rho, \varpi ) $ from being $0$). We will sometimes need refinements of such areas, namely similar sets localized both on charts of $ {\mathcal S} $ and in energy; if $ \kappa : U_{\kappa} \rightarrow V_{\kappa} $ is a chart of the atlas chosen in Section \ref{sectionNotation}, $ V \Subset V_{\kappa} $ and $ I \Subset (0,+\infty) $, we set
\begin{eqnarray}
 \widetilde{\Gamma}^{\pm}_{\epsilon} (R,V,I,\sigma)_{\kappa} := \{ (r,\theta,\rho,\eta) \in \Ra^{2n} \ | \ r > R , \ \theta \in V, \ p_{\epsilon,\kappa} \in I , \ \pm \rho > \sigma p_{\epsilon,\kappa}^{1/2} \} , \label{zonesortantelargecoordonnees}
\end{eqnarray}
where we recall that $ p_{\epsilon,\kappa} $ is defined in (\ref{symbolesprincipaux}). We will call such regions outgoing (+)/ incoming (-) regions according to a classical terminology. Note the difference with the strongly outgoing/incoming regions defined in (\ref{referencestrongarea})-(\ref{referencestrongarealow})  in the case when $ \sigma = 1 - \varepsilon^2 $ is close to $1$.

We record first non angularly localized estimates on the flow.


\begin{prop} \label{propagationclassiqueprop} For all $ \sigma \in (-1,1) $, there exists $ R \gg 1 $ such that
\begin{enumerate}
\item{  there exists $ c > 0 $ such that, for all $ \epsilon \in (0,1] $,
$$ \bar{r}^s_{\epsilon} \geq c \big(r + |s| p^{1/2}_{\epsilon} \big) , \qquad \mbox{for all}  \  \pm s \geq 0 \ \  \mbox{and} \ \ (r,\rho,\varpi) \in   \widetilde{\Gamma}^{\pm}_{\epsilon} (R ,\sigma). $$
  In particular, $  R_1$ can be chosen such that $ \phi^s_{\epsilon} $ is defined on $  \widetilde{\Gamma}^{\pm}_{\epsilon} (R ,\sigma) $ for all $ \pm s \geq 0 $.
  }
\item{ For all $ 0 < \varepsilon < 1 $, there exists $ T > 0 $ such that, for all $ \epsilon \in (0,1] $,
$$ \pm \bar{\varrho}^s_{\epsilon}  > (1-\varepsilon^2) p^{1/2}_{\epsilon} \qquad \mbox{provided that} \ \  \pm s \geq T r p^{- 1/2}_{\epsilon}    \ \  \mbox{and} \ \ (r,\rho,\varpi) \in   \widetilde{\Gamma}^{\pm}_{\epsilon} (R ,\sigma).    $$
}
\item{Let $ 0 < \varepsilon < 1 $ and $ t_0 > 0$ (as small as we want). One can find  $ \delta > 0 $ such that, for all $ \epsilon > 0 $ and all $ (r,\rho,\varpi) \in \widetilde{\Gamma}^{\pm}_{\epsilon} (R ,\sigma) $, we have
$$  \frac{|\rho|}{p^{1/2}_{\epsilon}} < 1 - \varepsilon^2 \qquad \mbox{and} \qquad \pm s \geq t_0 r p^{-1/2}_{\epsilon} \qquad \Longrightarrow  \qquad \pm \frac{\bar{\varrho}^s_{\epsilon}}{p^{1/2}_{\epsilon}} > \pm  \left( \frac{\rho}{p^{1/2}_{\epsilon}} + \delta \right) .$$
}
\end{enumerate}
\end{prop}


\noindent {\bf Remark.} As in previous sections, we will give all proofs in the case $ \epsilon = 1 $ and, when there is a sign condition, for $ s \geq 0$. This only simplifies the notation.

\bigskip

\noindent {\it Proof.}  We choose first $ R_1 \gg 1 $ such that $ - \frac{\partial p}{\partial r} \geq r^{-1} (p - \rho^2) $ for $ r > R_1 $. This implies that, as long as $ \bar{r}^s > R_1 $
$$ \frac{d^2}{d s^2}  ( \bar{r}^s )^2 = 2 \frac{d}{ds} (\bar{r}^s \bar{\varrho}^s) \geq 4 (\bar{\varrho}^s)^2 + 2 (p - (\bar{\varrho}^s)^2) \geq 2 p , $$
hence that
$$  (\bar{r}^s )^2 \geq r^2 + 2 s r \rho + s^2 p \geq r^2 - 2 |s \sigma| r p^{1/2} + s^2 p \geq (1-|\sigma|)  \big( r^2 + s^2 p \big). $$
By a simple bootstrap argument, using the above argument, one can see that  $ \bar{r}^s > R_1 $  for all $  s \geq 0 $ provided that $  r  > (1-|\sigma|)^{-1/2} R_1 $. This completes the proof of the item 1. For the item 2, we observe that $ \bar{r}^s \leq r + 2 s p^{1/2} $ hence, by integrating
$$ \frac{\dot{\bar{\varrho}}^s}{p - (\bar{\varrho}^s)^2} \geq \frac{1}{\bar{r}^s} \geq \frac{1}{r+ 2 s p^{1/2}} , $$
we get
$$ \mbox{artanh} \left( \frac{\bar{\varrho}^s}{p^{1/2}} \right) \geq \mbox{artanh} \left( \frac{\rho}{p^{1/2}} \right) + \frac{1}{2}  \ln \left( 1 + \frac{2 s p^{1/2}}{r} \right)$$
and for $s p^{1/2}/r$ large enough the right hand side is greater than $ \mbox{artanh} (1 - \varepsilon^2) $, yielding the result.
For the item 3, we observe that $ \bar{\varrho}^s $ is non decreasing in $s$ so if the estimate holds at some time before $   t_0 r p^{-1/2} $ then it holds for all larger times. By possibly increasing $R_1 $, we may assume that, for $r> R_1$, we have $  |\frac{\partial p}{\partial r}|  \leq 4  r^{-1} (p - \rho^2)   $. Therefore, using that $ \bar{r}^s > R_1 $ by the item 1, we have
$$ |\bar{\varrho}^s - \rho| \leq 4 p s /r , $$
so by assuming $ s p^{1/2}/r $ small enough, we have $ |\bar{\varrho}^s / p^{1/2} | \leq 1 - \frac{\varepsilon^2}{2} $. Thus, for such times, the first inequality in the proof of the item 1 yields
$$ \dot{\bar{\varrho}}^s \geq  ( 1 - (1- \varepsilon^2 /2)^2 ) \frac{p}{\bar{r}^s} . $$
On the other hand, using once more that $  \bar{r}^s  \leq r +  2 p^{1/2}s $, the above inequality yields
$$ \dot{\bar{\varrho}}^s \geq (\varepsilon^2 - \varepsilon^4 /4) \frac{p}{r\left( 1 + 2 s p^{1/2}/r \right)} \geq  \varepsilon^2 \frac{p}{2 r} ,  $$
provided $ s p^{1/2}/r $ is small enough, say not greater than $  4 \delta  $, where $ \delta > 0 $ can be chosen smaller than $  \varepsilon^2 t_0 / 2 $. By integration over such times, we get $ \bar{\varrho}^s \geq \rho + s \varepsilon^2 \frac{p}{2r} $ which yields $ \bar{\varrho}^s / p^{1/2} > \rho / p^{1/2} + \delta  $ if $ \varepsilon^2 s p^{1/2}/r > 2 \delta  $ hence in particular if $ s p^{1/2}/r > t_0 $.  \finpreuve

\bigskip

In the next proposition, we record estimates on the geodesic flow in a coordinate patch. We consider a chart $ \kappa : U_{\kappa}   \rightarrow V_{\kappa} $ on $ {\mathcal S} $ from the atlas chosen in Section \ref{sectionNotation}. We recall that $ \phi^s_{\epsilon,\kappa} (r,\theta,\rho,\eta) $ is the flow of $ p_{\epsilon,\kappa} $ on $ (R_{\mathcal M},\infty) \times V_{\kappa} $, the components of which we denote as in (\ref{componentsofphikappa}).


\begin{prop} \label{estimeeflotlocal}  Let $ V \Subset V_{\kappa} $ and $ I \Subset (0, + \infty) $. There exists $ t_1 > 0 $ and $ R_{V,I} \gg 1 $ such that, 
\begin{enumerate}
\item{ for  all $ \epsilon \in (0,1] $, $ \phi^s_{\epsilon,\kappa}(r,\theta,\rho,\eta) $ is defined for $ |s| \leq t_1 r $ and $  (\bar{r}^s_{\epsilon} , \bar{\vartheta}^s_{\epsilon} ) $ belongs to $ (R_{\mathcal M},\infty) \times V_{\kappa} $, provided that
\begin{eqnarray}
 r > R_{V,I}, \qquad \theta \in V , \qquad p_{\epsilon,\kappa}(r,\theta,\rho,\eta) \in I . 
  \label{pourlasecondecondition}
\end{eqnarray}}
\item{ For all $ ( j , \alpha ,k , \beta ) \in \Za_+^{2n}  $, there exists $ C > 0 $ such that, uniformly in $\epsilon \in (0,1] $,
\begin{eqnarray}
\big| \partial_r^j \partial_{\theta}^{\alpha} \partial_{\rho}^k \partial_{\eta}^{\beta}  \big( \bar{\vartheta}^s_{\epsilon} - \theta , \bar{\varrho}^s_{\epsilon} - \rho \big) \big| \leq C r^{-j-|\beta|}  \frac{|s|}{r}  \label{flotGronwall1} \\ 
\big| \partial_r^j \partial_{\theta}^{\alpha} \partial_{\rho}^k \partial_{\eta}^{\beta}  \big( \bar{r}^s_{\epsilon} - r , \bar{\eta}^s_{\epsilon} - \eta \big) \big| \leq C r^{1-j-|\beta|}   \frac{|s|}{r}  \label{flotGronwall2}
\end{eqnarray}
for all initial data satisfying (\ref{pourlasecondecondition}) and all $ |s| \leq t_1 r $.}
 \end{enumerate}
\end{prop}

\noindent {\it Proof.}  See \cite{MizutaniCPDE}.  \finpreuve

\bigskip

In Theorem \ref{Egorovquantifie} below, we will propagate observables which do not remain localized in a single chart. To handle this fact, the following coordinate invariance property will be useful. 


\begin{prop}[Normalizing the angular supports] \label{propnormalizingangular} Let $ \kappa_1 : U_1 \rightarrow V_1 $ be a chart on $ {\mathcal S} $ of the atlas chosen in Section \ref{sectionNotation} and $ \tilde{\psi}_{\kappa_1} $ as in (\ref{notationpartition}).  Let $ (a_R)_{R \gg 1} $  be a bounded family in $ \widetilde{S}^{-\infty,0} $  such that, 
$$ \emph{supp}(a_R) \subset  (R , \infty) \times K \times \Ra^n, \qquad \mbox{for \underline{some}} \ \ K \Subset V_{\kappa_1}  . $$
Then, for all given $ N \geq 0  $, one can write
\begin{eqnarray}
    O \! p^h_{\kappa_1} (a_R)  \tilde{\psi}_{\kappa_1}  = 
\left(  \sum_{\kappa_2}    O \! p^h_{\kappa_2} \big(a_{R,\kappa_2}(h) \big)  \tilde{\psi}_{\kappa_2} \right)  + O_{{\mathscr H}_{-N}^{-2N} \rightarrow {\mathscr H}_N^{2N}} (h^N R^{-N})
  \nonumber
\end{eqnarray}
and
 \begin{eqnarray}
   O \! p_{\epsilon,\kappa_1} (a_R)  \tilde{\psi}_{\kappa_1}(\epsilon r)    = 
 \left( \sum_{\kappa_2}    O \! p_{\epsilon,\kappa_2} \big(a_{R,\kappa_2,\epsilon} \big)   \tilde{\psi}_{\kappa_2} (\epsilon r)  \right) +  O_{{\mathscr L}_{-N}^{-2N} \rightarrow {\mathscr L}_N^{2N}} ( R^{-N}) ,
  \nonumber 
\end{eqnarray}
where  $ ( a_{R,\kappa_2}(h))_{R,h}$ and $ (a_{R,\kappa_2,\epsilon})_{R, \epsilon} $ belong to bounded subsets of $ \widetilde{S}^{-\infty,0} $ and, using the notation (\ref{notationtransition}) and (\ref{notationpartition}), are supported in
\begin{eqnarray}
  \left\{  \big(r , \tau_{12} (\theta),\rho,  \big( d \tau_{12} (\theta)^T \big)^{-1} \eta \big)  \ | \ (r,\theta,\rho,\eta) \in \emph{supp}(a_R)  \right\}  \cap  [ R , \infty ) \times \emph{supp}(\varphi_{\kappa_2}) \times \Ra^n . 
  \label{conditiondesupport}
\end{eqnarray}  
  If  $a_R$ depends in a bounded way on additional parameters, then so do the symbols $ a_{R,\kappa_2}(h),a_{R,\kappa_2, \epsilon} $ and the remainder terms.
\end{prop}

The meaning of this proposition is twofold:  it says first the natural fact that a (possibly rescaled) pseudodifferential operator  with symbol supported in $ (R,\infty) \times K \times \Ra^n $ with a compact set $ K $ contained in $ V_{\kappa_1} $ but possibly larger than the support of the angular cutoff $ \varphi_{\kappa_1} $, can be written  as a sum of operators with symbols angularly localized in the support of $ \varphi_{\kappa_2} $. The second point, which is technically important, is the control of the remainder terms with respect to $R$. This will be  useful to prove Theorem \ref{Egorovquantifie} below.

\bigskip

\noindent {\it Proof.} For definiteness, we consider rescaled operators, the other case is similar. 
By introducing the partition of unity (\ref{partitionfixee}), which is equal to $1$ near the range of the operator since its symbol is supported in $ r \geq R \gg 1 $, we have
\begin{eqnarray}
 O \! p_{\epsilon,\kappa}(	a_R) \tilde{\psi}_{\kappa_1}(\epsilon r) & = & \sum_{\kappa_2} \psi_{\kappa_2}(\epsilon r)
  O \! p_{\epsilon,\kappa_1} (a_R)  \tilde{\psi}_{\kappa_1}(\epsilon r)  \nonumber
 \end{eqnarray}
where we keep only those $ \kappa_2 $ such that $ U_{\kappa_1} \cap U_{\kappa_2} \ne \emptyset $  otherwise the corresponding operator vanishes by the support properties of $ \psi_{\kappa_2} $ and $a_R$.  In each term of the right hand side
   we write $ \tilde{\psi}_{\kappa_1}  = \tilde{\psi}_{\kappa_1} \tilde{\psi}_{\kappa_2}  + \tilde{\psi}_{\kappa_1} (1 - \tilde{\psi}_{\kappa_2}) $. The terms involving $ 1 - \tilde{\psi}_{\kappa_2} $ are of the form
 \begin{eqnarray}
 \Pi_{\kappa_1} {\mathscr D}_{\epsilon} \left(  \psi_{\kappa_2} (\breve{r}, \kappa_1^{-1}(\theta))   O \! p^1 (a_R)  \tilde{\psi}_{\kappa_1}  (1 - \tilde{\psi}_{\kappa_2}) (\breve{r},\kappa_1^{-1}(\theta)) \right) {\mathscr D}_{\epsilon}^{-1} \Pi_{\kappa_1}^{-1} = O_{{\mathscr L}_{-N}^{-N} \rightarrow {\mathscr L}_N^N} (R^{-N}) .  \label{restearappelertransition}
  \end{eqnarray}
 Indeed, since $ 1 - \tilde{\psi}_{\kappa_2} $ vanishes near the support of $ \psi_{\kappa_2} $, the composition rules of Proposition \ref{bonpseudodiff} show that the parenthese is a pseudodifferential operator with symbol $ O (R^{-\infty}) $ in $ \widetilde{S}^{-\infty,-\infty} $ which in turns show it is as  in the right hand side (for any $N$) by the third item of Proposition \ref{propositionpourestes}. Next,  using the notation  (\ref{notationtransition}) for the transition maps, the terms $ \psi_{\kappa_2}(\epsilon r)
  O \! p_{\epsilon,\kappa_1} (a_R) ( \tilde{\psi}_{\kappa_1} \tilde{\psi}_{\kappa_2} )(\epsilon r) $ can be written
$$  
  \Pi_{\kappa_2} {\mathscr D}_{\epsilon}\left( \Pi_{\tau_{12}}^{-1}   \psi_{\kappa_2} (\breve{r}, \kappa_1^{-1}(\theta)) O \! p^1 (a_R)  (\tilde{\psi}_{\kappa_1}\tilde{\psi}_{\kappa_2} ) (\breve{r}, \kappa_1^{-1}(\theta)) \Pi_{\tau_{12}}
\right) {\mathscr D}_{\epsilon}^{-1} \Pi_{\kappa_2}^{-1}  $$
by using that $( \Pi_{\kappa_1} {\mathscr D}_{\epsilon} u ) (r)  = \Pi_{\kappa_2} {\mathscr D}_{\epsilon} \big( \Pi_{\tau_{12}}^{-1} u (\breve{r}) \big)  $. We  then use the third item of Proposition \ref{bonpseudodiff} to write, for any $N$, the parenthese as the sum of  an operator with symbol supported in (\ref{conditiondesupport}) and  a remainder term with symbol $ O (R^{-N}) $ in $ \widetilde{S}^{-\infty,-2N} $ which produces a remainder as in (\ref{restearappelertransition}). This completes the proof. \finpreuve


\bigskip
We are now ready to prove the main result of this paragraph.
  We refer to  (\ref{definitionnotationdiffeo}) for the notation which is used extensively below.  We also refer to (\ref{zonesortantelargecoordonnees}) for $ \widetilde{\Gamma}_{\epsilon}^{\pm} (R,V,I,\sigma)_{\kappa} $.

 In the following theorem, given  a chart  $ \kappa : U_{\kappa} \rightarrow V_{\kappa} $ on $ {\mathcal S} $ of the atlas of Section \ref{sectionNotation}, as all charts below, we let 
 $$ C_{\kappa} : (R_{\mathcal M} , \infty) \times V_{\kappa} \times \Ra^n \rightarrow T^* ((R_{\mathcal M},\infty) \times U_{\kappa} ) $$ 
 be the inverse of the chart on $ T^* ((R_{\mathcal M},\infty) \times U_{\kappa} )  $ associated to $ \kappa $, namely that is defined by $ C_{\kappa} (r,\theta,\rho,\eta) = \rho dr + \sum_j \eta_j d \theta_j \in T^*_{(r , \kappa^{-1}(\theta))}  ((R_{\mathcal M},\infty) \times U_{\kappa} )  $. Notice in particular that
 $$ \phi_{\epsilon}^s \circ C_{\kappa} = C_k \circ \phi_{\epsilon,\kappa}^s $$
 on all initial data and times such that $  \phi_{\epsilon,\kappa}^s (r,\theta,\rho,\eta) $ remains localized inside $ (R_{\mathcal M},\infty) \times V_{\kappa} \times \Ra^n $.


\begin{theo} \label{Egorovquantifie} Let $ I \Subset (0,\infty) $, $ \sigma \in (-1,1) $ and $ V_0 \Subset V_{\kappa_0} $ for some given chart $ \kappa_0 $. There exists $ R_0 \gg 1  $ such that for all given $ T > 0 $, all $ N \geq 0 $   and all bounded family  $ ( b_{\epsilon,R} ) $  of  $ \widetilde{S}^{-\infty,0} $ (indexed by $R \geq R_0 $ and $ \epsilon \in (0,1]$) and satisfying
\begin{eqnarray}
 \emph{supp}(b_{\epsilon,R})  \subset \widetilde{\Gamma}^{\pm}_{\epsilon}(R,V_0,I,\sigma)_{ \kappa_0}
 \label{zonesortantesansledire}
\end{eqnarray} 
the following properties hold:
\begin{enumerate}
\item{ {\bf High frequency propagation} ($ \epsilon = 1 $ and $ h \in (0,1] $) {\bf :} as long as
$$  R \geq R_0 , \qquad h \in (0,1] , \qquad  0 \leq \pm  \frac{t}{h} \leq T R  , $$
   one can write
$$ e^{-itP} \left( \psi_{\kappa_0}   O \! p^h_{\kappa_0} (b_{1,R})  \tilde{\psi}_{\kappa_0} \right) e^{itP} = \sum_{\kappa} \psi_{\kappa}   O \! p^h_{\kappa} \big(b_{R}(t,h)_{\kappa} \big) \tilde{\psi}_{\kappa}  + O_{{\mathscr H}^{-2N}_{-N} \rightarrow {\mathscr H}^{2N}_{N}} \left( h^N R^{-N} \right) $$
with $  ( b_{R}(t,h)_{\kappa} )_{R,t,h}$ bounded in $  \widetilde{S}^{-\infty,0} $ and such that
$$  C_{\kappa} \left( \emph{supp} \big( b_{R}(t,h)_{\kappa} \big) \right) \subset \phi^{ t h^{-1}}_1 \big( C_{\kappa_0} ( \emph{supp} ( b_{1,R} ) ) \big)      . $$}
\item{{\bf Low frequency propagation} ($ h = 1 $ and $ \epsilon \in (0,1] $) {\bf :} as long as
$$  R \geq R_0 , \qquad \epsilon \in (0,1] , \qquad  0 \leq \pm t \epsilon^2 \leq T R  , $$
one can write
\begin{eqnarray*}
 e^{-itP} \left(  \psi_{\kappa_0} (\epsilon r) O \! p_{\epsilon,\kappa_0} (b_{\epsilon,R})  \tilde{\psi}_{\kappa_0} (\epsilon r)  \right) e^{itP} & = &  \left(  \sum_{\kappa}  \psi_{\kappa} (\epsilon r)  O \! p_{\epsilon,\kappa} \big(b_{\epsilon,R} (t)_{\kappa} \big) \tilde{\psi}_{\kappa} (\epsilon r) \right) \\ & & \ \ \ + \ O_{{\mathscr L}^{-2 N}_{-N} \rightarrow {\mathscr L}^{2 N}_{N}} \left(  R^{-N} \right) 
\end{eqnarray*} 
with $ ( b_{\epsilon,R} (t)_{\kappa} )_{\epsilon,R, t}  $ bounded in $ \widetilde{S}^{-\infty,0} $  and such that
$$ C_{\kappa}  \big( \emph{supp} \big( b_{\epsilon,R} (t)_{\kappa} \big) \big) \subset     \phi^{t \epsilon^2  }_{\epsilon} \big( C_{\kappa_0}( \emph{supp} ( b_{\epsilon,R})) \big) . $$
}
\end{enumerate}
\end{theo}

This is a quantitative version of the Egorov theorem. Its interests are to quantify (in terms of $R$) the range of times on which it holds, to  estimate  the remainder terms in suitable topologies and to include  a rescaled/low frequency version which is not  completely standard.


\bigskip

\noindent {\it Proof of Theorem \ref{Egorovquantifie}.} For definiteness, we consider the high frequency outgoing case (for which the notation is lighter since there is no $ \epsilon $ parameter). We use the general formula,
\begin{eqnarray}
 e^{-it P} A (0) e^{itP} = A (t) - \int_0^t e^{-i(t-\tau)P} \left( A^{\prime}(\tau) + i \big[ P , A (\tau) \big] \right) e^{i (t-\tau)P} d \tau  . \label{inthecasewemodify}
\end{eqnarray}
Choose  $ t_1 $ as in Proposition \ref{estimeeflotlocal} and consider first $ 0 \leq  s \leq t_1 R $ so that the flow remains localized in a single chart. We seek $ B (s) = A (sh) $, or equivalently $ A (t) = B (t/h) $, of the form
$$ B (0) = \psi_{\kappa_0}   O \! p^h_{\kappa_0} (  b_{R}) \tilde{\psi}_{\kappa_0}, \qquad B (s) = \sum_{j=0}^{J(N)} h^j O \! p^h_{\kappa_0} (b_j(s)) \tilde{\psi}_{\kappa_0} =: \Psi_N (s) \tilde{\psi}_{\kappa_0} , $$
for some $s$ dependent symbols $ b_j (\cdot) $ and some large enough order $ J (N) $ to be chosen.  Here and below we set $ b_R = b_{\epsilon,R} $ for $ \epsilon = 1 $. A simple calculation yields
\begin{eqnarray}
 h B^{\prime}(s) +h^2 [P,B(s)] = \left( h\Psi_N^{\prime} (s) + i \big[h^2 P , \Psi_N (s) \big]  \right) \tilde{\psi}_{\kappa_0} + i \Psi_N (s) [h^2 P,\tilde{\psi}_{\kappa_0}] . \label{pseudohorsdiagonale}
\end{eqnarray}
According to the usual procedure, we try first to make the first parentheses in the right hand side small. This is obtained
by constructing iteratively  the symbols $b_j $ as solutions to
\begin{eqnarray}
b_0^{\prime}(s) +  \{ b_0 (s) , p_{\kappa_0} \} & = & 0 , \qquad \ \ \ \ \ b_0 (0) = \psi_{\kappa_0}\big(r, \kappa_0^{-1}(\theta) \big) b_{R}  , \\
b_j^{\prime}(s) + \{ b_j(s) , p_{\kappa_0} \} & = & f_j(s), \qquad b_j (0) = 0 ,
 \label{bysolving}
\end{eqnarray}
with
$$ f_j (s) = - \sum_{j^{\prime} + k + l = j + 1 \atop
j^{\prime} < j } (  p_{\kappa_0,k} \# b_{j^{\prime}}(s) )_l - (   b_{j^{\prime}}(s) \# p_{\kappa_0,k} )_l , $$
where $  p_{\kappa_0, 0} = p_{\kappa_0} $ is the principal symbol of $ P $ in the chart associated to $ \kappa_0 $ and $ p_{\kappa_0,0} +p_{\kappa_0,1} $ is its full symbol, and where $ \{  a , b \} = \partial_x a \cdot \partial_{\xi} b - \partial_{\xi} a \cdot \partial_x b$ is the Poisson bracket. The solutions are given by
\begin{eqnarray}
 b_0 (s) = b_{R} \circ \phi^{-s}_{1,\kappa_0} , \qquad b_j (s) = \int_0^s f_j (u, \phi^{u-s}_{1 , \kappa_0}) d u ,  \label{transportEgorov}
\end{eqnarray}
with $ \phi_{1,\kappa_0}^s $  the Hamiltonian flow of $ p_{\kappa_0} $. According to the estimates (\ref{flotGronwall1}) and (\ref{flotGronwall2}), 
the formulas in (\ref{transportEgorov})  define symbols
$$ b_{j,R} (s) := b_{j}(s) \ \ \ \mbox{bounded in} \  \ \widetilde{S}^{-\infty, -j }(\Ra^{2n}) \ \ \mbox{for} \ \ R \geq R_0 , \ |s| \leq t_1 R . $$
Moreover, by choosing a relatively compact open subset $ K_0 \Subset V_{\kappa_0} $ and  $ R_1 \gg 1 $
so that
\begin{eqnarray}
   V_0  \Subset K_0, \qquad \tilde{\psi}_{\kappa_0} (r,\kappa_0^{-1}(\theta)) \equiv 1 \ \ \mbox{near} \ \  (R_1,\infty) \times K_0 , \label{pourlesurcutoff}
\end{eqnarray}     
we can ensure, by possibly taking a smaller $ t_1 > 0 $, that for  $ R \gg 1 $ and $ 0 \leq s \leq t_1 R $
$$ \mbox{supp}(b_{j,R}(s)) 
 \subset (R_1 , \infty) \times K_0 \times \Ra^n . $$
Hence, by the last condition in (\ref{pourlesurcutoff}),  $ \Psi_N (s)$ and $ [h^2 P,\tilde{\psi}_{\kappa_0}] $ in (\ref{pseudohorsdiagonale})  have disjoint supports. More precisely, the support of $ b_{j,R}(s) $ is contained in $\{ r \gtrsim R \} $ so the symbol of  $ \Psi_N (s) [h^2 P,\tilde{\psi}_{\kappa_0}]  $ is $ O (h^{\infty} R^{-\infty}) $ in $ \widetilde{S}^{-\infty,-\infty} $ which implies, by the third item of Proposition \ref{propositionpourestes}, that
$$ \Psi_N (s) [h^2 P,\tilde{\psi}_{\kappa_0}]  =  O_{{\mathscr H}^{-2N-2 \gamma ( \lceil N \rceil )  }_{-N} \rightarrow {\mathscr H}^{2N + 2 \gamma ( \lceil N \rceil ) }_{N}} \left( h^N R^{-N -2 \gamma ( \lceil N \rceil )} \right) . $$
Here $ \gamma ( \lceil N \rceil )  $ is as in Proposition \ref{proproughprop} (we will see the interest of this choice below).
  On the other hand, the construction of the $b_j (s) $ ensures that, for some $ \tilde{b}_{N,R} (s) $ bounded in $ \widetilde{S}^{-\infty,-J(N)} $ and supported in $ \{r \gtrsim R \}$,
\begin{eqnarray*}
   \left( h\Psi_N^{\prime} (s) + i \big[h^2 P , \Psi_N (s) \big]  \right) \tilde{\psi}_{\kappa_0} & = &  h^{J(N)} O \! p^h_{\kappa_0} (\tilde{b}_{N,R} (s)) \tilde{\psi}_{\kappa_0}  \\
   & = & O_{{\mathscr H}^{-2N-2 \gamma ( \lceil N \rceil )  }_{-N} \rightarrow {\mathscr H}^{2N + 2 \gamma ( \lceil N \rceil ) }_{N}} \left( h^N R^{-N -2 \gamma ( \lceil N \rceil )} \right)
 \end{eqnarray*}  
by choosing $ J (N) $ large enough and by using again the third item of Proposition \ref{propositionpourestes}. 
The interest of going to the order $\pm 2(N +  \gamma ( \lceil N \rceil ) ) $ in the remainder terms is that, by Proposition  \ref{proproughprop},
$$ e^{i (t-\tau) P} = O_{ {\mathscr H}^{-2N  }_{-N}    \rightarrow {\mathscr H}^{-2N-2 \gamma ( \lceil N \rceil )  }_{-N}  } (R^{ \gamma ( \lceil N \rceil )})  \qquad \mbox{for times} \  \ |t- \tau| \leq t_1 h R , $$
and similarly from $  {\mathscr H}^{2N + 2 \gamma ( \lceil N \rceil ) }_{N}  $ to $  {\mathscr H}^{2N }_{N} $. This allows to take into account the conjugation by propagators in the integral of  (\ref{inthecasewemodify}) and get, for our choice of $ A (t) = B (t/h) $,
$$  e^{-it P} A (0) e^{itP}  -  A (t) =  O_{{\mathscr H}^{-N}_{-N} \rightarrow {\mathscr H}^{N}_{N}} \left( h^{N-1} R^{1-N} \right) . $$
Here $N $ is arbitrary so getting $ h^{N-1}R^{1-N} $  rather than $ h^N R^{-N} $ is of course harmless.
Furthermore, one can rewrite $ A (t)  $ as a sum of $ \psi_{\kappa} O \! p_{\kappa}^h (b_R(t,h)_{\kappa} ) \tilde{\psi}_{\kappa} $ by mean of Proposition \ref{propnormalizingangular}, which yields the result for $ |t| \leq t_1 h R $. Then, by iterating this procedure a finite number ($ \approx O (T/t_1)$) of times, we get the result (note that along such an iteration, the symbols remain supported in $ r \gtrsim  R + |t/h| $ by Proposition \ref{propagationclassiqueprop}). 

The proof at low frequency is similar up to the replacement of pseudodifferential operators by rescaled ones and to the different time scaling $ s = \epsilon^2 t $. \finpreuve
 
 



\subsection{Resolvent estimates and their consequences} \label{soussectionresolvante}
In this short paragraph, we record some a priori decay estimates for $ e^{-itP} $ in weighted spaces, obtained as direct consequences of resolvent estimates. We consider both high and low frequency spectral localizations.
 
We recall first first well known consequences of the following Stone formula
$$  f (H) e^{-itH}  = \frac{1}{2 i \pi} \int_{\Ra} e^{-i t \lambda} f (\lambda) \left( (H - \lambda - i 0)^{-1} - (H- \lambda + i 0)^{-1} \right) d \lambda $$
valid for an arbitrary  self-adjoint operator $H$ and $ f \in C_0^{\infty} (\Ra) $. By integrations by part in $ \lambda $ together with the fact that 
$$ \partial_{\lambda}^k (H- \lambda \mp i 0)^{-1} = k! (H - \lambda \mp i 0)^{-1-k} $$
it allows to convert estimates on powers of the resolvent into time decay estimates for $ f (H) e^{-itH} $.

 Everywhere below, we let $ I \Subset (0,+\infty) $ and $ f \in C_0^{\infty} (I) $.  
  We consider first low frequency estimates for $ P $. Using the resolvent estimates of \cite{BoucletRoyer} namely
$$ \left| \left| \scal{\epsilon r}^{-k} \left(  \epsilon^{-2} P  - \lambda \pm i 0 \right)^{-k} \scal{\epsilon r}^{-k} \right| \right|_{L^2 \rightarrow L^2} \leq C_k , \qquad \lambda \in I , \ \epsilon \in (0,1] ,  $$
we obtain from the Stone formula, applied to $ H = P / \epsilon^2 $ and $ t$  replaced by $ \epsilon^2 t $, that for any $k \in \Na$
\begin{eqnarray}
 \left| \left| \scal{\epsilon r}^{-1-k}  f (P/\epsilon^2) e^{-itP}  \scal{ \epsilon r}^{-1-k} \right| \right|_{L^2 \rightarrow L^2} \lesssim \scal{\epsilon^2 t}^{-k}, \qquad t \in \Ra , \ \epsilon \in (0,1] .  \label{aprioriblackboxlowfrequency}
\end{eqnarray}
Another estimate from \cite{BoucletRoyer} that will be particularly useful is
\begin{eqnarray}
 \left| \left| \scal{ r}^{-1} \left(   P  - \lambda \pm i 0 \right)^{-1} \scal{ r}^{-1} \right| \right|_{L^2 \rightarrow L^2} \leq C , \qquad \lambda \in (0,1) ,  \label{remplacerrefvide0}
\end{eqnarray}
for it implies (see e.g. \cite[Thm XIII.25]{RS4}) that
\begin{eqnarray}
\left( \int_{\Ra } \left| \left| \scal{ r}^{-1} e^{-itP} f (P/\epsilon^2) u_0 \right| \right|_{L^2}^2 dt \right)^{1/2} \leq C || u_0 ||_{L^2} , \qquad \epsilon \in (0,1), \ u_0 \in L^2 . \label{remplacerefvide}
\end{eqnarray}

\bigskip

Getting similar estimates at high frequency, with polynomial growth in $ 1/h $, requires an assumption, for instance a non trapping condition. This is where the assumption (\ref{aprioriresolvente}) is useful since it allows to prove the following proposition.
\begin{prop}[Semiclassical power resolvent estimates]  \label{propositionpuissanceresolvante} Assume (\ref{aprioriresolvente}). Then for all $ k \geq 0 $
there exists $ N_k $ such that
\begin{eqnarray}
 \left| \left| \scal{ r}^{-1-k} \left( h^{2} P  - \lambda \pm i 0 \right)^{-1-k} \scal{ r}^{-1-k} \right| \right|_{L^2 \rightarrow L^2} \lesssim h^{-N_k} , \qquad \lambda \in I , \ h \in (0,1] . \label{resolvantsemiclassiquederivee}
\end{eqnarray}
\end{prop}

\noindent {\it Proof.}
  It is based on an argument in \cite[Prop. 1.3]{Vodev}. It consists in finding an operator $ P_0 $ defined on $ (0,+\infty) \times {\mathcal S} $ coinciding with $ P $ near infinity and satisfying nice resolvent estimates (as (\ref{polynomialgrowthinhminusone}) below) and then to use iterations of the resolvent identity.
  We explain schematically how to implement it in our context. We let  $ |D_{\mathcal S}| = (-\Delta_{\bar{g}})^{1/2} $ be the square root of the asymptotic Laplacian on $ {\mathcal S} $ and 
$$ h^j L_{j} , h^m L^{\prime}_m \in \{ 1 , h \partial_r , r^{-1} h |D_{\mathcal S}| \} $$
where $ j , m \in \{ 0 , 1 \} $ are the orders of the operators.
%
  Proceeding as in  \cite{BoucletRoyer}, one can find a second order differential operator $ P_0 $ on $ (0,+\infty) \times {\mathcal S} $ which is close everywhere to exact conical Laplacian $ - \partial_r^2 - r^{-2} \Delta_{\bar{g}} $ and equal to $ P $ near infinity in such a way that, letting $ P_{0,h} $ be the rescaled version of $ P_0 $, namely
$$ P_{0,h} = {\mathscr D}_{h} \big( h^2 P_0 \big) {\mathscr D}_{h}^{-1} , $$
we have, for any $ K \Subset \Ca \setminus 0 $ and $ k \in \Na $,
\begin{eqnarray}
 \left| \left|  \scal{ r}^{-k} L_j ( P_{0,h} - z )^{-k} L^{\prime}_m \scal{ r}^{-k}  \right| \right| \leq C, \qquad h \in (0,1],  \ z \in K \setminus \Ra ,  \label{estimeesemiclassiquerescalee}
 \end{eqnarray}
where for simplicity, $ || \cdot || $ is the operator norm on $ L^2 \big( (0,\infty) \times {\mathcal S} , r^{n-1} d r d {\rm vol}_{\bar{g}} \big) $. Such resolvent estimates follow from the techniques of  \cite{BoucletRoyer} (more precisely Proposition 3.13 and Lemma 4.2 there) which are based on a rescaling argument; they were used to prove low frequency estimates but work equally well at high frequency (one only uses that $ P_{0,h} $ is close to $  - \partial_r^2 - r^{-2} \Delta_{\bar{g}}  $ which satisfies a global positive commutator estimate at energy $1$).
%
Then, by unitarity of $ {\mathscr D}_{h}^{\pm 1} $ and (\ref{estimeesemiclassiquerescalee}), we find
\begin{eqnarray}
  \left| \left| \scal{r}^{-k} h^j L_j (h^2 P_0 - z )^{-k} h^j L^{\prime}_m \scal{r}^{-k} \right| \right| & = &  \left| \left| {\mathscr D}_{h}^{-1} \scal{h r}^{-k} L_j ( P_{0,h} - z )^{-k} L_m^{\prime} \scal{h r}^{-k} {\mathscr D}_{h}\right| \right| \nonumber \\
  & = &  \left| \left|  \scal{h r}^{-k} L_j ( P_{0,h} - z )^{-k} L^{\prime}_m \scal{h r}^{-k}  \right| \right| \nonumber \\
  & \lesssim & h^{-2 k}  \left| \left|  \scal{ r}^{-k} L_j ( P_{0,h} - z )^{-k} L^{\prime}_m \scal{ r}^{-k}  \right| \right| \nonumber
  \\ & \lesssim & h^{-2k} . \label{polynomialgrowthinhminusone}
\end{eqnarray}  
To illustrate the starting point of the method of \cite{Vodev}, we  check rapidly (\ref{resolvantsemiclassiquederivee}) for $k = 0 $, more precisely that
\begin{eqnarray}
  \left| \left| \scal{r}^{-1}  ( h^2 P - z )^{-1} \scal{r}^{-1} \right| \right| _{L^2 ({\mathcal M}) \rightarrow L^2 ({\mathcal M})} \lesssim h^{- 2- M} ,  \label{premierepolynomiale}
 \end{eqnarray} 
whose interest is to replace the compactly supported cutoffs $ \chi $ in (\ref{aprioriresolvente}) by the weight $ \scal{r}^{-1} $. By using the cutoffs $ \zeta , \tilde{\zeta} $ introduced in Section \ref{sectionNotation}  which are equal to $1$ near infinity and using the following resolvent identity
\begin{eqnarray}
 \zeta (r) (h^2 P - z)^{-1} = \zeta (r) (h^2 P_0 - z)^{-1}  \tilde{\zeta}(r)  - \zeta (r) (h^2 P_0 - z)^{-1}  \big[ \tilde{\zeta}(r) , h^2 P \big] (h^2 P - z)^{-1}  \label{resolvanteaiterer}
\end{eqnarray} 
together with (\ref{aprioriresolvente}) and (\ref{polynomialgrowthinhminusone}), we find that for any  $\chi \in C_c^{\infty}({\mathcal M} ) $
$$  \left| \left| \scal{r}^{-1}  ( h^2 P - z )^{-1} \chi \right| \right| _{L^2 ({\mathcal M}) \rightarrow L^2 ({\mathcal M})} \lesssim h^{- 2} + h^{-1 - M} .  $$
By using a second time the resolvent identity (\ref{resolvanteaiterer}) and using the above estimate, we obtain (\ref{premierepolynomiale}).
This leads to (\ref{resolvantsemiclassiquederivee}) for  $ k = 0 $. We get the result for higher $k$ by the same induction as in  \cite{Vodev}.
\finpreuve


\bigskip 

Using again the Stone formula with $ H = h^2 P $ and $ t $ replaced by $ t /h^2 $, Proposition \ref{propositionpuissanceresolvante} yields automatically
\begin{eqnarray}
 \left| \left| \scal{ r}^{-1-k}  f (h^2 P) e^{-itP}  \scal{ r}^{-1-k} \right| \right|_{L^2 \rightarrow L^2} \lesssim h^{-N_k} \scal{ t / h^2}^{-k}, \qquad t \in \Ra , \ h \in (0,1] , \nonumber
\end{eqnarray}
which in turn provides the weaker estimate
\begin{eqnarray}
 \left| \left| \scal{ r}^{-1-k}  f (h^2 P) e^{-itP}  \scal{ r}^{-1-k} \right| \right|_{L^2 \rightarrow L^2} \lesssim h^{-N_k} \scal{ t / h}^{-k}, \qquad t \in \Ra , \ h \in (0,1] ,  \label{aprioriblackboxhighfrequency}
\end{eqnarray}
which we record under this form to follow the natural semiclassical time scaling. Similarly to the estimate (\ref{remplacerefvide}), we also have the following consequence of (\ref{resolvantsemiclassiquederivee}) for $ k = 0 $,
\begin{eqnarray}
\left( \int_{\Ra } \left| \left| \scal{ r}^{-1} e^{-itP} f (h^2 P) u_0 \right| \right|_{L^2}^2 dt \right)^{1/2} \leq C h^{1-\frac{N_0}{2}} || u_0 ||_{L^2} , \qquad \epsilon \in (0,1), \ u_0 \in L^2 . \label{remplacerefvide2}
\end{eqnarray}
We recall that when the manifold is non trapping, one can take $ N_0 = 1 $, and the resulting $ h^{1/2} $ factor on the right hand side corresponds to the  $ H^{1/2} $ smoothing effect of the Schr\"odinger equation.

\subsection{Long time estimates} \label{soussectionpropagation}

In this paragraph, we prove several  $ L^2 $ propagation estimates on $ e^{-itP} $.   They will be used in Section \ref{SectionStrichartz} to control the remainder terms of the parametrices. However, their interest go beyond the applications to Strichartz inequalities. They generalize well known estimates (see e.g. \cite{Mourre,IsKi1}) in two ways: on one hand we consider the general geometric framework asymptotically conical manifolds and on the other hand we include a low frequency version of such inequalities which, to our knowledge, is an original result.

 We start with the following result on strongly outgoing/incoming microlocalizations (see (\ref{referencestrongarea}) and (\ref{referencestrongarealow}) for the related areas). This is a first application of Theorem \ref{Isozaki-Kitada-explicite}.

\begin{prop} \label{cellepourfortementsortant} Let $ k  \in \Na $, $ f \in C_0^{\infty} (0,+\infty) $, $ I_2 \Subset (0,+\infty) $ and $ V_2 \Subset V_{\kappa} $. Then, if $ R_2 \gg 1 $ and $ 0 <  \varepsilon_2 \ll 1 $,  we have the following estimates:
\begin{enumerate}
\item{{\bf High frequency:} Assume (\ref{aprioriresolvente}). If $ \chi_{\rm st}^{\pm} \in \widetilde{S}^{-\infty,0} (\Ra^{2n}) $ is supported in $ \widetilde{\Gamma}_{\rm st}^{\pm}(R_2,V_2,I_2,\varepsilon_2) $,
$$ \big| \big| \scal{r}^{-3k} e^{-i t P} f (h^2 P)  O \! p^h_{\kappa} (\chi_{\rm st}^{\pm}) \tilde{\psi}_{\kappa}  \scal{r}^{2k} \big| \big|_{L^2 \rightarrow L^2} \lesssim \scal{t/h}^{-k}, \qquad \pm t \geq 0 , \qquad h \in (0,1 ] . $$}
\item{{\bf Low frequency:} if $ (\chi_{\epsilon,{\rm st}}^{\pm})_{\epsilon} $ is a bounded family of $ \widetilde{S}^{-\infty,0} $ supported in $ \widetilde{\Gamma}_{{\rm st},\epsilon}^{\pm} (R_2,V_2,I_2,\varepsilon_2) $,
$$ \big| \big| \scal{\epsilon r}^{-3k} e^{-i t P} f ( P/ \epsilon^2)  O \! p_{\epsilon,\kappa} (\chi_{\epsilon,{\rm st}}^{\pm})  \tilde{\psi}_{\kappa}(\epsilon r ) \scal{\epsilon r}^{2k}  \big| \big|_{L^2 \rightarrow L^2} \lesssim \scal{\epsilon^2 t}^{-k}, \qquad \pm t \geq 0 , \qquad \epsilon \in (0,1 ] . $$ }
\end{enumerate}
\end{prop}

We point out that, in the high frequency estimate, we don't have any loss in $h$, as the $ h^{-N_k} $ in (\ref{aprioriblackboxhighfrequency}).

\bigskip

\noindent {\it Proof of Proposition \ref{cellepourfortementsortant}.} We may assume $ k \geq 1$. For definiteness, we consider the outgoing high frequency case. We use the notation of Theorem \ref{Isozaki-Kitada-explicite}, in particular (\ref{highfreqIKexpl}).  Note that, up to possibly decomposing $ \chi_{\rm st}^+ $ as a sum of symbols supported in balls with respect to $ \theta $, we assume that $ V_2 \Subset V $ for some convex open subset $ V \Subset V_{\kappa} $.
The contribution of the main term of the Isozaki-Kitada parametrix is
$$ \left( \scal{r}^{-3k} f (h^2 P) \scal{r}^{3k} \right) \scal{r}^{-3k} J^h_{\kappa} (a^h) e^{- i \tau D_x^2} J^h_{\kappa} (b^h)^{\dag} \scal{r}^{2k}. $$
Here the parenthese is a bounded operator on $ L^2 $ according to Theorem \ref{theoremcaclulfonctionnel} while the second factor decays as $ \scal{t/h}^{-k} $ by  Proposition \ref{propparam}. Next, the first term of the remainder $ R_N^h (t) $ of (\ref{highfreqIKexpl}) produces a term of the form
$$ \scal{r}^{-3k} f (h^2 P) e^{-itP} O_{{\mathscr H}_{-N}^{-2N} \rightarrow {\mathscr H}_N^{2N}} (h^N) \scal{r}^{2k}  $$
which is $ O \big( \scal{t/h}^{1-2k} h^{N-N_k} \big) $ in $ L^2  $ operator norm if $ N \geq 2k $ by (\ref{aprioriblackboxhighfrequency}) since one can write
\begin{eqnarray}
 O_{{\mathscr H}_{-N}^{-2N} \rightarrow {\mathscr H}_N^{2N}} (h^N) = \scal{r}^{-N} O_{L^2 \rightarrow L^2} (h^N) \scal{r}^{-N} . \label{utileaexpliciter}
\end{eqnarray} 
  By possibly increasing $N$ so that $ N \geq N_k $, we get an estimate by $ \scal{t/h}^{-k} $ (since $ 2 k -1 \geq k $). In the integral term of $ R_N^h(t) $, we consider first the contribution of $ J^h_{\kappa} (h^N r_N) $. By choosing  $ N $ large enough ($ N \geq 6k +1 $ and $ N \geq N_{k} $), Proposition \ref{propparam} and (\ref{aprioriblackboxhighfrequency}) imply that
 \begin{eqnarray} 
 \left| \left| \scal{r}^{-3k} f (h^2 P) e^{-i (t -\tau) P} J_{\kappa}^h (h^N r_{N}) e^{-i\tau D_x^2} J_{\kappa}^h (b^h)^{\dag} \scal{r}^{2k}\right| \right|_{L^2 \rightarrow L^2} \lesssim \scal{(t-\tau)/h}^{1-3k} \scal{\tau/h}^{-k-1} . 
 \nonumber
\end{eqnarray} 
After integration in $ \tau  $ between $0$ and $t$, we get an estimate by $ \scal{t/h}^{-k} $. It then remains to study the contributions of $ a_{\rm c}^h $ and $ \check{a}^h $. They follow as the one of $ h^N r_N $ once observed that we have the following estimates. By assuming $ R_2 $ large enough, the first item of Proposition \ref{nonstatprop} allows to write, for all $N$,
\begin{eqnarray}
 J^h_{\kappa} (a_{\rm c}^h) e^{- i \tau D_x^2} J^h_{\kappa} (b^h)^{\dag} =  O_{{\mathscr H}_{-N}^{-2N} \rightarrow {\mathscr H}_{N}^{2N}} ( h^{N} \scal{\tau/h}^{-N} ) , \qquad \pm \tau \geq 0 ,  \label{2emeN}
\end{eqnarray} 
since one has $ r \ll r^{\prime} $ on the support of the kernel of $  J^h (a_{\rm c}^h) e^{- i \tau D_x^2} J^h (b^h) $.  Using the second item of Proposition \ref{nonstatprop} and choosing $ \varepsilon_2 $ small enough (hence ensuring that $ |\theta-\vartheta| \gtrsim 1 $ and $ |\theta^{\prime} - \vartheta| \ll 1 $ on the support of the Schwartz kernel of $ J^h_{\kappa} (\check{a}^h) e^{- i \tau D_x^2} J^h_{\kappa} (b^h)^{\dag}  $), we obtain similarly.
\begin{eqnarray}
 J^h_{\kappa} (\check{a}^h) e^{- i \tau D_x^2} J^h_{\kappa} (b^h)^{\dag} =  O_{{\mathscr H}_{-N}^{-2N} \rightarrow {\mathscr H}_{N}^{2N}} ( h^{N} \scal{\tau/h}^{-N} ) , \qquad \pm \tau \geq 0 .  \label{3emeN}
\end{eqnarray} 
Using (\ref{utileaexpliciter}) with $ h^N \scal{\tau/h}^{-N}  $ instead of $h^N$,  we have the required spatial decay to use (\ref{aprioriblackboxhighfrequency}) and to control the growing weight $ \scal{r}^{2k} $. This completes the proof at high frequency. The proof is completely similar at low frequency by using (\ref{aprioriblackboxlowfrequency}) instead of (\ref{aprioriblackboxhighfrequency}). \finpreuve


\bigskip

In the next result, we partially relax the assumptions of Proposition \ref{cellepourfortementsortant} by replacing {\it strongly} outgoing (or incoming) microlocalizations by general outgoing (or incoming) ones, but at the expense of a stronger weight (which will eventually be harmless).
In the sequel, we denote
\begin{eqnarray}
\widetilde{\Gamma}^{\pm} (R,V,I,\sigma) & = & \{ (r,\theta,\rho,\eta) \ | \ r > R , \ \theta \in V, \ p_{\kappa} \in I, \ \pm \rho > \sigma p^{1/2}_{\kappa} \} \nonumber  \\
\widetilde{\Gamma}^{\pm}_{\epsilon} (R,V,I,\sigma) & = & \{ (r,\theta,\rho,\eta) \ | \ r > R , \ \theta \in V, \ p_{\epsilon,\kappa} \in I, \ \pm \rho > \sigma p^{1/2}_{\epsilon,\kappa} \} . \label{weaklowarea}
\end{eqnarray}
These regions  correspond to (\ref{zonesortantelargecoordonnees}) but we  now drop the index $ \kappa $ (unless it is necessary, {\it i.e.} in Proposition \ref{improvedoutgoingII}) and distinguish between the high and low frequency cases.
We recall that the difference with strongly outgoing/incoming regions considered in Proposition \ref{cellepourfortementsortant}  is that $ \sigma $ can be any real number $(-1,1) $, while  $ \sigma = 1 - \varepsilon^2 $ was close to $1$ in the previous proposition.

\bigskip

\begin{prop}[Half microlocalized propagation estimates] \label{halflocalized} Let $ k \in \Na $, $ I_2 \Subset (0,+\infty) $, $ V_2 \Subset V_{\kappa} $ and $ \sigma \in (-1,1) $. Then, if $ R_2 \gg 1 $, we have the following estimates:
%
\begin{enumerate}
\item{ {\bf High frequency estimates:} if $  \chi_{\pm} \in \widetilde{S}^{-\infty,0} $ is supported in $ \widetilde{\Gamma}^{\pm} (R_2,V_2,I_2,\sigma) $,
\begin{eqnarray}
 \left| \left| \scal{r}^{-4k} e^{-it P} f (h^2 P)  O \! p^h_{\kappa} (\chi_{\pm}) \tilde{\psi}_{\kappa} \scal{r}^k \right| \right|_{L^2 \rightarrow L^2} \lesssim  \scal{t/h}^{-k} , \qquad \pm t \geq 0  .   \label{weakoutgoingestimate}
\end{eqnarray} }
\item{{\bf Low frequency estimates:} if $(\chi_{\epsilon,\pm})_{\epsilon }$ is a bounded family of symbols in $\widetilde{S}^{-\infty,0} $ which are supported in $ \widetilde{\Gamma}^{\pm}_{\epsilon} (R_2,V_2,I_2,\sigma) $,
$$ \left| \left| \scal{\epsilon r}^{-4k} e^{-it P} f (P/ \epsilon^2)  O \! p_{\epsilon,\kappa} (\chi_{\epsilon,{\pm}}) \tilde{\psi}_{\kappa}(\epsilon r)   \scal{\epsilon r}^k \right| \right|_{L^2 \rightarrow L^2} \lesssim \scal{ \epsilon^2 t}^{-k} , \qquad \pm t \geq 0  .   $$}
\end{enumerate}
\end{prop}




 We will use here the results of paragraph \ref{paragraphetempsfini}.

\bigskip 

\noindent {\it Proof of Proposition \ref{halflocalized}.}
We consider in detail the high frequency outgoing case for $ t \geq 0 $. We can replace $ O \! p^h_{\kappa} (\chi_{+})  \scal{r}^k  $
by $ O \! p^h_{\kappa} (\chi_{+}^k) 
$ for some $ \chi_{+}^k \in \widetilde{S}^{-\infty,k} $ supported in the same set as $ \chi_{+} $; indeed, this is only at the expense of a remainder of the form $ \scal{r}^{-N} O_{L^2 \rightarrow L^2} (h^N) $ (for any fixed $N$) and whose contribution to the estimate  is a bound by $ \scal{t/h}^{-k} $ thanks to (\ref{aprioriblackboxhighfrequency}).
We then use a spatial dyadic partition of unity to split
\begin{eqnarray}
 \chi_{+}^k = \sum_{R = 2^l \atop l \geq l_0} \chi_{+,R}^k , \qquad \chi_{+,R}^k = \chi (r/R) \chi_{+}^k,  \label{decompositionspatiale}
\end{eqnarray} 
with some $ \chi \in C_0^{\infty} (0,+\infty) $ so that each $ \chi_{+,R}^k $ belongs to $ \widetilde{S}^{-\infty,0} $ with seminorms of order $R^k$.
For some small enough $ \varepsilon_2 > 0 $ to be chosen below, we pick $ T_+ > 0 $,  large enough such that for all $ l \geq l_0 $,
\begin{eqnarray}
 \phi^s \left( C_{\kappa} ( \mbox{supp} (\chi_{+,R}^k) ) \right) \subset \left\{ \rho > (1-\varepsilon^2_2) p^{1/2}, \ r > R_2 \right\} \qquad \mbox{for} \ \  s \geq R T_+ , \label{interetdeT+} 
\end{eqnarray} 
with $ \phi^s = \phi_1^s $ defined prior to Proposition \ref{propagationclassiqueprop} (in the low frequency case, we should consider $ \phi^s_{\epsilon} $). This is possible by the item 2 of Proposition \ref{propagationclassiqueprop} since, using the notation (\ref{notationpasterrible}) with $ \epsilon = 1 $,
$$ C_{\kappa} ( \mbox{supp} (\chi_{+,R}^k) ) \subset \widetilde{\Gamma}^{\pm}_1 (R_2,\sigma) . $$
  For each $R$, we then proceed as follows: 
  
\noindent \underline{If $ 0 \leq t \leq T_+ h R $}. We write $  \scal{r}^{-4k} e^{-it P} f (h^2 P)  O \! p^h_{\kappa} (\chi^k_{+,R}) \tilde{\psi}_{\kappa}  $  as 
$$\left( \scal{r}^{-4k} f (h^2 P) \scal{r}^{4k} \right) \scal{r}^{-4k} \left( e^{-it P}  O \! p^h_{\kappa} (\chi_{+,R}^k) \tilde{\psi}_{\kappa} e^{itP} \right) e^{-itP}  , $$
where, as in the proof of Proposition \ref{cellepourfortementsortant}, the first parenthese in the right hand side is bounded on $ L^2 $ thanks to Theorem \ref{theoremcaclulfonctionnel}.
 The second parenthese  can be computed  by mean of Theorem \ref{Egorovquantifie}. We get a sum of bounded pseudo-differential operators with symbols supported where $ r \sim R  $ (using the item 1 of Proposition \ref{propagationclassiqueprop} and that we propagate the support of $ \chi^k_{+,R} $ over a time $ t/h  \lesssim R $) plus a remainder which is, for any fixed $N$, of order $ h^N R^{-N} $, say in $ L^2 $ operator norm (here the stronger $ {\mathscr H}_{-N}^{-2N} \rightarrow {\mathscr H}_N^{2N} $ norm is not necessary). Since $ \scal{r}^{-4k} $ composed with  pseudo-differential operator localized in $ r \sim R $ has norm $ O (R^{-4k}) $  and since $ 0 \leq t / h \lesssim R $, we find
\begin{eqnarray}
 \left| \left| \scal{r}^{-4k} e^{-it P} f (h^2 P)  O \! p^h_{\kappa} (\chi_{+,R}^k) \tilde{\psi}_{\kappa} \right| \right|_{L^2 \rightarrow L^2} & \lesssim &  R^{-4k} R^k \nonumber \\
 & \lesssim & \scal{t/h}^{-k} R^{-2k} , \label{1ereasommerenR}
\end{eqnarray}
where the factor $ R^k $ takes into account that $ R^{-k} \chi_{+,R}^k $ is bounded in $ \widetilde{S}^{-\infty,0} $.

\noindent \underline{If $ t \geq T_+ h R $}. In this case, we write
$  \scal{r}^{-4k} e^{-it P} f (h^2 P)  O \! p^h_{\kappa} (\chi_{+,R}^k) \tilde{\psi}_{\kappa} $ as 
$$ \scal{r}^{-4k} f (h^2 P) e^{-i(t- T_+ h R) P}\left( e^{-i T_+ h R P}  O \! p^h_{\kappa} (\chi_{+,R}^k) \tilde{\psi}_{\kappa} e^{i T_+ h R P} \right) e^{-i T_+ h R P}  . $$
By (\ref{interetdeT+}), Theorem \ref{Egorovquantifie} and the seminorms estimates of $ \chi_{+,R}^k $, the parenthese is a sum of pseudo-differential operators with  symbols of size $ R^k $ in $ \widetilde{S}^{-\infty,0} $, supported in strongly outgoing areas,
\begin{eqnarray}
 \sum_{\kappa} O \! p_{\kappa}^h (\chi^k_{\kappa,R}(h)) \tilde{\psi}_{\kappa} , \qquad \mbox{supp}(\chi^k_{\kappa,R}(h)) \subset \widetilde{\Gamma}^+_{\rm st} (R / C, V_{\kappa} , I_2 , \varepsilon_2) \label{pseudo-differential-sum}
\end{eqnarray}
 with the additional property that $ r \sim R $ on their supports, and of a  remainder $ O _{ {\mathscr H}_{-2N}^{-N} \rightarrow {\mathscr H}_{2N}^N } ( h^N R^{-N})$ for any fixed $N$. In particular, if we take $ N \geq \max ( k + 1 , N_{k} ) $ (see (\ref{aprioriblackboxhighfrequency})), we get
\begin{eqnarray}
\left| \left| \scal{r}^{-4k} f (h^2 P) e^{-i(t- T_+ h R) P}  O _{ {\mathscr H}_{-N}^{-2N} \rightarrow {\mathscr H}_N^{2N} } ( h^N R^{-N}) \right| \right|_{L^2 \rightarrow L^2  } \qquad \qquad \qquad  \nonumber \\   \lesssim  h^{N-N_k} R^{- k-1}\left| \left| \scal{r}^{-k-1} f (h^2 P) e^{-i(t- T_+ h R) P} \scal{r}^{-k-1} \right| \right|_{L^2 \rightarrow L^2} \nonumber  \\
  \lesssim   \scal{t/h - T_+ R}^{-k}  R^{-k-1} \nonumber \\
  \lesssim   \scal{t/h}^{-k} R^{-1} . \label{2emeasommerenR}
\end{eqnarray}
To get the contribution of the pseudo-differential sum (\ref{pseudo-differential-sum}), we use Theorem \ref{Isozaki-Kitada-explicite}, which is why we need to choose $ \varepsilon_2 $ small enough. For any given $N$, we can write
$$  O \! p^h_{\kappa} (\chi^k_{\kappa,R}(h)) \tilde{\psi}_{\kappa}  =  J^h_{\kappa} ( a^h ) J^h_{\kappa} (b_R^h)^{\dag} + O _{ {\mathscr H}_{-N}^{-N} \rightarrow {\mathscr H}_N^N } ( h^N R^{-N}) $$ 
 where, by (\ref{supportestimate}),  $b_R^h$ is supported in $ r \sim R $ (this allows to get the additional factor $ R^{-N} $ in the remainder term) and belongs to $ S_0 $ with seminorms of order $R^k$ (uniformly $h$). The contribution of the remainder is estimated as above by choosing  $N$ large enough, while the contribution of the first term follows from Proposition \ref{cellepourfortementsortant} through
\begin{eqnarray}
 \left| \left| \scal{r}^{-4k} f (h^2 P) e^{-i(t- T_+ h R) P}   J^h_{\kappa} ( a^h ) J^h_{\kappa} (b_R^h)^{\dag}    \right| \right|_{L^2 \rightarrow L^2} \qquad \qquad \qquad \qquad \nonumber  \\  \lesssim   R^{-3k}\left| \left| \scal{r}^{-4k} f (h^2 P) e^{-i(t- T_+ h R) P}   J^h_{\kappa} ( a^h ) J^h_{\kappa} (b_R^h)^{\dag}   \scal{r}^{3k}  \right| \right|_{L^2 \rightarrow L^2} \nonumber \\ \lesssim   R^{-3k} \scal{t/h  - T_+ R }^{-k} R^k \nonumber \\
 \lesssim R^{-k} \scal{t/h}^{-k} , \label{3emeasommerenR}
\end{eqnarray}
where the factor $ R^k $ on the third line is the size of seminorms of $b_R^h $ in $ S_0 $.
Combining (\ref{1ereasommerenR}), (\ref{2emeasommerenR}) and (\ref{3emeasommerenR}), we get
$$  \left| \left| \scal{r}^{-4k} f (h^2 P) e^{-it P}  O \! p^h_{\kappa} (\chi^k_{+, R}) \tilde{\psi}_{\kappa}  \right| \right|_{L^2 \rightarrow L^2}    \lesssim R^{-1} \scal{t/h}^{-k} $$
which, once summed over $ R = 2^l $, provides the estimate  (\ref{weakoutgoingestimate}).
%
The low frequency case is obtained analogously by using the low frequency part of Theorem \ref{Egorovquantifie} together with (\ref{aprioriblackboxlowfrequency}). \finpreuve

\bigskip

Proposition \ref{halflocalized} provides time decay estimates with rate proportional to the  decay  rate of the weight. In the next two propositions, we get fast decay (and $ O (h^{\infty}) $ estimates at high frequency) for suitable microlocalizations.

\begin{prop}[Improved microlocal propagation estimates I] \label{fullyone} Let  $ I_2 \Subset (0,+\infty) $, $ V_2 \Subset V_{\kappa} $, $ \sigma \in (-1,1) $ and $ R_1 \geq 1 $. If $ R_2 \gg 1 $ then for each $ k \in \Na $ 
and $ \chi_{\pm} \in \widetilde{S}^{-\infty,0} $  supported in $ \widetilde{\Gamma}^{\pm} (R_2,V_2,I_2,\sigma) $,
one has
$$  \left| \left| {\mathds 1}_{[0,R_1]} (r)  f (h^2 P) e^{-it P}  O \! p^h_{\kappa} (\chi_{\pm}) \tilde{\psi}_{\kappa} (r)   \right| \right|_{L^2 \rightarrow L^2}    \lesssim h^k \scal{t/h}^{-k}, \qquad \pm t \geq 0, \ h \in (0,1] . $$
\end{prop}
This proposition reflects the intuitive fact that the forward (resp. backward) propagation of data localized in a far away outgoing (resp. incoming) area does not meet the region $ \{ r \leq R_1 \} $. Note that we consider  only the high frequency case, for which the estimate is improved by a factor $h^k $ compared to the  one of Proposition \ref{halflocalized}. At low frequency, 
 Proposition \ref{halflocalized} will be sufficient for us.

\bigskip

\noindent {\it Proof of Proposition \ref{fullyone}.} Here again we consider  the  outgoing case. We use the notation of the proof of Proposition \ref{halflocalized}, in particular  $ T_+$ and the decomposition (\ref{decompositionspatiale}). We distinguish two cases.

\noindent \underline{If $ 0 \leq t \leq T_+ h R $}. By Proposition \ref{Egorovquantifie}, we can write 
$$  {\mathds 1}_{[0,R_1]} (r)  f (h^2 P) e^{-it P}  O \! p^h_{\kappa} (\chi_{+,R}) \tilde{\psi}_{\kappa}  =
 \sum_{\kappa_1}  {\mathds 1}_{[0,R_1]} (r)  f (h^2 P) O \! p_{\kappa_1}^h (a^h_{R} (t)) e^{-itP} + O_{L^2 \rightarrow L^2} (h^N R^{-N}) $$
with symbols $ a_R^h (t) $ bounded in $ \widetilde{S}^{-\infty,0} $ as $ h,t,R $ vary and supported in $ r \sim R $ by the first item of Proposition \ref{propagationclassiqueprop}. In particular, they are supported in sets where $ r \gtrsim R_2 \gg 1$. Thus, using the pseudodifferential expansion of $ f (h^2 P) $  in Theorem \ref{theoremcaclulfonctionnel} (here the  localization $ \zeta (r) $ is implicit for we can write $ f (h^2 P) O \! p_{\kappa_1}^h (a^h_{R} (t)) = f (h^2 P) \zeta (r) O \! p_{\kappa_1}^h (a^h_{R} (t)) $), it follows from symbolic calculus and the form of the remainder terms in this theorem that
$$   {\mathds 1}_{[0,R_1]} (r)  f (h^2 P) O \! p_{\kappa_1}^h (a^h_{R} (t))  = O_{L^2 \rightarrow L^2} (h^N R^{-N}) $$
for any $N$. We thus conclude that, for any given $k$, 
\begin{eqnarray}
 \big| \big| {\mathds 1}_{[0,R_1]} (r)  f (h^2 P) e^{-it P}  O \! p^h_{\kappa} (\chi_{+,R}) \tilde{\psi}_{\kappa} \big| \big|_{L^2 \rightarrow L^2} = O (h^k R^{-1-k} ) = O \big(h^k \scal{t/h}^{-k} R^{-1} \big)  . \label{sommationcorrigee}
\end{eqnarray}

\smallskip

\noindent \underline{If $ t \geq T_+h R $}. In comparison to the proof of Proposition \ref{halflocalized}, it suffices to consider the terms
\begin{eqnarray*}
 \big| \big|  {\mathds 1}_{[0,R_1]} (r) f (h^2 P) e^{-i(t- T_+ h R) P}   J^h_{\kappa} ( a^h ) J^h_{\kappa} (b_R^h)^{\dag}  \big| \big|_{L^2 \rightarrow L^2} \lesssim \qquad \qquad \qquad
 \\
 \qquad \qquad \qquad  R^{-2k}  \big| \big|  {\mathds 1}_{[0,R_1]} (r) f (h^2 P) e^{-i(t- T_+ h R) P}   J^h_{\kappa} ( a^h ) J^h_{\kappa} (b_R^h)^{\dag} \scal{r}^{2k} \big| \big|_{L^2 \rightarrow L^2} 
\end{eqnarray*}
since all the other ones are remainder terms carrying an additional $ h^N $ factor with $ N $ arbitrarily large. To estimate the norm in the second line, we use the Isozaki-Kitada parametrix as in the proof of Proposition \ref{cellepourfortementsortant}. All remainders decay  as $ \scal{t / h -T_+R}^{-k} $ times $h^k$ (or even $h^N$) by pushing the expansion to a sufficiently high order exactly as in the proof of Proposition \ref{cellepourfortementsortant}. Thus, it remains to consider the main term which is
\begin{eqnarray}
  {\mathds 1}_{[0,R_1]} (r) f (h^2 P)   J^h_{\kappa} ( a^h ) e^{-i (t - T_+ h R)D_x^2) } J^h_{\kappa} (b_R^h)^{\dag} \scal{r}^{2k} . 
\label{contributionremainderin}
\end{eqnarray}
Using Theorem \ref{theoremcaclulfonctionnel}, one can write $$  {\mathds 1}_{[0,R_1]} (r) f (h^2 P)  =  {\mathds 1}_{[0,R_1]} (r) f (h^2 P) {\mathds 1}_{[0,\tilde{R}_1]} + O_{L^2 \rightarrow L^2 } (h^N) \scal{r}^{-N} $$
with $ \tilde{R}_1 > R_1 $. By choosing $ N \geq 3 k $, the contribution of the above remainder in (\ref{contributionremainderin}) is of the form $ O (h^k \scal{t/h - T_+ R}^{-k}) $ by Proposition \ref{propparam}. On the other hand, by choosing $ R_2 \gg \tilde{R}_1 $, the first item of Proposition \ref{nonstatprop} shows that (uniformly in $R$)
 $$  {\mathds 1}_{[0,R_1]} (r) f (h^2 P)   J^h_{\kappa} ( a^h ) e^{-i (t-T_+ h R) D_x^2 } J^h_{\kappa} (b_R^h)^{\dag} \scal{r}^{2k} = O (\scal{t/h - T_+ R}^{-\infty} h^{\infty}) . $$
 We thus get
 \begin{eqnarray*}
   R^{-2k}  \big| \big|  {\mathds 1}_{[0,R_1]} (r) f (h^2 P) e^{-i(t- T_+ h R) P}   J^h_{\kappa} ( a^h ) J^h_{\kappa} (b_R^h)^{\dag} \scal{r}^{2k} \big| \big|_{L^2 \rightarrow L^2} & \lesssim & h^k R^{-2k} \scal{t/h - T_+ R}^{-k} \\  & \lesssim &
    h^k R^{-k} \scal{t/h}^{-k} . 
\end{eqnarray*}    
Taking (\ref{sommationcorrigee}) into account, we conclude as in Proposition \ref{halflocalized} by summing all estimates over $R$. \finpreuve
 


\bigskip

In the next proposition, we use the notation (\ref{zonesortantelargecoordonnees}) (when $ \epsilon = 1 $ we do not indicate the dependence on $ \epsilon $). We let $ \chi_{\rm st}^{\pm} $ and $ \chi_{\epsilon,{\rm st}}^{\pm} $ be supported in an angular patch associated to a given chart $ \kappa $ (the same as in all previous propositions) but we allow the symbols $ \chi_{\mp} $ and $ \chi_{\epsilon,\mp} $ to be angularly supported in a possibly different patch associated to another chart $ \kappa^{\prime} $.

\begin{prop}[Improved microlocal propagation estimates II] \label{improvedoutgoingII}
Let $ V_2 \Subset V_0 \Subset V_{\kappa} $, $ I_2 \Subset (0,+\infty) $. Let also $ V^{\prime}_2 \Subset V_{\kappa^{\prime}} $, $ I_2^{\prime} \Subset (0,+\infty) $ and $ \sigma \in (-1,1) $. If $ \varepsilon_2 > 0 $ is small enough and $ R_2 > 0 $ is large enough, the following estimates hold for all $  k \in \Na$:
\begin{enumerate}
\item{{\bf High frequency case:} if $ \chi_{\pm} , \chi_{\rm st} \in \widetilde{S}^{-\infty,0} $ satisfy
$$ \emph{supp} (\chi_{\mp}) \subset \widetilde{\Gamma}^{\mp} (R_2,V^{\prime}_2,I^{\prime}_2,\sigma)_{\kappa^{\prime}}, \qquad \emph{supp}(\chi_{\rm st}^{\pm} ) \subset \widetilde{\Gamma}^{\pm}_{\rm st} \big( R_2 , V_2 , I_2 , \varepsilon_2 \big) $$
then
$$ \left| \left| \scal{r}^k  \tilde{\psi}_{\kappa^{\prime}} O \! p^h_{\kappa^{\prime}} (\chi_{\mp})^*   f (h^2 P) e^{-it P}  O \! p^h_{\kappa} (\chi_{{\rm st}}^{\pm}) \tilde{\psi}_{\kappa} \scal{r}^k  \right| \right|_{L^2 \rightarrow L^2}    \lesssim h^k \scal{t/h}^{-k}, \qquad \pm t \geq 0. $$
}
\item{{\bf Low frequency case:}
if $ (\chi_{\epsilon,\pm})_{\epsilon} , (\chi_{\epsilon,{\rm st}}^{\pm})_{\epsilon} $ are bounded families of $ \widetilde{S}^{-\infty,0} $ satisfying
$$ \emph{supp} (\chi_{\epsilon,\mp}) \subset \widetilde{\Gamma}^{\mp}_{\epsilon} (R_2,V_2^{\prime},I_2^{\prime},\sigma)_{\kappa^{\prime}}, \qquad \emph{supp}(\chi_{\epsilon,{\rm st}} ) \subset \widetilde{\Gamma}^{\pm}_{{\rm st},\epsilon} \big( R_2 , V_2 , I_2 , \varepsilon_2 \big) $$
then
$$ \left| \left| \scal{\epsilon r}^k \tilde{\psi}_{\kappa^{\prime}} (\epsilon r)  O \! p_{\epsilon,\kappa^{\prime}} (\chi_{\epsilon,\mp})^*    f ( P/\epsilon^2) e^{-it P}  O \! p_{\epsilon,\kappa} (\chi_{\epsilon,{\rm st}}^{\pm}) \tilde{\psi}_{\kappa}(\epsilon r) \scal{\epsilon r}^k  \right| \right|_{L^2 \rightarrow L^2}    \lesssim  \scal{\epsilon^2 t}^{-k}, \qquad \pm t \geq 0. $$
}
\end{enumerate}
 %

\end{prop}


\begin{lemm} \label{lemmedexpansionpseudoOIF} Let $ I_0 \Subset (0,+\infty) $. If $ \varepsilon_0 > 0 $ is small enough and $ R_0  > 0 $ is large enough, then for all $ a \in S_0 $ supported in $ \Theta^{\pm} (R_0,V_0,I_0,\varepsilon_0) $ and $ \chi \in S^{-\infty,k} $ supported in $ (R_2,\infty) \times V_2 \times \Ra^n $ with $ R_2 \geq R_0 $, one can write for any $ N $ 
$$  O \! p^h_{\kappa} (\chi) J^h_{\kappa} (a) =  J^h_{\kappa} (a_N (h)) + \ \mbox{finite sum of }  \  h^N \scal{r}^{-N} B_h \scal{r}^{-N } J_{\kappa}^h(r_N(h))  $$
with $ ||B_h||_{L^2 \rightarrow L^2} \lesssim 1 $, $ ( r_N (h) )_h $ bounded in $  \widetilde{S}^{-\infty,0} $ and supported in $ \Theta^{\pm} (R_0,V_0,I_0,2 \varepsilon_0) $, 
and with $( a_N (h) )_h $ bounded in $ \widetilde{S}^{-\infty,k} $ satisfying $$ \emph{supp} (a_N(h)) \subset \emph{supp} \big(\chi (.,.,\partial_r \varphi , \partial_{\theta} \varphi \big) \times a \big) . $$
More precisely, 
$$ a_N (h) = \chi (r,\theta,\partial_r \varphi , \partial_{\theta} \varphi) a (r,\theta,\varrho,\vartheta) + O (h) $$
where $ O (h) $ is a finite sum of products of derivatives of $ \chi $ (of order $ \geq 1 $) evaluated at $ (r,\theta,\partial_r \varphi , \partial_{\theta} \varphi)  $, of derivatives of $a$ and of rational fractions in derivatives of $ \varphi $.

\end{lemm}

\noindent {\it Proof.} It follows from the usual calculation of the action of a pseudo-differential operator on an oscillatory integral, see e.g. \cite{AlinhacGerard,Taylor0}. \finpreuve

\bigskip

\noindent {\bf Remark.} Of course, a completely parallel statement holds at low frequency but we do not quote it for we give the proof of Proposition \ref{improvedoutgoingII} only in the high frequency case. Also, the parameter $ 2 \varepsilon_0 $ in the support of the remainder terms (coming from technical considerations due to the non locality of $ O \! p_{\kappa}^h (\chi) $) could be replaced by any $ \tilde{\varepsilon}_0 > \varepsilon_0 $ but this is irrelevant for our purposes.

\bigskip

\noindent {\it Proof of Proposition \ref{improvedoutgoingII}.} We consider again the high frequency case for $ t \geq 0 $. Also, w.l.o.g. as in the proof of Proposition \ref{cellepourfortementsortant}, we may assume that $ V_0 $ is convex to be in position to use the expression of    $  e^{-it P}  O \! p^h_{\kappa} (\chi_{{\rm st}}^{+}) \tilde{\psi}_{\kappa} $ given by Theorem \ref{Isozaki-Kitada-explicite}. Proceeding
exactly as in the proof of Proposition \ref{cellepourfortementsortant}, up to the replacement of (\ref{aprioriblackboxhighfrequency}) by the new a priori estimate  
$$ \big| \big| \scal{r}^k \tilde{\psi}_{\kappa^{\prime}}  O \! p_{\kappa^{\prime}}^h (\chi_-)^*   f (h^2 P) e^{-i(t- \tau ) P} \scal{r}^{-N} \big| \big|_{L^2 \rightarrow L^2} \lesssim \scal{(t-\tau)/h}^{-k-1} , \qquad 0 \leq \tau \leq t, $$
which follows from the adjoint estimate to (\ref{weakoutgoingestimate}) for $N$ large enough, 
%
we  see that the contribution of the remainder $ R_N^h (t) $ is $ O (h^k \scal{t/h}^{-k}) $. Note that here, we do not have to care about the fact that $ \kappa $ and $ \kappa^{\prime} $ may be different. It then remains to consider the contribution of 
$$   \scal{r}^k \tilde{\psi}_{\kappa^{\prime}}  O \! p_{\kappa^{\prime}}^h (\chi_-)^*   f (h^2 P) J_{\kappa}^h (a^h) e^{-it D_x^2} J_{\kappa}^h (b^h)^{\dag} \scal{r}^k . $$
We  consider the case when $ \kappa = \kappa^{\prime} $ and explain at the end of the proof how to handle the general case.
 Using the expansion of Theorem \ref{theoremcaclulfonctionnel} and symbolic calculus, one can write for any $N$,
 $$ \scal{r}^k \tilde{\psi}_{\kappa}  O \! p_{\kappa}^h (\chi_-)^*  f (h^2 P) = O \! p_{\kappa}^h (\chi_-^k (h)) + O (h^N)_{L^2 \rightarrow L^2} \scal{r}^{-N}  $$
 with $ \chi_{-}^k (h) \in \widetilde{S}^{-\infty,k}$ with the same support as $ \chi_- $ and bounded with respect to $h$. Note that we do not need to put a localization $ \tilde{\psi}_{\kappa} $ on the right hand side since  $ J^h_{\kappa} (a^h) = \tilde{\psi}_{\kappa} J_{\kappa}^h (a^h) $ by the localization of  the support  of $a^h$. The contribution of the remainder follows from Proposition \ref{propparam}, provided we take $ N \geq k $. On the other hand, using Lemma \ref{lemmedexpansionpseudoOIF}, we can compute
 $$ O \! p_{\kappa}^h (\chi_-^k (h))  J_{\kappa}^h (a^h) e^{-it D_x^2} J_{\kappa}^h (b^h)^{\dag} =  J^h_{\kappa} (a_N (h)) e^{-it D_x^2} J_{\kappa}^h (b^h)^{\dag} + \mbox{remainder terms} . $$
 The contribution of the remainder terms follows from Proposition \ref{propparam}, using their fast decay in $r$ and $h$. On the other hand, 
on the support of $ a_N (h) $,  one must have
 $ (r,\theta,\partial_r \varphi, \partial_{\theta} \varphi) \in \widetilde{\Gamma}^- (R_2,V_2,I_2,\sigma) $ and $ (r,\theta,\varrho,\vartheta) \in \Theta^{+} (R_0,V_0,I_0,\varepsilon_0) $. This implies in particular that
 $$  - \partial_r \varphi > \sigma p_{\kappa}(r,\theta,\partial_r \varphi, \partial_{\theta} \varphi)^{1/2} = \sigma |\varrho| \qquad \mbox{and} \qquad   \varrho > 0.
 $$
 By (\ref{pourdesconditionsdesupport}), these conditions are incompatible  if $ \sigma \geq 0 $, so $ a_N (h) \equiv 0 $  in this case. On the other hand, if $ \sigma < 0 $, one has $ 0 < \partial_r \varphi < |\sigma| p (r,\theta,\partial_r \varphi, \partial_{\theta} \varphi)^{1/2} $, hence
 $$ \sigma^2 r^{-2} g^{jk}(r,\theta) \partial_{\theta_j} \varphi  \partial_{\theta_k} \varphi >  (1-\sigma^2) (\partial_r \varphi )^2$$
 which, together with  (\ref{pourdesconditionsdesupport}), implies that for some $ c_{\sigma} > 0 $
 $$ |\theta - \vartheta| > c_{\sigma} \varrho^2 . $$
 Thus, on the support of the kernel  of $  J^h (a_N (h) ) e^{-it D_x^2} J^h (b^h)^{\dag} $, we have
 $$ | \theta - \vartheta | > c_{\sigma} \varrho^2 \gg \varepsilon_2 \gtrsim | \theta^{\prime} - \vartheta| $$
 so we obtain the fast decay by mean of the item 2 of Proposition \ref{nonstatprop}, provided $ \varepsilon_2 $ is small enough and $R_2$ is large enough. This completes the proof (when $  \kappa = \kappa^{\prime}$).
When $ \kappa \ne \kappa^{\prime} $, we may split $ O \! p_{\kappa^{\prime}}^h (\chi_-)^* $ as $ O \! p_{\kappa^{\prime}}^h (\chi_-)^* \chi_{\kappa} + O \! p_{\kappa^{\prime}}^h (\chi_-)^*  (1-\chi_{\kappa}) $ with $ \chi_{\kappa} \equiv 1 $ near the spatial projection of the support of $ a^h $. The operator $  O \! p_{\kappa^{\prime}}^h (\chi_-)^* \chi_{\kappa}  $ can then be written in the chart $  \kappa $ as in Proposition \ref{propnormalizingangular} (and then treated as above), up to terms which decay fast in $h$ and $  r$. The contribution  of $ (1-\chi_{\kappa})  f (h^2 P) J_{\kappa}^h (a^h)  $ also produces terms which are $ O (h^{\infty}) $ and decay fast in $r$. All these decaying remainders can then be handled thanks to Proposition \ref{propparam}. \finpreuve

\section{Strichartz estimates} \label{SectionStrichartz}
\setcounter{equation}{0}

We focus on the low frequency case (in dimension $n \geq 3$), i.e. on $ u_{\rm low} $ defined in (\ref{defhighlow}). Indeed this one is slightly more technical than the high frequency case, for instance to handle the $ L^q \rightarrow L^q $ estimates of $ f (P/\epsilon^2) $. In paragraph \ref{minormodification}, we explain the minor modifications to handle high frequencies.

\subsection{Finite time estimates} \label{Finitetimeestimates}
In this paragraph, we use the well known geometric optics technique to derive propagator approximations  for finite times, but depending both on the frequency and spatial localizations. This follows previous similar arguments introduced in \cite{MizutaniCPDE} for high frequency localizations. Our main purpose  is  to give such an  approximation at low frequency, but we restate the high frequency case  both for completeness and for comparison with the low frequency regime.

For a given chart $ \kappa : U_{\kappa} \rightarrow V_{\kappa} $ on the angular manifold, $ V \subset V_{\kappa} $, $ I \Subset (0,+\infty) $, $ C \geq 1 $, $ \epsilon \in (0,1] $ and $ R \gg 1 $, we use the notation
\begin{eqnarray*}
\Omega_{ R}(V,I,C) &= & \{ (r,\theta,\rho,\eta) \in p_{\kappa}^{-1} (I) \ | \ r \in (R/C, C R), \ \theta \in V \}  \\
\Omega_{\epsilon,R} (V,I,C) & = & \{ (\breve{r},\theta,\breve{\rho},\eta) \in p_{\epsilon,\kappa}^{-1} (I) \ | \ \breve{r} \in (R/C , C R), \ \theta \in V \}. 
\end{eqnarray*}
Note that $ \Omega_{R} (V,I,C) = \Omega_{1,R} (V,I,C)  $ 

\begin{prop}[Existence of phase functions] Let $ V \Subset V_{\kappa} $ be a relatively compact open convex subset of $ V_{\kappa} $. Let $  V_0 \Subset V $, $ C_0  > 1 $ and $  I_0 \Subset (0,+\infty) $. There are $ 0 < t_0 \ll 1 $ and  $ R_0 \gg 1 $ such that one can find a family of smooth functions
$$  (\varphi_{\epsilon,R})_{ \epsilon \in (0,1],R \geq R_0}  $$
defined  on  $  (- t_0 R , t_0 R) \times \Omega_{\epsilon,R} (V_0,I_0,C_0) $, solving the eikonal equation
$$  \partial_s \varphi_{\epsilon,R} + p_{\epsilon,\kappa} (r,\theta,\partial_{r,\theta} \varphi_{\epsilon,R} ) =0, \qquad \varphi_{\epsilon,R} (0,r,\theta,\rho,\eta) = r \rho + \theta \cdot \eta ,$$
and satisfying the estimates
\begin{eqnarray}
\big| \partial_r^j \partial_{\theta}^{\alpha} \partial_{\rho}^k \partial_{\eta}^{\beta} \big( \varphi_{R,\epsilon}(s) -  \varphi_{\epsilon,R}(0) +   s p_{\epsilon,\kappa} \big) \big| \leq C_{\gamma} \frac{s^2}{R} R^{-j-|\beta|}, \label{approxphasstat}
\end{eqnarray}
for $ R \geq R_0 $, $ \epsilon \in (0,1] $, $ |s| < t_0 R $ and $ (r,\theta,\rho,\eta) \in \Omega_{\epsilon,R} (V_0,I_0,C_0) $. 
\end{prop}

\noindent {\it Proof.} It follows the usual local in time resolution of the Hamilton-Jacobi equation, by using the flow estimates given in Proposition \ref{estimeeflotlocal} which allow to show that the map $ (r,\theta,\rho,\eta) \mapsto (\bar{r}^s_{\epsilon}, \bar{\vartheta}_{\epsilon}^s,\rho,\eta) $ is a diffeomorphism if $ |s| \leq t_0 R $ with $ t_0 $ small enough. More precisely, to prove that this is a diffeomorphism, one can check that the map $ (x,\theta) \mapsto (R^{-1} \bar{r}^s_{\epsilon}, \bar{\vartheta}^s_{\epsilon}) (Rx,\theta,\rho,\eta) $ is close to the identity on $ (1/2C_0,2C_0) \times V $ provided that $ s/R $ is close enough, uniformly in $ \epsilon , \rho , \eta $. The convexity of $V$ allows to check that this map is injective while standard arguments show that the range will contain $ (1/C_0,C_0) \times V_0 $. \finpreuve


\bigskip

We can next consider the related Fourier integral operators
\begin{eqnarray}
W_{\epsilon,R} (s,A_{\epsilon}) u (\breve{r},\theta) = (2 \pi)^{-n} \int \int e^{i  \big( \varphi_{\epsilon,R}(s,\breve{r},\theta,\breve{\rho},\eta) - \breve{r}^{\prime} \breve{\rho} - \theta^{\prime} \cdot \eta \big)} A_{\epsilon} (s,\breve{r},\theta,\breve{\rho},\eta) u (\breve{r}^{\prime},\theta^{\prime}) d \breve{\rho} d \eta d \breve{r}^{\prime} d \theta^{\prime} \label{pourlapreuvepratique}
\end{eqnarray}
and, setting $ \varphi_R = \varphi_{1,R} $,
\begin{eqnarray}
W_{R}^h (s,A) u (r,\theta) = (2 \pi h)^{-n} \int \int e^{ \frac{i}{h} \big( \varphi_{R}(s,r,\theta,\rho,\eta) - r^{\prime} \rho - \theta^{\prime} \cdot \eta \big)} A (s,r,\theta,\rho,\eta) u (r^{\prime},\theta^{\prime}) d \rho d \eta d r^{\prime} d \theta^{\prime}
\nonumber
\end{eqnarray}
which are globally well defined on $ \Ra^n $ provided the amplitudes $ A_{\epsilon} $ and $A$ are supported  respectively in $ \Omega_{\epsilon,R} (V_0,I_0,C_0) $ and $ \Omega_{1,R} (V_0,I_0,C_0) $. Using the cutoffs $ \tilde{\psi}_{\kappa}(\epsilon r) $ and $ \tilde{\psi}_{\kappa}(r) $ chosen in (\ref{notationpartition}),
we can  pullback these operators on $ {\mathcal M} $, i.e.   define the operators
\begin{eqnarray}
W_{\epsilon,R,\kappa}(s,A_{\epsilon}) \tilde{\psi}_{\kappa}(\epsilon r) : = \Pi_{\kappa} \left( {\mathscr D}_{\epsilon} W_{\epsilon,R} (s,A_{\epsilon})  {\mathscr D}_{\epsilon}^{-1} \right) \Pi_{\kappa}^{-1} \tilde{\psi}_{\kappa}(\epsilon r)  \nonumber
\end{eqnarray}
and
\begin{eqnarray}
W_{R,\kappa}^h(s,A) \tilde{\psi}_{\kappa}( r) := \Pi_{\kappa}  W_{R}^h (s,A)   \Pi_{\kappa}^{-1} \tilde{\psi}_{\kappa}(r) . \nonumber
\end{eqnarray}
\begin{prop} \label{WKB} Let $ V \Subset V_{\kappa} $ be convex. Let  $ V_1 \Subset V_0 \Subset V $, $ C_0 > C_1 > 1 $ and $ I_1 \Subset I_0 \Subset (0,+\infty) $. There are $ 0 < t_0 \ll 1 $ and  $ R_0 \gg 1 $ such that for any $N \in \Na $ the following approximations hold.
\begin{enumerate}
\item{ {\bf Low energy WKB approximation:} Given a bounded family $ (a_{\epsilon,R})_{\epsilon,R} $ of $ \widetilde{S}^{-\infty,0} $ supported in $ \Omega_{\epsilon,R} (V_1,I_1,C_1) $, one can find a bounded family $  (A_{\epsilon,R} (\epsilon^2 t))_{\epsilon,R,t} $  of $ \widetilde{S}^{-\infty,0} $ supported in $ \Omega_{\epsilon,R} (V_0,I_0,C_0) $ and $ \chi \in C_0^{\infty} (0,+\infty) $ such that
$$ e^{-itP} O \! p_{\epsilon,\kappa} (a_{\epsilon,R})  \tilde{\psi}_{\kappa}(\epsilon r)  = W_{\epsilon,R,\kappa} \left( \epsilon^2 t  , A_{\epsilon,R} \right) \chi (\epsilon r / R )  \tilde{\psi}_{\kappa}(\epsilon r) + O_{L^1 \rightarrow L^2} ( \epsilon^{\frac{n}{2}} R^{-N}) $$
and
\begin{eqnarray}
 \big| \big|  W_{\epsilon,R,\kappa} \left( \epsilon^2 t  , A_{\epsilon,R} \right) \chi (\epsilon r / R) \tilde{\psi}_{\kappa}(\epsilon r) \big| \big|_{L^1 \rightarrow L^{\infty}} \lesssim \scal{t}^{-\frac{n}{2}} \label{displow1}
\end{eqnarray}
as long as
$$ \epsilon \in (0,1], \qquad R \geq R_0, \qquad |t| \leq t_0 \epsilon^{-2} R . $$
} 
\item{{\bf High energy WKB approximation:} Given a bounded family $ (a_{R})_{R} $ of $ \widetilde{S}^{-\infty,0} $ supported in $ \Omega_{R} (V_1,I_1,C_1) $, one can find a bounded family $  (A_{R}^h ( \frac{t}{h}))_{R,h,t} $  of $ \widetilde{S}^{-\infty,0} $ supported in $ \Omega_{R} (V_0,I_0,C_0) $ and $ \chi \in C_0^{\infty}(0,+\infty) $ such that
$$ e^{-itP} O \! p_{\kappa}^h (a_R) \tilde{\psi}_{\kappa}( r)  = W_{R,\kappa}^h \left( t / h , A^h_{R} \right) \chi (r/R) \tilde{\psi}_{\kappa}( r) + O_{L^1 \rightarrow L^2} (h^N R^{-N}) $$
and
\begin{eqnarray}
 \big| \big| W_{R,\kappa}^h \left( t / h , A^h_{R} \right) \chi ( r / R) \tilde{\psi}_{\kappa}( r) \big| \big|_{L^1 \rightarrow L^{\infty}} \lesssim |t|^{-\frac{n}{2}} \label{disphig2}
\end{eqnarray}
as long as
$$ h \in (0,1], \qquad R \geq R_0, \qquad |t| \leq t_0 h R . $$}
\end{enumerate}
\end{prop}

\begin{lemm} \label{conversionnoyaux} Let $ K  (\breve{r},\theta,\breve{r}^{\prime},\theta^{\prime}) $ be the kernel of an operator $ W $ with respect to the Lebesgue measure $ d \breve{r} d \theta $. Assume that $ K $ is supported in $ \big( (R_0,\infty) \times V \big)^2 $ for some $ V \Subset V_{\kappa} $. Then, the Schwartz kernel $ {\mathcal K}_{\epsilon} $ of $ \Pi_{\kappa} \left( {\mathscr D}_{\epsilon} W   {\mathscr D}_{\epsilon}^{-1} \right) \Pi_{\kappa}^{-1} $ with respect to the Riemannian measure satisfies
$$ \big| {\mathcal K}_{\epsilon} (r,\omega ,r^{\prime},\omega^{\prime}) \big| \leq C \epsilon^n \big| K (\epsilon r , \theta , \epsilon r^{\prime},\theta^{\prime}) (\epsilon r^{\prime})^{1-n} \big|, \qquad \omega = \kappa^{-1}(\theta), \ \omega^{\prime} = \kappa^{-1}(\theta^{\prime}), $$
for some constant $ C $ depending on $V$ but not on $K$ nor $ \epsilon $.
\end{lemm}

\noindent {\it Proof.}  We omit the conjugation by $ \Pi_{\kappa} $ whose role is irrelevant here. Then
\begin{eqnarray*}
 {\mathscr D}_{\epsilon} W u (r,\theta) &= & \epsilon^{\frac{n}{2}}  \int \int K (\epsilon r , \theta , \breve{r}^{\prime} , \theta^{\prime}) u (\breve{r}^{\prime},\theta^{\prime}) d \breve{r}^{\prime} d \theta^{\prime}  \\
  &= & \epsilon^{n}  \int \int K (\epsilon r , \theta , \epsilon r^{\prime} , \theta^{\prime}) (\epsilon r^{\prime})^{1-n}  ({\mathscr D}_{\epsilon}u) (r^{\prime},\theta^{\prime}) (r^{\prime})^{n-1} d r^{\prime} d \theta^{\prime}
\end{eqnarray*} 
so that the kernel of $ {\mathscr D}_{\epsilon} W {\mathscr D}_{\epsilon}^{-1} $ with respect to $ (r^{\prime})^{n-1} d r^{\prime} d \theta^{\prime} $ is $ \epsilon^n K (\epsilon r , \theta, \epsilon r^{\prime},\theta^{\prime}) (\epsilon r^{\prime})^{1-n} $. Since $ (r^{\prime})^{n-1} dr^{\prime} d \theta^{\prime} $ is comparable to the Riemannian density
$ (r^{\prime})^{n-1} \mbox{det}(g(r^{\prime},\theta))^{1/2} d r^{\prime} d\theta^{\prime} $, we get the result. \finpreuve
 

\bigskip

\noindent {\it Proof of Proposition \ref{WKB}.} We consider the low energy case.  Dropping the spatial cutoff for simplicity, one has the  identity
\begin{eqnarray}
 e^{-itP}  W_{\epsilon,R,\kappa} (0,A_{\epsilon,R}) = W_{\epsilon,R,\kappa} (\epsilon^2 t,A_{\epsilon,R}) - \int_0^{\epsilon^2 t} e^{-i \left( t - \frac{s}{\epsilon^2} \right) P} W_{\epsilon,R,\kappa} (s,b_{\epsilon,R}) ds  \label{integralin}
\end{eqnarray}
where
$$ W_{\epsilon,R} (s,b_{\epsilon,R}) =  \partial_s W_{\epsilon,R} (s,A_{\epsilon})  + i P_{\epsilon, \kappa} W_{\epsilon,R} (s,A_{\epsilon})   . $$
By the usual geometric optics construction, we can find, for any $N$, symbols $ A_{\epsilon,R}(s,r,\theta,\rho,\eta) $ in a bounded set of $ \widetilde{S}^{-\infty,0} $ (as $ \epsilon \in (0,1] $, $R \geq R_0$ and $ |s| \leq t_0 R $  vary), supported in $ \Omega_{\epsilon,R} (V_0,I_0,C_0)  $
and such that
$$ A_{\epsilon,R}|_{s=0}= a_{\epsilon,R} ,\qquad  b_{\epsilon,R}  (s) \ \mbox{in a bounded subset of} \ \widetilde{S}^{-\infty,-N} \qquad \mbox{supp} \big( b_{\epsilon,R} (s) \big) \subset \Omega_{\epsilon,R} (V_0,I_0,C_0) .  $$
This follows by solving iteratively transport equations in the usual manner and by observing that in the iterative construction of the amplitude $ A_{\epsilon,R} $ the symbols decay faster and faster in $\breve{r} $; in other words, the scale of classes $ \widetilde{S}^{-\infty,-j} $ replaces here the scale of powers $h^j$ in the usual semiclassical framework.  The boundedness in $s$ of the solutions to the transport equations follow from the flow estimates of Proposition \ref{estimeeflotlocal}. To get the remainder estimate and (\ref{displow1}), we proceed as follows. Since $ a_{\epsilon,R} $ is supported in a region where $ \breve{r} \sim R $, we can write
$$ O \! p^1 (a_{\epsilon,R}) =  O \! p^1 (a_{\epsilon,R}) \chi (\breve{r}/R) + O \!p^1 (a_{\infty,\epsilon,R}) $$
with $ a_{\infty,\epsilon,R} = O (R^{-N}) $ in $ \widetilde{S}^{-N,-N} $ for any $N$. In particular, using Lemma \ref{conversionnoyaux}, it is not hard to check that
$$ ||  O \!p_{\epsilon} (a_{\infty,\epsilon,R}) \tilde{\psi}_{\kappa} (\epsilon r) ||_{L^1 \rightarrow L^2} \lesssim_N R^{-N} \epsilon^{n/2}. $$
This allows to replace $  e^{-itP} O \! p_{\kappa}^h (a_R) \tilde{\psi}_{\kappa}( r)  $ by   $ e^{-itP} O \! p_{\kappa}^h (a_R) \chi (\epsilon r / R) \tilde{\psi}_{\kappa}( \epsilon r) $ and we are left with two types of terms: the main term of the expansion $ W_{\epsilon,R,\kappa} (s,A_{\epsilon,R}) $, which will produce (\ref{displow1}), and the remainder involving $ W_{\epsilon,R,\kappa} (s,b_{\epsilon,R}) \chi (\epsilon r / R) \tilde{\psi}_{\kappa} (\epsilon r) $ coming from the integral in (\ref{integralin}). We start with this remainder. Using (\ref{pourlapreuvepratique}), with $ b_{\epsilon,r} $ instead of $A_{\epsilon}$, and using the decay in $ \breve{r} $ together with the fact that we integrate over a fixed bounded in region in $ \eta / \breve{r} $, the Schwartz kernel of $ W_{\epsilon,R} (s,b_{\epsilon,R}) \chi (\breve{r} / R) $ with respect to $ d \breve{r} d \theta $ is bounded by $ C \breve{r}^{-N+(n-1)} $ and is supported in a region where both $\breve{r}$ and $ \breve{r}^{\prime} $ are of size $R$. Note that the power $ \breve{r}^{n-1} $ comes from the fact that the kernel is given by an integral where $ \eta $ belongs to a region of volume $ \breve{r}^{n-1} $. Then, by Lemma \ref{conversionnoyaux}, the kernel of $ W_{\epsilon,R,\kappa} (s,b_{\epsilon,R}) $ with respect to the Riemannian measure is bounded by
$$ \epsilon^n \scal{\epsilon r}^{-N/3} \scal{\epsilon r^{\prime}}^{-N/3} R^{-N/3} . $$
The corresponding operator has an $ L^1 \rightarrow L^2 $ norm of order $ \epsilon^{n/2} R^{-N/3} $ (if $ N/3 > n/2 $).
 Since $N$ is arbitrary,  $ | \epsilon^2 t |\lesssim R $ and the propagator is unitary on $ L^2 $, we get the control on the remainder of (\ref{integralin}) in $  L^1 \rightarrow L^2 $ operator norm. 
   Finally, the dispersion  estimate (\ref{displow1})  follow from the fact that
the $ L^1  \rightarrow L^{\infty} $ norm of $  W_{\epsilon,R,\kappa} \left( s  , A_{\epsilon,R} \right) \tilde{\psi}_{\kappa}(\epsilon r) \chi (\epsilon r / R) $ is controled by
$$ \epsilon^n \sup_{\breve{r},\theta,\breve{r}^{\prime},\theta^{\prime},\epsilon} \left| \left( \int e^{i  \big( \varphi_{\epsilon,R}(s,\breve{r},\theta,\breve{\rho},\eta) - \breve{r}^{\prime} \breve{\rho} - \theta^{\prime} \cdot \eta \big)} A_{\epsilon} (s,\breve{r},\theta,\breve{\rho},\eta)  d \breve{\rho} d \eta \right)  \scal{\breve{r}^{\prime}}^{1-n} \chi (\breve{r^{\prime}}/R)\right| \lesssim \epsilon^n \scal{s}^{-n/2}  $$
where 
the estimate by $ \scal{s}^{-n/2} $ follows from a  standard non stationary phase argument  by  exploiting that 
 $$ \varphi_{\epsilon}(s,\breve{r},\theta,\breve{\rho},\eta) = (\breve{r}-\breve{r}^{\prime}) \breve{\rho} + (\theta-\theta^{\prime} )\cdot \eta - s p_{\kappa,\varepsilon}(\breve{r},\theta,\breve{\rho},\eta) + O (s^2 / R) , $$
 (by (\ref{approxphasstat})). Note that the weight $ \scal{\breve{r}^{\prime}}^{1-n} $  is crucial to compensate that we integrate over a region of volume $O (\breve{r}^{n-1})$ in $ \eta $ (recall that both $ \breve{r} $ and $\breve{r^{\prime}}$ are of order $ R $ here). With $ s = \epsilon^2 t $ we find that $ \epsilon^n \scal{\epsilon^2 t}^{-n/2} \lesssim \scal{t}^{-n/2} $. 
  We refer to \cite{MizutaniCPDE} for more details on the stationary phase.
  The proof is similar at high energy. Up to the scaling in time, the main differences are that 
 we drop the scaling operators $ {\mathscr D}_{\epsilon}^{\pm} $ and that in the iterative construction of the amplitude we gain both decay in $h$ and in $ r $. 
   \finpreuve

\subsection{Proof of Theorem \ref{theoremlow}} \label{proofTheoremlow}
 It suffices to prove the result for the endpoint pair $ (p,q) = (2,2^*) = \big( 2 , \frac{2n}{n-2}\big)$, the other ones following by interpolation with the trivial estimate for $ (p,q) = (\infty,2) $.

For $ u_0 \in L^2 $, we use the notation (\ref{defhighlow}).
The starting point is the estimate 
\begin{eqnarray}
 ||  u_{\rm low} ||_{L^2 (\Ra ; L^{2^*})} \lesssim \left(  \sum_{\epsilon^2 = 2^{-k}} || (1-\chi(\epsilon r)) f (P/\epsilon^2) u ||_{L^2 (\Ra ; L^{2^*} ) }^2 + || \scal{r}^{-1} f (P/\epsilon^2) u ||_{L^2 (\Ra ; L^2)}^2 \right)^{1/2} \label{quiserarecitee}
\end{eqnarray} 
which follows from Theorem \ref{thlowfreqLq}.  By the integrated $L^2$ decay estimate (\ref{remplacerefvide}), we have
$$  || \scal{r}^{-1} f (P/\epsilon^2) u ||_{L^2 (\Ra ; L^2)} \lesssim || u_0 ||_{L^2} $$
where, in the right hand side, we may replace $ u_0 $ by $ \tilde{f} (P/\epsilon^2) u_0  $ with $ \tilde{f} \in C_0^{\infty} (0,+\infty) $ equal to $1$ near the support of $f$. We thus only have to prove 
\begin{eqnarray}  
|| (1-\chi(\epsilon r)) f (P/\epsilon^2) u ||_{L^2 (\Ra ; L^{2^*} ) }  \lesssim ||u_0||_{L^2}, \qquad \epsilon \in (0,1], \ u_0 \in L^2 .  \label{lestimeenonL2}
\end{eqnarray}
Indeed, with (\ref{lestimeenonL2}) (whose right hand side can be replaced by $ || \tilde{f} (P/\epsilon^2) u_0 ||_{L^2} $) at hand, (\ref{quiserarecitee}) yields 
%
\begin{eqnarray*}
 || u_{\rm low} ||_{L^2(\Ra, L^{2^*} )} & \lesssim & \left( \sum_{\epsilon^2 = 2^{-k}} || \tilde{f} (P/\epsilon^2) u_0 ||_{L^2}^2 \right)^{1/2} \\
 & \lesssim & || u_0 ||_{L^2}
 \end{eqnarray*}
by quasi-orthogonality in the second line and  which completes the proof of Theorem \ref{theoremlow} .

\medskip

The rest of this  paragraph is thus devoted to the proof of (\ref{lestimeenonL2}).

\medskip

We write $ (1-\chi)(\epsilon r) f (P/\epsilon^2) = (1-\chi)(\epsilon r) \tilde{f}(P/\epsilon^2) f (P/\epsilon^2)  $ with $ \tilde{f} \in C_0^{\infty}(0,+\infty) $ equal to $1$ on the support of $f$. Then, using Theorem \ref{theoremcaclulfonctionnel}, we can decompose 
\begin{eqnarray}
  (1-\chi)(\epsilon r) \tilde{f}(P/\epsilon^2) = \sum_{\kappa} \tilde{\psi}_{\kappa} (\epsilon r) O \! p_{\epsilon,\kappa} (\chi_{\epsilon,\kappa})^* + {\mathscr R}_{\epsilon} \label{pourchaqueadjoint}
\end{eqnarray}
where, for some $N$ as large as we wish and some bounded family $ (B_{\epsilon})_{\epsilon \in (0,1]} $ of bounded operators on $L^2$, 
$$ {\mathscr R}_{\epsilon}= \zeta (\epsilon r) (P/\epsilon^2 + 1)^{-N} B_{\epsilon} \scal{\epsilon r}^{-N} .
$$
Each $ \chi_{\epsilon,\kappa} = \chi_{\epsilon,\kappa}(\breve{r},\theta,\breve{\rho},\eta) $ belongs to $ \widetilde{S}^{-\infty,0} $, has uniform bounds in $ \epsilon $ and is supported in a way that $ (\breve{r},\theta ) \in \mbox{supp}(1-\chi) \times V_{\kappa} $ and $ p_{\epsilon,\kappa }(\breve{r},\theta,\breve{\rho},\eta) \in \mbox{supp}(f)  $. Furthermore, $  \tilde{\psi}_{\kappa} \equiv 1$ near the support of $ \chi_{\epsilon,\kappa} $.
Note that we use adjoint pseudo-differential operators $ O \! p_{\epsilon,\kappa} (\chi_{\epsilon,\kappa})^* $ (this is possible by (\ref{adjointgloballf})), which is not essential but will be more convenient.

\begin{prop} \label{TTstar1} If $ N \geq n/2 + 1 $, one has
$$ \left( \int_{\Ra} || {\mathscr R}_{\epsilon} f (P/\epsilon^2) e^{-itP} u_0 ||^2_{L^{2^*}} dt \right)^{1/2} \lesssim ||u_0||_{L^2}, \qquad \epsilon \in (0,1], \ u_0 \in L^2 .  $$
\end{prop}

\noindent {\it Proof.}  It follows from Proposition \ref{injectionsdeSobolevfinales} and (\ref{aprioriblackboxlowfrequency}) that
\begin{eqnarray*}
\big| \big| {\mathscr R}_{\epsilon} f (P/\epsilon^2) e^{-i (t-t^{\prime})P} f (P/\epsilon^2) {\mathscr R}_{\epsilon}^* \big| \big|_{L^1 \rightarrow L^{\infty}} & \lesssim & \epsilon^n \big| \big| \scal{\epsilon r}^{-N} e^{-i(t-t^{\prime})P} f^2 (P/\epsilon^2) \scal{\epsilon r}^{-N} \big| \big|_{L^2 \rightarrow L^2} \\ & \lesssim & \epsilon^n \scal{\epsilon^2 (t-t^{\prime})}^{1-N}  \\
& \lesssim & \scal{t-t^{\prime}}^{-n/2} .
\end{eqnarray*}
The result follows then from the $ TT^* $ criterion of \cite{KeelTao} since $  {\mathscr R}_{\epsilon} f (P/\epsilon^2) e^{-itP} $ is bounded on $L^2$ (uniformly in $ \epsilon $ and $t$).
\finpreuve

\bigskip

We are  left with the (rescaled) pseudodifferential terms in (\ref{pourchaqueadjoint}).
For each $ \kappa $ (which we omit in the notation below), we split
\begin{eqnarray}
 \chi_{\epsilon,\kappa} = \chi_{\epsilon,{\rm st}}^+ + \chi_{\epsilon,{\rm int}} + \chi_{\epsilon,{\rm st}}^- ,  \label{decompositiondedepart}
\end{eqnarray} 
with $ \chi_{\epsilon,{\rm st}}^{\pm} ,  \chi_{\epsilon,{\rm int}} \in \widetilde{S}^{-\infty,0} $ (with uniform bounds in $\epsilon$) and supported in  strongly outgoing/incoming areas (see (\ref{referencestrongarealow})), {\it i.e.}
\begin{eqnarray}
 \mbox{supp} \big(  \chi_{\epsilon,{\rm st}}^{\pm} \big) \subset \widetilde{\Gamma}^{\pm}_{\epsilon,{\rm st}} (R,V,I,\varepsilon) 
 \label{fixelevarepsilon}
 \end{eqnarray}
for some $  R \gg 1 $ and $ 0< \varepsilon \ll 1 $ to be chosen below independently of $ \epsilon $, and $ V \Subset V_{\kappa} $, $ I \Subset (0,+\infty) $. Note that to be able to choose $ R $ large, we have to assume that $ (1-\chi) (\breve{r}) $ is supported in $ \breve{r} \geq R $ which is not a restriction since, in (\ref{quiserarecitee}) and Theorem \ref{thlowfreqLq}, we may choose $ \chi \equiv 1 $ on a set as large as we wish.
The third symbol $ \chi_{\epsilon,{\rm int}} $ satisfies
\begin{eqnarray}
 \mbox{supp} \big( \chi_{\epsilon,{\rm int}} \big) \subset \widetilde{\Gamma}^+_{\epsilon} (R,V,I,\sigma) \cap \widetilde{\Gamma}^-_{\epsilon} (R,V,I,\sigma) \label{incomingoutgoing}
\end{eqnarray}
for some $ \sigma $ independent of $ \epsilon $ (see (\ref{weaklowarea}) for the notation of the areas). The decomposition (\ref{decompositiondedepart}) follows easily by applying a partition of unity to $ \breve{\rho} / p_{\epsilon,\kappa} (\breve{r},\theta,\breve{\rho},\eta) $ adapted to regions where this quotient is either lower than $ -1 + \varepsilon^2 $, greater than $ 1-  \varepsilon^2 $ or between $ -1 + \varepsilon^2 /2 $ and $ 1 - \varepsilon^2 /2 $.

\begin{prop} \label{proofpropIK} If $ \varepsilon $ is small enough and $ R  $ is large enough, one has
$$ \left( \int_{\Ra} || \tilde{\psi}_{\kappa} (\epsilon r) O \! p_{\epsilon,\kappa} (\chi_{\epsilon,{\rm st}}^{\pm})^* f (P/\epsilon^2) e^{-itP} u_0 ||^2_{L^{2^*}} dt \right)^{1/2} \lesssim ||u_0||_{L^2}, \qquad \epsilon \in (0,1], \ u_0 \in L^2 .  $$
\end{prop}


\noindent {\it Proof.} We consider the $+$ case. We use again the $TT^*$ criterion and show that
\begin{eqnarray}
\big| \big| \tilde{\psi}_{\kappa} (\epsilon r) O \! p_{\epsilon,\kappa} (\chi_{\epsilon,{\rm st}}^{+})^* f^2 (P/\epsilon^2) e^{-i t P} O \! p_{\epsilon,\kappa} (\chi_{\epsilon,{\rm st}}^{+}) \tilde{\psi}_{\kappa} (\epsilon r) \big| \big|_{L^1 \rightarrow L^{\infty}} \lesssim |t|^{-n/2}, 
\end{eqnarray}
for $t \ne 0 $ and $ \epsilon \in (0,1] $.  Up to taking the adjoint, it suffices to consider  $ t \leq 0 $ (following a trick of \cite{BoucletTzvetkov1}). For simplicity, we let 
$$ K_{\epsilon}^* =  \tilde{\psi}_{\kappa} (\epsilon r) O \! p_{\epsilon,\kappa} (\chi_{\epsilon,{\rm st}}^{+})^* f^2 (P/\epsilon^2)  . $$  We then use Theorem \ref{Isozaki-Kitada-explicite} to expand $  e^{-i t P} O \! p_{\epsilon,\kappa} (\chi_{\epsilon,{\rm st}}^{+}) \tilde{\psi}_{\kappa} (\epsilon r)   $. Consider first the main term $ J_{\epsilon,\kappa} (a_{\epsilon}) e^{-it\epsilon^2 D_x^2} J_{\epsilon,\kappa}(b_{\epsilon,\kappa})^{\dag}$ of this expansion. Using Proposition \ref{statphaseIK} (with $h=1$ and $ s = \epsilon^2 t $) together with Proposition \ref{borneL1Linfiniepsilon} to handle the contribution of the scaling operators, 
we find
\begin{eqnarray}
 \big| \big| J_{\epsilon,\kappa} (a_{\epsilon}) e^{-i \epsilon^2 tD_x^2} J_{\epsilon,\kappa} (b_{\epsilon})^{\dag} \big| \big|_{L^1 \rightarrow L^{\infty}} & \lesssim & \epsilon^{\frac{n}{2}} \scal{\epsilon^2 t}^{-\frac{n}{2}} \epsilon^{\frac{n}{2}} \nonumber \\
 & \lesssim & \scal{t}^{-\frac{n}{2}} . \nonumber
\end{eqnarray}
Note that no sign condition on $t$ is required here. Observing that the support of $ \chi_{\epsilon,{\rm st}}^+ $ allows to write $ O \! p_{\epsilon,\kappa} (\chi^+_{\epsilon,{\rm st}})^* =  O \! p_{\epsilon,\kappa} (\chi^+_{\epsilon,{\rm st}})^* \zeta (\epsilon r)  $, we see that $ || K_{\epsilon}^* ||_{L^{\infty} \rightarrow L^{\infty}} \lesssim 1 $ by Propositions \ref{bornepseudoLq} and \ref{clarificationLinfini},  hence that
\begin{eqnarray}
 \big| \big| K_{\epsilon}^* J_{\epsilon,\kappa} (a_{\epsilon}) e^{-i \epsilon^2 tD_x^2} J_{\epsilon,\kappa} (b_{\epsilon})^{\dag} \big| \big|_{L^1 \rightarrow L^{\infty}} 
 & \lesssim & \scal{t}^{-\frac{n}{2}} . \label{dispersionprequeusuelle}
\end{eqnarray}
We  next consider the first term  of the remainder $ R_{\epsilon,N}(t) $ of (\ref{lowfreqIKexpl}), where $N$ is as large as we wish. It is of the form
$$ e^{-itP}  O_{{\mathscr L}_{-N}^{-2N} \rightarrow {\mathscr L}_{N}^{2N}} \big( 1 \big) =  e^{-itP}\scal{\epsilon r}^{-N} B^{ }_{\epsilon} (P/\epsilon^2 + 1)^{-N} \zeta (\epsilon r), $$
with $ || B_{\epsilon} ||_{L^2 \rightarrow L^2} \lesssim 1 $. To get the time decay, we exploit that this operator is composed to the left with $ K_{\epsilon}^* $ which we can rewrite as
\begin{eqnarray}
 K_{\epsilon}^* = \zeta (\epsilon r) (P/\epsilon^2 + 1)^{-N} \left(  \tilde{\psi}_{\kappa} (\epsilon r) O \! p_{\epsilon,\kappa} (\tilde{\chi}_{\epsilon,{\rm st}}^{+})^* + B^{\prime }_{\epsilon} \scal{\epsilon r}^{-N} \right) f^2 (P/\epsilon^2)  \label{decompositionK*}
\end{eqnarray} 
with $ \tilde{\chi}^{+}_{\epsilon, {\rm st}} \in \widetilde{S}^{-\infty,0} $ with the same support as $ \chi_{\epsilon,{\rm st}}^{+} $ and $ B_{\epsilon}^{\prime} $ bounded on $ L^2 $. This follows simply by expanding $ (P/\epsilon^2 + 1)^N \tilde{\psi}_{\kappa} (\epsilon r) O \! p_{\epsilon,\kappa} (\chi_{\epsilon,{\rm st}}^{+})^*  $.  Then, as in the proof of Proposition \ref{TTstar1},
\begin{eqnarray}
\big| \big| \zeta (\epsilon r) (P/\epsilon^2 + 1)^{-N} B^{\prime }_{\epsilon} \scal{\epsilon r}^{-N}   f^2 (P/\epsilon^2) e^{-itP}\scal{\epsilon r}^{-N} B^{ }_{\epsilon} (P/\epsilon^2 + 1)^{-N} \zeta (\epsilon r) \big| \big|_{L^1 \rightarrow L^{\infty}} 
\lesssim \scal{t}^{-\frac{n}{2}} . \nonumber
\end{eqnarray}
On the other hand, the adjoint estimates of Proposition \ref{cellepourfortementsortant} together with Proposition \ref{injectionsdeSobolevfinales} yield
\begin{eqnarray}
\big| \big| \zeta (\epsilon r) (P/\epsilon^2 + 1)^{-N} \tilde{\psi}_{\kappa} (\epsilon r) O \! p_{\epsilon,\kappa} (\tilde{\chi}_{\epsilon,{\rm st}}^{+})^*  f^2 (P/\epsilon^2) e^{-itP}\scal{\epsilon r}^{-N} B^{}_{\epsilon} (P/\epsilon^2 + 1)^{-N} \zeta (\epsilon r) \big| \big|_{L^1 \rightarrow L^{\infty}} \nonumber \\
  \lesssim  \epsilon^n \big| \big|  \tilde{\psi}_{\kappa} (\epsilon r) O \! p_{\epsilon,\kappa} (\tilde{\chi}_{\epsilon,{\rm st}}^{+})^*  f^2 (P/\epsilon^2) e^{-itP}\scal{\epsilon r}^{-N}  \big| \big|_{L^2 \rightarrow L^{2}}   \lesssim \epsilon^n \scal{\epsilon^2 t}^{-N/3}  \nonumber
\end{eqnarray}
for $ t \leq 0 $. 
Therefore, if $N$ is large enough,
\begin{eqnarray}
 \big| \big|  K_{\epsilon}^* e^{-itP} O_{{\mathscr L}_{-N}^{-2N} \rightarrow {\mathscr L}_{N}^{2N}} \big( 1 \big)  \big| \big|_{L^1 \rightarrow L^{\infty}} \lesssim \scal{t}^{-\frac{n}{2}}, \qquad t \leq 0 .
\end{eqnarray}
It remains to treat the integral terms of $ R_{N,\epsilon}(t) $, involving the operator $ J_{\epsilon,\kappa} (a_{\epsilon,{\rm c}} + r_{\epsilon,N} + \check{a}_{\epsilon}) $. At low frequency, the contribution of $ a_{\epsilon,{\rm c}}  + r_{\epsilon,N} $ follows only from its spatial decay (see the slight difference with the high frequency case in paragraph \ref{minormodification}). We thus only exploit that
$$ J_{\epsilon,\kappa} (a_{\epsilon,{\rm c}}) + J_{\epsilon,\kappa} (r_{\epsilon,N}) := \scal{\epsilon r}^{-N} J_{\epsilon,\kappa} (\tilde{a}_{\epsilon,N}), $$
for some  bounded family of symbols $ (\tilde{a}_{\epsilon,N})_{\epsilon} $ in $ S_{0} $, supported in $ \Theta^+ (R_0,V_0,I_0,\varepsilon_0) $ with $ R_0  $ as large as we wish by taking $ R $  large enough. To estimate the contribution of this term in $ R_{N , \epsilon}(t) $, we use the estimate
\begin{eqnarray*} \big| \big|  K_{\epsilon}^* e^{-i (t - \frac{s}{\epsilon^2})P} \scal{\epsilon r}^{-N} J_{\epsilon,\kappa} (\tilde{a}_{\epsilon,N}) e^{-isD_x^2} J_{\epsilon,\kappa} (b_{\epsilon})^{\dag} \big| \big|_{L^1 \rightarrow L^{\infty}}  \lesssim \epsilon^n \scal{\epsilon^2 t - s}^{-\frac{N}{6}} \scal{s}^{-\frac{N}{2}} 
\end{eqnarray*}
for $ t \leq t - \frac{s}{\epsilon^2} \leq 0 $ and which, after integration in $s$, provides an upper bound by $ \scal{t}^{-\frac{n}{2}} $ if $  N$ is chosen large enough. To get the above estimate, we use on one hand that
$$  \big| \big|  K_{\epsilon}^* e^{-i (t - \frac{s}{\epsilon^2})P} \scal{\epsilon r}^{-N/2} \big| \big|_{L^2 \rightarrow L^{\infty}}  \lesssim \epsilon^{\frac{n}{2}} \scal{\epsilon^2 t - s}^{-\frac{N}{6}}  $$
by using the decomposition (\ref{decompositionK*}) together with the propagation estimates given by (\ref{aprioriblackboxlowfrequency}) and (the adjoint estimates of) Proposition \ref{cellepourfortementsortant}.  On the other hand, we use 
$$ \big| \big|   \scal{\epsilon r}^{-N/2} J_{\epsilon,\kappa} (\tilde{a}_{\epsilon,N}) e^{-isD_x^2} J_{\epsilon,\kappa} (b_{\epsilon})^{\dag} \big| \big|_{L^1 \rightarrow L^{2}}  \lesssim \epsilon^{\frac{n}{2}}  \scal{s}^{-\frac{N}{2}}  $$
which  comes from Proposition \ref{propparam} for the time decay, up to the replacement of the source space $ L^2 $ by $ L^1 $  which provides the additional $ \epsilon^{n/2} $ factor. This replacement is  possible by writing
$  J_{\epsilon,\kappa} (b_{\epsilon})^{\dag} = J_{\epsilon,\kappa} (\tilde{b}_{\epsilon})^{\dag}  (P/\epsilon^2 + 1)^{-N}  \zeta (\epsilon r) $ for some $ \tilde{b}_{\epsilon} $ with the same properties as $ b_{\epsilon} $ (it is obtained by computing $  J_{\epsilon,\kappa} (b_{\epsilon})^{\dag} (P/\epsilon^2 + 1)^{N} = J_{\epsilon,\kappa} (\tilde{b}_{\epsilon})^{\dag}  $) and by using  Proposition \ref{injectionsdeSobolevfinales}.

The last term of $ R_{\epsilon,N}(t) $ to consider is the one containing $ J_{\epsilon,\kappa} (\check{a}_{\epsilon}) $. Here the crucial observation is that $ |\theta - \vartheta| $ is bounded below on the support of $ \check{a}_{\epsilon} $. In particular, using (\ref{pourdesconditionsdesupport}) we see that  $ |\partial_{\theta} \varphi_{\epsilon}| / r \partial_r \varphi_{\epsilon} $ is bounded from below on the support of $ \check{a}_{\epsilon} $, which implies that $ (\breve{r},\theta,\partial_{\breve{r}} \varphi_{\epsilon},\partial_{\theta} \varphi_{\epsilon}) $ must belong to an incoming area. More precisely, according to (\ref{pourphasenonstatbassefreq}) and (\ref{pourdesconditionsdesupport}), we must have $ \partial_{\breve{r}} \varphi_{\epsilon} < \sigma_1 p_{\epsilon,\kappa}(\breve{r},\theta,\breve{r} \varphi_{\epsilon},\partial_{\theta} \varphi_{\epsilon})^{1/2} $  on the support of $ \check{a}_{\epsilon} $ with $ \sigma_1 = 1 - \varepsilon_1^2 / C $ independent of $ \varepsilon $ ({\it i.e.} of $ \varepsilon_2 $ in Theorem \ref{Isozaki-Kitada-explicite}).    Thus, using
 Lemma \ref{lemmedexpansionpseudoOIF} we can replace $ J_{\epsilon,\kappa} (\check{a}_{\epsilon})  $ by $ O\!p_{\epsilon,\kappa} (\tilde{\chi}_{\epsilon}^-) J_{\epsilon,\kappa} (\check{a}_{\epsilon}) $ with $ \tilde{\chi}_{\epsilon}^- $ supported in an incoming region,  up to decaying remainders that can be treated as before. We can then proceed as before by using the adjoint a priori estimate of Proposition \ref{improvedoutgoingII} (since one can choose $ \varepsilon $ as small as we want, without affecting the value of $ \sigma_1 $ above) which provides the estimate
$$ \big| \big|  K_{\epsilon}^* e^{-i (t - \frac{s}{\epsilon^2})P} O \! p_{\epsilon,\kappa} (\tilde{\chi}^-_{\epsilon}) J_{\epsilon,\kappa} (\check{a}_{\epsilon}) e^{-isD_x^2} J_{\epsilon,\kappa} (b_{\epsilon})^{\dag} \big| \big|_{L^1 \rightarrow L^{\infty}} \lesssim \epsilon^n \scal{s}^{-\frac{n}{2}} \scal{\epsilon^2 t - s}^{-N} $$
and then the final estimate by $ \scal{t}^{-n/2} $ after integration in $s$. The result follows.
\finpreuve

\bigskip

To complete the proof of (\ref{lestimeenonL2}), it remains to study the contribution of $ \chi_{\epsilon,{\rm int}} $ in (\ref{decompositiondedepart}). We follow the idea of \cite{BoucletAPDE,MizutaniCPDE}, by adapting it to the low frequency and global in time case.

 Everywhere below, we choose $ t_0 > 0 $ small enough as in Proposition \ref{WKB}. Also, the parameter $ \varepsilon $ used in (\ref{fixelevarepsilon}) (and hence the parameter $ \sigma$ in (\ref{incomingoutgoing})) is chosen according to Proposition \ref{proofpropIK}. We then choose $ \delta > 0 $ small enough, according to the third item of  Proposition \ref{propagationclassiqueprop}, and we split $ \chi_{\epsilon,{\rm int}} $ as a sum
 $$ \chi_{\epsilon,{\rm int}} = \sum_{j \in J} \chi_{\epsilon,j}, \qquad \mbox{supp} (\chi_{\epsilon,j}) \subset \mbox{supp} (\chi_{\epsilon,{\rm int}}) \cap  \left\{ \frac{\breve{\rho}}{p_{\epsilon}^{1/2}} \in (j \delta , (j+1) \delta)  \right\} , $$
  where $ J $ is a finite subset of $ \Za $ (depending on $ \delta $) and $ (\chi_{\epsilon,j})_{\epsilon} $ is a bounded family of $ \tilde{S}^{-\infty,0} $.  It now suffices to prove global in time dispersion estimates, say for $ t \geq 0 $, for the operators
\begin{eqnarray}  
 \tilde{\psi}_{\kappa} (\epsilon r) O \! p_{\epsilon,\kappa} (\chi_{\epsilon,j})^* f^2 (P/\epsilon^2) e^{-itP} O \! p_{\epsilon,\kappa}(\chi_{\epsilon,j}) \tilde{\psi}_{\kappa} (\epsilon r) ,  \label{operateuradisperser}
\end{eqnarray} 
uniformly in $ \epsilon $. To do so, we introduce a spatial partition of unity on the support of the symbols, 
$$ 1 = \sum_{\ell \geq \ell_0} \phi  (\breve{r}/ R_\ell), \qquad R_\ell = 2^\ell, \qquad \phi \in C_0^{\infty} (0,\infty) $$
and define
$$ \chi_{\epsilon,j}^{(\ell)} (\breve{r},\theta,\breve{\rho},\eta) = \phi (\breve{r}/R_\ell) \chi_{\epsilon,j} (\breve{r},\theta,\breve{\rho},\eta) . $$
Picking $ \tilde{\phi} \in C_0^{\infty} (0,\infty) $ equal to $1$ near the support of $ \phi $ and using that $ 1 - \tilde{\phi} (\breve{r}/R_\ell) $ vanishes near the support of $ \chi_{\epsilon,j}^{(\ell)} $, we obtain by symbolic calculus that, for any given $N$,
\begin{eqnarray}
  O \! p_{\epsilon,\kappa} \big(\chi_{\epsilon,j}^{(\ell)} \big) \tilde{\psi}_{\kappa} (\epsilon r) =  O \! p_{\epsilon,\kappa}\big(\chi_{\epsilon,\delta}^{(\ell)} \big) \tilde{\psi}_{\kappa} (\epsilon r) \chi (\epsilon r /R_\ell) + \scal{\epsilon r}^{-N} B (\epsilon,R_\ell) (P/\epsilon^2 + 1)^{-N} \zeta (\epsilon r) , \label{pourreferenceaureste}
\end{eqnarray}
where, uniformly in $ \epsilon $,
$$ || B (\epsilon,R_\ell) ||_{L^2 \rightarrow L^2}  \lesssim R_\ell^{-N} .$$
The contribution of the remainder term of (\ref{pourreferenceaureste}) can be treated as the remainders in the above proof of Proposition \ref{proofpropIK} by propagation estimates and we get
 $$ \big| \big| \tilde{\psi}_{\kappa} (\epsilon r) O \! p_{\epsilon,\kappa} (\chi_{\epsilon,j})^* f^2 (P/\epsilon^2) e^{-itP}   \scal{\epsilon r}^{-N} B (\epsilon,R_\ell) (P/\epsilon^2 + 1)^{-N} \zeta (\epsilon r)  \big| \big|_{L^1 \rightarrow L^{\infty}} \lesssim \scal{t}^{-n/2} R_\ell^{-N} $$
 for all $ t \geq 0 $ (actualy this holds for all $ t \in \Ra $ since $ \chi_{\epsilon,j} $ is both incoming and outgoing by (\ref{incomingoutgoing})). These estimates can be easily summed over $k$. On the other hand, using the general fact that
$$ \left| \left|\sum_\ell A_\ell \tilde{\phi}(\epsilon r / R_\ell) v \right| \right|_{L^{\infty}}  \leq \left(  \sup_\ell || A_\ell \tilde{\phi}(\epsilon r / R_\ell) ||_{L^1 \rightarrow L^{\infty}} \right) \sum_\ell \int_{  \frac{\epsilon r}{R_\ell} \in {\rm supp} \tilde{\phi} } |v|$$ 
  where the last sum is bounded above by $ C || v ||_{L^1} $ (with $ C $ independent of $ \epsilon $ and $k$), we see that the dispersion estimate for (\ref{operateuradisperser}) is a  consequence of the following uniform estimates.
\begin{prop} There exists $ C > 0 $ such that for all $ \ell \geq \ell_0 $, all $ \epsilon \in (0,1] $ and all $ t \geq 0 $,
\begin{eqnarray}
 \left| \left| \tilde{\psi}_{\kappa} (\epsilon r) O \! p_{\epsilon,\kappa} (\chi_{\epsilon,j})^* f^2 (P/\epsilon^2) e^{-itP} O \! p_{\epsilon,\kappa}\big(\chi_{\epsilon,j}^{(\ell)} \big) \tilde{\psi}_{\kappa} (\epsilon r) \tilde{\phi} (\epsilon r / R_\ell) \right| \right|_{L^1 \rightarrow L^{\infty}} \leq C \scal{t}^{-\frac{n}{2}} .  \label{pourallegernotation}
\end{eqnarray}
\end{prop}



\noindent {\it Proof.} For $ 0 \leq \epsilon^2 t \leq t_0 R_\ell $, the estimate follows from Proposition \ref{WKB} together with the fact that
$$ ||  \tilde{\psi}_{\kappa} (\epsilon r) O \! p_{\epsilon,\kappa} (\chi_{\epsilon,j})^* f^2 (P/\epsilon^2) ||_{L^{\infty} \rightarrow L^{\infty}} \lesssim 1, \qquad ||  \tilde{\psi}_{\kappa} (\epsilon r) O \! p_{\epsilon,\kappa} (\chi_{\epsilon,j})^* f^2 (P/\epsilon^2) ||_{L^2 \rightarrow L^{\infty}} \lesssim \epsilon^{n/2} , $$
the second estimate being used to treat the remainder term of the parametrix of Proposition \ref{WKB}, which provides an $ L^1 \rightarrow L^{\infty} $ estimate by $ \epsilon^n R_\ell^{-n} \lesssim \scal{t}^{-n/2} $.
Then, for $ t \geq \epsilon^{-2} t_0 R$, we use $L^2$ propagation estimates as follows. First, we  write for an arbitrary $N>0$,
$$ O \! p_{\epsilon,\kappa}\big(\chi_{\epsilon,j}^{(\ell)} \big) \tilde{\psi}_{\kappa} (\epsilon r) \tilde{\phi} (\epsilon r / R_\ell)  =    \left( O \! p_{\epsilon,\kappa} ( \tilde{\chi}_{\epsilon,j}^{(\ell)}) \tilde{\psi}_{\kappa} (\epsilon r)  + \scal{\epsilon r}^{-N} \tilde{B}(\epsilon,R_\ell ) \right) (P/\epsilon^2 + 1)^{-N} \zeta (\epsilon r)  $$
with $ ||  \tilde{B}(\epsilon,R_\ell ) ||_{L^2 \rightarrow L^2} \lesssim R_\ell^{-N} $ and $ (\tilde{\chi}_{\epsilon,j}^{(\ell)} )_{\epsilon,\ell}  $ bounded in $ \tilde{S}^{-\infty,0} $ with the same support as $ \chi_{\epsilon,j}^{(\ell)} $. This is obtained by expanding
$ O \! p_{\epsilon,\kappa}\big(\chi_{\epsilon,j}^{(\ell)} \big) \tilde{\psi}_{\kappa} (\epsilon r) \tilde{\phi} (\epsilon r / R_\ell) (P/\epsilon^2 + 1)^N  $. Then the contribution of  the term involving $ \tilde{B}(\epsilon,R_\ell) $ is similar to the one of the remainder of (\ref{pourreferenceaureste}) and provides a $ L^1 \rightarrow L^{\infty} $ estimate by $ R_\ell^{-N} \scal{t}^{-n/2} $. We are thus left with the contribution of $ \tilde{\chi}^{(\ell)}_{\epsilon,j} $.
For this term, we distinguish between two cases
$$ t_0 R_\ell \leq \epsilon^2 t \leq T R_\ell , \qquad \epsilon^2 t > T R_\ell $$
with $ T > 0 $ large enough (independent of $ \epsilon $ and $\ell$)  chosen according to the item 2 of Proposition \ref{propagationclassiqueprop}, namely such that the support of $ \tilde{\chi}_{\epsilon,j}^{(\ell)} $ is mapped into a stronly outgoing region by the classical flow at time $T$.
Indeed, for $ \epsilon^2 t > T R_\ell $, we can write the contribution of $ \tilde{\chi}^{(\ell)}_{\epsilon,j}  $ to the estimate (\ref{pourallegernotation}), as the one of
$$   \tilde{\psi}_{\kappa} (\epsilon r) O \! p_{\epsilon,\kappa} \big(\chi_{\epsilon,j} \big)^* f^2 (P/\epsilon^2) e^{-i \big( t - \frac{T R_{\ell}}{\epsilon^2} \big) P} \left( e^{-i\frac{T R_{\ell}}{\epsilon^2} P}  O \! p_{\epsilon,\kappa} \big( \tilde{\chi}_{\epsilon,j}^{(\ell)} \big) \tilde{\psi}_{\kappa} (\epsilon r)  e^{i \frac{T R_{\ell}}{\epsilon^2} P} \right) e^{- i \frac{T R_{\ell}}{\epsilon^2} P} (P/\epsilon^2 + 1)^{-N} \zeta (\epsilon r) .  $$
Using Proposition \ref{Egorovquantifie}, we can write for any given $N$ the parenthese as a sum (over angular charts $ \kappa_2 $) of operators of the form
$$R_\ell^{-N} O \! p_{\epsilon,\kappa_2} ( \hat{\chi}_{\epsilon,{\rm st},\kappa_2}^{(\ell)}) \tilde{\psi}_{\kappa_2} (\epsilon r) \scal{\epsilon r}^\ell + \scal{\epsilon r}^{-N} O_{L^2 \rightarrow L^2} (R^{-N}_\ell)  $$
with $ (\hat{\chi}_{\epsilon,{\rm st},\kappa_2}^{(\ell)} )_{\epsilon,\ell} $ bounded in $ \widetilde{S}^{-\infty,0} $ and supported in a an outgoing region with parameter $ \sigma^{\prime} $ as close to $1$ as we wish, hence in particular disjoint from the support of $ \chi_{\epsilon,j} $.  Using Propositions \ref{halflocalized} and \ref{improvedoutgoingII}, we get a dispersion estimate of order $  \epsilon^n R^{-N}_\ell \scal{\epsilon^2 t - T R_\ell}^{-N} \lesssim \scal{t}^{-\frac{n}{2}}  $.
Finally, for $ t_0 R_\ell \leq \epsilon^2 t  \leq T R_\ell $, we  write the contribution of $ \tilde{\chi}^{(\ell)}_{\epsilon,j}  $ to the estimate (\ref{pourallegernotation}), as the one of
$$   \tilde{\psi}_{\kappa} (\epsilon r) O \! p_{\epsilon,\kappa} (\chi_{\epsilon,j})^* f^2 (P/\epsilon^2)  \left( e^{-i tP}  O \! p_{\epsilon,\kappa} ( \tilde{\chi}_{\epsilon,j}^{(\ell)}) \tilde{\psi}_{\kappa} (\epsilon r)  e^{i t P} \right) e^{- i t P} (P/\epsilon^2 + 1)^{-N} \zeta (\epsilon r)  . $$
By Theorem \ref{Egorovquantifie} together with the third item Proposition \ref{propagationclassiqueprop} and our choice of $ \delta $, the parenthese is microlocalized in a set where $ \breve{\rho} / p_{\epsilon,\kappa}^{1/2} > (j+1) \delta $, hence disjoint from the support of $ \chi_{\epsilon,j} $. Thus, only residual terms contribute and they produce a norm of order $ \epsilon^n R_\ell^{-\infty} = O( \scal{t}^{-n/2} ) $ since $ \epsilon^2 t $ is of order $ R_\ell $ in this case. This completes the proof. \finpreuve




\subsection{Proof of Theorem \ref{theohighinfty}} \label{minormodification}

Here the analysis is very similar to the one of \cite{MizutaniCPDE}, the main difference  being that we control the remainder terms globally in time. This is done  as in the low frequency case, by using the high frequency propagations estimates of paragraph \ref{soussectionpropagation}, so the techniques are the same. We only record here that the estimates of paragraph \ref{soussectionpropagation} are not sensitive to a possible trapping since the moderate growth in $ \lambda \sim h^{-2} $ in (\ref{aprioriresolvente}) is controlled by the large powers of $h$ provided by the remainders in the expansions (the a priori resolvent estimates are only used to control the remainders).


We also mention the following minor technical point in the transposition of the proof of Proposition \ref{proofpropIK} to high frequencies. In the remainder $ R_N^h (t) $ of the high frequency Isozaki-Kitada parametrix (see after (\ref{highfreqIKexpl})) neither $ a_{\rm c}^h $ nor $ \check{a}^h $ decay  in $h$, so it is not clear that they will have a negligible contribution in the end. To make sure they are negligible in the derivation of dispersion bounds, we have to make sure that these terms have a $ O (h^{\infty}) $ contribution. For $ a_{\rm c}^h $ this follows by using Proposition \ref{fullyone}. The contribution of $ \check{a}^h $ is handled by   the propagation estimates of Proposition \ref{improvedoutgoingII} which provide the fast decay in $h$.

\subsection{Proof of Theorem \ref{globalestimates}}
Thanks to Theorem \ref{theohighinfty}, it suffices to prove that for any given $ \chi \in C_0^{\infty} ({\mathcal M}) $, one has the global Strichartz estimates
$$ || \chi u_{\rm hi} ||_{L^p (\Ra,L^q ({\mathcal M}))} \lesssim || u_0 ||_{L^2} . $$ 
This follows from the technique of \cite{Burq,StTa}. In the non trapping case, this is a classical fact. For hyperbolic trapping, the analysis is detailed in \cite{BGH} for local in time estimates, but it  holds also globally in time.
Note also that we are allowed to use  Theorem \ref{theohighinfty} since, under the assumptions of Theorem \ref{globalestimates}, the resolvent has high energy bounds growing at worst like $ \lambda^{-1/2} \log \lambda $ (see \cite{NZ}).

 \section{Nonlinear equations} \label{sectionnonlineaire}
 \setcounter{equation}{0}
 
 In this paragraph, we use the global Strichartz inequalities of Theorem \ref{globalestimates} to study  the $ L^2 $ critical nonlinear Schr\"odinger equation 
 \begin{equation}
 i \partial_t u - P u = \sigma |u|^{\frac{4}{n}} u \tag{NLS} 
\end{equation}  
 where $ n \geq 3 $ is the space dimension and $ \sigma $ is a sign; $ \sigma = 1 $ corresponds to the defocusing case and $ \sigma = - 1 $ to the focusing case. Here the sign will not  matter since we are going to consider small data. 
 We will solve (NLS) in
$$ X := L^{2 + \frac{4}{n}} (\Ra \times {\mathcal M}) \cap C_{\rm scat}(\Ra, L^2 ({\mathcal M}) ) $$ 
 where
$$ C_{\rm scat}(\Ra, L^2 ({\mathcal M}) ) = \left\{ u \in C (\Ra,L^2 ({\mathcal M})) \ | \ \mbox{the limits} \ \lim_{t \rightarrow \pm \infty} e^{itP} u (t) \ \mbox{exist in } \  L^2 ({\mathcal M})  \right\}  $$
is a Banach space for the norm $ ||u||_{L^{\infty}L^2} := \sup_{t \in \Ra}||u(t)||_{L^2({\mathcal M})} $  (it is a closed subspace of the space of bounded uniformly continuous  functions $ u : \Ra \rightarrow L^2 ({\mathcal M}) $). We then equip $X$ with the norm
$$ ||u||_X = ||u||_{L^{2 + \frac{4}{n}}(\Ra \times {\mathcal M})} + ||u||_{L^{\infty} L^2} , $$
which makes it a Banach space. 
\begin{theo} \label{theoremenonlineaire} Let $ \sigma = 1 $ or $ -1 $. Under the assumptions of Theorem \ref{globalestimates}, there exists $ \varepsilon > 0 $ such that, for all $ u_0 \in L^2 ({\mathcal M}) $ satisfying $ || u_0 ||_{L^2} < \varepsilon $, there exists a unique $ u \in X $ such that
$$ u (0) = u_0 \qquad \mbox{and} \qquad u \ \mbox{solves (NLS) in the distributions sense}. $$
In particular, since it belongs to $  C_{\rm scat}(\Ra,L^2 ({\mathcal M})) $, this solution scatters as $ t \rightarrow \pm \infty $, i.e. there are $ u_{\pm} \in L^2 ({\mathcal M}) $ such that
$$ || u (t) - e^{-itP} u_{\pm} ||_{L^2} \rightarrow 0, \qquad t \rightarrow \pm \infty . $$
\end{theo}

This theorem is of course  similar to the well known result for (NLS) on $ \Ra^n $. Its novelty stems in the fact that we work on an asymptotically conical manifold and that a possible hyperbolic trapping on $ {\mathcal M} $ will not change the usual picture, namely the global well posedness and the existence of scattering for small data.
 
The proof follows the usual scheme, the main tool being the global Strichartz estimates. We record the main lines below to point out the where one has to be careful in the transposition of the proof on $ \Ra^n $.

\bigskip

\noindent {\it Proof of Theorem \ref{theoremenonlineaire}.} The principle is to solve (NLS) in the Duhamel form
\begin{equation}
 u (t) = e^{-itP} u_0 + \frac{\sigma}{i} \int_0^t e^{-i(t-s)P} |u(s)|^{\frac{4}{n}} u(s)ds ,  \tag{Duh}
\end{equation} 
by a fixed point argument on a ball $ \overline{B}_X (0,r) $ with $r$ small enough.
We note first that the pair $ (p,q) $ defined by  $ p = q = 2 + \frac{4}{n} $
is Schr\"odinger admissible, so the homogeneous Strichartz inequalities of Theorem \ref{globalestimates} show that the map
$$ U: L^2 ( {\mathcal M} ) \ni u_0 \mapsto [ t \mapsto e^{-itP} u_0 ] \in X $$
is well defined and that one has
$$ U \left( B_{L^2} (0,\varepsilon) \right) \subset \overline{B}_X (0, C \varepsilon) .$$
 Also, since $ (p,q) $ is not an endpoint pair (i.e. $ p \ne 2 $), the homogeneous inequalities provide inhomogeneous Strichartz inequalities thanks to the Christ-Kiselev lemma \cite{ChKi}. This means that, if we set
\begin{eqnarray} 
 (D f)(t) := \int_0^t e^{-i (t-s)P} f (s) ds ,  \label{integralform}
\end{eqnarray} 
we have
\begin{equation}
 || Df ||_{L^{2 + \frac{4}{n}}(\Ra \times {\mathcal M})} \leq C || f ||_{L^{\frac{2n+4}{n+4}} (\Ra \times {\mathcal M})} , \label{inStr}
\end{equation} 
where $ \frac{2n+4}{n+4} $ is the conjugate exponent to $ 2 + \frac{4}{n} $. More precisely, the integral defining $Df$ has a clear sense if $ f \in C (\Ra,L^2 ({\mathcal M})) $ so the precise meaning of (\ref{inStr}) is that it holds on the dense subset  $C (\Ra,L^2 ({\mathcal M})) \cap L^{\frac{2n+4}{n+4}} (\Ra \times {\mathcal M}) $ and that  $ D $ can then be extended by density to $  L^{\frac{2n+4}{n+4}} (\Ra \times {\mathcal M}) $. The adjoint estimates to the the homogeneous Strichartz estimates also imply that
$$  || Df ||_{L^{\infty} L^2} \leq C || f ||_{L^{\frac{2n+4}{n+4}} (\Ra \times {\mathcal M})} , $$
and  that
$$ \big| \big| e^{it P}  (Df) (t) - e^{ i t^{\prime} P}  (Df) (t^{\prime}) \big| \big|_{L^2({\mathcal M})} = \left| \left| \int_{t^{\prime}}^t e^{is P} f (s) ds \right| \right|_{L^2 ({\mathcal M})} \leq C ||f||_{L^{\frac{2n+4}{n+4}}([t^{\prime},t] \times {\mathcal M})} , $$
for all $f \in  C (\Ra,L^2 ({\mathcal M}))  \cap L^{\frac{2n+4}{n+4}} (\Ra \times {\mathcal M})  $. This last inequality
implies that $ e^{itP} (Df)(t) $ has  limits as $ t \rightarrow \pm \infty $ hence that $ D f $ belongs to $ C_{\rm scat} (\Ra,L^2({\mathcal M})) $.
 Thus
$$  D :   L^{\frac{2n+4}{n+4}} (\Ra \times {\mathcal M}) \rightarrow X $$
is well defined and continuous, by taking the closure of $ D :C (\Ra,L^2({\mathcal M})) \cap L^{\frac{2n+4}{n+4}} (\Ra \times {\mathcal M}) \rightarrow X$. One has however to be careful that the closure of $D$ is no longer clearly given by the explicit integral form (\ref{integralform}).


To handle the nonlinearity $ u \mapsto N(u):= |u|^{\frac{4}{n}} u $,
we use the estimate on complex numbers
\begin{eqnarray}
 \big||z|^{\frac{4}{n}}z - |\zeta|^{\frac{4}{n}} \zeta\big| \leq C_n |z - \zeta| \big( |z|^{\frac{4}{n}} + |\zeta|^{\frac{4}{n}} \big) , 
 \label{forzzeta}
\end{eqnarray}
to derive the estimate
$$ \big| \big| N (u) - N (v) \big| \big|_{ L^{\frac{2n+4}{n+4}} (\Ra \times {\mathcal M}) } \leq C_n || u - v ||_{L^{2 + \frac{4}{n}}(\Ra \times {\mathcal M})} \left( || u ||_{L^{2 + \frac{4}{n}}(\Ra \times {\mathcal M})}^{\frac{4}{n}} + ||v||_{L^{2 + \frac{4}{n}}(\Ra \times {\mathcal M})}^{\frac{4}{n}} \right)$$
which implies in particular that 
$$ N :  X \subset  L^{2 + \frac{4}{n}} (\Ra \times {\mathcal M})  \rightarrow L^{\frac{2n+4}{n+4}} (\Ra \times {\mathcal M}) $$
is well defined and Lipschitz on balls of $X$. The above estimate with $v=0$ also implies that
$$ N \big( \overline{B}_X (0,r) \big) \subset B_{L^{ \frac{2n+4}{n+4}} (\Ra \times {\mathcal M})} \left(0, C_n r^{1 + \frac{4}{n}} \right) . $$
We can thus define the map $F_{u_0} : X \rightarrow X $ by
$$ F_{u_0} (u)  = U( u_0 ) + \frac{\sigma}{i} D (N(u))  $$
which gives a precise sense to the right hand side of (Duh). Furthermore, for $ u , v \in \overline{B}_X (0,r) $ and $ u_0 \in B_{L^2}(0,\varepsilon) $, one has
$$ || F_{u_0} (u) ||_X \leq || U (u_0) ||_{X} + ||D(N(u))||_X \leq C \varepsilon + C r^{1 + \frac{4}{n}}$$
and
\begin{eqnarray}
 || F_{u_0} (u) - F_{u_0} (v) ||_X = || D (N(u) - N(v)) ||_X \leq C r^{\frac{4}{n}} || u - v ||_X  , \label{inegcontinuite}
\end{eqnarray}
so, if $r$ is small enough and $ \varepsilon \ll r $,  the ball $ \overline{B}_X (0,r) $ is stable by $F_{u_0} $ on which it is a contraction. This provides a solution to the equation $ u = F_{u_0}(u) $. To complete the proof, one has to observe that this solution is a solution in the distributions sense and, conversely, that if we have a distributional solution which belongs to $X$ then it satisfies $ F_{u_0} (u) = u $.

To prove these two facts, we will use that,
 if $ \chi \in C_0^{\infty} (\Ra) $ is equal to $1$ near $0$,  then for every given $u \in X $
\begin{eqnarray}
 \chi (2^{-j} P) u \rightarrow u \qquad \mbox{in} \ X \ \mbox{as} \ j \rightarrow \infty. \label{convergenceutilenonlineaire}
\end{eqnarray}
Here $ \chi (2^{-j} P) u = [t \mapsto \chi (2^{-j} P) u(t)] $.
The convergence (\ref{convergenceutilenonlineaire}) follows from the  strong convergence of $ \chi (2^{-j} P) $ to the identity on both $ L^2 ({\mathcal M}) $ and $ L^{2+\frac{4}{n}} ({\mathcal M}) $, which can be proved as on $ \Ra^n $ for Fourier multipliers by using the pseudo-differential description of $ \chi (2^{-j} P) $. We omit the details of the proof but only record that to prove
$$ \sup_{t \in \Ra} || u (t) - \chi (2^{-j}P ) u (t) ||_{L^2 ({\mathcal M})} \rightarrow 0, \qquad j \rightarrow \infty $$
we may replace the norm by $  || e^{itP }u (t) - \chi (2^{-j}P ) e^{itP} u (t) ||_{L^2 ({\mathcal M})}   $ and exploit that $ t \mapsto e^{itP}u(t) $ is uniformly continuous with limits at $ \pm \infty $ to get the uniform convergence as $ j \rightarrow \infty $. 
Thus, given a solution $u$ to $ u = F_{u_0} (u) $ and letting $u_j =\chi(2^{-j}P)$, one has $ F_{u_0} (u_j) \rightarrow F_{u_0} (u) = u$ by   (\ref{inegcontinuite}) and (\ref{convergenceutilenonlineaire}).
Since  $ |u_j|^{\frac{4}{n}} u_j  $ belongs to $ C (\Ra,L^2 ({\mathcal M}))$ (this can be checked by using (\ref{forzzeta}) and  that $ \chi (2^{-j} P) $ maps $ L^2 ({\mathcal M}) $ into $ L^{\infty}({\mathcal M}) $), we can write
$$ F_{u_0} (u_j) (t)  = e^{-itP} u_0 + \frac{\sigma}{i} \int_0^t e^{-i(t-s)P} |u_j(s)|^{\frac{4}{n}} u_j (s) ds $$
(i.e. the integral has a clear sense) and, from this expression, we easily infer that
$$ ( i \partial_t - P ) F (u_j) = \sigma |u_j|^{\frac{4}{n}} u_j $$
in the distributions sense on $ \Ra \times {\mathcal M} $. Letting $j \rightarrow \infty $, we conclude that $u$ solves (NLS) in the distributions sense. 

Conversely, if $ u \in X $ solves (NLS) in the distributions sense, it remains to prove that $ u = F_{u_0}(u) $. By definition, we have
\begin{eqnarray}
 \int_{\Ra} \int_{\mathcal M}  \overline{ (i \partial_t - P ) \phi (t,x)} u (t,x) d {\rm vol}_g dt = \sigma \int_{\Ra} \int_{\mathcal M} \overline{\phi(t,x)} |u(t,x)|^{\frac{4}{n}} u (t,x)  d {\rm vol}_g dt  \label{pourdistfort}
\end{eqnarray}
for all $ \phi \in C_0^{\infty} (\Ra \times {\mathcal M}) $ and then for all $ \phi \in C_0^{\infty} (\Ra, {\mathscr S}({\mathcal M})) $ by a simple limiting argument (see (\ref{definitiondeSdeM}) for $ {\mathscr S}({\mathcal M}) $). The interest of allowing $ \phi (t) = \phi (t,.) $ to belong to $ {\mathscr S}( {\mathcal M}) $, is that we can write the left hand side of (\ref{pourdistfort}) as
$$ \int_{\Ra}  \left( i\partial_t (e^{itP} \phi (t)) , e^{itP} u (t)\right)_{L^2 ({\mathcal M})} dt , $$
since $ e^{itP} $ leaves $ {\mathscr S} ({\mathcal M}) $ stable but not $ C_0^{\infty} ({\mathcal M}) $. On the other hand, by approximating $ u $ by $ u_j = \chi (2^{-j} P) u $ using (\ref{convergenceutilenonlineaire}), the right hand side of (\ref{pourdistfort}) reads
$$  \sigma \int_{\Ra} \int_{\mathcal M} \overline{\phi(t,x)} |u_j(t,x)|^{\frac{4}{n}} u_j (t,x)  d {\rm vol}_g dt + O \big( ||u- u_j||_{X} \big)  . $$
Using that $ t \mapsto |u_j(t)|^{\frac{4}{n}} u_j (t) $ is continuous with values in $ L^2 ({\mathcal M}) $, one can write
$$  \sigma  |u_j(t)|^{\frac{4}{n}} u_j (t) = e^{-it P}  i \partial_t \left( \frac{\sigma}{i} \int_0^{t} e^{is P} |u_j(s)|^{\frac{4}{n}} u_j(s) ds  \right). $$
Then, by integration by part, (\ref{pourdistfort}) yields
$$  \int_{\Ra}  \left( i\partial_t (e^{itP} \phi (t)) , e^{itP} u (t) -  G (N(u_j)) (t) \right)_{L^2 ({\mathcal M})} dt =  O \big( || u - u_j ||_X \big)  , $$
where
$$ G (f) (t) :=  \frac{\sigma}{i} \int_0^t e^{isP} f (s) ds  $$
is well defined for $ f \in C (\Ra,L^2 ({\mathcal M})) $ with values on $ C (\Ra,L^2 ({\mathcal M})) $  but can be extended to all $ f \in L^{\frac{2n+4}{n+4}} (\Ra \times {\mathcal M}) $ by the adjoint of homogeneous Strichartz estimates. Letting $ j \rightarrow \infty $ and choosing $ \phi (t) = e^{-itP} \psi (t) $ with $ \psi \in C_0^{\infty} (\Ra \times {\mathcal M}) $, we find that
$$ \int_{\Ra} \int_{\mathcal M} \overline{i \partial_t \psi (t,x)} \left\{ e^{it P} u (t,x) - G (N(u))(t,x) \right\} d {\rm vol}_g dt = 0 $$
hence that $  e^{it P} u (t,x) - G (N(u))(t,x)  $ is independent of $t$. By evaluation at $ t = 0 $, we find
$$ e^{itP} \left( u(t) - \frac{\sigma}{i} D (N(u))(t) \right) = u_0 , \qquad t \in \Ra ,$$
since $ e^{-itP} G (N(u))(t) = \frac{\sigma}{i} D (N(u))(t) $. This proves
 that $ u = F_{u_0}(u) $ and completes the proof. \finpreuve


\appendix

\section{Putting the metric in normal form} \label{appendiceformenormale}
\setcounter{equation}{0}

\begin{prop} If $ ({\mathcal M}, G) $ is asymptotically conic, $G$ can be put in normal form.
\end{prop}

\noindent {\it Proof.} The main steps are described in \cite{HTW}, but locally with respect to the angular variable. We briefly describe here how to globalize the construction on $ {\mathcal S} $.
It is sufficient to prove the existence of sequences of compact subsets  $ {\mathcal K}_k \Subset {\mathcal M} $, real numbers $ R_k > 0 $ and diffeomorphisms
$$ \Omega_k : {\mathcal M} \setminus {\mathcal K}_k \ni m \mapsto \big( r_k (m) , \omega_k (m) \big) \in (R_k,\infty) \times {\mathcal S} , $$
 with $ r_k / r $  bounded from above and below on $ {\mathcal M} \setminus {\mathcal K}_k $ (so that preimages of  bounded intervals by $ r_k $ are relatively compact in $ {\mathcal M} $),
  through which $ G = \Omega^*_k \big(  A_k(r_k)dr_k^2 + 2 r_k B_k (r_k) dr_k + r^2_k g_k (r_k) \big)  $ with
\begin{eqnarray}
  A_k(\cdot)-1 \in S^{-k \nu} , \qquad  B_k (\cdot) \in S^{- k \nu}
\qquad  g_k (\cdot) -  \bar{g}  \in  S^{-\nu } . \label{decayk}
\end{eqnarray}
If we achieve this,  then in a finite number of steps we have $ k \nu > 1 $ and can put the metric in normal form by using \cite{BoucletOsaka}. We  proceed by induction by setting first $ \Omega_1 = \Omega $. We seek  $ \Omega_ k = D_k^{-1} \circ \Omega_{k-1} $, between suitable open subsets of $ \Ra_x \times {\mathcal S} $, by constructing  a diffeomorphism of the form
$$ D_k (x , \omega) = \big( x +x  \sigma_k (x,\omega) , \exp_{\omega} ( V_k (x) ) \big) $$
for some symbol $  \sigma_k$ and some $x$ dependent vector field $ V_k (x) $ on $ {\mathcal S} $.
For $ R_k $ large enough, we define  $ \sigma_k $ and then $ V_k  $ on $ (R_k,\infty) \times {\mathcal S} $ as the unique solutions in $ S^{(1-k) \nu} $ to 
\begin{eqnarray}
 2 \big( x \partial_x \sigma_k + \sigma_k \big) = 1- A_{k-1} (x), \qquad  x \partial_x V_k (x) = - \bar{g}^{-1} \big( d_{\omega} \sigma_{k}(x) + B_{k-1}(x) \big) ,  \label{balance}
\end{eqnarray}
where $ \bar{g}^{-1} $ stands for the isomorphism $ T^* {\mathcal S} \rightarrow T {\mathcal S} $ induced by $\bar{g}$, and $ d_{\omega} $ is the differential on $ {\mathcal S} $.  These objects are globally defined with respect to the angular variable on $ {\mathcal S} $. Note in particular that, since $ V (x) \rightarrow 0 $ as $ x \rightarrow \infty $, $ \exp_{\omega} (V (x)) $ is close to the identity on $ {\mathcal S} $. It is then not hard to check that, for $ R_k $ large enough, $ D_k $ is a diffeomorphim between $ (R_k , \infty) \times {\mathcal S} $ and an open subset of $ (R_{k-1}, \infty) \times {\mathcal S} $ which contains $ (\tilde{R}_{k-1} , \infty) \times {\mathcal S} $ for some $ \tilde{R}_{k-1} $ large enough. We find that
$$ D_k^* \big( G_{k-1}(r_{k-1})  \big) = A_k(r_k)dr_k^2 + 2 r_k B_k (r_k) dr_k + r^2_k g_k (r_k)  $$
with
\begin{eqnarray*}
A_k (r_k) & = & 1 + 2 \sigma_k (r_k) + r_k \partial_{r_k} \sigma_k (r_k) + ( A_{k-1} - 1 )(r_k) + S^{- k \nu} \\
B_k (r_k)& = &  \bar{g} \big( r_k \partial_{r_k} V_k (r_k) \big) +  B_{k-1} (r_k) + d_{\omega } \sigma_k (r_k) + S^{- k \nu}  \\
g_k (r_k)& = & \bar{g} + S^{-\nu} .
\end{eqnarray*}
By (\ref{balance}), we see that (\ref{decayk}) is satisfied. Furthermore, the form of $ D_k $ implies that $ r_k / r_{k-1} $ is bounded from above and below, so by the induction assumption on $ r_{k-1} $ the same holds for $ r_k / r $.  The result follows. \finpreuve



\section{Weak type $ (1,1) $ estimates} \label{sectionCalderonZygmund}
\setcounter{equation}{0}
In this section, we explain how to reduce the proof of  weak type $ (1,1) $ estimates  on $ L^1 ({\mathcal M}) $ for the operators of Propositions \ref{singularintegral} and \ref{squarehautefrequence}
 to the  standard theory of Calder\'on-Zygmund operators on $ \Ra^n $ (Theorem \ref{thCalderonZygmundusuel} below).
 
 We first recall  some general and elementary facts. Assume that $ {\mathcal X} $ is a manifold equipped with a measure $ \mu $ which is a positive smooth density. We recall that a linear map $ T $ on $ L^1 ({\mathcal X},\mu) $ (with values on measurable functions on $ {\mathcal X} $) is said to be of weak type $(1,1)$ with bound $ C $ if
$$ \mu  \big( \left\{ |Tf| > \lambda  \right\} \big) \leq \frac{C}{\lambda} ||f||_{L^1({\mathcal X}, \mu)} $$
for all $ \lambda > 0 $ and $ f \in L^1 ({\mathcal X},\mu) $.
\begin{prop} \label{propositionsimpleweak11}  Let $ T $ be of weak type $ (1,1) $ on $ L^1 ({\mathcal X},\mu) $ with bound $C$.
 \begin{enumerate} \item{Let $ b : {\mathcal X} \rightarrow [ m,M ] $, with $ 0 < m < M $, be  measurable  and let $ \mu_b $ be the measure defined by
$$ \mu_b (B) := \int_B b  d \mu . $$
Then $T$ is of weak type $ (1,1) $ on $ L^1 ({\mathcal X},\mu_b) $ with bound $CM/m$. }
\item{Let  $ \Phi  : {\mathcal X} \rightarrow {\mathcal Y} $ be a diffeomorphism between $ {\mathcal X} $ and another manifold $ {\mathcal Y} $. 
\begin{enumerate} \item{Then $ \Phi_* T \Phi^* $ is of weak type $ (1,1) $ on $ L^1 ({\mathcal Y}, \Phi_* \mu) $ with bound $ C $.}
\item{If $ T $ is bounded on $ L^2 ({\mathcal X},\mu)$  (but not necessarily of weak type $ (1,1) $), then  $ \Phi_* T \Phi^* $ is bounded on $ L^2 ({\mathcal Y},\Phi_* \mu) $ with the same operator norm.}
\end{enumerate}
}
\end{enumerate}
\end{prop}
In this proposition, $ \Phi_* \mu $ is the usual pushforward measure ({\it i.e.} $ \Phi_* \mu (B) = \mu \big( \Phi^{-1}(B) \big) $) and $ \Phi_*,\Phi^* $ are respectively the pushforward and pullback operators ({\it i.e.} $ \Phi_* v = v \circ \Phi^{-1} $ and $ \Phi^* f = f \circ \Phi $).

We will apply Proposition \ref{propositionsimpleweak11} to prove the weak type $ (1,1) $ bounds  stated in the proofs of Propositions \ref{singularintegral} and \ref{squarehautefrequence}, that is for operators of the form
$$ T_{\rm low} (M,t) : =\sum_{ \ell = 0}^M \varrho_\ell (t)  {\mathscr D}_{\varepsilon} \Pi_{\kappa}   O \! p_{1}  \big( a_{\epsilon} \big) \psi  \Pi_{\kappa}^{-1}  {\mathscr D}_{\varepsilon}^{-1}  , \qquad \epsilon^2 = 2^{-\ell}, $$
and
$$ T_{\rm high}(M,t) :=   \sum_{ \ell = 1}^M \varrho_\ell (t) \Pi_{\kappa}   O \! p_{h}  \big( a_h \big) \psi \Pi_{\kappa}^{-1} , \qquad h^2 = 2^{-\ell} . $$
We recall that $ \Pi_{\kappa} $ is  associated to the angular chart $ \kappa : U \rightarrow V $ by (\ref{notationpullpush}), $ \psi $ is a smooth cutoff supported in $ (R_0,\infty) \times V $ and that $ a_{\epsilon} , a_h$ are symbols of the form
$$  b \left( r , \theta , \rho , \frac{\eta}{r} \right), $$
with $ b (r,\theta,\xi) \in S^0 $ (possibly depending on $ \epsilon $ or $h$ in a bounded fashion) supported in $ (R_0,\infty) \times K \times \{ c \leq |\xi| \leq C \} $ for some $ K \Subset V $ and $ C > c > 0 $ independent of $ \epsilon $ or $h$.

We proceed as follows.
When $ {\mathcal X} = {\mathcal M} $ and $ \mu $ is the Riemannian measure $ |g(r,\theta)|r^{n-1} dr d\theta $, the item 2 (a) with $ \Phi = \Pi_{\gamma} $ allows to transfer the analysis from $ {\mathcal M} $ to a chart $ (R,\infty) \times V $ equipped with the measure $ |g(r,\theta)| r^{n-1} dr d d \theta $. The item 1 allows to drop the factor $ |g(r,\theta)| $. We next introduce the diffeomorphism 
$$ \Phi (r,\theta) := (r, r \theta) $$
between $ \Ra_+ \times \Ra^{n-1}_{\theta} $ and $ \Ra_+ \times \Ra^{n-1}_z $, whose interest is that
$$ \Phi_* \big( r^{n-1} d r d \theta \big) = dr dz . $$
Then another application of the item 2 (a) shows that it suffices  that
$$  A_{\rm low} (M,t) := \Phi_* \Pi_{\kappa}^{-1} \left( T_{\rm low} (M,t)\right) \Pi_{\kappa} \Phi^* , \qquad A_{\rm high} (M,t) := \Phi_* \Pi_{\kappa}^{-1} \left( T_{\rm high} (M,t)\right) \Pi_{\kappa} \Phi^* ,  $$
satisfy weak type $ (1,1) $ estimates on $ L^1 (\Ra^n,drdz) $. To prove the latter, it suffices to check they
satisfy the assumptions of the following theorem.

\begin{theo}[Calder\'on-Zygmund operators] \label{thCalderonZygmundusuel} Let $( A_M ) $ be a sequence of operators on $  \Ra_r \times \Ra^{n-1}_{z} $ with Schwartz kernel $ K_M $ such that, for some $ C > 0 $ and all $ M $,
$$ || A_M ||_{L^2(\Ra^n ,drdz) \rightarrow L^2 (\Ra^n,dr dz)} \leq C, \qquad M \geq 0, $$
and, for any $j,\alpha$ such that $ j + |\alpha| \leq 1 $,
$$ | \partial_{r^{\prime}}^j \partial^{\alpha}_{z^{\prime}} K_M (r,z,r^{\prime},z^{\prime})| \leq C (|r-r^{\prime}|+|z-z^{\prime}|)^{-n-j-|\alpha|}, \qquad (r,z,r^{\prime},z^{\prime}) \in \Ra^{2n}, \ \ M \geq 0 . $$
Then $ A_M $ is of weak type $ (1,1) $ on $ L^1 (\Ra^n,drdz) $ with bound uniform in $ M $.
\end{theo}

We refer for instance to  \cite{Taylor} for a proof of this theorem.

The uniform $ L^2 (dr dz) $ boundedness of $ A_{\rm low} (M,t) $ follows from the item 2 (b) of Proposition \ref{propositionsimpleweak11} together with the Cotlar-Stein argument described in the proof of Proposition \ref{singularintegral}. For $ A_{\rm high}(M,t) $, it suffices to observe that 
$$ \sum_{\ell=1}^M \varrho_\ell (t) a (h) \in \widetilde{S}^{0,0} , $$
uniformly in $ M $ and $t$. This follows from the form of $ a (h) $. Therefore $ T_{\rm high}(M,t) $ is uniformly bounded on $ L^2 ({\mathcal M}) $ so $ A_{\rm high}(M,t) $ is uniformly bounded on $ L^2 (drdz) $ by the item 2 (b) of Proposition \ref{propositionsimpleweak11}.

We next consider the kernel estimates.
To put both cases under a single form, we compute the Schwartz kernel of
$$  A_{\epsilon}^h := \Phi_*  {\mathscr D}_{\epsilon}  O \! p_{h}  \big( a \big) \psi {\mathscr D}_{\epsilon}^{-1}   \Phi^* $$
with respect to $ dr dz $, with
$$ a (r,\theta,\rho,\eta) = b \left( r , \theta , \rho , \frac{\eta}{r} \right), \qquad b \in S^{-\infty} . $$
The Schwartz kernel of $ O \! p_{h} (a) $ with respect to $ d r d \theta $ is of the form
$$ (2 \pi h)^{-n} \hat{b} \left( r , \theta , \frac{r-r^{\prime}}{h} , \frac{r (\theta - \theta^{\prime})}{h}  \right) $$
where $ \hat{b} $ is the Fourier transform with respect to $ (\rho,\eta) $. After elementary calculations, we find that the Schwartz kernel of $ A_{\epsilon}^h $ reads (up to the irrelevant factor $ (2 \pi)^{-n} $)
$$ K_{\epsilon}^h (r,z,r^{\prime},z^{\prime}) = \left( \frac{\epsilon}{h} \right)^n  \left( \frac{r}{r^{\prime}} \right)^{n-1}  \hat{b} \left(\epsilon  r ,  \frac{z}{r} , \frac{\epsilon}{h}(r-r^{\prime}) , \frac{\epsilon}{h} \big( z - (r/r^{\prime}) z^{\prime} \big)  \right) \psi \left( \epsilon r^{\prime} , \frac{z^{\prime}}
{r^{\prime}} \right) . $$  
We want to show that $ \sum \varrho_\ell (t) K_1^h $ and $ \sum \varrho_\ell (t) K_{\varepsilon}^1 $ satisfy the second assumption of Theorem \ref{thCalderonZygmundusuel}. By exploiting that $ z^{\prime} / r^{\prime} $ belongs to a compact set, as well and the fact that $ \epsilon r^{\prime} $  is bounded below by some $ R \gg 1 $, these kernel estimates follow from the following lemma which we use either with $ \lambda = h $ or $ \lambda = \epsilon^{-1} $.

\begin{lemm} \label{superposition} \begin{enumerate} \item{ For all $ N \geq 0 $, there exists $ C > 0 $ such that
$$ \left( 1 + \frac{|r-r^{\prime}|}{\lambda} + \frac{\big| z - \frac{r}{r^{\prime}} z^{\prime} \big|}{\lambda} \right)^{-3N} \left(1+\frac{|z^{\prime}|}{r^{\prime}} \right)^{-N} \leq C \left( 1 + \frac{|r-r^{\prime}|}{\lambda} + \frac{| z - z^{\prime} |}{\lambda} \right)^{-N} $$
for all $ \lambda > 0 $, all $ r ,r^{\prime} > 0 $ and all $ z , z^{\prime} \in \Ra^{n-1}. $ }
\item{Let $ c > 0 $. There exists $ C > 0 $ such that, for all $ r , r^{\prime} > 0 $ and $  \lambda > 0$, we have
$$ \frac{r}{r^{\prime}} \leq C \left( 1 + \frac{|r-r^{\prime}|}{\lambda} \right) $$ provided that
  \begin{eqnarray}
  \frac{r^{\prime}}{\lambda} \geq c . 
  \nonumber
 \end{eqnarray} 
  }
\item{ Let $ s\in [0,1]  $ and $ N > n+1 $. Then
\begin{eqnarray*}
  \sum_{\lambda = 2^\ell  \atop \ell \in \Za} \lambda^{-n-s} \left(1+ \frac{|x-y|}{\lambda} \right)^{-N} & \lesssim &   |x-y|^{-n-s} 
 \end{eqnarray*} 
for all $ x , y \in \Ra^n$ such that $ x \ne y $. }
\end{enumerate}
\end{lemm}

\noindent {\it Proof.}  In the item 1, the left hand side is not greater than
$$  \left( 1 + \frac{|r-r^{\prime}|}{\lambda}  \right)^{-2N}   \left( 1 + \frac{\big| z - \frac{r}{r^{\prime}} z^{\prime} \big|}{\lambda} \right)^{-N} \left(1+ \frac{|z^{\prime}|}{r^{\prime}} \right)^{-N} . $$ 
 Writing $ z - \frac{r}{r^{\prime}} z^{\prime} = z - z^{\prime} + \frac{r^{\prime}-r}{r^{\prime}} z^{\prime} $ and using the Peetre inequality for the term in the middle, we obtain an upper bound of the form
$$ C  \left( 1 + \frac{|r-r^{\prime}|}{\lambda}  \right)^{-2N}   \left( 1 + \frac{\big| z -  z^{\prime} \big|}{\lambda} \right)^{-N} \left(1+ \frac{|r-r^{\prime}|}{\lambda} \frac{|z^{\prime}|}{r^{\prime}} \right)^N \left(1+ \frac{|z^{\prime}|}{r^{\prime}} \right)^{-N} , $$  
which in turn is bounded by
$$ C  \left( 1 + \frac{|r-r^{\prime}|}{\lambda}  \right)^{-2N}   \left( 1 + \frac{\big| z -  z^{\prime} \big|}{\lambda} \right)^{-N} \left(1+ \frac{|r-r^{\prime}|}{\lambda} \right)^N  . $$
 This yields the result once observed that
$$ \left( 1 + \frac{|r-r^{\prime}|}{\lambda}  \right)^{-N}   \left( 1 + \frac{\big| z -  z^{\prime} \big|}{\lambda} \right)^{-N} \leq \left( 1 + \frac{|r-r^{\prime}|}{\lambda}  + \frac{\big| z -  z^{\prime} \big|}{\lambda} \right)^{-N} . $$
The item 2 follows simply from the fact that $ \frac{r}{r^{\prime}} = 1 + \frac{r-r^{\prime}}{\lambda} \frac{\lambda}{r^{\prime}} $.
  The item 3 is standard. \finpreuve

\section{Sobolev estimate} \label{sectionSobolevhomogene}
\setcounter{equation}{0}
In this section we provide a short proof of the homogeneous Sobolev estimate (\ref{Sobolevestimate}).

Using the same cutoff $ f_0 $ as in (\ref{defhighlow}), we have
$$  || (1 - f _0) (P) v ||_{L^{2^*}({\mathcal M})} \lesssim  || (P+1)^{1/2} (1 - f _0) (P) v ||_{L^{2}({\mathcal M})} \lesssim || P^{1/2} v ||_{L^2 ({\mathcal M})}  $$
thanks to  the inhomogeneous Sobolev estimate (see {\it e.g.}  \cite{BoucletFourier}) 
\begin{eqnarray}
 || u ||_{L^{2^*}({\mathcal M})} \lesssim || (P+1)^{1/2} u ||_{L^2 ({\mathcal M})}   \label{Sobolevinhomogenestandard}
 \end{eqnarray}
and  the spectral theorem.
 Thus we have to show that
$$ || f _0 (P) v ||_{L^{2^*}({\mathcal M})} \lesssim || P^{1/2} v ||_{L^2 ({\mathcal M})} . $$
To do so, we choose $ \chi \in C_c^{\infty} ({\mathcal M}) $ which is equal to $1$ on a large enough compact set and 
observe  that
$$ || \chi f_0 (P) v ||_{L^{2^*}({\mathcal M})} \lesssim || \scal{r}^{-1} v ||_{L^2 ({\mathcal M})} \lesssim || P^{1/2} v ||_{L^2 ({\mathcal M})}  $$
using first that $ \chi f_0 (P) \scal{r} $ is bounded from $ L^2 $ to $ L^{2^*} $ (which follows from (\ref{Sobolevinhomogenestandard}) and a standard commutator argument) and then the Hardy inequality (see {\it e.g.} \cite[Prop. 2.2]{BoucletRoyer}). Using a partition of unity $ \sum_{\kappa} \varphi_{\kappa}(\omega) = 1 $ on $ {\mathcal S} $ with functions supported in coordinates patches, we can see that
$$  || (1-\chi) \varphi_{\kappa} (\omega) f_0 (P) v ||_{L^{2^*} ({\mathcal M})}  \lesssim    \big| \big| \nabla_g \big( (1-\chi) \varphi_{\kappa}(\omega) f (P_0) v \big) \big| \big|_{L^2 ({\mathcal M})}  $$
using  the usual proof of the Sobolev inequality on $ \Ra^n $ since the cutoff $ (1-\chi) \varphi_{\kappa}(\omega) $ localizes in the product of a half line  and a patch. From this estimate, we then obtain
\begin{eqnarray*}
 || (1-\chi) f_0 (P) v ||_{L^{2^*} ({\mathcal M})} 
 & \lesssim &  || \nabla_g  f (P_0) v ||_{L^2 ({\mathcal M})} + || \scal{r}^{-1} f (P_0) v ||_{L^2 ({\mathcal M})} \\
 & \lesssim & || P^{1/2} v ||_{L^2 ({\mathcal M})} 
\end{eqnarray*}
using again the Hardy inequality.

\bigskip

 {\sc Institut de Math\'ematiques de Toulouse (UMR CNRS 5219),  Universit\'e Paul Sabatier, 118 route de Narbonne, F-31062 Toulouse FRANCE }

{\it E-mail address:} \verb+ jean-marc.bouclet@math.univ-toulouse.fr  +

\medskip

 {\sc Department of Mathematics, Osaka University Toyonaka, Osaka 560-0043, JAPAN}

{\it E-mail address:} \verb+ haruya@math.sci.osaka-u.ac.jp  +
\end{document}